\documentclass[preprint]{elsarticle}

\usepackage{a4}
\usepackage[all]{xy}

\usepackage{amsmath}
\usepackage{amsfonts}
\usepackage{amssymb}
\usepackage{latexsym}
\usepackage{stmaryrd}
\usepackage{array}
\usepackage{exscale}
\newcommand{\nc}{\newcommand}
\newcommand{\ol}{\overline}
\newcommand{\ul}{\underline}
\newcommand{\es}{\emptyset}
\newcommand{\sm}{\setminus}
\newcommand{\ve}{\varepsilon}
\newcommand{\vp}{\varphi}

\newcommand{\bc}{\bigcup}
\newcommand{\bca}{\bigcap}
\newcommand{\Lra}{\Leftrightarrow}

\newcommand{\Lla}{\Longleftarrow}

\newcommand{\Ra}{\Rightarrow}
\newcommand{\La}{\Leftarrow}
\newcommand{\ra}{\rightarrow}

\newcommand{\sse}{\subseteq}

\newcommand{\spe}{\supseteq}
\newcommand{\fa}{\forall}
\newcommand{\ex}{\exists}
\newcommand{\mr}{\mathrm}
\newcommand{\mc}{\mathcal}
\newcommand{\mf}{\mathfrak}

\newcommand{\DMO}{\DeclareMathOperator}
\newcommand{\DST}{\displaystyle}

\newcommand{\NN}{\mathbb{N}}
\newcommand{\NNZ}{\NN_0}
\newcommand{\ZZ}{\mathbb{Z}}
\newcommand{\QQ}{\mathbb{Q}}
\newcommand{\RR}{\mathbb{R}}

\newcommand{\UU}{\mathbb{U}}

\mathchardef\breakingcomma\mathcode`\,
{\catcode`,=\active
  \gdef,{\breakingcomma\discretionary{}{}{}}
}

\usepackage{listings}
\newcommand{\inl}[1]{\lstinline$#1$}
\newcommand{\und}{{\:\wedge\:}} 
\newcommand{\oder}{{\:\vee\:}} 
\newcommand{\mb}{{\:|\:}} 
\newcommand{\set}[1]{\{ #1 \}}
\newcommand{\setb}[1]{\big \{ \, #1 \, \big \}}
\DeclareMathOperator{\dom}{dom}

\nc{\simlvi}[1]{\!\sim_{#1}}
\newcommand{\addcup}{\mathbin{\mathaccent\cdot\cup}} 
\nc{\apprel}[3]{{#1}(#2)_{(#3)}} 
\nc{\cmpli}[1]{\complement^1_{#1}} 
\nc{\cmplzi}[1]{\complement^0_{#1}} 
\nc{\cmplzoi}[1]{\complement^*_{#1}} 

\nc{\zf}{\mr{ZF}}
\nc{\zfmf}{\zf^0} 
\nc{\zfc}{\mr{ZFC}}
\nc{\zfcmf}{\zfc^0} 
\nc{\bst}{\mr{BST}} 
\newcommand{\nni}{\NNZ \cup \{+\infty\}} 
\newcommand{\nnpi}{\NN \cup \{+\infty\}} 
\newcommand{\tb}[2]{\set{#1, \dots, #2}} 
\providecommand{\abs}[1]{\lvert #1 \rvert} 
\providecommand{\norm}[1]{\lVert #1 \rVert} 
\DeclareMathOperator{\ld}{ld} 
\DeclareMathOperator{\fld}{fld} 
\DeclareMathOperator{\rank}{rank} 
\DeclareMathOperator{\Q}{\mf{Q}} 
\makeatletter
\DeclareRobustCommand{\genericinterval}[2]{%
  \@ifstar{\genericinterval@star{#1}{#2}}{\genericinterval@nostar{#1}{#2}}}
\newcommand{\genericinterval@star}[4]{\mathopen{}\mathclose{\left#1#3,#4\right#2}}
\newcommand{\genericinterval@nostar}[4]{\mathopen{#1}#3,#4\mathclose{#2}}

\makeatother
\nc{\untit}[2]{{#1}^{#2 \downarrow}} 
\nc{\obit}[2]{{#1}^{#2 \uparrow}} 
\nc{\inzEKi}[1]{\mc{I}^{\mr{V}}_{#1}}

\nc{\inzKEi}[1]{\mc{I}^{\mr{E}}_{#1}}

\nc{\adjEi}[1]{\mc{A}^{\mr{V}}_{#1}}

\DeclareMathOperator{\Tr}{Tr} 
\nc{\BD}[1]{{#1}\text{-}\mr{BD}}

\nc{\konv}[2]{{#1}[{#2}]} 
\nc{\actpres}[1]{\Phi_{#1}} 
\newcommand{\floor}[1]{\lfloor #1 \rfloor}

\nc{\Prim}{\mc{PR}} 

\nc{\sselr}{\sse^{\mapsto}}
\nc{\sserl}{\sse^{\mapsfrom}}
\nc{\spelr}{\spe^{\mapsto}}
\nc{\sperl}{\spe^{\mapsfrom}}
\nc{\ball}[1]{\mr{B}^{#1}} 
\nc{\oball}[1]{\breve{\mr{B}}^{#1}} 
\nc{\pball}[1]{\dot{\mr{B}}^{#1}} 
\nc{\prr}[1]{\dot{\RR}^{#1}} 
\nc{\sph}[1]{\mr{S}^{#1}} 
\nc{\ssim}[1]{s\sigma_{#1}} 
\nc{\koerper}[1]{\norm{#1}}
\nc{\Ccovdim}{\mc{CD}}
\nc{\Cinddim}{\mc{SID}}

\nc{\CInddim}{\mc{LID}}

\DeclareMathOperator{\diffop}{D} 
\DeclareMathOperator*{\diffoplimit}{D} 
\nc{\diffopc}[1]{\sideset{_{#1}}{}\diffoplimit} 
\nc{\diffopp}[1]{\diffop_{#1}} 
\nc{\diffopcp}[2]{\sideset{_{#2}}{_{#1}}\diffoplimit} 
\nc{\meanH}[2]{\mf{M}_{#1,#2}} 
\nc{\emean}[2]{\mf{M}_{\exp_{#1},#2}} 
\DeclareMathOperator{\mor}{Mor}
\DeclareMathOperator{\Hom}{Hom} 
\nc{\autoerw}[1]{\mr{Aut}^{#1}} 
\nc{\komma}[2]{(#1 \downarrow #2)} 
\nc{\Kmat}{\mf{MAT}} 
\nc{\Khmat}{\mf{HMAT}} 
\nc{\homfun}[1]{\mor_{#1}(-_1,-_2)} 
\nc{\homfunae}[1]{\mor_{#1}(-_1)} 
\nc{\homfunbe}[1]{\mor_{#1}(-_2)} 
\nc{\homfunxy}[3]{\mor_{#1}(#2(-_1), #3(-_2))}
\nc{\homfunx}[2]{\mor_{#1}(#2(-_1), -_2)}
\nc{\homfuny}[2]{\mor_{#1}(-_1, #2(-_2))}
\nc{\homfuna}[2]{\mor_{#1}(#2, -)} 
\nc{\homfunb}[2]{\mor_{#1}(-, #2)} 
\nc{\hhomfuna}[2]{\Hom_{#1}(#2, -)} 
\nc{\hhomfunb}[2]{\Hom_{#1}(-, #2)} 
\newcommand{\Va}{\mc{V\hspace{-0.1em}A}}

\newcommand{\Lit}{\mc{LIT}}
\newcommand{\Cl}{\mc{CL}}
\newcommand{\Cls}{\mc{CLS}}

\newcommand{\Pcls}[1]{#1\mbox{--}\Cls}

\newcommand{\Pass}{\mc{P\hspace{-0.32em}ASS}}
\newcommand{\epa}{\pab{}} 

\newcommand{\Sat}{\mc{SAT}}

\newcommand{\Usat}{\mc{USAT}}

\newcommand{\Musat}{\mc{M\hspace{0.8pt}U}} 
\newcommand{\Musati}[1]{\Musat_{\!#1}} 
\newcommand{\Smusat}{\mc{S}\Musat} 
\newcommand{\Smusati}[1]{\Smusat_{\!#1}}

\newcommand{\Vmusat}{\mc{V}\Musat} 
\newcommand{\Vmusati}[1]{\Vmusat_{\!#1}}
\newcommand{\symu}{\bar{\mc{U}}} 
\newcommand{\Unicls}{\symu\mc{CLS}} 
\newcommand{\Punicls}[1]{{#1}\mbox{-}\Unicls} 
\newcommand{\Unimusat}{\symu\mc{MU}} 
\newcommand{\Punimusat}[1]{{#1}\mbox{-}\symu\mc{MU}} 
\nc{\Clsoo}{\Cls^{1,1}} 
\DeclareMathOperator{\lit}{lit}
\DeclareMathOperator{\var}{var}

\DMO{\dos}{ds} 
\DMO{\mdos}{mds} 
\newcommand{\Clash}{\mc{HIT}} 
\newcommand{\Clashi}[1]{\Clash_{\!\!#1}}
\newcommand{\Uclash}{\mc{U}\Clash} 
\newcommand{\Uclashi}[1]{\Uclash_{\!\!#1}}

\DeclareMathOperator{\hdef}{\delta_{\mr{h}}} 

\newcommand{\Lsat}{\mc{L}\Sat}

\newcommand{\Lean}{\mc{LEAN}}
\newcommand{\Leani}[1]{\Lean_{\!#1}}
\newcommand{\Llean}{\mc{L}\Lean}

\newcommand{\Mlean}{\mc{M}\Lean}
\newcommand{\Mleani}[1]{\Mlean_{\!#1}}
\newcommand{\Msat}{\mc{M}\Sat}

\newcommand{\Mcls}{\mc{M}\Cls} 

\newcommand{\Sed}{\mc{SED}} 
\DeclareMathOperator{\res}{\diamond} 
\DeclareMathOperator{\dpl}{DP} 
\newcommand{\dpi}[1]{\dpl_{\!#1}}
\DMO{\premr}{ax} 
\DMO{\concr}{C} 
\DMO{\allcr}{cl} 

\DMO{\thardness}{thd} 
\DMO{\phardness}{phd} 
\DMO{\whardness}{awid} 
\DMO{\dep}{dep} 
\DMO{\hts}{hs} 
\DMO{\semspace}{css} 
\DMO{\resspace}{crs} 
\DMO{\treespace}{cts} 
\newcommand{\pab}[1]{\langle #1 \rangle}
\newcommand{\pao}[2]{\langle #1 \ra #2 \rangle}

\nc{\bth}[1]{\langle{#1}\rangle} 
\newcommand{\A}{\mc{A}} 

\DMO{\rsub}{r_S} 
\DMO{\rk}{r} 
\DMO{\ro}{\rk_1} 
\DMO{\rki}{\rk_{\infty}} 
\DMO{\rpl}{r^{pl}} 
\DMO{\ropl}{\rk_1^{pl}} 
\nc{\rslur}{\xrightarrow{\text{SLUR}}} 
\nc{\rslurs}{\rslur_{\!*}} 
\DMO{\slur}{slur} 
\nc{\Slur}{\mc{SLUR}} 
\nc{\rkslur}[1]{\xrightarrow{\text{SLUR}_{#1}}} 
\nc{\rkslurs}[1]{\rkslur{#1}_{\!*}} 
\nc{\Altsluri}[1]{\Slur(#1)}
\nc{\Altslurstari}[1]{\Slur\text{\textasteriskcentered}(#1)}
\nc{\Canoni}[1]{\mr{CANON}(#1)}
\nc{\rkslurstar}[1]{\xrightarrow{\text{SLUR\textasteriskcentered}#1}} 
\nc{\rkslursstar}[1]{\rkslurstar{#1}_{\!*}} 
\DMO{\slurstar}{\slur\!\text{\textasteriskcentered}}
\nc{\Urefc}{\mc{UC}}
\nc{\Propc}{\mc{PC}}
\nc{\Wrefc}{\mc{WC}} 
\DeclareMathOperator{\scf}{CM} 
\DeclareMathOperator{\cmg}{cmg} 
\DeclareMathOperator{\bcp}{bcp} 
\DeclareMathOperator{\ldeg}{ld} 
\DeclareMathOperator{\minldeg}{\mu\!\ldeg} 
\DeclareMathOperator{\vdeg}{vd} 
\DeclareMathOperator{\minvdeg}{\mu\!\vdeg} 
\DeclareMathOperator{\maxvdeg}{\nu\!\vdeg} 
\DMO{\varmvd}{\var_{\minvdeg}} 
\DMO{\nfc}{fc} 
\DMO{\maxnfc}{\nu\!\nfc} 
\nc{\Dt}[1]{\mc{F}_{#1}} 

\DeclareMathOperator{\surp}{\sigma} 
\DeclareMathOperator{\nonmer}{nM} 
\newcommand{\inonmer}{\operatorname{i}_{\mathrm{nM}}} 
\nc{\svbf}{\mc{VB}} 
\nc{\svbfs}{\mc{VB}^*} 
\DMO{\potp}{pp} 
\DMO{\potprec}{NM} 
\DMO{\minnonmer}{VDM} 
\DMO{\minnonmerh}{VDH} 
\DMO{\maxsmar}{FCM} 
\DMO{\maxsmarh}{FCH} 
\DMO{\varsing}{\var_s} 
\DMO{\varosing}{\var_{1s}} 
\DMO{\varnosing}{\var_{\neg1s}} 
\DMO{\nsv}{\mathit{n}_s} 
\DMO{\nosv}{\mathit{n}_{1s}}
\DMO{\nnosv}{\mathit{n}_{\neg1s}}
\nc{\Musatns}{\Musat'} 
\nc{\Musatnsi}[1]{\Musati{#1}'}
\nc{\Smusatns}{\Smusat'} 
\nc{\Smusatnsi}[1]{\Smusati{#1}'}
\nc{\Uclashns}{\Uclash'} 
\nc{\Uclashnsi}[1]{\Uclashi{#1}'}
\nc{\tsdp}{\xrightarrow{\text{sDP}}}
\nc{\tsdps}{\tsdp_{\!*}}
\nc{\tosdp}{\xrightarrow{\text{1sDP}}}
\nc{\tosdps}{\tosdp_{\!*}}
\DMO{\sdp}{sDP} 
\DMO{\osdp}{sDP_1} 
\nc{\cflmusat}{\mc{CF}\Musat} 
\nc{\cflmusati}[1]{\mc{CF}\Musati{#1}}
\nc{\cflimusat}{\mc{CFI}\Musat} 
\DMO{\sNF}{sNF} 
\DMO{\eqp}{eqp} 
\DMO{\sgp}{sp} 
\DMO{\singind}{si} 
\DMO{\osingind}{si_1} 
\DMO{\shyp}{svh} 
\DMO{\sdph}{ssh} 
\DMO{\msdph}{mss} 
\DMO{\osdph}{ssh_1} 
\DMO{\mosdph}{mss_1} 
\DMO{\mps}{mps} 
\DMO{\purec}{puc} 
\DMO{\doping}{D}
\nc{\glue}[4]{\operatorname{glue}((#1,#2), (#3,#4))} 
\nc{\gluea}[3]{#1 \mathbin{\boxplus}_{#3} #2} 
\DMO{\saturate}{S}
\DMO{\marginalise}{M}
\DMO{\frl}{fl} 
\nc{\Con}{\mr{Con}}
\nc{\Log}{\mr{Log}}
\nc{\Lin}{\mr{Lin}}
\nc{\Pol}{\mr{Pol}}
\nc{\ExL}{\mr{ExL}}
\nc{\ExP}{\mr{ExP}}
\nc{\CTime}{\mr{CTime}}
\nc{\CSpace}{\mr{CSpace}}
\nc{\LTime}{\mr{LTime}}
\nc{\LSpace}{\mr{L}}
\nc{\NLSpace}{\mr{NL}}
\nc{\LinTime}{\mr{LinTime}}
\nc{\LinSpace}{\mr{LinSpace}}
\nc{\PTime}{\mr{P}}
\nc{\PSpace}{\mr{PSpace}}
\nc{\Np}{\mr{NP}}
\nc{\Conp}{\text{coNP}}
\nc{\NPSpace}{\mr{NPSpace}}
\nc{\CoNPSpace}{\mr{coNPSpace}}
\nc{\ELTime}{\mr{ELTime}}
\nc{\ELSpace}{\mr{ELSpace}}
\nc{\EPTime}{\mr{EPTime}}
\nc{\EPSpace}{\mr{EPSpace}}
\nc{\NEPTime}{\mr{NEPTime}}
\nc{\polydelta}[1]{\Delta_{#1}^{\mr P}}
\nc{\polypi}[1]{\Pi_{#1}^{\mr P}}
\nc{\polysigma}[1]{\Sigma_{#1}^{\mr P}}
\nc{\Ph}{\mr{PH}}

\nc{\Dp}{D^P}
\nc{\PllC}[2]{{\text{$\mr{PT}$/$\mr{WK}$}(#1, #2)}} 
\nc{\Nc}{\mr{NC}}
\nc{\Nci}[1]{\Nc^{#1}}
\nc{\Ac}{\mr{AC}}
\nc{\Aci}[1]{\Ac^{#1}}
\nc{\pmodpoly}{P / \mathrm{poly}}
\nc{\Wh}[1]{\mr{W}[#1]} 
\nc{\Rl}{\mr{RL}}
\nc{\coRl}{\mr{coRL}}
\nc{\Rp}{\mr{RP}}
\nc{\coRp}{\mr{coRP}}
\nc{\Zpp}{\mr{ZPP}}
\nc{\Bpp}{\mr{BPP}}
\nc{\Pp}{\mr{PP}}
\nc{\Reach}{\mr{STCON}} 
\nc{\Undreach}{\mr{USTCON}} 
\nc{\Pcol}[2]{\mr{COL}(#1,#2)} 
\nc{\Pscol}[2]{\mr{SCOL}(#1,#2)} 
\nc{\Psorcol}[2]{\mr{SORCOL}(#1,#2)} 
\DMO{\slp}{slp}
\nc{\Mss}{\mr{MSS}}
\nc{\Key}{\mr{KEY}}
\nc{\Keyi}[1]{\Key_{\!#1}}
\nc{\Nbmss}{N_{\mr{bm}}} 
\nc{\Nbkey}{N_{\mr{bk}}} 
\nc{\Rnb}{N_{\mr{b}}}
\nc{\Rnk}{N_{\mr{k}}}
\nc{\Rnr}{N_{\mr{r}}}

\nc{\Byte}{\mr{B}[8]}
\nc{\QByte}{\mr{B}[4,8]}
\nc{\KByte}{\mc{B}} 
\nc{\RQByte}{\mc{QB}} 

\nc{\ramz}[3]{\mr{ram}_{#1}^{#2}(#3)} 
\nc{\waez}[2]{\mr{vdw}_{#1}(#2)} 
\nc{\gtz}[2]{\mr{grt}_{#1}(#2)} 
\nc{\pdwaez}[2]{\mr{vdw}_{#1}^{\mr{pd}}(#2)} 
\nc{\absfeh}[1]{\delta_{#1}} 
\nc{\relfeh}[1]{\ve_{#1}} 
\usepackage{theorem} 
\usepackage[breaklinks]{hyperref}
\newtheorem{defi}{Definition}[section]
\newtheorem{lem}[defi]{Lemma}
\newtheorem{thm}[defi]{Theorem}
\newtheorem{corol}[defi]{Corollary}

\newtheorem{conj}[defi]{Conjecture}

\newtheorem{examp}[defi]{Example}

\newtheorem{quest}[defi]{Question}

\theorembodyfont{\rmfamily}

\theorembodyfont{}
\newenvironment{prf}{\noindent\textbf{Proof:}\;}{\par\noindent\ignorespacesafterend}

\newcommand{\Qed}{\hfill $\square$}

\newcounter{dDef} 

\newcounter{dLem} 

\newcounter{dThm} 

\newcounter{dPro} 

\newcounter{Beispielzaehler}

\nc{\bm}{\boldmath}
\nc{\bmm}[1]{\mbox{\bm$\DST #1$}}
\nc{\mi}[1]{\bmm{\mathrm{(#1):}} \quad}
\usepackage{enumerate}
\usepackage[active]{srcltx}

\newcommand{\Poscls}{\mc{PCLS}}
\newcommand{\Hyp}{\mc{HYP}}
\newcommand{\Mnc}{\mc{MNC}}
\DMO{\defh}{\delta_H}
\newcommand{\Ihyp}{\mc{IHYP}}

\DMO{\odef}{nA}

\newcommand{\Schrift}{report}
\newcommand{\Zusammenfassung}{We investigate connections between SAT (the propositional satisfiability problem) and combinatorics, around the minimum degree of variables in various forms of redundancy-free boolean conjunctive normal forms (clause-sets). Extensive introductions, overviews, conclusions, examples and open problems are provided.

 Let $\minvdeg(F) \in \NN$ for a clause-set $F$ denote the \emph{minimum variable-degree}, the minimal number of occurrences of a variable. A central result is the upper bound $\surp(F) + 1 \le \minvdeg(F) \le \nonmer(\surp(F)) \le \surp(F) + 1 + \log_2(\surp(F))$ for \emph{lean clause-sets} $F \in \Lean$ in dependency on the \emph{surplus} $\surp(F) \in \NN$. Lean clause-sets are defined by having no non-trivial autarkies (partial assignments satisfying some clauses and not touching the other clauses), and generalise \emph{minimally unsatisfiable clause-sets} (which become satisfiable iff any single clause is removed), i.e., $\Lean \supset \Musat$. For the surplus we have $\surp(F) \le \delta(F) = c(F) - n(F)$, using the deficiency $\delta(F)$ of clause-sets, the difference between the number $c(F)$ of clauses and the number $n(F)$ of variables. And $\nonmer(k) \in \NN$ is the $k$-th ``non-Mersenne'' number, skipping in the sequence of natural numbers all numbers of the form $2^n - 1$. As an application of the upper bound we obtain, that clause-sets $F$ violating $\minvdeg(F) \le \nonmer(\surp(F))$ must have a non-trivial autarky, and thus clauses can be removed satisfiability-equivalently. We obtain a polynomial time autarky reduction (removing the clauses), but where it is open whether the autarky itself can be found in polynomial time.

 We show that the upper bound on $\minvdeg$ is sharp, indeed already on $\Musat \subset \Vmusat \subset \Lean$, the class of \emph{variable-minimally unsatisfiable clause-sets} (these become satisfiable iff all clauses containing any single variable are removed).  That is, $\minvdeg(\Vmusati{\delta=k}) = \nonmer(k)$ for all deficiencies $k \in \NN$, where $\minvdeg(\Vmusati{\delta=k})$ is the maximum of $\minvdeg(F)$ over $F \in \Vmusat$ with $\delta(F) = k$. The determination of $\minvdeg(\Musati{\delta=k}) =: \minnonmer(k)$ seems to be a much more involved question. We show that for $k \le 5$ we have $\minnonmer(k) = \nonmer(k)$, but for $k = 6$ we have $\minnonmer(k) = \nonmer(k) - 1$. Moreover this correction by $-1$ causes further corrections by $-1$ for infinitely many other deficiencies, resulting in the upper-bound function $\nonmer_1: \NN \ra \NN$, an instance of a \emph{generalised non-Mersenne function}, found by a novel recursion scheme.

 We also show various auxiliary results, especially concerning $\Vmusat$.}
\newcommand{\Liste}{conjunctive normal form \sep deficiency \sep minimally unsatisfiable \sep variable minimal unsatisfiable \sep saturated minimal unsatisfiable \sep hitting clause-set \sep disjoint tautology \sep orthogonal tautology \sep lean clause-set \sep autarky \sep surplus \sep matching-lean \sep minimally non-2-colourable hypergraph \sep L-matrix \sep minimally sign-central matrix \sep polynomial time \sep SAT decision \sep non-Mersenne number \sep minimum variable degree \sep singular DP-reduction (variable elimination) \sep subsumption resolution}
\newcommand{\correspauta}{}
\newcommand{\correspautb}{}

\DMO{\ulcls}{cls}
\DMO{\procls}{\ul{cls}}

\newcommand{\Mlcr}{\mc{MLCR}} 

\nc{\tsubres}{\xrightarrow{\text{sfsR}}}
\nc{\tsubresk}[1]{\xrightarrow{\text{sfsR}}_{\!#1}}
\nc{\tsubress}{\tsubresk{*}}

\DMO{\stabpar}{sir} 

\begin{document}

\title{Bounds for variables with few occurrences\\ in conjunctive normal forms}

\author[sw]{Oliver Kullmann\correspauta}
\ead[url]{http://cs.swan.ac.uk/~csoliver}

\author[gh]{Xishun Zhao\fnref{cor2}}
\ead[url]{http://logic.sysu.edu.cn/faculty/zhaoxishun/en/}

\correspautb
\fntext[cor2]{Supported by Grants NSFC 61272059, NSSFC 13\&ZD186 and MOE 11JJD7200020.}

\address[sw]{Computer Science Department, Swansea University, UK}
\address[gh]{Institute of Logic and Cognition, Sun Yat-sen University, Guangzhou, P.R.C.}

\begin{abstract}
  \Zusammenfassung
\end{abstract}

\begin{keyword}
  \Liste
\end{keyword}

\maketitle

\tableofcontents

\section{Introduction}
\label{sec:intro}

In this work we aim at bringing together some aspects of combinatorics with the developing theory of SAT. We concentrate on degree considerations in ``clause-sets'' (conjunctive normal forms as set-systems), which can be considered as generalised hypergraphs, namely hypergraphs with ``polarities''. The general goal is to develop an \emph{understanding} of propositional (un)satisfiability, which corresponds for hypergraphs to an understanding of (non-)2-colourability. The considerations of this \Schrift{}, belonging to the intersection between \href{https://en.wikipedia.org/wiki/Propositional_calculus}{Propositional Logic} and \href{https://en.wikipedia.org/wiki/Extremal_combinatorics}{Extremal Combinatorics}, yield first steps, through the study of basic numerical parameters.

SAT, the prototypical NP-complete problem (\cite{Co71}), took a strong development in the past two decades also regarding (industrial) applications (see the handbook \cite{2008HandbuchSAT} for a recent overview). It is often mainly considered as belonging to complexity theory, algorithms and heuristics (with \cite{DP60,DLL62} the basic papers here), and finally implementations and experimentation (``SAT solvers''). ``Understanding'' SAT in a precise sense is considered to be impossible, and only various investigations on random and approximation structures (including ``islands of tractability'') in general are deemed fruitful. We want to challenge this view, starting to build a new bridge, towards an \emph{understanding of unsatisfiability}. We note here that understanding unsatisfiability seems easier than to understand satisfiability, since unsatisfiability means a form of completion, \emph{all} assignments have been excluded as potential satisfying assignments (``models''), while satisfiability means a lack of such completion. More precisely, we aim at understanding \emph{minimal unsatisfiability}, the building blocks of unsatisfiability --- similar to critical colourability, here removal of any clause renders the clause-set satisfiable. The main first goal of ``understanding'' here is the `Finitely many patterns'' conjecture on the classification of minimally unsatisfiable clause-sets, as discussed in Subsection \ref{sec:concclassmu} --- for a fixed deficiency, there shall be only finitely many basic ``ideas'' to establish unsatisfiability, and the rest is some form of ``trivial embellishment''.

A fundamental question, the subject of this study, is the existence of ``simple'' variables in clause-sets. ``Simple'' here means a variable occurring not very often (i.e., with a low ``degree''). A major use of the existence of such variables is in inductive proofs of properties of minimally unsatisfiable clause-sets, using \emph{splitting} on a variable to reduce $n$, the number of variables, to $n-1$: here it is vital that we have control over the changes imposed by the substitution, and so we want to split on a variable occurring as few times as possible. ``Splitting'' of a clause-set $F$ on variable $v$ means the consideration of the two clause-sets $\pao v0 * F$, $\pao v1 * F$, that is, instantiating variable $v$ by both truth values $0, 1$. A feature of clause-sets is the closure under splitting, and splitting is a major tool for investigations into minimal unsatisfiability. In the remainder of the introduction we give an overview on the results of this \Schrift{} and their context.

In Subsection \ref{sec:introdef} we introduce minimal unsatisfiability and the main complexity parameter, the deficiency, and discuss the first main result, the upper bound for the minimum variable-degree for minimally unsatisfiable clause-sets in the deficiency. Its basic importance for the (beginning of) understanding of minimal unsatisfiability can be seen for example by the fact, that the precise knowledge for small deficiencies is central for determining the structure, as \cite{DDK98} for deficiency $1$ (degree $2$), \cite{KleineBuening2000SubclassesMU} for deficiency $2$ (degree $4$), and \cite{KullmannZhao2016UHitSAT} for deficiency $3$ (degree $5$) apply. The easiest case is deficiency $1$, and indeed, once we know that there must be a variable of degree $2$, then the simple (tree-)structure is fairly easy to deduce (as done in \cite{DDK98,Ku99dKo}).

In Subsection \ref{sec:introsurp} we discuss the generalisation to ``lean'' clause-sets and their ``surplus''. After having given the basic definitions and main results, in Subsection \ref{sec:introfour} we then come back to the situation for minimal unsatisfiability, and present some basic methods and arguments, together with some basic intuition regarding the shape of the number-theoretical upper-bound function. This concludes the introduction into the results of this \Schrift, and in the rest of the introduction, we discuss the background and context. Related work on minimal unsatisfiability is reviewed in Subsection \ref{sec:prelimMU}. The generalisation to lean clause-sets is based on ``autarkies'', which are reviewed in Subsection \ref{sec:introautgen}. In Subsection \ref{sec:relwork} then the connections to combinatorics are discussed: hypergraph colouring, hypergraph transversals, combinatorial matrix theory, and biclique partitions of (multi)graphs, where always SAT is treated in disguise, and autarkies play an important role. The introduction is concluded by going through all main results in Subsection \ref{sec:overv}.

According to the goal of bringing different communities together, we provide and explain much of the relevant background, so that this \Schrift{} is mostly self-contained, and the results cited from the literature can be treated as black-boxes.

\subsection{Deficiency as the main structural parameter}
\label{sec:introdef}

The definition of the class $\Cls$ of ``clause-sets'', and of the class $\Musat \subset \Cls$ of ``minimally unsatisfiable clause-sets'', can be quickly (and precisely) given as follows, using (just) natural numbers as ``variables''. A ``literal'' $x$ is an element of $\ZZ \sm \set{0}$, i.e., a non-zero integer. A  ``clause'' $C$ is a finite set of literals, such that there is no $x \in C$ with $-x \in C$. Using $-L := \set{-x : x \in L}$ for sets $L$ of literals, the ``clash-freeness'' condition for $C$ becomes $C \cap -C = \es$. A ``clause-set'' $F$ is a finite set of clauses, the set of all clause-sets is denoted by $\Cls$. The set $\var(F)$ of variables of $F$ is the set of $v \in \NN = \ZZ_{\ge 1}$ with $\set{v,-v} \cap \bc F \ne \es$. The basic measurements for $F \in \Cls$ are:
\begin{itemize}\setlength{\itemsep}{0pt}
\item the number $c(F) := \abs{F} \in \NNZ$ of clauses of $F$;
\item the number $n(F) := \abs{\var(F)} \in \NNZ$ of variables of $F$;
\item The ``deficiency'' $\delta(F) := c(F) - n(F) \in \ZZ$.
\end{itemize}
The deficiency is only informative when certain (weak) assumptions are made for $F$, and for general $F$ the ``maximal deficiency'' $\delta^*(F) := \max_{F' \sse F} \delta(F') \in \NNZ$ is to be used.
 A clause-set $F$ is ``satisfiable'' if there exists a partial assignment $\vp$, which here in this introduction is just a clause, such that $\vp \cap D \ne \es$ for all $D \in F$.\footnote{The clause $\vp$ is the set of satisfied literals of the corresponding ``partial assignment''. This definition of ``satisfying assignments'', via clauses intersecting every clause of $F$, generalises transversals of hypergraphs, by taking complementation into account ($\vp$ does not contain clashes).} The set of all satisfiable clause-sets is $\Sat \subset \Cls$, the set of all unsatisfiable clause-sets is $\Usat := \Cls \sm \Sat$. Finally $\Musat \subset \Usat$ is the set of $F \in \Usat$ such that for all $C \in F$ we have $F \sm \set{C} \in \Sat$. The background for the investigations of this \Schrift{} is the enterprise of classifying $F \in \Musat$ in dependency on $\delta(F)$. The basic facts are $\delta^*(F) = \delta(F)$ (as will be discussed in Subsection \ref{sec:introautgen}), and the well-known $\delta(F) \ge 1$, as first shown in \cite{AhLi86} (``Tarsi's Lemma''); we give the simple proof we learned from \cite{FlRe94} in Subsection \ref{sec:introTarsi}. For $\delta(F) = 1$ the structure is best understood (\cite{AhLi86,DDK98,Ku99dKo}; see Example \ref{exp:MU1}), for $\delta(F) = 2$ the structure after reduction of singular variables (occurring in one sign only once) is known (\cite{KleineBuening2000SubclassesMU}; see Example \ref{exp:MU2}), while for $\delta(F) \in \set{3,4}$ only basic cases have been classified (\cite{XD99}).

The starting point of our investigation is \cite[Lemma C.2]{Ku99dKo}, where it is shown that $F \in \Musat$ with $n(F) > 0$ must have a variable $v \in \var(F)$ with at most $\delta(F)$ positive and at most $\delta(F)$ negative occurrences (we will give the short proof in Subsection \ref{sec:introbasicdegree}); we write this as $\ldeg_F(v) \le \delta(F)$ and $\ldeg_F(-v) \le \delta(F)$, using the notion of \emph{literal-degrees} (the number of occurrences of the literal), where for a literal $x$ its degree is $\ldeg_F(x) := \abs{\set{C \in F : x \in C}} \in \NNZ$.
Thus we have $\vdeg_F(v) \le 2 \delta(F)$, using the \emph{variable-degree}
\begin{displaymath}
  \vdeg_F(v) := \ldeg_F(v) + \ldeg_F(-v) \in \NNZ.
\end{displaymath}
Using the \emph{minimum variable-degree} (min-var-degree)
\begin{displaymath}
  \minvdeg(F) := \min_{v \in \var(F)} \vdeg_F(v) \in \NN
\end{displaymath}
of $F$ with $n(F) > 0$, the upper bounds becomes $\minvdeg(F) \le 2 \delta(F)$. A main theme of this \Schrift{} is the consideration of $\minvdeg(\Musati{\delta=k}) \in \NN$ for $k \in \NN$, the maximum of $\minvdeg(F)$ for $F \in \Musat$ with $\delta(F) = k$. The upper bound now becomes $\minnonmer(k) := \minvdeg(\Musati{\delta=k}) \le 2 k$. We show a sharper bound on $\minvdeg(F)$, namely we show that the worst-cases $\ldeg_F(v), \ldeg_F(-v) \le \delta(F)$ can not occur at the same time (for a suitable variable; see Subsection \ref{sec:introdef3} for the basic example), but actually $\ldeg_F(v) + \ldeg_F(\ol{v}) - \delta(F)$ only grows logarithmically in $\delta(F)$. The really interesting aspect here is the precise determination of $\minnonmer(k)$, and we investigate the (elementary) number-theoretic function $\nonmer(k)$, which yields the upper bound $\minnonmer(k) \le \nonmer(k)$ for all $k \in \NN$, where the function $\nonmer: \NN \ra \NN$ fulfils $k + \floor{\log_2(k+1)} \le \nonmer(k) \le k + 1 + \floor{\log_2(k)}$ for $k \in \NN$.

\subsection{Refining deficiency by surplus}
\label{sec:introsurp}

After having settled this basic min-var-degree upper bound for $\Musati{\delta=k}$, we show a sharper bound on $\minvdeg(F)$ for a larger class of clause-sets $F$:
\begin{itemize}
\item The larger class of clause-sets considered is the class $\Lean$ of \emph{lean clause-sets} (introduced in \cite{Ku98e}), which are clause-sets having no non-trivial autarky. For an overview on the theory of minimally unsatisfiable clause-sets and on the theory of autarkies see \cite{Kullmann2007HandbuchMU}. $\Lean \subset \Cls$ is the set of $F \in \Cls$ such that there is no partial assignment $\vp$ (a ``non-trivial autarky'') with the properties
  \begin{itemize}\setlength{\itemsep}{0pt}
  \item for every clause $D \in F$ with $-\vp \cap D \ne \es$ we have $\vp \cap D \ne \es$ (note that this generalises the satisfaction criterion);
  \item there exists $v \in \var(F)$ with $\set{v,-v} \cap \vp \ne \es$.
  \end{itemize}
  Note $\Lean \cap \Sat = \set{\top}$, where $\top := \es \in \Cls$ is the empty clause-set (the standard satisfiable clause-set).
\item The deficiency $\delta(F) \in \ZZ$ is strengthened by the \emph{surplus} $\surp(F) \in \ZZ$, defined in case of $n(F) > 0$ as follows.

  Consider the bipartite clause-variable graph of $F$ (generalising the incidence graph of a hypergraph), with the clauses $C \in F$ on one side of the bipartition, and the variables $v \in \var(F)$ on the other side, and an edge between $v$ and $C$ if $\set{v,-v} \cap C \ne \es$. The ``expansion'' of a set $\es \ne V \sse \var(F)$ of variables is $\abs{\Gamma(V)} - \abs{V}$, where $\Gamma(V)$ is the set of neighbours of $V$ (incident clauses), and the surplus then is the minimum expansion, i.e., $\surp(F) = \min_{\es \ne V \sse \var(F)} \abs{\Gamma(V)} - \abs{V}$.

  In the terminology of \cite[Section 1.3]{LP86neu}, $\delta^*(F)$ is the deficiency of the bipartite clause-variable graph (with bipartition $(F,\var(F))$), while $\surp(F)$ is the surplus of the bipartite variable-clause graph (with bipartition $(\var(F),F)$).

  Note that by considering $V = \var(F)$ we have $\surp(F) \le \delta(F)$, and by considering $V = \set{v}$ for $v \in \var(F)$ we get $\surp(F) \le \minvdeg(F) - 1$.

  We have $\surp(F) \ge 1$ for $F \in \Lean$ with $n(F) > 0$ (\cite[Lemma 7.7]{Ku00f}), generalising the basic fact $\delta(F) \ge 1$ for $F \in \Musat$ (``Tarsi's Lemma'').
\end{itemize}
Now a central result of this \Schrift{} (Theorem \ref{thm:leanminvardeg}) is
\begin{displaymath}
  \minvdeg(F) \le \nonmer(\surp(F))
\end{displaymath}
for $F \in \Lean$ with $n(F) > 0$. As an application we obtain (Theorem \ref{thm:genautred}), that by removing some clauses (satisfiability-equivalently, implicitly using some autarky), we can reduce every (multi-)clause-set $F$ in polynomial time to a clause-set $F'$ containing a variable with degree at most $\surp(F') + 1 + \log_2(\surp(F'))$. It is an open problem whether the witnessing-autarky can be found in polynomial time; we conjecture (Conjecture \ref{con:findauthard}) that this is possible. A central tool here are ``variable-minimally unsatisfiable clause-sets'' (class $\Vmusat$; introduced in \cite{ChenDing2006VMU}), where more generally than for MU, removal of clauses might not destroy unsatisfiability, but as soon as all occurrences of any variable disappear, then satisfiability is guaranteed. The basic insight is, that the surplus is based on sub-instances of $F$ of \emph{deficiency} $\surp(F)$, and for $F \in \Lean$ these sub-instances are indeed variable-minimally unsatisfiable, and thus we can use the bound for MU. We also show sharpness of the upper bound, i.e., $\minvdeg(\Leani{\delta=k}) = \nonmer(k)$ for all $k \in \NN$, in Corollary \ref{cor:leansharp} (proving Conjecture 23 from the conference version \cite{KullmannZhao2011Bounds}), which indeed holds for every class of clause-sets between $\Vmusat$ and $\Lean$.

We then come back to the special case of minimal unsatisfiability. Here things are much more complicated, and the numbers $\minnonmer(k)$, the guaranteed minimum variable-degrees for minimally unsatisfiable clause-sets of deficiency $k$, are very interesting quantities. We prove the sharpened bound $\minnonmer(k) \le \nonmer_1(k)$, which improves on $\nonmer(k)$ for infinitely many $k$.

\subsection{Four case studies}
\label{sec:introfour}

We outline now some typical basic arguments, to develop an intuition about the subjects of this \Schrift. Consider $F \in \Musat$. We want to gain information on $F$, by taking some $v \in \var(F)$, consider the splitting results $F_{\ve} := \pao v{\ve} * F$ for $\ve = 0,1$, i.e., setting $v$ to false and true, and ``reconstructing'' $F$ from these two pieces. How is $\delta(F) = c(F) - n(F)$ related to $\delta(F_{\ve})$ ? The deficiency goes down by the number of satisfied clauses (exercise: why can more clauses be lost than that, if considering general $F \in \Cls$ ?), and goes up by the number of eliminated variables. The first main observation here is that if $v$ has \emph{minimum degree}, then no further variable can be lost in $F_{\ve}$, that is, $\var(F_{\ve}) = \var(F) \sm \set{v}$ (exercise: again, here we use some aspect of $F \in \Musat$). So there is a direct relation between the literal-degrees $m_0, m_1$ of $\ol{v}$ resp.\ $v$, and the deficiencies of $F_0, F_1$, namely $\delta(F_{\ve}) = \delta(F) - m_{\ve} + 1$. But we have the problem that in general $F_{\ve} \notin \Musat$. To the rescue comes \emph{saturation}, which adds literals (in the given variables) to some clauses of $F$, obtaining some (non-unique) \emph{saturated} $F' \in \Musat$, with $n(F') = n(F)$ and $c(F') = c(F)$, thus $\delta(F') = \delta(F)$. Adding the literals in the saturation repairs the non-minimality of the $F_{\ve}$, by making sure that the superfluous clauses get eliminated by the assignment; see Subsection \ref{sec:prelimMU} for more on saturation. If we can argue now that the process of saturation is ``harmless'' (as it is regarding the deficiency), then w.l.o.g.\ we can assume that $F$ is saturated, and then we have $F_{\ve} \in \Musat$ for both $\ve$.

\subsubsection{A simple proof of Tarsi's Lemma}
\label{sec:introTarsi}

We are now ready to show $\fa\, F \in \Musat : \delta(F) \ge 1$, by induction on $n(F)$. We use $\bot := \es$ for the empty clause. If $n(F) = 0$, then $F = \set{\bot}$, and $\delta(F) = 1 - 0 = 1$; assume $n(F) \ge 1$. W.l.o.g.\ we can assume $F$ is saturated, since saturation does not change the deficiency. Consider $v \in \var(F)$ of minimum degree. By the above discussion we have $\delta(F_{0}) \le \delta(F) - 1 + 1 = \delta(F)$, where by induction hypothesis we have $\delta(F_0) \ge 1$. QED

\subsubsection{The basic degree bound}
\label{sec:introbasicdegree}

Also the proof, that for $F \in \Musat$ with $n(F) \ge 1$ there is a variable $v \in \var(F)$ with $\ldeg_F(v), \ldeg_F(\ol{v}) \le \delta(F)$, is easy now (details in Lemma \ref{lem:uppbdldg}): Again, w.l.o.g.\ $F$ is saturated (saturation can not decrease literal-degrees). Choose $v \in \var(F)$ of minimum degree. Assume w.l.o.g.\ that $\ldeg_F(v) \ge \delta(F)+1$. Then by the above discussion $\delta(F_1) \le \delta(F) - (\delta(F) + 1) - 1 \le 0$, contradicting Tarsi's Lemma. QED

\subsubsection{Improving the basic bound for deficiency $3$}
\label{sec:introdef3}

We get $\fa\, F \in \Musat \sm \set{\set{\bot}} : \minvdeg(F) \le 2 \delta(F)$.
Consider $\delta(F) = 3$; we show $\minvdeg(F) \le 5$. W.l.o.g.\ $F$ is saturated. Take $v \in \var(F)$ with minimum degree, and assume $\vdeg_F(v) = 6$ (note $n(F) \ge 2$). Now $1 \le \delta(F_0) \le \delta(F) - 3 + 1$, so $\delta(F_0) = 1$, and there is $w \in \var(F_0)$ with $\vdeg_{F_0}(w) = 2$, and we get $\vdeg_F(w) \le 2 + 3 = 5$. QED

\subsubsection{Some basic intuitions about the upper bound $\nonmer$}
\label{sec:introbasicint}

The function $\nonmer: \NN \ra \NN$ is strictly increasing with range
\begin{displaymath}
  \nonmer(\NN) = \NN \sm \set{2^n-1: n \in \NN} = \set{\ul{2},\quad 4,5,\ul{6}, \quad 8,\dots,\ul{14},\quad 16,17,\dots}.
\end{displaymath}
We show $\minvdeg(\Leani{\delta=k}) = \nonmer(k)$ for deficiencies $k \in \NN$, that is, every lean clause-set $F$ with $n(F) > 0$ contains a variable $v \in \var(F)$ with $\vdeg_F(v) \le \nonmer(\delta(F))$, and for every deficiency $k \ge 1$ there are lean clause-sets $F$ with $\minvdeg(F) = \nonmer(\delta(F))$.

The underlined values $2, 6, 14, \dots$, which have the form $2^n-2$ for $n \ge 2$, are the function values at the ``jump positions'' $1, 4, 11, \dots$, which are of the form $2^n - n - 1$ for $n \ge 2$ (where the function values changes by $+2$, while otherwise it changes by $+1$ for an increment of the argument). This basic structure of $\nonmer$ can be motivated by the following constructions of $F \in \Musat$ with ``high'' min-var-degree. Indeed these considerations only concern the lower bounds, given by appropriate constructions, while the arithmetic nature of the upper bound $\nonmer(k)$ rests on different considerations, but for the deficiencies considered here, lower and upper bounds are equal, and the lower bounds are easier to understand here.

The basic clause-sets are the $A_n$ for $n \in \NNZ$, which consist of all $2^n$ sets (clauses) of numbers $\pm 1, \dots, \pm n$, using the natural numbers $1,\dots,n$ as variables. So $A_0 = \set{\bot}$, $A_1 = \set{\set{-1},\set{1}}$, $A_2 = \set{\set{1,2},\set{-1,2},\set{1,-2},\set{-1,-2}}$ and so on. It is easy to see that we have $A_n \in \Musat$ with $n(A_n) = n$, $c(A_n) = 2^n = \minvdeg(A_n)$, and $\delta(A_n) = 2^n - n$. We will see that the $A_n$ have the largest possible min-var-degree $2^n$ for given deficiency $2^n - n$, and we also have $\nonmer(2^n - n) = 2^n$ for $n \in \NN$. These deficiencies $k = 2^n - n$ (numerical values are $1, 2, 5, 12, \dots$) are the positions directly after the jump positions (excluding deficiency $k=1$ as a special case).

How can we obtain from that more clause-sets in $\Musat$ with high min-var-degree? Consider $A_3$: we have e.g.\ $\set{1,2,3}, \set{1,2,-3} \in A_3$; now logically these two clauses are equivalent to the single clause $\set{1,2}$ (i.e., we have the same satisfying assignments; technically, a ``strict full subsumption resolution'' is performed), and we obtain $A_3' := (A_3 \sm \set{\set{1,2,3}, \set{1,2,-3}}) \cup \set{\set{1,2}} \in \Musat$. Performing this process in general, using $\set{1,\dots,n}, \set{1,\dots,n-1,-n} \in A_n$, yields $A_n' \in \Musat$ for $n \ge 2$, with $n(A_n') = n$, $c(A_n') = 2^n - 1$, $\delta(A_n') = 2^n - n - 1$ , and $\minvdeg(A_n') = 2^n - 2$ (the (single) variable with minimum occurrences is $n$). These deficiencies are precisely the jump positions $2^n - n - 1$, and accordingly we have $\nonmer(2^n - n - 1) = 2^n - 2$.

Performing the same trick again to $A_3'$, we can replace $\set{-1,2,3}, \set{-1,-2,3} \in \A_n'$ by $\set{-1,3}$, obtaining $A_3'' \in \Musat$. Again for general $n \ge 3$ we get $A_n'' \in \Musat$, $n(A_n'') = n$, $c(A_n'') = 2^n - 2$, $\delta(A_n'') = 2^n - n - 2$, and $\minvdeg(A_n'') = 2^n - 3$; note here the crucial difference, that the min-var-degree has only been changed by $-1$. The reason is that there are two variables now with minimum occurrences, namely $n-1, n$, where the degree of variable $n$ changed first by $-2$, then by $-1$, while for variable $n-1$ the degree first changed by $-1$, and then by $-2$ (and for the other variables $1, \dots, n-2$ we have degree changes by $-1$, $-1$).

We will apply this lower bound method to simple cases in Subsection \ref{sec:preciseval}. One might imagine this process of ``strict full subsumption resolutions'' continuing until deficiency $2^{n-1} - (n - 1) + 1$, always with change of the min-var-degree by $-1$, just before the deficiency of the previous $A_{n-1}$ --- this would yield the function $\nonmer$. However the combinatorial reality is more complicated, and as we prove in this \Schrift{} (Section \ref{sec:strbMU}), at least we can not get until $2^{n-1} - (n - 1) + 1 = 2^{n-1}-n+2$ for $n \ge 4$ (in effect), that is, at these deficiencies $k = 6, 13, 28, \dots$ we have $\minnonmer(k) \le \nonmer_1(k) = \nonmer(k) - 1$. In general the determination of $\minnonmer(k)$ seems complicated, however the construction of examples showing the weaker $\minvdeg(\Leani{\delta=k}) \ge \nonmer(k)$ is rather easy (see Lemma \ref{lem:leansharpaux}).

After having developed now some intuition on the main results of this \Schrift, we turn to the discussion of the background and context.

\subsection{Related work on $\Musat$}
\label{sec:prelimMU}

A general overview on minimally unsatisfiable clause-sets (also ``minimal unsatisfiable clause-sets/formulas'', or ``MU'') is \cite{Kullmann2007HandbuchMU}; later developments are in \cite{Kullmann2007ClausalFormZI,Kullmann2007ClausalFormZII} (generalisations to non-boolean clause-sets) and in \cite{KullmannZhao2012ConfluenceC,KullmannZhao2012ConfluenceJ} (studying ``singular DP-reduction'', the elimination of variables which occur in one sign only once).

Two early papers on the complexity aspects are \cite{PY84,PW88}, who introduced the complexity class $D^P$ and showed that the decision ``$F \in \Musat$ ?'' for input $F \in \Cls$ is complete for this class. Another important early paper is \cite{AhLi86}, which showed $\delta(F) \ge 1$ for $F \in \Musat$, where the notion of ``deficiency'' was introduced by \cite{FrGe98}. Furthermore \cite{AhLi86} showed polytime-decision of the sub-class $\Smusati{\delta=1} \subset \Musati{\delta=1}$ (called ``strongly minimal unsatisfiable'' there), where $\Smusat \subset \Musat$ is the set of $F \in \Usat$ such that for all $C \in F$ and all $x \in \ZZ \sm \set{0}$ with $\set{x,-x} \cap C = \es$ holds $(F \sm \set{C}) \cup \set{C \cup \set{x}} \in \Sat$, that is, adding any literal to any clause renders the clause-set satisfiable. We use the terminology ``saturated minimally unsatisfiable'' introduced in \cite{FlRe94}, where the important connection to splitting was introduced, and a simpler proof of $\delta(F) \ge 1$ for $F \in \Musat$ was given (as presented in Subsection \ref{sec:introTarsi}). Recall that in this introduction we handle ``partial assignments'' via clauses $\vp$ (containing the satisfied literals; thus $-\vp$ is the set of falsified literals), so for a literal $x$ the partial assignment $\pao x0$ is given by $\set{-x}$, while $\pao x1$ is given by $\set{x}$. The application of $\vp$ to $F \in \Cls$ is defined as
\begin{displaymath}
  \vp * F := \set{C \sm -\vp : C \in F \und \vp \cap C = \es} \in \Cls,
\end{displaymath}
that is, removing first the satisfied clauses from $F$, and then the falsified literals from the remaining clauses. For $F \in \Cls$ holds $F \in \Smusat$ iff for all $x \in \ZZ \sm \set{0}$ holds $\pao x1 * F \in \Musat$ (the ``only if''-direction was shown in \cite{FlRe94}, the ``if''-direction in \cite{Ku99dKo}; see Lemma \ref{lem:auxSMUSAT} for details). Since every $F \in \Musat$ can be ``saturated'' by adding literals to clauses (as explained in Subsection \ref{sec:introfour}), the class $\Smusat$ is thus an important helping class for investigations into $\Musat$ via the splitting method, splitting up $F \in \Musat$ into $\pao v0 * F$ and $\pao v1 * F$ for selected variables $v$.

We have already mentioned the literature concerned with characterising the classes $\Musati{\delta=k}$ (and subclasses) for small deficiencies $k \le 4$. Less ambitious is the goal of polytime decision of these classes: the problem was raised in \cite{Kl98}, and has been solved via two independent approaches in \cite{Ku99dKo} and \cite{FKS00} (indeed establishing polytime SAT decision for inputs $F \in \Cls$ and fixed $\delta^*(F)$), later strengthened in \cite{Szei2002FixedParam} (showing that SAT decision is even fixed-parameter tractable in $\delta^*(F)$; see also \cite{Kullmann2007ClausalFormZI} for generalisations and simplifications).

We now consider the three other main areas in the literature on MU, generalisations (Subsection \ref{sec:prelimbeyond}), minimally unsatisfiable \emph{sub-}clause-sets (Subsection \ref{sec:prelimMUS}), and \emph{maximum} variable-degrees (Subsection \ref{sec:prelimTovey's}).

\subsubsection{Beyond boolean CNFs}
\label{sec:prelimbeyond}

For our overview on generalisations of boolean $\Musat$, we restrict attention to areas which have some form of ``Tarsi's Lemma'', that is, where for (generalised) ``minimal unsatisfiable formulas'' the number of variables is upper-bounded by some function of the number of some form of (generalised) ``clauses''.

The first generalisation, in the Schaefer framework (\cite{Sch78}), can be understood as restricting clauses to a bounded size $k$, but allowing arbitrary boolean functions with $k$ variables as constraint-templates (instead of clauses); it is furthermore assumed that each template does not have forced assignments (no variable is fixed to some value in all satisfying assignments). \cite[Proposition 3.6]{CreignouDaude2002ThresholdConstraints} shows that for every instantiation of such (boolean) constraint-templates, the number of variables in a minimally unsatisfiable constraint set of size $s$ is at most $(k-1) \cdot s$ (instead of the trivial $k \cdot s$). A further sharpening to ``$(k-1) \cdot s - 1$'' in case none of the constraint templates ``depends strongly on a 2XOR relation'' is \cite[Theorem 4.3]{CreignouDaude2002ThresholdConstraints}.

A different generalisation is considered in \cite{BenSassonNordstroem2011Substitutions,NordstroemRazborov2011MU}, called ``$k$-DNF'', which we can understand as to allow for ``super-literals'' in a CNF (a clause-set), which are conjunctions of up to $k$ ordinary literals (so ordinary CNFs are obtained for $k=1$); in other words, we consider conjunctions of DNFs, where in each DNF the ``terms'' (the conjunctions) have size at most $k$. The notion of ``minimal unsatisfiability'' means here that removing \emph{any} literal from any super-literal (the innermost conjunctions) makes the whole formula satisfiable; this generalises the case $k=1$, since empty conjunctions are constant true, and thus lead to the removal of the whole containing clause. In \cite[Theorem 15]{BenSassonNordstroem2011Substitutions} the upper bound $n(F) \le (k \cdot c(F))^{k+1}$ for such formulas $F$ is shown, with $n(F)$ the number of variables and $c(F)$ the number of ``super-clauses'' (or DNFs), while \cite[Theorem 4]{NordstroemRazborov2011MU} shows the lower bound $\Omega(c(F)^k)$. We now turn to generalised notions of ``deficiency'', which yield a ``precise Tarsi's Lemma'', i.e., the generalised deficiency for ``MU'' is at least $1$, and at least the bottom layer is polytime-decidable (and thus we obtain layers of complexity, with deficiency $1,2,\dots$).

For general propositional formulas in NNF (``negation normal form'', arbitrary \texttt{or}'s and \texttt{and}'s of literals), the deficiency has been generalised to ``cohesion'' in \cite{KleineBueningZhao2007DeficiencyQBF}, which is $1$ plus the number of \texttt{and}'s minus the number of variables. ``Minimal unsatisfiability'' here means that replacing any \texttt{or}-term (including trivial ones) with \texttt{true} yields a satisfiable clause-set. In this way the basic facts about MU and deficiency for CNFs are properly generalised, and especially the cohesion of an MU-NNF is at least $1$. We note that compared to ``$k$-DNF'' as above, on the one hand the notion of ``MU'' for NNF is more general than the (special) notion of ``MU'' for (the subclass) $k$-DNF, since for NNF we only consider to replace whole disjunctions by \texttt{true}. On the other hand the cohesion for every additional super-literal (an innermost conjunction) of size $m$ is increased by $m-1$, while for ``$k$-DNF'' as above the super-clauses are counted simply as $1$, whatever they contain (note that a single super-clause can contain an unbounded number of \texttt{and}'s). A further generalisation to arbitrary circuits is given in \cite{BelovMarquesSilvaMUCircuits}.

The natural generalisation of deficiency from CNF to QCNF (quantified boolean CNF), just ignoring universal variables, has been introduced and studied in \cite{BueningZhao2008QBFFixDef}, and generalised to quantified formulas in NNF in \cite{KleineBueningZhao2007DeficiencyQBF}. It is open whether we have here polytime decision for bounded (maximal) deficiency, only deficiency $1$ has been resolved. In a different direction, the generalisation to hypergraph-2-colouring (involving a translation) will be considered in Subsection \ref{sec:prelimcolouring}. Another environment in which many questions regarding $\Musat$ (and $\Lean$) are reformulated (and generalised) is ``qualitative matrix analysis'', discussed in Subsection \ref{sec:prelimQCA}.

All above generalisations use \emph{boolean} variables, and go beyond CNF --- staying with CNF, but admitting variables with arbitrary finite domains is studied in \cite{Kullmann2007ClausalFormZI,Kullmann2007ClausalFormZII}. This generalisation is closer to this \Schrift, and indeed we think it is an interesting, challenging and important endeavour to generalise our results to such non-boolean CNFs; we will discuss this further in Subsection \ref{sec:concgennb}, while we refer to existing generalisations where appropriate. Roughly the relations between \cite{Kullmann2007ClausalFormZI,Kullmann2007ClausalFormZII} and this \Schrift{} consist in two points:
\begin{enumerate}
\item The fundamental ``Tarsi's Lemma'' is generalised in \cite[Corollaries 9.8, 9.9, Lemma 11.1]{Kullmann2007ClausalFormZI} to ``non-boolean clause-sets'.
\item The classification of $\Musat$ layered by deficiency is started in \cite[Chapter 5]{Kullmann2007ClausalFormZII}, yielding also first generalisations related to the minimum variable-degree .
\end{enumerate}
Finally, a different field of investigations, this time in number theory, turns out also to be related to our field, namely the study of covers of the integers by congruence relations, and covers of lattice parallelotopes by certain types of cells. As remarked in \cite[Section 6]{BergerFelzenbaumFraenkel1990Covers}, the main result of \cite{BergerFelzenbaumFraenkel1988Covers} generalises the (original) Tarsi's Lemma. And indeed, this result, generalising \cite{Znam1975Covers}, is closely related to the above ``Tarsi's Lemma'' for non-boolean clause-sets (\cite[Corollary 9.9]{Kullmann2007ClausalFormZI}); especially \cite[Corollary 4a]{BergerFelzenbaumFraenkel1988Covers} (concerning parallelotope covers) can be seen as equivalent to it (where the proof of \cite[Corollary 9.9]{Kullmann2007ClausalFormZI} is much easier). A measure for congruence covers corresponding to deficiency is called the ``Mycielski-Zn\'{a}m abundance'' in \cite[Remark 2.11]{Korec1987Covers}.

\subsubsection{MUSs}
\label{sec:prelimMUS}

As we have already mentioned, we consider $\Musat$ as the ``primal'' building block for understanding unsatisfiability. In general an unsatisfiable clause-set can contain many minimally unsatisfiable \emph{sub}-clause-sets, called ``MUSs''. The task of enumerating all of them or at least some ``good'' ones is also of practical importance, to extract more information on the ``causes'' of unsatisfiability. A recent overview is \cite{MarquesSilva2012MUS}, while a clean approach to enumerate all MUSs, via hypergraph transversals, is in \cite{LiffitonSakallah2005AllMUS} (the earliest appearance of the underlying observation seems \cite[Theorem 2]{Bruni2003ApprMU}; compare also \cite[Subsection 4.3]{Kullmann2007ClausalFormZII} for generalisations of the fundamental approach). See also \cite{KullmannLynceSilva2005Autarkies} for a reflection on various types of such sub-clause-sets, and on the connection to autarky theory (compare Subsection \ref{sec:introLEAN}). For non-boolean variables and arbitrary constraints, the problem of finding good or all minimally unsatisfiable sub-constraint-sets is also of importance, and an influential paper is \cite{BaileyStuckey2005MUS} (another source of the above mentioned approach via hypergraph transversals).

\subsubsection{Tovey's problem (uniform clause-sets)}
\label{sec:prelimTovey's}

This \Schrift{} appears to be the first systematic study of the problem of \emph{minimum} variable occurrences / degrees in minimally unsatisfiable clause-sets and generalisations, in dependency on the deficiency --- asking for the existence of a variable occurring ``infrequently'' in general, or for extremal examples where all variables occur not infrequently. The ``dual'' problem is to consider \emph{maximum} variable occurrences / degrees --- asking for the existence of a variable occurring frequently in general, or for extremal examples where all variables occur not frequently. More precisely, the \emph{maximum variable-degree} is
\begin{displaymath}
  \maxvdeg(F) := \max_{v \in \var(F)} \vdeg_F(v) \in \NN,
\end{displaymath}
for $n(F) > 0$, while for a class $\mc{C} \sse \Cls$ of clause-sets, the quantity $\maxvdeg(\mc{C})$ is the minimum of $\maxvdeg(F)$ for $F \in \mc{C}$. This problem has been well-studied for \emph{$p$-uniform} minimally unsatisfiable clause-sets, starting with \cite{Tovey1984NPcomplete,Dubois1990Tovey,KratochvilSavickyTuza1993Jump}.\footnote{We remark that typically in the literature the connections to \emph{minimally} unsatisfiable clause-sets are not emphasised, but it is clear that when considering (uniform) unsatisfiable clause-sets with a maximum variable-degree as small as possible, then one can restrict attention w.l.o.g.\ to (uniform) minimally unsatisfiable clause-sets (as worst-cases).} We denote by $\Pcls{p} \subset \Cls$ for $p \in \NNZ$ the set of all $F \in \Cls$ with $\fa\, C \in \Cls : \abs{C} \le p$, while by $\Unicls \subset \Cls$ we denote the set of all uniform clause-sets, i.e., those $F \in \Cls$ such that for $C, D \in F$ holds $\abs{C} = \abs{D}$. Finally $\Punicls{p} := \Pcls{p} \cap \Unicls$ and $\Punimusat{p} := \Punicls{p} \cap \Musat$. Now the basic fact is
\begin{displaymath}
  \maxvdeg(\Punimusat{p}) \ge p+1
\end{displaymath}
for $p \in \NN$ (\cite{Tovey1984NPcomplete}, generalised in \cite[Corollary 7.3]{Kullmann2007ClausalFormZII}). Trivially $\maxvdeg(\Punimusat{1}) = 2$, and easily one sees $\maxvdeg(\Punimusat{2}) = 3$, while by \cite{Tovey1984NPcomplete} holds $\maxvdeg(\Punimusat{3}) = 4$. As reported in \cite{HoorySzeider2005FewOcc}, we have $\maxvdeg(\Punimusat{4}) = 5$, and these are all known precise values of $\maxvdeg(\Punimusat{p})$ (where the notation $f(p) := \maxvdeg(\Punimusat{p}) - 1$ was introduced in \cite{KratochvilSavickyTuza1993Jump}). In \cite{HoorySzeider2005FewOcc} it was observed, that extremal examples might be found in $\Punimusat{p}_{\delta=1}$: the intuition here is that the elements of $\Musati{\delta=1}$ contain the maximal number of variables for a given number of clauses, so the average variable-degree is lowest, and so minimising the maximal var-degree seems easiest. This work was recently extended in \cite{GebauerSzaboTardos2010LocalLemma}, establishing the asymptotically tight bound $\lim_{p \ra \infty} \frac 2e \frac{2^p}{p} / \maxvdeg(\Punimusat{p}) = 1$ (where indeed $\Punimusat{p}_{\delta=1}$ is considered).

An open question is the computability of $\maxvdeg(\Punimusat{p})$. If $\maxvdeg(\Punimusat{p}) = \maxvdeg(\Punimusat{p}_{\delta=1})$ holds, then via simple search we have computability; this is somewhat similar to our study of the min-var-degree in the deficiency, where also computability is open, while we obtain it (by simple search) when assuming, that the class of clause-sets to be considered can be restricted (Lemma \ref{lem:conjimplcomp}).

\paragraph{Uniformity versus non-uniformity}

For studying the classes $\Musati{\delta=k}$, the max-var-degree is not very relevant, since we have $\maxvdeg(\Musati{\delta=1}) = 2$, while $\maxvdeg(\Musati{\delta=k}) = 3$ for $k \ge 2$. This can be seen as follows: As already noticed in \cite{Tovey1984NPcomplete}, there is a poly-time transformation from $\Cls$ to the class $\Cls(1,2) \subset \Cls$, consisting of those $F \in \Cls$ where for every variable $v \in \var(F)$ we have $\ldeg_F(-v) = 1$ and $\ldeg_F(v) \le 2$. Namely if there is a literal $x$ and two clauses $C, D \in F$ with $x \in C \cap D$, then we can introduce a new variable $v$, replace $x$ in $C,D$ by $v$, and add the new clause $\set{-v,x}$, obtaining $F'$. We study such extensions under the name of ``singular DP-extension'', but it is also easy to see directly that $F'$ is satisfiable iff $F$ is, that $F'$ is minimally unsatisfiable iff $F$ is, and that $\delta(F') = \delta(F)$. By repeating this transformation, we obtain $t^{1,2}: \Cls \ra \Cls(1,2)$. So for $F \in \Musat$ we get $t^{1,2}(F) \in \Musat \cap \Cls(1,2)$ with $\delta(t(F)) = \delta(F)$. Whence for all $k \in \NN$ we have $\maxvdeg(\Musati{\delta=k}) \le 3$. Now trivially $\maxvdeg(\Musati{\delta=1}) = 2$ due to $\set{\set{1}, \set{-1}} \in \Musati{\delta=1}$. On the other hand, if for $F \in \Musat$ holds $\maxvdeg(F) \le 2$ (thus $\maxvdeg(F) = 2$), then via so-called singular DP-reduction this clause-set can be reduced to $\set{\bot}$, whence $F \in \Musati{\delta=1}$ (this is well-known; compare Example \ref{exp:MU1} later). The above transformation $t^{1,2}$ relied on the introduction of singular variables; now denote by $\Musatns$ the set of nonsingular elements of $\Musat$, i.e., where every literal occurs at least twice:
\begin{quest}\label{que:maxvdegns}
  Is $k \mapsto \maxvdeg(\Musatnsi{\delta=k})$ strictly increasing?
\end{quest}

So for the study of the max-var-degree, the uniformity restriction seems important. This is similar to investigations into (colour-)critical hypergraphs (discussed in Subsection \ref{sec:prelimcolouring} below), where uniformity is often a crucial assumption, and the hyperedge-length $p$ is the main parameter. For investigations into the case of uniform (general) clause-sets, where clauses share at most one variable, see \cite{PorschenSpeckenmeyerZhao2009LinearCNF,Scheder2010Linear}. The maximal number of clauses in $F \in \Punimusat{p}$ has been studied in \cite{Lee2009MU3CNF}, showing that for $p=2$ holds $c(F) \le 4 n - 2$, while for $p \ge 3$ there are $F$ with $c(F) = \Omega(n(F)^p)$. Finally, the number of conflicts (clashes) in $F \in \Punimusat{p}$ is considered in \cite{Scheder2013Conflicts}, and for a review of the use of the Lov{\'{a}}sz Local Lemma in this context see \cite{GebauerMoserSchederWelzl2009LocalLemma}.

In contrast, for the study of the minimum variable-degree as in this \Schrift, in dependency on the deficiency, the restriction to uniformity seems not interesting, and is also not needed, but unrestricted clause-sets are considered. We remark that for every $p \in \NN$, $p \ge 3$, there is a polytime translation $t_p: \Cls \ra \Punicls{p}$, such that $t_p(F)$ is satisfiable iff $F$ is, $t_p(F)$ is minimally unsatisfiable iff $F$ is, and $\delta(t_p(F)) = \delta(F)$. This works by replacing clauses $C$ with $\abs{C} < p$ by clauses $C \cup \set{v}, C \cup \set{-v}$ for some new variable $v$ (in the MU-case we will call this a ``non-strict full subsumption extension''), and by replacing clauses $C$ with $\abs{C} > p$ by clauses $C' \cup \set{v}$, $C'' \cup \set{-v}$ for some new variable $v$ and choosing clauses $C', C''$ with $C = C' \cup C''$ and $\abs{C'} = p-1$, $\abs{C''} \ge p-1$ (in the MU-case again we have a singular DP-extension). But the transformation $t_p$ appears to be useless for structural investigations. Supposing nonsingularity makes uniformity more interesting --- now there seem to be very few choices for fixed deficiency:
\begin{quest}\label{que:uniformns}
  Is $n(F)$ bounded for fixed $k$ for  $F \in \Unimusat \cap \Musatnsi{\delta=k}$?
\end{quest}

\subsection{Autarkies}
\label{sec:introautgen}

An important tool, used in this \Schrift{} to go beyond $\Musat$, is the theory of autarkies, which also provides links to various areas of combinatorics; the relations to hypergraph colouring will be discussed in Subsection \ref{sec:introLEAN}. Recall that a partial assignment $\vp$ is an autarky for $F \in \Cls$ iff every clause $C \in F$ touched by $\vp$ (i.e., $\vp \cap (C \cup -C) \ne \es$) satisfies $C$ (i.e., $\vp \cap C \ne \es$), which is equivalent to $\fa\, F' \sse F : \vp * F' \sse F'$. Autarkies were introduced in \cite{MoSp85} for improved worst-case upper bounds for SAT decision, applying that obviously $\vp * F$ is sat-equivalent to $F$ ($\vp * F \in \Sat \Lra F \in \Sat$) for an autarky $\vp$. For a recent overview see \cite{Kullmann2007HandbuchMU}.

\paragraph{Autarky reduction}

The reduction of $F \in \Cls$ to $\vp * F \in \Cls$ for a non-trivial autarky $\vp$ is an essential concept, algorithmically as well as for theoretical understanding; see \cite[Subsection 11.10]{Kullmann2007HandbuchMU} for an overview on finding autarkies. If we reduce all autarkies, then we obtain the (unique) \emph{lean kernel} of $F$. If there are no non-trivial autarkies, then we have a lean clause-sets, i.e., $F \in \Lean$, as already mentioned in Subsection \ref{sec:introsurp}; this concept was introduced in \cite{Ku98e}, and \cite[Subsection 11.8.3]{Kullmann2007HandbuchMU} contains more information. The lean kernel of $F$ is the largest lean sub-clause-set of a clause-set, that is, $\bc \set{F' \sse F : F' \in \Lean}$; for recent work on the computation of the lean kernel see \cite{SIMML2014EfficientAutarkies,KullmannSilva2015OracleSATC}. The decision of leanness is coNP-complete, and so consideration of special autarkies is of interest; actually, these considerations are not just ``algorithmic hacks'', but in a sense represent various areas of combinatorics (for example matching theory) via ``autarky systems''.

\paragraph{Autarky systems}

The notion of an ``autarky system'', as a selection of special autarkies with similar good properties as general autarkies, was introduced in \cite{Ku00f}, partially further expanded in \cite{Kullmann2007Balanciert}, and overviewed in \cite[Subsection 11.11]{Kullmann2007HandbuchMU}. The starting point for an autarky system is to single out a restricted notion of autarky. This restricted autarky notion implies a restricted satisfiability notion, namely clause-sets satisfiable via (iterated) autarky reduction, using only these special autarkies. This is indeed equivalent for ``normal autarky systems'' to satisfiability by a single such special autarky.\footnote{``Normal autarky systems'' were called ``strong autarky systems'' in \cite[Section 8]{Ku00f}.} Furthermore, for such normal systems reversely the general autarkies of the system can be derived from those special autarkies, just using the satisfying assignments amongst them. For arbitrary autarky systems also the notions ``minimal unsatisfiability'' and ``lean'' are defined, and are central properties.

\emph{Balanced autarkies} yield an example of a rather general autarky system, the basis of autarkies for hypergraph colouring; here for an autarky, touched clauses need not only have some satisfied literal, but also some \emph{falsified} literal. The corresponding satisfiability notion is ``NAE-satisfiability'', and will be further discussed in Subsection \ref{sec:introLEAN}.

\paragraph{Matching autarkies}

The autarky system especially of importance in this \Schrift, besides the full system, is that of \emph{matching autarkies}; for a short introduction see \cite[Subsection 11.11.2]{Kullmann2007HandbuchMU}. They yield the set $\Mlean \supset \Lean$ of \emph{matching-lean clause-sets}, and the set $\Msat \subset \Sat$ of \emph{matching-satisfiable clause-sets} (called ``matched clause-sets'' in \cite{FrGe98}):
\begin{itemize}
\item A matching autarky for $F \in \Cls$ is an autarky $\vp$ for $F$ such that for all $C \in F$ touched by $\vp$ one can select $x_C \in C \cap \vp$ such that the underlying variables $\var(x_C)$ are pairwise different.
\item We have $F \in \Msat \Lra \fa\, F' \sse F : \delta(F) \le 0$, i.e., $\delta^*(F) = 0$.
\item And $F \in \Mlean \Lra \fa\, F' \subset F : \delta(F') < \delta(F)$.
\item Thus for $F \in \Mlean$ holds $\delta^*(F) = \delta(F)$, and for $F \ne \top$ holds $\delta(F) \ge 1$ (note $\delta(\top) = 0$), a vast generalisation of this fact for $\Musat$.
\item Stronger we have $F \in \Mlean \Lra \surp(F) \ge 1$ for $F \ne \top$ (recall the surplus).
\item Every $F \in \Cls$ has a largest matching-lean sub-clause-set, the \emph{matching-lean kernel}, namely $\bc \set{F' \sse F : F' \in \Mlean}$, computable in polynomial time (for example via reduction by matching autarkies).
\end{itemize}

\paragraph{Linear autarkies}

A stronger autarky system than matching autarkies is given by ``linear autarkies''; we will not use them for the results of this \Schrift, but they are an important link to combinatorics, and so we discuss them here; see \cite[Subsection 11.11.3]{Kullmann2007HandbuchMU} for a more elaborated introduction. ``Simple linear autarkies'' for $F \in \Cls$ have been introduced in \cite{Ku98e}, based on linear programming. For $F \in \Cls$ we consider the clause-variable matrix $M(F)$, which is a $c(F) \times n(F)$ matrix over $\RR$ (or over $\QQ$ for computational purposes), which encodes in the rows the clauses and in the columns the variables, by using $0$ for absence of the variable, and $\pm 1$ for positive resp.\ negative sign. Now the \emph{simple linear autarkies} $\vp$ are obtained from solutions $\vec{x} \in \RR^{n(F)}$ of $M(F) \cdot \vec{x} \ge 0$, by translating the values $\vec{x}_i$, where the indices $i$ correspond to the variables of $F$, into ``unassigned'' for $\vec{x}_i = 0$, ``true'' (i.e., $1$) for $\vec{x}_i > 0$, and ``false'' (i.e., $0$) for $\vec{x}_i < 0$. It is an easy exercise to see that this yields indeed autarkies. We have a non-trivial simple linear autarky iff $M(F) \cdot \vec{x} \ge 0$ has a non-trivial solution. We obtain the classes $\Llean$ of ``linearly lean clause-sets'' (not having a non-trivial simple linear autarky), with $\Lean \subset \Llean \subset \Mlean$, and $\Lsat$ of ``linearly satisfiable clause-sets'' (satisfiable by a sequence of simple linear autarkies), with $\Msat \subset \Lsat \subset \Sat$.

\emph{Linear autarkies}, as introduced in \cite{Ku00f}, are obtained from simple linear autarkies by composition, corresponding to iterated reduction by simple linear autarkies; simple linear autarkies yield an autarky system, while linear autarkies yield a normal autarky system. The point here is, that the reduction to the linearly-lean kernel can now be done by a \emph{single} linear autarky, and linearly satisfiable clause-sets are now satisfiable by a \emph{single} linear autarky. In Subsection \ref{sec:introLEAN} we discuss the special case of ``balanced linear autarkies''. For recent developments see \cite{KimuraMakino2014LinearForms}.

\subsection{Connections to combinatorics}
\label{sec:relwork}

We now discuss the connections between SAT and combinatorics in a wider context than the degree considerations of this \Schrift, concentrating on aspects related to minimal unsatisfiability and autarkies (if one is only interested in the results of this \Schrift, then these discussions may be ignored). A general source on SAT is the handbook \cite{2008HandbuchSAT}; a classical connection to combinatorics, random satisfiability, is discussed in Chapter 8 (\cite{Ach09HBSAT}) there, and of further general interest to combinatorics is Chapter 10 (\cite{Sak09HBSAT}) on symmetry (group theory), Chapter 13 (\cite{SS09HBSAT}) on fixed-parameter tractable problems (for example treewidth and related notions), and Chapter 17 (\cite{Zha09HBSAT}) on the handling of various combinatorial designs by SAT solving, for example from Ramsey theory. Ramsey theory has strong connections to hypergraph colouring, which we discuss next. Indeed, applying SAT solving to hypergraph colouring problems is a powerful tool, and a recent overview can be found in \cite{AhmedKullmannSnevily2011VdW3kArt}, where especially van-der-Waerden numbers are discussed, while a recent success concerning colouring Pythagorean triples can be found in \cite{HeuleKullmannMarek2016Pythagorean}.

\subsubsection{Hypergraph colouring}
\label{sec:prelimcolouring}

Hypergraph-colouring, especially 2-colouring, and SAT are closely connected; see \cite[Section 5]{Duchet1995Hypergraphs} for a general introduction and overview on hypergraph colouring (from the combinatorial point of view), while a monograph is given with \cite{JensenToft1995Graphenfaerbung}. An overview especially on the question of the minimum number of hyperedges for a given number of vertices in non-$k$-colourable hypergraphs is given in \cite{Kostochka2006UebersichtKritischeFaerbung}.

\paragraph{Hypergraphs}

For this introduction, a hypergraph $G$ is a finite set of finite subsets of $\ZZ$; so $G$ itself is the set of hyperedges, i.e., $E(G) := G$, while $\bc G$ is the set of vertices, i.e., $V(G) := \bc G$. The set of all hypergraphs is denoted by $\Hyp$. Let the deficiency be $\defh(G) := \abs{E(G)} - \abs{V(G)} \in \ZZ$. Note that clause-sets are special hypergraphs ($\Cls \subset \Hyp$), but their deficiency is defined differently. Hypergraphs $G$ with $\defh(G) = 0$ are called \emph{square hypergraphs}. Special hypergraphs are the positive clause-sets, and the set of all positive clause-sets is denoted by $\Poscls := \set{F \in \Cls : \bc F \subset \NN} = \set{G \in \Hyp : V(G) \subset \NN}$. For $F \in \Poscls$ we have $\delta(F) = \defh(F)$; obviously every hypergraph can be renamed to a positive clause-set. From general clause-sets $F \in \Cls$ we obtain (directly) two hypergraphs:
\begin{enumerate}
\item $F$ itself is a hypergraph (breaking the link between positive and negative literals, which are now just unrelated vertices).

  We note that we could have allowed $\Cls = \Hyp$, by allowing tautological clauses (i.e., clauses containing clashing literals) and self-complementary literals ($-0 = 0$). In certain contexts allowing such degenerations has advantages, but in our context is seems best to ban them (for example so we have a direct correspondence between clauses and partial assignments).
\item The ``variable-hypergraph'' of $F$ is $\set{\var(C) : C \in F} \in \Poscls$. This formation for example is important to apply methods from matching theory.
\end{enumerate}
For positive clause-sets both formations collapse to the identity, and we treat positive clause-sets as representing (general) hypergraphs by (special) clause-sets.

\paragraph{Colouring}

A $k$-colouring for $k \in \NNZ$ of $G$ is a map $f: V(G) \ra \tb 1k$ such that for all $H \in G$ there are $x, y \in H$ with $f(x) \ne f(y)$. $G$ is called $k$-colourable if there exists a $k$-colouring of $G$; instead of saying that $G$ is 2-colourable, one also says that $G$ ``has Property B''. Note that if there are $H \in G$ with $\abs{H} \le 1$, then $G$ is not $k$-colourable for any $k$. A hypergraph $G$ is \emph{critically $k$-colourable} if $G$ is $k$-colourable and not $(k-1)$-colourable, while for each $H \in G$ the hypergraph $G \sm \set{H}$ is $(k-1)$-colourable. In the SAT-context there is no need to discard hyperedges containing at most one vertex, and then \emph{minimally non-$k$-colourability} is more appropriate, that is $G$ is not $k$-colourable (possibly not colourable at all), while after removal of any hyperedge $G$ becomes $k$-colourable. The set of all minimally non-$k$-colourable hypergraphs is denoted by $\Mnc[k] \subset \Hyp$ for $k \in \NNZ$. We have $\set{\es}, \set{\set{x}} \in \Mnc[k]$ for all $k \in \NNZ$ and $x \in \ZZ$.

We are especially interested in $\Mnc[2]$. For $G \in \Mnc[2]$ holds $\defh(G) \ge 0$, as shown in \cite{Seymour19742Faerbung}, and so we can consider the sets $\Mnc[2]_{\defh=k}$ for deficiencies $k \in \NNZ$ (all minimally non-$2$-colourable hypergraphs of deficiency (exactly) $k$).\footnote{Indeed in \cite[Corollary 8.2]{Kullmann2005b2} it is shown $\defh(G) \ge 0$ for all $G \in \Mnc[k]$ and all $k \ge 2$, as a simple application of the autarky method; note that for $G := \set{\set{1,\dots,n}} \in \Mnc[k]$ for $k \le 1$ and $n \ge 2$ holds $\defh(G) = 1 - n < 0$.}  The famous problem of deciding in polynomial time, whether a directed graph contains an even cycle, is equivalent to the problem of deciding ``$F \in \Mnc[2]_{\defh=0}$ ?'' for $F \in \Hyp$ (via simple transformations), and this problem was finally solved in \cite{RobertsonSeymourThomas1999GeradeKreise,McCuaig2004PolyasProblem}. It was conjectured in \cite{Kullmann2007Balanciert}, that for all $k \in \NNZ$ the classes $\Mnc[2]_{\defh=k}$ are decidable in polynomial time (see also \cite[Conjecture 11.12.1]{Kullmann2007HandbuchMU}). More on this in Subsection \ref{sec:prelimQCA}. In \cite{AbbottHare1999SquareColor3critical} one finds more information on vertex degrees in uniform elements of $\Mnc[2]_{\defh=0}$ (i.e., where all hyperedges have the same length).

\paragraph{Translating hypergraphs into clause-sets}

For a positive clause-set $G \in \Poscls \subset \Hyp$ we obtain the translation of $2$-colouring to satisfiability via $F_2: \Poscls \ra \Cls$
\begin{displaymath}
  F_2(G) := G \cup \set{-H : H \in G} \in \Cls.
\end{displaymath}
For a general discussion of such translations, also considering more colours, see \cite[Subsection 1.2]{Kullmann2007ClausalFormZII}. A hypergraph $G \in \Poscls$ is 2-colourable iff $F_2(G)$ is satisfiable, and $G$ is minimally non-$2$-colourable iff $F_2(G)$ is minimally unsatisfiable, i.e., $F_2(G) \in \Musat \Lra G \in \Mnc[2]$ (this is easy to prove, and a special case of \cite[Lemma 8.1]{Kullmann2005b2}). Regarding the deficiency we have $\delta(F_2(G)) = \defh(G) + \abs{E(G)}$ for $\es \notin G$, and thus e.g.\ $F_2(\Mnc[2]_{\defh=0} \cap \Poscls)$ is not contained in any $\Musati{\delta=k}$ for some $k \in \NN$.

A slight generalisation of the image $F_2(\Poscls)$ under this translation is the class of complementation-invariant clause-sets $F \in \Cls$, characterised by $C \in F \Lra -C \in F$ for clauses $C$, as introduced in \cite{Kullmann2007Balanciert} (see also \cite[Subsection 11.4.5]{Kullmann2007HandbuchMU}): The image $F_2(\Poscls)$ is the set of complementation-invariant PN-clause-sets, that is, clause-sets $F$ where every clause $C \in F$ is positive (i.e., $C \subset \NN$) or negative ($-C \subset \NN$). See Subsection \ref{sec:prelimQCA} for how autarkies, considered for $F_2(G)$, can help understanding $G$.

We mention here the study of the quantity $m(p)$ for $p \in \NN$, which is the minimal number of hyperedges in $p$-uniform elements of $\Mnc[2]$, where the known precise values are $m(1) = 1$, $m(2) = 3$, $m(3) = 7$, $m(4) = 23$; see \cite{MathewsPandaShannigrahi2015MNC2} for a recent article. By the above remarks we see that $m(p)$ equals the minimal $c(F) / 2$ for complementation-invariant PN-clause-sets $F \in \Punimusat{p}$ ($p$-uniform MUs). As an aside, the minimal $c(F)$ for general $F \in \Punimusat{p}$ is $2^p$, realised by the $A_p$ (proof is a simple exercise for the reader). Back to hypergraph colouring, for the more general $m_k(p)$, the minimal number of hyperedges in $p$-uniform elements of $\Mnc[k]$, see \cite{Shabanov2012MNCr} for a recent article (also discussing $m^*(p,k)$, which considers only linear hypergraphs (hyperedges share at most one vertex) and further variations). Generalisations to non-uniform hypergraphs are discussed in \cite{Shabanov2014NonUniform}, considering various hyperedge-weights with exponential decay, similar to the use of clause-weights in heuristics for SAT solvers (see \cite{Kullmann2007HandbuchTau} for an overview). The argumentation of \cite{Shabanov2012MNCr,Shabanov2014NonUniform} is probabilistic (lower bounds), while \cite{MathewsPandaShannigrahi2015MNC2} works constructively (upper bounds).

\paragraph{Translating clause-sets into hypergraphs}

In the other direction a translation $e: \Cls \ra \Hyp$ was provided in \cite{LinialTarsi1985HResolution}. For $F \in \Cls$ let
\begin{displaymath}
  e(F) := \set{C \addcup \set{0} : C \in F} \addcup \set{\set{v,-v} : v \in \var(F)} \in \Hyp
\end{displaymath}
(where ``$\addcup$'' is just union, but highlighting disjointness); e.g.\ $e(\set{\bot}) = \set{\set{0}}$. The hypergraph $e(F)$ is 2-colourable iff $F$ is satisfiable, and $F$ is minimally unsatisfiable iff $e(F)$ is minimally non-$2$-colourable, i.e., $e(F) \in \Mnc[2] \Lra F \in \Musat$ (the direction ``$\La$'' of the latter statement is stated in the proof of Theorem 3 in \cite{AhLi86}, the other direction is (also) very easy). Furthermore $\defh(e(F)) = \delta(F) - 1$. Thus $e$ embeds the classes $\Musati{\delta=k}$ into $\Mnc[2]_{\defh=k-1}$, which motivates the conjecture, that all $\Mnc[2]_{\defh=k}$ for $k \in \NNZ$ are polytime decidable, as a strengthening of the polytime decision of the $\Musati{\delta=k}$ for $k \in \NN$ (recall Subsection \ref{sec:prelimMU}).

We remark that via this embedding $e$ we obtain a proof of $\delta(F) \ge 1$ for $F \in \Musat$ from $\defh(G) \ge 0$ for $G \in \Mnc$ (this is one of the proofs given in \cite{AhLi86}). In \cite{AhLi86} also an alternative proof of $\defh(G) \ge 0$ is given, based on matching theory, plus one further proof of $\delta(F) \ge 1$, using linear algebra, as in \cite{Seymour19742Faerbung}. In Subsection \ref{sec:introLEAN} we will further comment on these proofs, as they are unfolded in the theory of autarkies.

We also remark, that the hypergraph class $e(\Musati{\delta=1}) \subset \Mnc[2]_{\defh=0}$ has the property, that every hypergraph in it different from $e(\set{\bot}) = \set{\set{0}}$ has a vertex of degree $2$ (since every $F \in \Musati{\delta=1}$ different from $\set{\bot}$ has a variable of degree $2$). More generally, for all $k \in \NN$ every hypergraph in $e(\Musati{\delta=k}) \sm \set{\set{\set{0}}} \subset \Mnc[2]_{\defh=k-1}$ has a vertex of degree at most $k+1$ (by the basic degree bound as shown in Subsection \ref{sec:introbasicdegree}).

\begin{quest}\label{que:vertdeghyp}
  Are the minimum vertex-degrees of general $G \in \Mnc[2]_{\defh=k}$ for (fixed) $k \in \NNZ$ bounded?
\end{quest}
In \cite[Proposition 1]{Kostochka2006UebersichtKritischeFaerbung} we find the easy proof that for every $G \in \Mnc[k]$, where every hyperedge has size at least $2$, every vertex has degree at least $k$ (exercise: show that hyperedges containing vertices of degree at most $k-1$ can be removed -- this leads to the ``$k$-core'' of a hypergraph $G$, the largest $G' \sse G$ such that every $v \in V(G')$ occurs at least $k$ times in $G'$). This has been algorithmically exploited for $k=2$ in \cite[Subsection 5.2]{HeuleKullmannMarek2016Pythagorean}, in the form of removal of ``blocked clauses''. Question \ref{que:vertdeghyp} is related to literal-degrees in $\Musat$, while the apparently more interesting question about \emph{variable-}degrees in $\Musat$ does not seem to have a natural equivalent for hypergraph colouring.

\subsubsection{Hypergraph transversals}
\label{sec:introtrans}

For $G \in \Hyp$ let $\Tr(G) \in \Hyp$, the \emph{transversal hypergraph} of $G$, be defined as the set of all minimal $T \sse V(G)$ such that $T \cap H \ne \es$ for all $H \in G$. The \emph{Transversal Hypergraph Problem} is the computational problem, given $G, G' \in \Hyp$, to decide whether $\Tr(G) = G'$ holds. Equivalently, the input is $G \in \Hyp$, and it is to be decided whether $G = \Tr(G)$ holds (obviously this is a special case of the Transversal Hypergraph Problem, and by a polynomial-time translation the general case can be reduced to it). For an overview on this important problem and its many guises see \cite{EiterMakinoGottlob2008MonotoneDual}. It is known that the problem is solvable in quasi-polynomial time, and the long outstanding problem is whether it can be solved in polynomial time.

An \emph{intersecting hypergraph} is a hypergraph $G \in \Hyp$ with $G \sse \Tr(G)$, and the class of all intersecting hypergraphs is denoted by $\Ihyp \subset \Hyp$. It is not hard to see that for $G \in \Ihyp$ holds $G \in \Mnc[2]$ iff $\Tr(G) = G$. Since from $G = \Tr(G)$ follows $G \in \Ihyp$, thus the Transversal Hypergraph Problem is equivalent to the problem of deciding, whether an intersecting hypergraph is minimally non-$2$-colourable. The natural question arises for the decision of the classes $(\Mnc \cap \Ihyp)_{\defh=k}$ for $k \in \NNZ$. The case $k=0$ has been handled in \cite{Seymour19742Faerbung}, indeed not just deciding the class in polynomial time, but efficiently classifying the elements. The cases $k \ge 1$ appear to be open, and whether decision is possible in polynomial time for fixed $k$, or is even fixed-parameter tractable (fpt) in $k$, is an interesting test case for the general Hypergraph Transversal Problem, as well as it is relevant for the understanding of minimally non-$2$-colourable hypergraphs.

The translation of intersecting hypergraphs $G \in \Ihyp$ into clause-sets $F_2(G) \in \Cls$ yields also a natural and interesting class of clause-sets. \emph{Bihitting clause-sets}, introduced in \cite[Subsection 4.2]{GalesiKullmann2003bHermitian}, are those $F \in \Cls$ where $F', F'' \sse F$ with $F' \cup F'' = F$, $F' \cap F'' = \es$ exist, such that for all $C' \in F', C'' \in F''$ holds $C' \cap -C'' \ne \es$, while $F', F''$ itself are clash-free (i.e., $(\bc F') \cap - (\bc F') = \es$, and $(\bc F'') \cap - (\bc F'') = \es$). Obviously, the images under $F_2$ of intersecting hypergraphs are precisely the bihitting complementation-invariant PN-clause-sets (i.e., the set of bihitting clause-sets in the image of $F_2$) different from $\set{\bot}$, and deciding their minimal unsatisfiability is thus another manifestation of the Hypergraph Transversal Problem, directly related to the decision ``$\Tr(G) = G $?''. And another one is to decide SAT for general bihitting clause-sets (as can be easily seen, and is discussed in \cite[Subsection 4.3]{GalesiKullmann2003bHermitian}), this time directly related to the decision ``$\Tr(G) = G'$ ?''.

In \cite[Theorem 8.14]{Kullmann2005b2} (the first 6 sections are covered by \cite{Kullmann2007ClausalFormZI,Kullmann2007ClausalFormZII}) the characterisation of \cite{Seymour19742Faerbung} (the intersecting $G \in \Mnc[2]_{\defh=0}$) is translated into $\Cls$-language.

\subsubsection{Autarkies for hypergraphs}
\label{sec:introLEAN}

We discuss here now two autarky systems (recall Subsection \ref{sec:introautgen} for a general introduction), which are especially relevant for hypergraph colouring.

\paragraph{Balanced autarkies}

\cite{Kullmann2007Balanciert} introduced \emph{balanced autarkies} for $F \in \Cls$, partial assignments $\vp$, which in every clause of $F$ they touch satisfy as well as falsify at least one literal (that is, for $C \in F$ with $C \cap (\vp \cup -\vp) \ne \es$ holds $C \cap \vp \ne \es$ as well as $C \cap -\vp \ne \es$); \cite[Subsection 11.11.4]{Kullmann2007HandbuchMU} provides an introduction. This is a normal autarky system, and thus we basically have all the good property general autarkies have. Balanced autarkies are closely related to hypergraph colouring. The balanced autarkies for $F$ are precisely the autarkies of $F \cup \set{-C : C \in F}$, and every autarky for a complementation-invariant clause-set is automatically balanced. A clause-set is balanced-satisfiable, i.e., can be satisfied by a balanced autarky, iff it is NAE-satisfiable (``not-all-equal''; see \cite{PorschenRanderathSpeckenmeyer2003NAESAT} for basic results).

Balanced autarkies provide the general autarky form for $\Poscls$ (whose elements are all trivially satisfiable, and thus unrestricted autarkies are not of interest here), which represents hypergraphs for the $2$-colouring problem: an $F \in \Poscls$ is $2$-colourable iff it is balanced-satisfiable, and $F$ is minimally non-$2$-colourable iff it is minimally balanced-unsatisfiable.  Finally we have \emph{balanced lean} clause-sets (i.e., having no non-trivial balanced autarkies), and this is the appropriate notion of ``leanness'' for hypergraphs, as represented by the class $\Poscls$; more precisely, a hypergraph $G$ is lean iff for an isomorphic $F \in \Poscls$ (isomorphic as hypergraph) we have that $F$ is balanced lean. For lean hypergraphs $G$ we have $\defh(G) \ge 0$, and this is indeed more generally treated by ``balanced \emph{linear} autarkies''.

\paragraph{Balanced linear autarkies}

The special case of ``balanced linear autarkies'' was introduced in \cite[Section 6]{Ku00f}; these are the simple linear autarkies for $F \cup \set{-C : C \in F}$ (recall Subsection \ref{sec:introautgen}).\footnote{More precisely one should speak of ``balanced simple linear autarkies'', but for convenience ``simple'' is dropped. We note that ``balanced linear autarkies'' are balanced and linear autarkies, but in general a balanced and linear autarky need not be a balanced linear autarky, and thus one should speak of ``balanced-linear autarkies''; again we abuse language, motivated by the fact that linear autarkies which are also balanced are apparently too general a concept to be useful.} Equivalently, the balanced linear autarkies $\vp$ for $F \in \Cls$ are obtained from solutions $\vec{x} \in \RR^{n(F)}$ of $M(F) \cdot \vec{x} = 0$, by translating the values $\vec{x}_i$ as discussed before (it is an easy exercise to see that this yields indeed balanced autarkies). We have a non-trivial balanced linear autarky iff $M(F) \cdot \vec{x} = 0$ has a non-trivial solution, and so, in other words, $F$ is balanced linearly lean iff the columns of $M(F)$ are linearly independent (iff $\rank(M(F)) = n(F)$). Thus if $F \in \Cls$ is balanced linearly-lean, then $\delta(F) \ge 0$ holds; furthermore, as shown in \cite[Lemma 7.2]{Kullmann2005b2}, there is then a matching in the clause-variable graph covering all variable nodes, and thus even $\delta^*(F) = \delta(F)$ holds. By noting that $F \in \Cls$ is balanced linearly lean iff $F \cup \set{-C : C \in F}$ is linearly lean, and considering $\Poscls$, we obtain that for lean hypergraphs $G$ (especially, minimally non-2-colourable) we have $\defh(G) \ge 0$. To say the argument again explicitly: Consider a hypergraph $G \in \Poscls$; then $G$ (as a clause-set) is balanced linearly lean iff the variable-clause matrix has linearly independent rows, iff $F_2(G)$ is linearly lean (again, as a clause-set), which is implied by $F_2(G)$ being minimally unsatisfiable (or weaker, being lean), which in turn is equivalent to $G$ (as a hypergraph) being minimally-non-$2$-colourable. This conclusion ``The rows of the incidence matrix [our variable-clause matrix] of a minimally-non-$2$-colourable hypergraph are linearly independent over $\RR$.'' is shown in \cite{Seymour19742Faerbung}; see \cite[Lemma 4.7]{Godsil1995LinearAlgebra} for this and related results, while the conclusion ``$\defh(G) \ge 0$'' is discussed as Principle 2.1 in \cite{Godsil1995LinearAlgebra}. For properties of minimally balanced linearly unsatisfiable clause-sets see \cite[Section 4]{Kullmann2007Balanciert}.

\paragraph{Fundamental inequalities}

We have yet seen two fundamental inequalities, namely $\delta(F) \ge 1$ for $F \in \Mlean$, as first shown in \cite{AhLi86} for minimally unsatisfiable clause-sets (``Tarsi's Lemma''), and $\delta(F) \ge 0$ for balanced linearly lean clause-sets, first shown in \cite{Seymour19742Faerbung} (as $\defh(G) \ge 0$ for minimally non-2-colourable hypergraphs).\footnote{An application yielding Fisher's inequality (design theory) is discussed in  \cite[Subsection 7.4]{Kullmann2005b2} (while Seymour's inequality is discussed there in Subsection 8.2).}Autarky theory shows the general structure of the arguments: We find ``obstructions'', which prevent these bounds from holding, where such obstructions are given by a subset $F' \sse F$ where there is a partial assignment $\vp$ with $\vp * F' = \top$, while $\var(\vp) \cap \var(F' \sm F) = \es$. Now minimally unsatisfiable $F$ do not have such $F'$, and thus the envisaged bound holds for them, and this is the argumentation in e.g.\ \cite{Seymour19742Faerbung,AhLi86}.

But one can go beyond this, exploiting autarky \emph{reduction} $F \leadsto F \sm F'$ (as long as possible). Note that the above $\vp$ is precisely an autarky, and furthermore possibly one of a special structure. If we just look at general autarkies, then we obtain the first generalisation, to lean clause-sets (having no non-trivial autarkies) or balanced lean clause-sets (having no non-trivial balanced autarkies, covering the hypergraph cases). However often, due to the special structure, these special autarkies can be found in polynomial time, and their application yields some $F_0 \sse F$, such that the bound holds for $F_0$ (while for $F \in \Musat$ we just have $F_0 = F$). If we have even an ``autarky system'', then $F_0$ is uniquely determined, that is, does not rely on the choice of the autarkies in the reduction process. The case of main importance for this \Schrift{} is $\delta(F) \ge 1$, where the autarkies are matching autarkies, and the reduced $F_0$ is the matching-lean kernel of $F$, while those $F$ with $F_0= F$ are precisely the $F \in \Mlean$. On the other hand, for hypergraph colouring the fundamental fact is $\delta(F) \ge 0$ for balanced linearly-lean clause-sets, where the autarkies are balanced linear autarkies, and the reduced $F_0$ is the balanced-linearly-lean kernel of $F$. In fact, via autarky reduction we obtain a general method to study \emph{decompositions}, which we will discuss next, in the context of matrix analysis.

\subsubsection{Qualitative matrix analysis (QMA)}
\label{sec:prelimQCA}

QMA can be understood as the analysis of matrices $M$ over the real numbers in abstraction of the absolute value of the entries, but only their signs count. More precisely, one considers the \emph{qualitative class} $\Q(M)$, which consists of all matrices with the same dimensions as $M$, which have entry-wise the same signs as $M$ (positive, zero, negative), and investigates when a property of $M$ holds for all $M' \in \Q(M)$. For example, a matrix $M$, such that all $M' \in \Q(M)$ have linearly independent rows, is called an \emph{$L$-matrix}. The monograph \cite{BS95} is an excellent source on QMA until the 1990's, while a more recent overview is given in \cite{HallLi2007SignPatternMatrices}.\footnote{For us the original notion of $L$-matrix, before \cite{BS95}, with (qualitatively) independent \emph{columns}, would be more convenient, but we stick to the (important) monograph \cite{BS95}.} An early demonstration of the close relations to SAT is \cite{KleeLadner1981WeakSAT}, where ``weak satisfiability'' was introduced, which is \emph{precisely} the existence of a non-trivial autarky. \cite[Theorem 5]{KleeLadner1981WeakSAT} shows that weak satisfiability is NP-complete, which is the earliest known proof of $\Lean$ being coNP-complete. However these connections to SAT have not been pursued further.

Starting from \cite{DD92}, which exploits Farkas' lemma to understand (un)satisfiability, the connections to QMA have been explored in \cite[Sections 3, 5]{Ku00f}; see \cite[Subsection 11.12.1]{Kullmann2007HandbuchMU} for a more substantial introduction. It is shown in \cite[Remark 5 in Section 5]{Ku00f}, that $L$-matrices correspond (nearly) precisely (up to transposition and handling of zero-rows/columns and repeated rows/columns) to balanced lean clause-sets, while lean clause-sets correspond (nearly) precisely to so-called \emph{$L^+$-matrices} (as investigated in \cite{LS1998}). The square $L$-matrices are called \emph{SNS-matrices}; SNS-matrices are at the heart of the poly-time decision for $\Mnc[2]_{\defh=0}$ (recall Subsection \ref{sec:prelimcolouring}), and the connections to autarky theory are explored in \cite{Kullmann2007Balanciert}; see \cite[Subsection 11.12.2]{Kullmann2007HandbuchMU} for an overview.

Further with the translation of terminology, now regarding unsatisfiability: unsatisfiable clause-sets correspond to \emph{sign-central matrices}, minimally unsatisfiable clause-sets correspond to \emph{minimal sign-central matrices}, and unsatisfiable hitting clause-sets (MUs, where every pair of clauses has a clash) correspond to \emph{tight sign-central matrices} (\cite{HKKZ2003Tight}). So \cite[Theorem 5.4.3]{BS95} is yet another proof of $\delta(F) \ge 1$ for $F \in \Musat$ (``Tarsi's Lemma''). The variable-degree, as studied in the current \Schrift, corresponds to the number of non-zero entries in the rows of the matrices (while the deficiency is the difference of the number of columns and the number of rows). The elements of $\Musati{\delta=1}$ correspond to \emph{$S$-matrices} (thus \cite[Corollary 2.2]{KLM1984} yields perhaps the first proof of $\minnonmer(1) = 2$), while the elements of $\Smusati{\delta=1}$ correspond to \emph{maximal $S$-matrices}.

So our Theorem ``$F \in \Lean \Ra \minvdeg(F) \le \nonmer(\delta(F))$'' is equivalent to the property of $L^+$-matrices of dimension $m \times n$ (thus $n > m$), that there always exists a row with at most $\nonmer(n-m)$ many non-zero entries.

As mentioned, autarky systems $\mc{A}$ (like balanced autarkies, matching autarkies, etc.) also yield a framework for decomposition theorems. The basic decomposition is into $\mc{A}$-lean and $\mc{A}$-satisfiable sub-clause-sets, as given in \cite[Theorem 8.5]{Ku00f} for normal autarky systems. This corresponds to a unique representation of the clause-variable matrix as a $2 \times 2$ triangular block matrix, and generalises various matrix decompositions in QMA, as discussed in \cite[Footnote 7, Page 246]{Ku00f}. Furthermore, $\mc{A}$-lean clause-sets itself can be decomposed, considering a triangular decomposition into $\mc{A}$-lean blocks for the clause-variable matrix. The main result is \cite[Lemma 6]{Kullmann2007Balanciert}, reviewed in \cite[Subsection 11.11.5]{Kullmann2007HandbuchMU} and generalising \cite[Theorem 2.2.5]{BS95}: A clause-set $F \in \Cls$ is minimally $\mc{A}$-unsatisfiable iff $F$ is barely $\mc{A}$-lean (it is lean, but removal of any single clause destroys this) and $\mc{A}$-indecomposable.

\subsubsection{Biclique partitions of (multi-)graphs, and algebraic graph theory}
\label{sec:introbicliques}

We finish this overview on related themes in combinatorics by a field of graph theory, which, like QMA, can be understood as a study of clause-sets from a special angle, focusing on the conflict-structure of clauses.

\paragraph{Certain aspects of algebraic graph theory}

The starting point is \cite{GH1971}, where the problem of ``addressing a graph'' was introduced. One considers a symmetric matrix $D$ of dimension $m \in \NN$ over $\NNZ$, with a zero-diagonal, where the entries are interpreted as ``distances'' (in \cite{GH1971} the $D_{i,j}$ are the distances between the nodes of some graph), and asks for the smallest $N \in \NNZ$, such that there are $m$ codewords $c_1,\dots,c_m \in \set{0,1,*}^N$ with the property, that the modified Hamming distance between $c_i$ and $c_j$, which simply ignores positions with $*$, is $D_{i,j}$. See \cite[Chapter 9]{LW1992b} for an introduction.

A basic result of \cite{GH1971} is that if $D$ has all entries outside the diagonal equal to $1$, then $N = m-1$ (see also \cite[Lemma 6.6]{Godsil1995LinearAlgebra} for a discussion in the context of eigenvalue methods; for a direct combinatorial proof see \cite{Vishwanathan2013GrahamPollak}). This follows from the general result $N \ge \max(n_+(D), n_-(D))$ of \cite{GH1971} (the ``Lemma of Witsenhausen''), where $n_+(D)$ resp.\ $n_-(D)$ is the number of positive/negative eigenvalues of $D$. For the general case in \cite{Wink1983Squashed} it is shown, that if the distances $D_{i,j}$ are indeed the distances between the nodes of some graph, then we have $N \le m-1$. 

Actually, the Lemma of Witsenhausen works for arbitrary symmetric matrices $D$ over $\NNZ$ with zero diagonal, and a shift of perspective is useful. A codeword over $\set{0,1,*}^N$ is nothing else than a clause over the variables $1,\dots,N$, while the modified Hamming distance is the number of clashes (conflicts). So the question is about the existence of clauses $C_i$ for $i \in \tb 1m$ over variables $1,\dots,N$, such that $D_{i,j} = \abs{C_i \cap -C_j}$ for $i, j \in \tb 1m$. The above parameter $N = N(D)$ is the minimal number $N$ of variables in a (multi-)clause-set representing $D$ (it is an easy exercise to see that $N$ is finite, i.e., a representation is always possible). Considering $D$ as the adjacency matrix of some multigraph (where parallel edges are allowed), we see that $N$ is also equal to the minimum number of bicliques into which the edge-set of that multigraph can be partitioned, and $N$ is therefore denoted by $\bcp(A) \in \NNZ$ (the ``biclique partition number'' of $A$ resp.\ the corresponding multigraph).\footnote{In \cite{Ku2003e} we used ``$n_{\mathrm{s}}(A)$'' instead, the ``symmetric conflict number''.}

\paragraph{Clause-sets as biclique partitions}

The essential observation is now that we can go back and forth between biclique partitions of multigraphs and clause-sets. In one direction we can understand clause-sets $F$ as representations of biclique partitions of multigraphs, where for each vertex we get a clause, and from each biclique we obtain a variable, where the two sides of the biclique are the positive and negative occurrences of the variable. So we can understand a multigraph together with a biclique partition as a clause-set, and we can use tools from clause-set-logic to analyse the pair multigraph with biclique-partition. The deficiency then becomes the difference between the number of nodes and the number of bicliques. Satisfiability means that it is possible to select from each biclique one side such that all vertices are covered, and minimal unsatisfiability means that such a covering is not possible, but becomes possible as soon as any vertex is removed. The upper bound $\minnonmer(k) \le \nonmer(k)$ thus is equivalent to the theorem, that for every biclique partition with $n$ bicliques of a multigraph with $m$ vertices and with this minimal non-coverability property, there exists a biclique with at most $\nonmer(m - n)$ vertices.

For example $A_0 = \set{\bot}$ represents the $K_1$, the graph with one vertex, and $A_1 = \set{\set{1},\set{-1}}$ represents $K_2$, consisting of one edge (and two vertices), with this edge constituting the biclique partition. More generally, $\set{\set{1,\dots,n},\set{-1},\dots,\set{-n}}$ represents the star with $n$ satellites, where each edge is taken as a biclique; if all the edges shall be taken as one biclique, then we have to use the multi-clause-set $\set{\set{1}, n * \set{-1}}$. To conclude these simple examples, $\set{\set{1,,\dots,n},\set{-1,\dots,-n}}$ represents the dipole $D_n$, consisting of two vertices and $n$ (parallel) edges, together with the biclique partition with $n$ bicliques.

In the other direction we can understand biclique partitions of multigraphs as representations of clause-sets $F$, namely the nodes of the \emph{conflict multigraph} $\cmg(F)$ are given by the clauses, while the edges are the conflicts (clashing literal occurrences $x, -x$), and the bicliques are given by the variables (their positive and negative occurrences).  In this way we can analyse the influence of the ``conflict structure'' on properties of clause-sets; the basic notions, as introduced in \cite{Ku2003e} with underlying report \cite{Ku2003c}, are as follows.

\paragraph{Conflict analysis}

The notion ``hermitian rank'' has been introduced and studied in \cite{GWS1999} for arbitrary hermitian matrices $A$ (square matrices with complex numbers as entries, such that transposing the matrix and taking the complex conjugate of each entry yields back the original matrix), denoted by $h(A) := \max(n_+(A), n_-(A)) \in \NNZ$. So the Lemma of Witsenhausen takes the form, that for symmetric matrices $A$ over $\NNZ$ with zero diagonal holds $\bcp(A) \ge h(A)$.

For $F \in \Cls$ let $\scf(F)$ (the \emph{conflict matrix}) be the square matrix of dimension $c(F)$ over $\NNZ$, with entries $\abs{C \cap -D}$ for $C, D \in F$ (thus with zero diagonal), i.e., $\scf(F)$ is the adjacency matrix of $\cmg(F)$. So we can use the hermitian rank as a measure $h: \Cls \ra \NNZ$ (as first done in \cite[Subsection 3.2]{Ku2003e}), namely
\begin{displaymath}
  h(F) := h(\scf(F));
\end{displaymath}
see Points 1, 3 in \cite[Section 2]{GalesiKullmann2003bHermitian} for various equivalent characterisations.\footnote{For Point 1(c) there it must be a ``\ul{diagonal} matrix $A'$''.} By definition we have $\bcp(F) := \bcp(\scf(F)) \le n(F)$, and thus $h(F) \le \bcp(F) \le n(F)$. Since for a principal submatrix $A'$ of a hermitian matrix $A$ holds $h(A') \le h(A)$ (this follows by ``interlacing''; see \cite[Theorem 9.1.1]{GR2001}), we get $h(\vp * F) \le h(F)$ for all partial assignments $\vp$, and also $h(F') \le h(F)$ for all $F' \sse F$, which gives motivation to consider $h(F)$ as a \emph{complexity measure} for $F \in \Cls$.

In \cite{Ku2003e} also the \emph{hermitian defect} $\hdef: \Cls \ra \NNZ$ has been introduced as
\begin{displaymath}
  \hdef(F) := c(F) - h(F),
\end{displaymath}
and thus $\delta(F) \le \hdef(F)$; see Point 2 in \cite[Section 2]{GalesiKullmann2003bHermitian} for a geometric characterisation (as the ``Witt index'' of the quadratic form associated with $\scf(F)$). Actually $\delta^*(F) \le \hdef(F)$ holds and even stronger properties (see \cite[Subsection 3.3]{Ku2003e}). An important property is (again) $\hdef(\vp * F) \le \hdef(F)$ for all $F \in \Cls$ and partial assignments $\vp$, together with $\hdef(F') \le \hdef(F)$ for all $F' \sse F$, by \cite[Corollary 9]{Ku2003e}, and so we might consider the hermitian defect as a stabilised version of the maximal defect (both are also complexity measures; recall that we have fixed-parameter tractable SAT decision for input $F \in \Cls$ in the parameter $\delta^*(F)$). Note that in general we can have $\delta^*(\vp * F) > \delta^*(F)$, for example $F := \set{\set{1}}$ has $\delta^*(F) = \delta(F) = 0$, while for $F' := \pao 10 * F$ we get $F' = \set{\bot}$, and thus $\delta^*(F') = \delta(F') = 1$. See \cite[Subsection 3.3]{Ku99dKo} and \cite[Subsection 11.2]{Kullmann2007ClausalFormZI} for more information on $\delta^*(\vp * F)$; splitting on a single(!) variable is very important for this \Schrift, with the basic fact $\delta^*(\pao x1 * F) \le \delta(F)$ for $F \in \Musat$ and any literal $x$.

\paragraph{Applications}

The first direct application applied the fact $\delta(F) \le \hdef(F)$ for $F \in \Cls$, namely that for a hitting clause-set $F \in \Clash$ (equivalently, all entries of $\scf(F)$ outside the diagonal are non-zero) with a regular conflict multigraph (i.e., all entries of $\scf(F)$ outside the diagonal are equal) we have $\hdef(F) \le 1$, and thus $\delta(F) \le 1$ (\cite[Theorem 33]{Ku2003e}).\footnote{In \cite{Ku2003e} unfortunately the term ``uniform'' was (mis)used instead of ``regular''.} We get that $\Smusati{\delta=1} = \Uclashi{\delta=1}$ is (precisely) the class of unsatisfiable hitting clause-sets with regular conflict multigraph (\cite[Corollary 34]{Ku2003e}; a combinatorial proof of this was independently found in \cite[Lemma 11]{SloanSzoerenyiTuran2005Primimplikanten_1}), and is also (precisely) the class of unsatisfiable clause-sets $F$ with $\hdef(F) \le 1$ (\cite[Theorem 26]{Ku2003e}).

A clause-set $F \in \Cls$ is called \emph{exact} (\cite[Subsection 3.4]{Ku2003e}) if $\bcp(F) = n(F)$, that is, $F$ is optimal in realising $\cmg(F)$ with respect to the number of variables. Deciding exactness is coNP-complete, while the special class of \emph{eigensharp clause-sets}, defined by $h(F) = n(F)$, or, equivalently, $\hdef(F) = \delta(F)$, is decidable in polynomial time. With \cite[Theorem 14]{Ku2003e} every eigensharp clause-set is matching lean. This leads to \cite[Conjecture 16]{Ku2003e}, ``Every exact clause-set, whose conflict-matrix is the distance matrix of some connected graph, is matching lean.'', which generalises the already mentioned main result of \cite{Wink1983Squashed} (the proof of the ``squashed cube conjecture'').

As already mentioned, we consider $h(F)$ for $F \in \Cls$ as some form of complexity measure, measuring the complexity of representing the conflicts of $F$ via simple matrices. In \cite{GalesiKullmann2003bHermitian} polytime SAT decision in case $h(F) \le 1$ was shown, while the cases $h(F) \le k$ for fixed $k \ge 2$ are open; an interesting stepping stone would be to show polytime SAT decision for $F \in \Cls$ with $\bcp(F) \le k$ (recall $\Cls_{\bcp \le k} \sse \Cls_{h \le k}$). The notion of \emph{blocked clauses}, a special type of clauses which can be removed sat-equivalently, introduced in \cite{Ku97b}, is important here, and \cite[Theorem 3]{GalesiKullmann2003bHermitian} shows, that from $F \in \Cls_{h \le 1}$ after elimination of all blocked clauses (which yields a unique sub-clause-set) we obtain $F' \sse F$ with $\bcp(F') \le 1$. We recall from Subsection \ref{sec:introtrans}, that SAT-decision for $F'$ is now a special case of the Transversal Hypergraph Problem, namely, as shown in \cite[Lemma 11]{GalesiKullmann2003bHermitian}, the problem is exactly the \emph{Exact Transversal Hypergraph Problem}, where every transversal must be ``exact'', that is, must intersect every hyperedge in exactly one vertex; this problem is decidable in polynomial time by \cite{Eiter1994ExactTransversal}, and thus we get SAT-decision for $\Cls_{h \le 1}$ in polynomial time. The characterisation of $F \in \Musat$ with $\bcp(F) \le 1$ is an open problem (while we have polytime membership decision for $\Musati{\bcp \le 1}$), and by \cite[Conjecture 16]{GalesiKullmann2003bHermitian} they would have a very simple structure.

In \cite[Section 6]{Kullmann2007ClausalFormZII} the above basic facts are generalised to non-boolean clause-sets, and that by extending the reduction of multiclique partitions to biclique partitions from \cite{GregoryMeulen1996RPartiteDecompositions} a new and interesting translation from non-boolean to boolean clause-sets was obtained (applied in \cite{Kullmann2010GreenTao} to Ramsey theory).

\subsection{Overview on results}
\label{sec:overv}

Sections \ref{sec:prelim} to \ref{sec:2subr} provide foundations for the main results in the later sections. In Section \ref{sec:prelim} basic notions and concepts regarding clause-sets and autarkies are reviewed. In Section \ref{sec:MUprelim} we discuss minimal unsatisfiability, with some auxiliary results on splitting and saturation (adding literal occurrences to clauses, to make minimal unsatisfiability robust against splitting). Section \ref{sec:vmusat} reviews ``variable-minimal unsatisfiability'', as introduced in \cite{ChenDing2006VMU}, i.e., the class $\Musat \subset \Vmusat \subset \Usat$. There are mistakes in this paper, and we rectify them:
\begin{itemize}
\item we show that $\Vmusat \subset \Lean$ holds (Lemma \ref{lem:vmusatlean});
\item we provide a corrected characterisation of $\Vmusat$ (Lemma \ref{lem:crtivmusat});
\item and we give a corrected proof of polytime decision of $\Vmusati{\delta=k}$ for fixed $k$ in Theorem \ref{thm:vmusatk}, where we also obtain the stronger result, that decision of $\Vmusat$ is fixed-parameter tractable in the deficiency $\delta(F)$.
\end{itemize}
Section \ref{sec:elimcreatesing} is then concerned with singular variables, eliminating them via ``singular DP-reduction'', and creating them via ``singular extensions''. An important auxiliary result is Lemma \ref{lem:DPminvdeg}, showing that eliminating singular variables is harmless for upper bounds on the minimum variable-degree; we also show various auxiliary results on unit-clauses in minimally unsatisfiable clause-sets. This block of preparatory sections is concluded by Section \ref{sec:2subr} on ``full subsumption resolution'', the ubiquitous reduction $ C \oder x, C \oder \ol{x} \leadsto C$, which becomes an ``extension'' in the other direction; as an application, in Theorem \ref{thm:minnumvarmu} we can determine precisely the possible $n(F)$ and $c(F)$ for $F \in \Musati{\delta=k}$.

The first main results (but still on the preparation side) one finds in Section \ref{sec:nonmer}, which introduces the numbers $\nonmer(k) \in \NN$ and proves exact formulas and sharp lower and upper bounds. The point here is that the introduction of $\nonmer(k)$ happens via a recursion which is tailor-made for our application in Section \ref{sec:specialcaseMU}, but which makes it somewhat difficult to determine the numbers in a global way. The analysis of the recursion yields as the main result Theorem \ref{thm:characjumpc}, which shows how the sequence is made up of ``blocks'' of a simple recursive structure. From this in Theorem \ref{thm:solveN} we obtain the general formula.

In Section \ref{sec:specialcaseMU} then we find a basic central result of this \Schrift{}, the upper bound $\minnonmer(k) = \minvdeg(\Musati{\delta=k}) \le \nonmer(k)$ (Theorem \ref{thm:MUminvdegdef}). Section \ref{sec:leansurp} is concerned with generalising this upper bound. An interesting auxiliary class $\Sed \subset \Cls$, clause-sets where deficiency and surplus coincide, is introduced in Subsection \ref{sec:classSED}. The main characterisation of $\Sed$ is given in Theorem \ref{thm:characSEDMSAT} (removal of all clauses containing any given variable yields a matching-satisfiable sub-clause-set), and we obtain Corollary \ref{cor:auxminvardegsigma}, that unsatisfiable elements of $\Sed$ are in fact in $\Vmusat$. In Subsection \ref{sec:proofgencase} the upper bound for $\Musat$ then is lifted to lean clause-sets in Theorem \ref{thm:leanminvardeg}, and also sharpened via replacing the deficiency $\delta$ by the surplus $\surp$. Theorem \ref{thm:vmusharp} shows that the upper bound is sharp for any class between $\Vmusat \cap \Sed$ and $\Lean$.

Section \ref{sec:algoimpl} concerns algorithmic applications. A corollary of Theorem \ref{thm:leanminvardeg} is, that if the asserted upper bound on the minimum variable-degree is not fulfilled, then a non-trivial autarky must exist (Lemma \ref{lem:charakappcor}). Since the variable-set of such a non-trivial autarky is polytime computable, we show in Theorem \ref{thm:genautred} that we can indeed establish the upper bound shown for lean clause-sets also for general clause-sets, after a polytime autarky-reduction. In Subsection \ref{sec:findaut} then the problem of finding a witnessing autarky is discussed, with Conjecture \ref{con:findauthard} making precise our believe that one can find such autarkies efficiently. Theorem \ref{thm:HardMlcr} pinpoints the ``critical'' class $\Mlcr \subset \Sat \cap \Sed$, which is polytime decidable, and where we know that these clause-sets are satisfiable, but where we even don't know how to find any non-trivial autarky efficiently. This block on generalisations of the min-var-degree upper bound is concluded by Section \ref{sec:genboundml}, where we discuss the possibilities to generalise it to matching-lean clause-sets (only the absence of very special (non-trivial) autarkies is guaranteed).

In Section \ref{sec:strengtheningbound} we then turn to the study of the numbers $\minnonmer(k)$, looking now for improved upper bounds and matching lower bounds. We present two infinite classes of deficiencies $k$ with $\minnonmer(k) = \nonmer(k)$, and present a general method of obtaining lower bounds for $\minnonmer(k)$, via counting full clauses (clauses containing all variables --- these clause are strong structural anchors). In Section \ref{sec:improveupbdMU} we introduce a general recursive method to obtain upper bounds like $\nonmer(k)$, via the ``non-Mersenne operator'' $\potprec(f)$, which takes a ``valid bounds function'' $f$, that is, some partial information on $\minnonmer(k)$, and improves it (Definition \ref{def:potprec}). Theorem \ref{thm:basicproppp} shows that this indeed yields a valid method for improving upper bounds on $\minnonmer(k)$, while in Theorem \ref{thm:altrecnonmersenne} we demonstrate how this method recovers $\nonmer(k)$, by just starting with the information $\minnonmer(1) = 2$. In Section \ref{sec:strbMU} we harvest (first) fruits of these methods. First in Theorem \ref{thm:def15} we show $\minnonmer(k) = \nonmer(k)$ for $k \le 5$. Then in Theorem \ref{thm:supminvdegk6} we prove $\minnonmer(6) = \nonmer(k) - 1$ (using a variety of structural results on $\Musat$ provided in this \Schrift{}). Plugging this information on $\minnonmer$ into our machinery, we obtain the improved upper bound $\minnonmer \le \nonmer_1$ in Theorem \ref{thm:imprupperbound}, while in Theorem \ref{thm:imprupperbound} we determine $\nonmer_1(k)$ numerically.

Finally, in Section \ref{sec:open} open problems are stated, thoroughly discussing research perspectives, including (altogether) nine conjectures. Subsection \ref{sec:concmnM} discusses improved upper bounds for $\minnonmer(k)$ from the forthcoming work \cite{KullmannZhao2014Sharper}. Subsection \ref{sec:concSmar} is about improved lower bounds, via counting full clauses. In Lemma \ref{lem:smarandache} we present the lower bound via the ``Smarandache primitive function'' $S_2(k)$ from \cite{KullmannZhao2015FullClauses}, yielding the first-order asymptotic determination of $\minnonmer(k) \sim k$ (Corollary \ref{cor:asympdet}), where now the open question is about the asymptotic determination of $\minnonmer(k) - k$. In Subsection \ref{sec:concgennb} we discuss generalisations to non-boolean clause-sets.

The central Conjecture \ref{con:classmu} of the project of ``understanding MU'', on the finitely many ``characteristic patterns'' for each $\Musati{\delta=k}$, is discussed in Subsection \ref{sec:concclassmu}. An important special case is Conjecture \ref{con:finhit} (now a fully precise statement), about the classification of unsatisfiable hitting clause-sets (or ``disjoint/orthogonal tautologies'' in the terminology of DNFs). In Lemma \ref{lem:conjimplcomp} we show how two of the conjectures together yield computability of $\minnonmer(k)$.

This \Schrift{} is a substantial extension of the conference paper \cite{KullmannZhao2011Bounds}: Section 3 there has been extended to Section \ref{sec:nonmer} here, with considerable more details and examples on non-Mersenne numbers. Section 4 there is covered by Sections \ref{sec:specialcaseMU}, \ref{sec:leansurp} and \ref{sec:algoimpl}, with various additional results (for example showing sharpness of the upper bound for $\Lean$). And the results for Section 5 there are contained in Subsection \ref{sec:genboundml} here. All other sections in this \Schrift{} are new.

Before starting with the work, a few words on dependencies between the sections: It is possible to jump directly to Section \ref{sec:specialcaseMU}, the degree upper-bound for $\Musat$, and only to look up the results of previous sections if needed. Sections \ref{sec:leansurp}, \ref{sec:algoimpl}, \ref{sec:genboundml} on the generalisation to $\Lean$ might also be considered independently, or skipped if only interested in $\Musat$ (which is taken up again in Sections \ref{sec:strengtheningbound}, \ref{sec:improveupbdMU}, \ref{sec:strbMU}). The earlier Section \ref{sec:nonmer} on $\nonmer(k)$ is free-standing. Section \ref{sec:vmusat} on $\Vmusat$ can also be read independently. In general reading of Sections \ref{sec:prelim}, \ref{sec:MUprelim} on clause-sets and on minimal unsatisfiability is recommended to firmly establish the basic notions. The two sections on reductions, Section \ref{sec:elimcreatesing} on singular DP-reduction and Section \ref{sec:2subr} on full subsumption resolution, can again be considered independently, or only looked up if needed.

Finally, in \ref{sec:appOverview} an overview is given on all notations of this \Schrift.

\section{Preliminaries}
\label{sec:prelim}

We follow the general notations and definitions as outlined in \cite{Kullmann2007HandbuchMU}, where also further background on autarkies and minimal unsatisfiability can be found. We use $\NN = \set{n \in \ZZ : n \ge 1}$ and $\NNZ = \NN \cup \set{0}$. For the \href{https://en.wikipedia.org/wiki/Binary_logarithm}{binary logarithm} we use $\ld(x) := \log_2(x) \in \RR$ (``logarithm dualis'') for $x \in \RR_{>0}$, and $\fld(x) := \floor{\ld(x)} \in \ZZ$ for the integral part (e.g., $\fld(1) = 0$, $\fld(2) = \fld(3) = 1$, $\fld(4) = \fld(5) = 2$; this is sequence \url{http://oeis.org/A000523} in the ``On-Line Encyclopedia of Integer Sequences'' (\cite{Sloane2008OEIS})).

We apply standard set-theoretic concepts, like that of a map as a set of pairs, and standard set-theoretic notations, like $f(S) = \set{f(x) : x \in S}$ for maps $f$ and $S \sse \dom(f)$, and ``$\subset$'' for the strict subset-relation. We use also the less-known notation ``$A \addcup B$'' for union in case $A, B$ are disjoint, that is,  $A \addcup B := A \cup B$ is only defined for $A \cap B = \es$. For maps $f, g$ with the same domain $X$ we use $f \le g :\Lra \fa\, x \in X : f(x) \le g(x)$ (i.e., pointwise comparison), while $f < g :\Lra \fa\, x \in X : f(x) < g(x)$.

\subsection{Clause-sets}
\label{sec:prelimcls}

The basic structure is a set $\Lit$, the elements called ``literals'', together with a fixed-point free involution called ``complementation'', written $x \in \Lit \mapsto \ol{x} \in \Lit$; so the laws are $\ol{x} \ne x$ and $\ol{\ol{x}} = x$, which imply $x \ne y \Ra \ol{x} \ne \ol{y}$, for $x, y \in \Lit$. We assume $\ZZ \sm \set{0} \sse \Lit$, with $\ol{x} = -x$ for $x \in \ZZ \sm \set{0}$. For a set $L$ of literals we define $\ol{L} := \set{\ol{x} : x \in L}$. Furthermore a set $\NN \sse \Va \subset \Lit$, the elements called ``variables'', is given, with $\Lit = \Va \addcup \ol{\Va}$. Variables are also called ``positive literals'', while literals $\ol{v}$ for $v \in \Va$ are called ``negative literals''. The ``underlying variable'' of a literal is given by the operation $\var: \Lit \ra \Va$ (``forgetting complementation''), with $\var(v) := v$ and $\var(\ol{v}) := v$ for $v \in \Va$.

\begin{examp}\label{exp:lits}
  We can thus write e.g.\ $1, 6$ for two (different) variables, and $1, 5, -1$ for three (different) literals. In examples we will also use $v, w, a, b, c$ and such letters for variables (as it is customary), and accordingly $\ol{v}$ etc.\ for literals, and in this context (only) it is then understood that these variables are pairwise different. So $\set{v,w,x, \ol{x}}$, when given in an example (without further specification), denotes a set of literals with $\abs{\set{v,w,x, \ol{x}}} = 4$ and $\abs{\set{v,w,x, \ol{x}} \cap \Va} = 3$.

  Without restriction we could assume $\Lit = \ZZ \sm \set{0}$ (as we did in the Introduction), but it is often convenient to use arbitrary mathematical objects as variables. All our objects built from literals are finite, and thus, because of the infinite supply of variables, there will always be ``new variables'' (that's the mathematical point of having natural numbers as variables --- we won't use the arithmetical structure).
\end{examp}

A \textbf{clause} $C$ is a finite and clash-free set of literals (i.e., $C \cap \ol{C} = \es$), the set of all clauses is $\Cl$. A \textbf{clause-set} is a finite set of clauses, the set of all clause-sets is $\Cls$. The simplest clause is the \emph{empty clause} $\bmm{\bot} := \es \in \Cl$, the simplest clause-set is the \emph{empty clause-set} $\top := \es \in \Cls$. The set of all \textbf{hitting clause-sets} is denoted by $\Clash \subset \Cls$, those $F \in \Cls$ such that two different clauses $C, D \in F$, $C \ne D$, have at least one clash, i.e., $C \cap \ol{D} \ne \es$. In the language of DNF, hitting clause-sets are known as ``orthogonal'' or ``disjoint'' DNF's; see \cite[Chapter 7]{CramaHammer2011BooleanFunctions}.

\begin{examp}\label{exp:cls}
  We have e.g.\ $\set{1,2,-3} \in \Cl$, while $\set{-1,1} \notin \Cl$. The only clause-set in $\Clash$ containing the empty clause is $\set{\bot} \in \Clash$. An example of a non-hitting clause-set is $\set{\set{1,2},\set{-1,2},\set{3}} \in \Cls \sm \Clash$, where we obtain an element of $\Clash$ if we add literal $-2$ to the third clause.
\end{examp}

We use $\var(F) := \bc_{C \in F} \var(C) \subset \Va$ for the set of variables of $F \in \Cls$, where $\var(C) := \set{\var(x) : x \in C} \subset \Va$ is the set of variables of clause $C \in \Cl$. The possible literals for a clause-set $F$ are given by $\lit(F) := \var(F) \addcup \ol{\var(F)} \subset \Lit$, while the actually occurring literals are just given by $\bc F \subset \Lit$ (the union of the clauses of $F$). A literal $x$ is \textbf{pure} for $F$ if $\ol{x} \notin \bc F$. For a clause-set $F$ we use the following measurements:
\begin{itemize}
\item $n(F) := \abs{\var(F)} \in \NNZ$ is the number of variables,
\item $c(F) := \abs{F} \in \NNZ$ is the number of clauses,
\item $\delta(F) := c(F) - n(F) \in \ZZ$ is the \textbf{deficiency} (the difference of the number of clauses and the number of variables),
\item $\ell(F) := \sum_{C \in F} \abs{C} \in \NNZ$ is the number of literal occurrences.
\end{itemize}

We call a clause $C$ \textbf{full} for a clause-set $F$ if $C \in F$ and $\var(C) = \var(F)$, while a clause-set $F$ is called full if every clause is full. For a finite set $V$ of variables let
\begin{displaymath}
  \bmm{A(V)} := \set{C \in \Cl : \var(C) = V} \in \Cls.
\end{displaymath}
Obviously $A(V) \in \Clash$ is the set of all $2^{\abs{V}}$ full clauses over $V$, and $F \in \Cls$ is full iff $F \sse A(\var(F))$. We use $\bmm{A_n} := A(\tb 1n)$ for $n \in \NNZ$. Dually, a variable $v \in \Va$ is called \textbf{full} for a clause-set $F$ if for all $C \in F$ holds $v \in \var(C)$. A clause-set is full iff every $v \in \var(F)$ is full.

\begin{examp}\label{exp:measurecls}
  For $F := \set{\bot, \set{1},\set{-1,2}}$ we have:
  \begin{enumerate}
  \item $\var(F) = \set{1,2}$, $\lit(F) = \set{-1,1,-2,2}$, $\bc F = \set{-1,1,2}$.
  \item Literal $2$ is pure for $F$ (the other literals in $\lit(F)$ are not pure).
  \item $n(F) = 2$, $c(F) = 3$, $\delta(F) = 1$, $\ell(F) = 3$.
  \item $\set{-1,2}$ is a full clause of $F$, while the two other clauses are not full.
  \item $F$ has no full variable, while $F \sm \set{\bot}$ has the (single) full variable $1$.
  \end{enumerate}
  The standard ``complete'' full clause-sets are $A_0 = \set{\bot}$, $A_1 = \set{\set{-1},\set{1}}$, and $A_2 = \set{\set{-1,-2},\set{-1,2},\set{1,-2},\set{1,2}}$, and so on.
\end{examp}

We often define a class of clause-sets via some measure $\mu$ as follows:
\begin{defi}\label{def:classescls}
  Consider a class $\mc{C} \sse \Cls$ and a measure $\mu: \Cls \ra \RR$. For $a \in \RR$ we use $\bmm{\mc{C}_{\mu=a}} := \set{F \in \mc{C} : \mu(F) = a}$, and similarly we use \bmm{\mc{C}_{\mu<a}} etc.
\end{defi}
When we use the form ``$\mc{C}_{\mu \Box a}$'', then $\mu$ stands for a measure (e.g., $\mu = \delta$ or $\mu = n$).
\begin{examp}\label{exp:classescls}
  $\Cls_{n=0} = \Cls_{\ell=0} = \set{\top, \set{\bot}}$, $\Cls_{c=0} = \set{\top}$, and $\Cls_{n<0} = \es$.
\end{examp}

\subsection{Semantics}
\label{sec:prelimsem}

A \textbf{partial assignment} is a map $\vp: V \ra \set{0,1}$ for some finite (possibly empty) set $V \subset \Va$ of variables, where $\var(\vp) := V$ and $\lit(\vp) := \lit(\var(\vp)) = \var(\vp) \addcup \ol{\var(\vp)}$. The set of all partial assignments is denoted by $\Pass$. For a literal $x \in \lit(\vp)$ we also define $\vp(x) \in \set{0,1}$, via $\vp(\ol{v}) := 1- \vp(v)$ for $v \in \var(\vp)$. Via a small abuse of language we define $\vp^{-1}(\ve) := \set{x \in \lit(\vp) : \vp(x) = \ve} \in \Cl$ for $\ve \in \set{0,1}$. Special partial assignments are the \emph{empty partial assignment} $\epa := \es$, and for literals $x \in \Lit$ and $\ve \in \set{0,1}$ the partial assignment $\bmm{\pao{x}{\ve}} \in \Pass$, with $\var(\pao{x}{\ve}) = \set{\var(x)}$ and $\pao{x}{\ve}(x) = \ve$.

 The application of a partial assignment $\vp \in \Pass$ to a clause-set $F$ is denoted by $\bmm{\vp * F}$, which yields the clause-set obtained from $F$ by removing all satisfied clauses (which have at least one literal set to $1$), and removing all falsified literals from the remaining clauses:
 \begin{displaymath}
   \vp * F := \set{C \sm \vp^{-1}(0) : C \in F \und C \cap \vp^{-1}(1) = \es} \in \Cls.
 \end{displaymath}
This definition is motivated by the default interpretation of a clause-set as a ``conjunctive normal form'' (CNF), where a clause is understood as a disjunction of literals (thus is satisfied iff at least one literal is satisfied), while a clause-set is understood as a conjunction of its clauses (thus is satisfied iff all clauses are satisfied). A clause-set $F$ is \emph{satisfiable} iff there is a partial assignment $\vp$ with $\vp * F = \top$, otherwise $F$ is \emph{unsatisfiable}. The set of satisfiable clause-sets is denoted by $\bmm{\Sat} \subset \Cls$, while $\bmm{\Usat} := \Cls \sm \Sat$ denotes the set of all unsatisfiable clause-sets.

\begin{examp}\label{exp:sat}
  If $F \in \Usat$ and for $F' \in \Cls$ holds $F \sse F'$, then also $F' \in \Usat$ (satisfying a clause-sets gets harder the more clauses there are).

  By definition we have $\vp * F = \top$ iff $\fa\, D \in F : \vp^{-1}(1) \cap D \ne \es$; thus $F \in \Sat$ iff there is a clause $C \in \Cl$ with $C \cap D \ne \es$ for all $D \in F$. (We could write ``$C \cap D \ne \bot$'' here, but it appears somewhat more natural to use ``$\es$'' here.)
\end{examp}

 The unsatisfiable hitting clause-sets are denoted by $\bmm{\Uclash} := \Usat \cap \Clash$.

\begin{examp}\label{exp:hit}
  $\top \in \Sat \cap \Clash$ and $\set{\bot} \in \Uclash$. In general a full clause-set $F$ is unsatisfiable iff $F = A(\var(F))$, and thus $A(V) \in \Uclash$ for all finite $V \subset \Va$.

  The fundamental property for $F \in \Clash$ is that two different clauses do not have a common falsifying assignment. More precisely, consider $\vp, \psi \in \Pass$, such that there are $C, D \in F$, $C \ne D$, with $\vp * \set{C} = \psi * \set{D} = \set{\bot}$ (that is, $\bot \in \vp * F \cap \psi * F$, where there are different falsified clauses for these two partial assignments). Then $\vp, \psi$ are incompatible, i.e., there is $v \in \var(\vp) \cap \var(\psi)$ with $\vp(v) \ne \psi(v)$.

  It follows easily that for $F \in \Clash$ holds $F \in \Usat \Lra \sum_{C \in F} 2^{-\abs{C}} = 1$.

   A nice exercise is to show $\Uclashi{\delta \le 0} = \es$ (in Section \ref{sec:prelimAut} a more general result is stated).
\end{examp}

Finally, the \emph{semantical implication} $F \models C$ for $F \in \Cls$ and clauses $C \in \Cl$ holds iff $\fa\, \vp \in \Pass : \vp * F = \top \Ra \vp * \set{C} = \top$. We have $F \in \Usat \Lra F \models \bot$.

\subsection{Resolution}
\label{sec:prelimResolution}

Two clauses $C, D \in \Cl$ are \textbf{resolvable} if $\abs{C \cap \ol{D}} = 1$, i.e., they clash in exactly one variable (called the resolution variable $\var(x)$, while $x$ is called the resolution literal). For two resolvable clauses $C$ and $D$, the \textbf{resolvent} $C \res D := (C \cup D) \sm \set{x,\ol{x}} \in \Cl$ for $C \cap \ol{D} = \set{x}$ is the union of the two clauses minus the resolution literal and its complement. As it is well-known (the earliest source is \cite{B37,Ki38}), a clause-set $F$ is unsatisfiable iff via resolution (i.e., closing $F$ under addition of resolvents) we can derive $\bot$, and, more generally, we have $F \models C$ iff from $F$ via resolution a clause $C' \sse C$ is derivable.

An important reduction for clause-sets $F \in \Cls$ and variables $v \in \Va$, resulting in a clause-set satisfiability-equivalent to $F$ (satisfiable iff $F$ is; sometimes called ``equi-satisfiable'') and with variable $v$ eliminated, is \textbf{DP-reduction}
\begin{displaymath}
\bmm{\dpi{v}(F)} := \set{C \in F : v \notin \var(C)} \cup \set{C \res D : C, D \in F \und C \cap \ol{D} = \set{v}} \in \Cls
\end{displaymath}
(also called ``variable elimination''), obtained from $F$ by removing all clauses containing variable $v$ from $F$, and replacing them by their resolvents on $v$. See \cite{KullmannZhao2012ConfluenceJ} for a fundamental study of DP-reduction. The satisfying assignments $\vp$ of $\dpi{v}(F)$ (i.e., $\vp * \dpi{v}(F) = \top$) with $\var(\vp) = \var(F) \sm \set{v}$ are precisely the satisfying assignments $\vp$ of $F$ with $\var(\vp) = \var(F)$, when restricting $\vp$ to $\var(F) \sm \set{v}$. Logically, $\dpi{v}(F)$ is equivalent to $\ex\, v : F$, the existential quantification of $v$ for $F$ (but we do not use quantifiers in this \Schrift{}, so this remark might be ignored here).

\subsection{Multi-clause-sets}
\label{sec:prelimmulti}

These notions are generalised to \textbf{multi-clause-sets}, which are maps $F: \Cl \ra \NNZ$, such that the underlying set of clauses $\set{C \in \Cl : F(C) \ne 0}$ is finite, and so we speak of the \textbf{underlying clause-set}; the set of all multi-clause-sets is denoted by $\ul{\Cls} := \set{F: \Cl \ra \NNZ \mb \Cl \sm F^{-1}(0) \text{ is finite}}$ (in earlier papers we used ``$\Mcls$'' instead of ``$\ul{\Cls}$''). Clause-sets are implicitly promoted to multi-clause-sets, if needed, by using their characteristic functions, and multi-clause-sets are implicitly cast down, if needed, to clause-sets by considering the underlying clause-set; ``if needed'' refers to operations which either require multi-clause-sets or clause-sets. If however we want to make explicit these operations, we use $\ulcls: \ul{\Cls} \ra \Cls$ (with $\ulcls(F) := \Cl \sm F^{-1}(0)$) and $\procls: \Cls \ra \ul{\Cls}$ (with $\procls(F)(C) := 1$ if $C \in F$, and $\procls(F)(C) := 0$ otherwise). For $F \in \ul{\Cls}$ we extend the basic operations:
\begin{itemize}
\item $\var(F) := \var(\ulcls(F))$, $\lit(F) := \lit(\ulcls(F))$, $\bc F := \bc \ulcls(F)$.
\item $n(F) := n(\ulcls(F)) \in \NNZ$, $c(F) := \sum_{C \in \Cl} F(C) \in \NNZ$, $\delta(F) := c(F) - n(F) \in \ZZ$, $\ell(F) := \sum_{C \in F} F(C) \cdot \abs{C} \in \NNZ$.
\end{itemize}
The application of partial assignments $\vp \in \Pass$ to a multi-clause-set $F \in \ul{\Cls}$ yields a \emph{multi-}clause-set $\vp * F \in \ul{\Cls}$, where the multiplicity of a clause $C \in \Cl$ in $\vp * F$ is the sum of all multiplicities of clauses $D \in F$ (i.e., $D \in \ulcls(F)$) which are shortened to $C$ by $\vp$:
\begin{displaymath}
  (\vp * F)(C) := \sum_{D \in F,\, D \cap \vp^{-1}(1) = \es,\, D \sm \vp^{-1}(0) = C} F(D).
\end{displaymath}

\begin{examp}\label{exp:apppamcls}
  If $\vp$ is a total assignment for $F$ (assigns all variables of $F$, that is, $\var(\vp) = \var(F)$), then $\vp * F$ is $\set{m*\bot}$, denoting the multiplicity of a clause by a (formal) factor, with $m = \sum_{C \in F,\, C \cap \vp^{-1}(1) = \es} F(C) \in \NNZ$ (so $m = 0 \Lra \vp * F = \top$).
\end{examp}

For us, typically clause-sets are the objects of interests, while multi-clause-sets are mostly auxiliary devices, created by the operation of ``restriction'' defined below (Definition \ref{def:restriction}). However we have to take care of the details, and thus together with introducing a class $\mc{C} \sse \Cls$ we also introduce the corresponding class $\ul{\mc{C}} \sse \ul{\Cls}$ of multi-clause-sets, using the generalised definition of $\mc{C}$, and where $\mc{C} = \set{\ulcls(F) : F \in \ul{\mc{C}}} = \ulcls(\ul{\mc{C}})$ holds for the following two main forms of relations between $\mc{C}$ and $\ul{\mc{C}}$:
\begin{itemize}
\item For example, the classes $\ul{\Sat}$ and $\ul{\Usat}$ are \textbf{invariant under multiplicities}, that is, a multi-clause-set is in it iff the underlying clause-set is in the underlying class of clause-sets ($\Sat$ resp.\ $\Usat$). In general $\ul{\mc{C}}$ is invariant under multiplicities iff for all $F, F' \in \ul{\Cls}$ with $\ulcls(F) = \ulcls(F')$ either $F, F' \in \ul{\mc{C}}$ or $F, F' \notin \ul{\mc{C}}$ holds.
\item The other extreme we have with the class $\ul{\Clash}$ of multi-hitting-clause-sets, which \textbf{disallows multiplicities}, that is, all multiplicities must be $1$ (since clauses can not clash with themselves, by definition of clauses), and thus up to the canonical identification, the classes $\Clash$ and $\ul{\Clash}$ are identical. In general $\ul{\mc{C}}$ disallows multiplicities iff for all $F \in \ul{\mc{C}}$ holds $F(\Cl) \sse \set{0,1}$.
\end{itemize}
In these cases, the class $\ul{\mc{C}}$ does not carry more information than the class $\mc{C}$. However in general we only have $\mc{C} = \ul{\mc{C}} \cap \Cls$ (using the implicit conversions between clause-sets and multi-clause-sets), and the class $\ul{\mc{C}}$ can not be derived from $\mc{C}$. For classes $\ul{\mc{C}}$ sensitive to multiplicities, in this \Schrift{} we have examples for:
\begin{itemize}
\item \textbf{downward closure for multiplicities}, i.e., reducing a multiplicity does not leave the class: $\ul{\Msat}$ (``matching satisfiable'', Subsection \ref{sec:prelimAut}) and $\ul{\Sed}$ (``surplus equal deficiency'', Definition \ref{def:sed});
\item \textbf{upward closure for multiplicities} (increasing any non-zero multiplicity does not leave the class): $\ul{\Mlean}$ (``matching lean'', Subsection \ref{sec:prelimAut});
\item neither: $\ul{\Mlcr}$ (``matching lean critical'', Definition \ref{def:Mlcr}).
\end{itemize}

Clause-sets $F, G$ are called \textbf{isomorphic}, if the variables of $F$ can be renamed and potentially flipped so that $F$ is turned into $G$. More precisely, an isomorphism $\alpha$ from $F$ to $G$ is a bijection $\alpha: \lit(F) \ra \lit(G)$ which preservers complementation ($\alpha(\ol{x}) = \ol{\alpha(x)}$), and which maps the clauses of $F$ precisely to the clauses of $G$; when considering multi-clause-sets, then the isomorphism must preserve the multiplicity of clauses (that is, $G(\alpha(C)) = F(C)$ for all $C \in \Cl$). All classes of (multi-)clause-sets we consider in this \Schrift{} are closed under isomorphisms (besides some special cases only mentioned in the Introduction).

\subsection{Restrictions}
\label{sec:prelimrestr}

An important operation with multi-clause-set is the ``restriction'' to a set of variables (see \cite[Subsection 3.5]{Kullmann2007ClausalFormZI} for more information):
\begin{defi}\label{def:restriction}
  For a set $V \sse \Va$ of variables and a multi-clause-set $F \in \ul{\Cls}$ by $\bmm{F[V]} \in \ul{\Cls}$ the \textbf{restriction} of $F$ to $V$ is denoted, which is the multi-clause-set obtained by removing clauses from $F$ which have no variables in common with $V$, and removing from the remaining clauses all literals where the underlying variable is not in $V$; so  $F[V](\bot) := 0$, while for $C \in \Cl \sm \set{\bot}$
  \begin{displaymath}
    F[V](C) := \sum_{D \in F,\, D \cap \lit(V) = C} F(D).
  \end{displaymath}
\end{defi}
It is essential that $F[V]$ is a multi-clause-set (when considering classes of multi-clause-sets sensitive to multiplicities), even when $F$ is just a clause-set, and if previously unequal clauses become equal, then accordingly their multiplicity is increased.

\begin{examp}\label{exp:restr}
  $\set{\set{a},\set{a,b},\set{b},\set{\ol{a},\ol{b}}}[\set{a}] = \set{2*\set{a}, \set{\ol{a}}}$.
\end{examp}

Simple properties of this operation are (for multi-clause-sets $F$ and $V, V' \sse \Va$):
\begin{enumerate}
\item $F[\es] = \top$, $F[V] = F \sm \set{\bot}$ for $\var(F) \sse V$ (where $F \sm F'$ for a clause-set $F'$ means that all occurrences of clauses from $F'$ are removed from $F$).
\item $(F[V])[V'] = F[V \cap V']$.
\item $c(F[V]) = \sum_{C \in F, \var(C) \cap V \ne \es} F(C)$ (the number of clauses containing some variable from $V$).
\end{enumerate}

\subsection{Degrees}
\label{sec:prelimDegrees}

For the number of occurrences of a literal $x \in \Lit$ in a (multi-)clause-set $F \in \ul{\Cls}$ we write
\begin{displaymath}
  \bmm{\ldeg_F(x)} := \sum_{C \in F,\, x \in C} F(C) \in \NNZ,
\end{displaymath}
called the \textbf{literal-degree}; note $\ldeg_F(x) = 0 \Lra x \notin \bc F$. The \textbf{variable-degree} of a variable $v$ is defined as $\bmm{\vdeg_F(v)} := \ldeg_F(v) + \ldeg_F(\ol{v}) \in \NNZ$; note $\vdeg_F(v) = 0 \Lra v \notin \var(F)$ and $\vdeg_F(v) = c(F[\set{v}])$. We remark that the literal-degree $\ldeg_F(x)$ (for $x \in \Lit$) always has the index $F$, the (multi-)clause-set in which the literal-degree is counted, and thus there is little danger of confusion with the binary logarithm $\ld(x)$ (for $x \in \RR_{>0}$), which does not have an index. A (multi-)clause-set $F$ is called \textbf{variable-regular} if all variables $v \in \var(F)$ have the same degree, or, stronger, \textbf{literal-regular}, if all literals $x \in \lit(F)$ have the same degree. A \textbf{singular variable} in a (multi-)clause-set $F$ is a non-pure variable occurring in one sign only once, that is, $\min(\ldeg_F(v), \ldeg_F(\ol{v})) = 1$, while $F$ is called \textbf{nonsingular} if it does not have singular variables, i.e., iff $\fa\, x \in \bc F : \ldeg_F(\ol{x}) \ne 0 \Ra \ldeg_F(x) \ge 2$. The central concept for this \Schrift{} is the minimum degree of a variable in a clause-set:
\begin{defi}\label{def:minvdeg}
  We define the \textbf{minimum variable-degree} $\bmm{\minvdeg} : \ul{\Cls} \ra \nnpi$ (``min-var-degree'' for short) as follows: For $F \in \ul{\Cls}$ with $n(F) \ne 0$ we let $\minvdeg(F) := \min_{v \in \var(F)} \vdeg_F(v) \in \NN$, while for $n(F) = 0$ we set $\minvdeg(F) := +\infty$.

  For a class $\mc{C} \sse \ul{\Cls}$ of (multi-)clause-sets let $\bmm{\minvdeg(\mc{C})} \in \nni$ be the supremum of $\minvdeg(F)$ for $F \in \mc{C}$ with $n(F) > 0$, where we set $\minvdeg(\mc{C}) := 0$ if there is no such $F$ (while otherwise we have $\minvdeg(\mc{C}) \ge 1$).
\end{defi}
We have $\minvdeg(F) = +\infty$ iff $n(F) = 0$ (otherwise $\minvdeg(F) \in \NN$), and the motivation for this setting is, that in this way for all $F \in \ul{\Cls}$ and $K \in \RR$ holds $\minvdeg(F) \ge K$ iff $\fa\, v \in \var(F) : \vdeg_F(v) \ge K$. By definition we have $\minvdeg(\mc{C}) \le \minvdeg(\mc{C}')$ for $\mc{C} \sse \mc{C}' \sse \ul{\Cls}$, and furthermore for $K \in \RR_{>0}$ we have $\minvdeg(\mc{C}) \ge K$ iff there is $F \in \mc{C}$ with $n(F) > 0$ and $\fa\, v \in \var(F) : \vdeg_F(v) \ge K$, while for $K \in \RR_{\ge 0}$ we have $\minvdeg(\mc{C}) \le K$ iff for all $F \in \mc{C}$ with $n(F) > 0$ there is $v \in \var(F)$ with $\vdeg_F(v) \le K$.

\begin{examp}\label{exp:vardegree}
  For $F := \set{2*\set{a,b},\set{\ol{a},b},\set{\ol{b},c}} \in \ul{\Cls}$ we have:
  \begin{itemize}
  \item $\ldeg_F(a) = 2$, $\ldeg_F(\ol{a}) = 1$, $\ldeg_F(b) = 3$, $\ldeg_F(\ol{b}) = 1$, $\ldeg_F(c) = 1$, $\ldeg_F(\ol{c}) = 0$.
  \item $\vdeg_F(a) = 3$, $\vdeg_F(b) = 4$, $\vdeg_F(c) = 1$, thus $\minvdeg(F) = 1$.
  \end{itemize}
  Every full clause-set is variable-regular (but in general not literal-regular). Examples for $\minvdeg(\mc{C})$ are $\minvdeg(\es) = 0$, $\minvdeg(\Cls) = +\infty$, $\minvdeg(\set{\top,\set{\bot},\set{\set{v},\set{\ol{v}}}, F}) = 2$.
\end{examp}

The simplest but relevant class of clause-sets for us is given by the $A(V)$, the unsatisfiable full clause-sets; these are the simplest unsatisfiable clause-sets:
\begin{lem}\label{lem:defAn}
  For $n \in \NNZ$ we have
  \begin{enumerate}
  \item $n(A_n) = n$, $c(A_n) = 2^n$, $\delta(A_n) = 2^n - n$.
  \item $A_n$ is full and unsatisfiable, and thus $A_n \in \Uclashi{\delta=2^n-n}$.
  \item $A_n$ is literal-regular (thus variable-regular).
  \item $\minvdeg(A_n) = 2^n$ for $n \ge 1$.
  \end{enumerate}
\end{lem}
Further properties of unsatisfiable full clause-sets one finds in Example \ref{exp:surp}, Lemma \ref{lem:charakmarsat}, Corollary \ref{cor:charakmarsat}, Lemmas \ref{lem:inv2subrescan}, \ref{lem:inv2satusubrescan}, Corollaries \ref{cor:upbdefn}, \ref{cor:maxdefmun}, and Examples \ref{exp:sed}, \ref{exp:sharpbmvd}. Properties of satisfiable full clause-sets are found in Example \ref{exp:mlcr}.

\subsection{Autarkies}
\label{sec:prelimAut}

Besides algorithmic considerations, which were present since the introduction of the notion of an ``autarky'' in \cite{MoSp85}, also a kind of a ``combinatorial SAT theory'' has been developed around this notion of generalised satisfying assignments. A general overview is given in \cite{Kullmann2007HandbuchMU}, with recent additions and generalisations in \cite{Kullmann2007ClausalFormZI}. An \textbf{autarky} (see \cite[Section 11.8]{Kullmann2007HandbuchMU}) for a clause-set $F \in \Cls$ is a partial assignment $\vp \in \Pass$ which satisfies every clause $C \in F$ it touches, i.e., for all $C \in F$ with $\var(\vp) \cap \var(C) \not= \es$ holds $\vp * \set{C} = \top$; equivalently, for all $C \in F$ holds $C \cap \vp^{-1}(0) \ne \es \Ra C \cap \vp^{-1}(1) \ne \es$. The simplest examples for autarkies are:

\begin{examp}\label{exp:autarkies}
  The empty partial assignment $\epa$ is an autarky for every $F \in \Cls$ (no clause is touched), and more generally all $\vp \in \Pass$ with $\var(\vp) \cap \var(F) = \es$ are autarkies for $F$, the \textbf{trivial autarkies}. On the other end of the spectrum every satisfying assignment for $F$ (i.e., $\vp * F = \top$) is an autarky for $F$ (every clause is satisfied). A literal $x \in \Lit$ is a pure literal for $F$ iff $\pao x1$ is an autarky for $F$.
\end{examp}

If $\vp$ is an autarky for $F$, then $\vp * F \sse F$ holds, and thus $\vp * F$ is satisfiability-equivalent to $F$. Autarkies mark redundancies, and the corresponding notion of clause-sets without such redundancies was introduced in \cite{Ku98e}, namely a clause-set $F$ is \textbf{lean} if there is no non-trivial autarky for $F$, and the set of all lean clause-sets is denoted by $\Lean \subset \Usat \addcup \set{\top}$. The class $\ul{\Lean}$ of lean multi-clause-sets is invariant under multiplicities.
\begin{examp}\label{exp:LEAN}
  Some simple examples:
  \begin{enumerate}
  \item $\top, \set{\bot}, \set{\set{v},\set{\ol{v}}}, \set{\set{v},\set{\ol{v}},\set{w},\set{\ol{w}}} \in \Lean$.
  \item If $F, F' \in \Lean$, then $F \cup F' \in \Lean$.
  \item If $F \in \Lean$ and $F' \in \Cls$ with $\var(F') \sse \var(F)$, then $F \cup F' \in \Lean$.
  \item $\set{\set{v},\set{\ol{v}},\set{w}} \notin \Lean$.
  \end{enumerate}
\end{examp}
A weakening is the notion of a \textbf{matching-lean} multi-clause-set $F$ (introduced in \cite[Section 7]{Ku00f}; see \cite[Section 11.11]{Kullmann2007HandbuchMU} for an overview), which has no non-trivial \textbf{matching autarky}, special autarkies given by a matching condition: for every clause touched (taking, of course, multiplicities into account), a  satisfied literal with unique underlying variable must be selectable; the class of all matching-lean multi-clause-sets is denoted by $\ul{\Mlean} \supset \ul{\Lean}$, while the class of all multi-clause-sets satisfiable by some matching autarky is denoted by $\ul{\Msat} \subset \ul{\Sat}$.

\begin{examp}\label{exp:MLEAN}
  $\ul{\Mlean} \cap \ul{\Msat} = \set{\top}$. $F := \set{\set{1,3},\set{2,-3},\set{3},\set{-3}}$ has the matching autarky $\pab{1 \ra 1, 2 \ra 1}$ (but is not satisfiable), while for $F' := F \cup \set{\set{1,2}}$ we have $F' \in \Mlean$. Note $\delta(F) = 1 = \delta(\set{\set{3},\set{-3}})$, while $\delta(F') = 2$.
\end{examp}

The class $\ul{\Mlean}$ is upward closed for multiplicities, $\ul{\Msat}$ is downward closed for multiplicities. A multi-clause-set $F \in \ul{\Cls}$ is matching-lean iff for all $F' \le F$, $F' \ne F$, holds $\delta(F') < \delta(F)$ (\cite[Theorem 7.5]{Ku00f}). Thus for every matching-lean multi-clause-set $F \not= \top$ we have $\delta(F) \ge 1$ (\cite{Ku00f}, generalising \cite{AhLi86}). It is decidable in polynomial time whether $F \in \ul{\Mlean}$ holds (which follows for example by the characterisation of $\ul{\Mlean}$ via the surplus below).

\begin{examp}\label{exp:ulMLEAN}
  $\set{\set{v}} \in \Msat$, but $\set{2*\set{v}} \in \ul{\Mlean}$, and more generally $\set{\set{v_1,\dots,v_n}} \in \Msat$ for $n \ge 1$, while $\set{(n+1) * \set{v_1,\dots,v_n}} \in \ul{\Mlean}$. Indeed it is easy to see that for every $F \in \Cls$ there is $F' \in \ul{\Cls}$ with $\ulcls(F') = F$ and $F' \in \ul{\Mlean}$. So matching autarkies can be destroyed by increasing multiplicities --- they stay \emph{autarkies}, but the matching criteria is made to fail.
\end{examp}

The process of applying autarkies as long as possible to a clause-set $F \in \Cls$ is confluent, yielding the \textbf{lean kernel} of $F$, the largest lean sub-clause-set of $F$, that is, $\bc \set{F' \sse F : F' \in \Lean}$; see \cite[Section 3]{Ku98e}. Computation of the lean kernel is NP-hard, since the lean kernel of satisfiable clause-sets is $\top$. But the \textbf{matching-lean kernel} of $F$, the largest matching-lean sub-clause-set of $F$ (that is, $\bc \set{F' \sse F : F' \in \Mlean}$; see \cite[Section 3]{Ku00f}), now obtained by applying \emph{matching} autarkies as long as possible (again a confluent process), is computable in polynomial time. Note that a clause-set $F$ is lean resp.\ matching lean iff the lean resp.\ matching-lean kernel is $F$ itself. Due to the polytime computability of the matching-lean kernel, a sub-clause-set obtained by removing clauses redundant in a strong sense, ``w.l.o.g.'' for SAT-decision one might consider the inputs as matching-lean. The matching-lean kernel is applied in this \Schrift{} only to clause-sets, the main inputs of the algorithms we consider, possibly after throwing away clause-multiplicities --- dedicated matching considerations only play a role for \emph{auxiliary} multi-clause-sets, created by restriction and considered in a combinatorial way.

\begin{examp}\label{exp:MLEANfpt}
  For inputs $F \in \Mlean$ by \cite[Theorem 4]{Szei2002FixedParam} we have SAT-decision in time $O(2^{\delta(F)} \cdot n(F)^3)$ (see \cite{Kullmann2007ClausalFormZI} for generalisations), and thus SAT-decision for inputs $F \in \Mlean$ is fixed-parameter tractability (fpt) in the parameter $\delta(F)$.

  We note here (though we won't use it in this \Schrift{}), that for inputs $F \in \Mleani{\delta=k}$ the computation of the lean kernel can be done in polynomial time for fixed $k$ (\cite[Theorem 10.3]{Kullmann2007ClausalFormZI}; this computational problem appears not to be fpt).
\end{examp}

\subsection{The surplus}
\label{exp:prelimsurplus}

The maximal generalisation of ``Tarsi's Lemma'' for (boolean) CNF is: A multi-clause-set $F \not= \top$ is matching lean (has no non-trivial matching autarky) iff we have $\surp(F) \ge 1$ for the \textbf{surplus}, the deficiency minimised over all non-empty restrictions (\cite[Lemma 7.7]{Ku00f}). The precise definition is as follows (see \cite[Subsection 11.1]{Kullmann2007ClausalFormZI} for more information; in \cite{Szei2002FixedParam} a clause-set has ``$q$-expansion'' iff $\surp(F) \ge q$):
\begin{defi}\label{def:surp}
  For a multi-clause-set $F$ let $\bmm{\surp(F)} \in \ZZ$ be defined as the minimum of $\delta(F[V])$ (recall Definition \ref{def:restriction}) over all $\es \not= V \sse \var(F)$ if $n(F) > 0$, while $\surp(F) := 0$ in case of $n(F) = 0$.
\end{defi}
So for $F \in \ul{\Cls}$, $n(F) > 0$, and $K \in \RR$ holds $\surp(F) \ge K$ iff for all $\es \ne V \sse \var(F)$ the number of clauses of $F$ containing some variable of $V$ is at least $K + \abs{V}$. The special case of $n(F)=0$ is handled so that we always have $- n(F) \le \surp(F) \le c(F)$. Note that for $\es \not= V \sse \var(F)$ we have
\begin{displaymath}
  \delta(F[V]) = c(F[V]) - \abs{V} = \big(\sum_{C \in F,\, \var(C) \cap V \ne \es} \hspace{-1.5em}F(C) \big) \: - \abs{V}.
\end{displaymath}
The surplus is computable in polynomial time. Some basic properties are:
\begin{enumerate}
\item $\surp(F)$ is independent of $F(\bot)$ (the number of occurrences of the empty clause).
\item $\surp(F) \le \delta(F[\var(F)]) = \delta(F \sm \set{\bot}) \le \delta(F) \le c(F)$.
\item $- n(F) \le \surp(F)$, and for $n(F) > 0$ holds $1 - n(F) \le \surp(F) \le c(F) - 1$.
\item For $v \in \var(F)$: $\surp(F) \le \delta(F[\set{v}]) = \vdeg_F(v) - 1$.
\item Thus $\surp(F) \le \minvdeg(F) - 1$.
\item For every $\es \subset V \sse \var(F)$ holds $\surp(F[V]) \ge \surp(F)$.
\item For $F' \le F$ with $\var(F') = \var(F)$ we have $\surp(F') \le \surp(F)$.
\end{enumerate}

\begin{examp}\label{exp:surp}
  $\surp(A_0) = 0$, while $\surp(A_n) = 2^n - n = \delta(A_n)$ for $n \in \NN$. If we take $F \in \Cls$ and some $v \in \Va \sm \var(F)$, then $\surp(F \addcup \set{\set{v}}) \le 0$.
\end{examp}

The basic intuition is that $\surp(F)$ gives us the ``easiest subinstances'' of $F$, ``easiest'' in the sense of deficiency, ``subinstance'' in the sense of restriction. The theory of the surplus of clause-sets is further developed in this \Schrift{} in Section \ref{sec:leansurp}, especially Subsection \ref{sec:classSED}, and Section \ref{sec:algoimpl}, especially Subsection \ref{sec:remsurp}.

\section{Minimally unsatisfiable clause-sets}
\label{sec:MUprelim}

In this section we review minimally unsatisfiable clause-sets; see \cite{Kullmann2007HandbuchMU} for an overview, while \cite{Kullmann2007ClausalFormZII,KullmannZhao2012ConfluenceJ} contain recent developments. First the basic definitions and examples are given in Subsection \ref{sec:musubcl}. In Subsection \ref{sec:Saturation} we consider in some detail the fundamental process of ``saturation'', which is about adding ``missing literal occurrences'' to minimally unsatisfiable clause-sets. The dual problem of deleting ``superfluous literal occurrences'' (``marginalisation'') is considered in Subsection \ref{sec:marginal}. Saturation repairs the problem that splitting of $F \in \Musat$ into $\pao v0 * F$ and $\pao v1 * F$ may destroy minimal unsatisfiable, i.e., $\pao v0 * F \notin \Musat$ or $\pao v1 * F \notin \Musat$ might hold, due to some clauses missed to be deleted by the partial assignment, and this process of splitting is considered in Subsection \ref{sec:Splitting}.

\subsection{$\Musat$ and subclasses}
\label{sec:musubcl}

An unsatisfiable clause-set $F$ is called \textbf{minimally unsatisfiable}, if for every clause $C \in F$ the clause-set $F \sm \set{C}$ is satisfiable, and the set of minimally unsatisfiable clause-sets is denoted by $\bmm{\Musat} \subset \Usat$. A clause-set $F \in \Musat$ is called \textbf{saturated}, if replacing any $C \in F$ by any super-clause $C' \supset C$ yields a satisfiable clause-set, and the set of saturated minimally unsatisfiable clause-sets is denoted by $\bmm{\Smusat} \subset \Musat$.

\begin{examp}\label{exp:simpleMU}
  The simplest element of $\Usat \sm \Musat$ is $\set{\bot,\set{1}}$, while the simplest element of $\Musat \sm \Smusat$ is $\set{\set{1,2},\set{-1},\set{-2}}$ (see Example \ref{exp:additionlo} for a ``saturation'').
\end{examp}

Unsatisfiable hitting clause-sets fulfil $\Uclash \subset \Smusat$ (see \cite[Lemma 2]{KullmannZhao2012ConfluenceJ} for the proof). The subsets of nonsingular elements (i.e., there is no literal occurring only once) are denoted by $\bmm{\Musatns} \subset \Musat$, $\bmm{\Smusatns} \subset \Smusat$, and $\bmm{\Uclashns} \subset \Uclash$.

\begin{examp}\label{exp:MU1}
  By \cite{DDK98} holds $\Musatnsi{\delta=1} = \Smusatnsi{\delta=1} = \Uclashnsi{\delta=1} = \set{\set{\bot}}$ , while for the characterisation of $\Musati{\delta=1} \supset \Smusati{\delta=1} = \Uclashi{\delta=1}$ see also \cite{AhLi86,Ku99dKo}. As shown in \cite{DDK98}, for $F \in \Musati{\delta=1}$ with $n(F) > 0$ holds $\minvdeg(F) = 2$.
\end{examp}

We consider the ``reasons'' for unsatisfiability as given by the elements of $\Musati{\delta=1}$ as ``noise'', only ``masking'' the pure contradiction of the only element of $\Musatnsi{\delta=1} = \set{\set{\bot}}$ (in Section \ref{sec:elimcreatesing} the elimination of singular variables will be discussed). ``Real reasoning'' starts with deficiency $2$:

\begin{examp}\label{exp:MU2}
 By \cite{KleineBuening2000SubclassesMU}, the elements of $\Musatnsi{\delta=2}$ are up to isomorphism precisely the $\bmm{\Dt{n}} := \set{\set{1,\dots,n},\set{-1,\dots,-n},\set{-1,2}, \dots, \set{-(n-1),n}, \set{-n,1}}$ for $n \ge 2$,

  All $\Dt{n}$ are literal-regular, with $\minvdeg(\Dt{n}) = 4$. It is easy to see that all $\Dt{n}$ are saturated, and thus $\Musatnsi{\delta=2} = \Smusatnsi{\delta=2}$. The only hitting clause-sets amongst the $\Dt{n}$ are for $n=2,3$, and thus up to isomorphism the elements of $\Uclashnsi{\delta=2}$ are $\Dt{2}, \Dt{3}$, with $\Dt{2} = A_2$ and $\Dt{3} = \set{\set{1,2,3},\set{-1,-2,-3},\set{-1,2},\set{-2,3},\set{-3,1}}$.

  We have $\surp(\Dt{n}) = 2 = \delta(\Dt{n})$, since any $m \le n$ variables occur at least in $m$ different binary clauses plus in the two ``long clauses'' (this is also obtained from the later Lemma \ref{lem:charsurp1}, Part \ref{lem:charsurp1a}). Further properties of the $\Dt{n}$ we have in Examples \ref{exp:addlitcluhit}, \ref{exp:s2subsr}, \ref{exp:splitting}, \ref{exp:sed}. See Section 7 in \cite{KullmannZhao2012ConfluenceJ} for more information.
\end{examp}
  As shown in \cite[Theorem 74]{KullmannZhao2012ConfluenceJ}, for every $F \in \Musati{\delta=2}$ there is a unique $n \ge 2$ such that $\Dt{n}$ ``embeds'' into $F$, and this $n$ is called the ``nonsingularity type'' of $F$. So for $\Musati{\delta=2}$ we have identified the (in a sense) unique reason of unsatisfiability, the (possibly hidden) presence of a cycle $v_1 \ra \dots \ra v_n \ra v_1$ together with the assertions, that one $v_i$ must be true and one must be false (only the $n$ is unique in general, not the $v_i$). We will come back to the theme of classifying $\Musati{\delta=k}$ in the Conclusion, Subsection \ref{sec:concclassmu}.

By definition, $\ul{\Musat}$ disallows multiplicities (since a duplicated clause is the trivial logical redundancy), and this also holds for the subclasses $\ul{\Smusat}$ and $\ul{\Uclash}$ (as well as for all other subclasses of $\Musat$ considered here). A fundamental fact is $\delta(F) \ge 1$ for all $F \in \Musat$ (note that every minimally unsatisfiable clause-set is lean), which motivates the investigation of the layers $\Musati{\delta=1}, \Musati{\delta=2}, \dots$. Special elements of $\Uclash$ are the $A(V)$ for finite sets $V$ of variables (recall Lemma \ref{lem:defAn}), which are the minimally unsatisfiable clause-sets with maximal deficiency for a given number of variables, as we will see in Corollary \ref{cor:upbdefn}.

Finally, the two main quantities studied in this \Schrift{} need dedicated names:
\begin{defi}\label{def:nonmerstar}
  For $k \in \NN$ let $\bmm{\minnonmer(k)} := \minvdeg(\Musati{\delta=k}) \in \NN$ and $\bmm{\minnonmerh(k)} := \minvdeg(\Uclashi{\delta=k}) \in \NN$ (note $\minnonmer(k) \le 2 k$ by \cite[Lemma C.2]{Ku99dKo}).
\end{defi}
By definition holds $\minnonmerh \le \minnonmer$.

\subsection{Saturation}
\label{sec:Saturation}

We recall the fact (\cite{FlRe94,KullmannZhao2011Bounds}) that every minimally unsatisfiable clause-set $F \in \Musat$ can be \textbf{saturated}, i.e., by adding literal occurrences to $F$ we can obtain $F' \in \Smusat$ with $\var(F') = \var(F)$ (there is then a bijection $\alpha: F \ra F'$ with $C \sse \alpha(C)$ for all $C \in F$). Since we consider saturation in many situations, we introduce some special notations for it from \cite[Subsection 2.2]{KullmannZhao2012ConfluenceJ}. First we introduce the notation $\saturate(F,C,x)$ for adding a literal $x$ to a clause $C$ in a clause-set $F$:
\begin{defi}[\cite{KullmannZhao2012ConfluenceJ}]\label{def:saturation}
  The operation $\bmm{\saturate(F,C,x)} := (F \sm \set{C}) \cup (C \addcup \set{x}) \in \Cls$ (adding literal $x$ to clause $C$ in $F$) is defined if $F \in \Cls$, $C \in F$, and $x$ is a literal with $\var(x) \in \var(F) \sm \var(C)$.
\end{defi}
Some technical remarks:
\begin{enumerate}
\item $\var(\saturate(F,C,x)) = \var(F)$.
\item If $C \addcup \set{x} \notin F$, then $c(\saturate(F,C,x)) = c(F)$, and thus also $\delta(\saturate(F,C,x)) = \delta(F)$.
\item For $F \in \Musat$ we have:
  \begin{enumerate}
  \item $\saturate(F,C,x) \in \Musat$ iff $\saturate(F,C,x)$ is unsatisfiable (since all what happened is that a clause has been weakened, i.e., extended).
  \item If $\saturate(F,C,x) \in \Musat$, then $c(\saturate(F,C,x)) = c(F)$ (no subsumption here).
  \item $F$ is saturated iff there are no $C, x$ such that $\saturate(F,C,x) \in \Usat$.
  \end{enumerate}
\end{enumerate}

\begin{examp}\label{exp:additionlo}
  For $F := \set{\set{a,b},\set{\ol{a}},\set{\ol{b}}} \in \Musat \sm \Smusat$ we have $\saturate(F,\set{\ol{a}},b) = \set{\set{a,b},\set{\ol{a},b},\set{\ol{b}}} \in \Smusat$.
\end{examp}

A ``saturation'' of a minimally unsatisfiable clause-set is obtained by adding literals to clauses as long as possible while maintaining unsatisfiability (which is the same as maintaining minimal unsatisfiability):
\begin{defi}[\cite{KullmannZhao2012ConfluenceJ}]\label{def:partsaturation}
  A \textbf{saturation} $F' \in \Smusat$ of $F \in \Musat$ is obtained by a \textbf{saturation sequence} $F = F_0, \dots, F_m = F'$, $m \in \NNZ$, such that
  \begin{enumerate}
  \item[(i)] for $0 \le i < m$ there are $C_i, x_i$ with $F_{i+1} = \saturate(F_i,C_i,x_i)$,
  \item[(ii)] for all $1 \le i \le m$ we have $F_i \in \Usat$,
  \item[(iii)] the sequence cannot be extended (without violating conditions (i) or (ii)).
  \end{enumerate}
  Note that $n(F') = n(F)$, $c(F') = c(F)$, $\delta(F') = \delta(F)$, and $\ell(F') = \ell(F_0) + m$ holds. If we drop requirement (iii), then we speak of a \textbf{partial saturation sequence}, while $F' \in \Musat$ is a \textbf{partial saturation} of $F \in \Musat$.
\end{defi}

Obviously every $F \in \Musat$ has a saturation $F' \in \Smusat$. Also by definition follows that $F' \in \Musat$ is a partial saturation of $F$ iff $\var(F') = \var(F)$ and there is a bijection $\alpha: F \ra F'$ such that for all $C \in F$ we have $C \sse \alpha(C)$. And $F'$ is a saturation of $F$ iff $F'$ is a partial saturation of $F$ with $F' \in \Smusat$.

\begin{examp}\label{exp:saturation}
  A saturation sequence for $F := \set{\set{a,b,c},\set{\ol{a}},\set{\ol{b}},\set{\ol{c}}}$ with $m=3$ is obtained by adding literals $b, c$ to clause $\set{\ol{a}}$, and adding literal $c$ to clause $\set{\ol{b}}$.
\end{examp}

We can perform a partial saturation $F \leadsto \saturate(F,C,x)$ iff $F$ without $C$ implies (logically) $C \addcup \set{\ol{x}}$ (note that $C \addcup \set{x}, C \addcup \set{\ol{x}}$ implies $C$):
\begin{lem}\label{lem:basicharacpsat}
  Consider $F \in \Musat$, $C \in F$, and a literal $x$ with $\var(x) \in \var(F) \sm \var(C)$. Then $\saturate(F,C,x)$ is a partial saturation of $F$ if and only if $F \sm \set{C} \models C \addcup \set{\ol{x}}$.
\end{lem}
\begin{prf}
First assume that $\saturate(F,C,x)$ is a partial saturation of $F$, but $F \sm \set{C} \not\models C \addcup \set{\ol{x}}$. So there is a partial assignment $\vp$ with $\vp * (F \sm \set{C}) = \top$ but $\vp * \set{C \addcup \set{\ol{x}}} = \set{\bot}$ (whence $\vp(x) = 1$). But then we have $\vp * \saturate(F,C,x) = \top$. Reversely assume $F \sm \set{C} \models C \addcup \set{\ol{x}}$, but that  $\saturate(F,C,x)$ is not a partial saturation of $F$. So $\saturate(F,C,x)$ has a satisfying assignment $\vp$; due to $F \in \Usat$ we have $\vp(x) = 1$ and $\vp * \set{C} = \set{\bot}$. But this yields $\vp * (F \sm \set{C}) = \top$ and $\vp * \set{C \addcup \set{\ol{x}}} = \set{\bot}$. \Qed
\end{prf}

See Lemma \ref{lem:aux2subsumpR}, Part \ref{lem:aux2subsumpR3ab}, for another characterisation of partial saturations.

\subsection{Marginal MUs}
\label{sec:marginal}

The dual notion of ``saturated'' is ``marginal'': $F \in \Musat$ is \textbf{marginal} iff replacing any clause by a strict subclause yields a clause-set not in $\Musat$. The decision ``$F$ marginal minimally unsatisfiable ?'' for inputs $F \in \Cls$ is $D^P$-complete (\cite[Theorem 2]{KleineBueningZhao2007ComplexitySomeSubclassesMU}). By \cite[Theorem 8]{KleineBueningZhao2003StructureSomeClassesMU} however this decision is easy for inputs $F \in \Smusat$, namely there is the following characterisation of minimally unsatisfiable clause-sets which are marginal and saturated at the same time:
\begin{lem}[\cite{KleineBueningZhao2003StructureSomeClassesMU}]\label{lem:charakmarsat}
  $F \in \Musat$ is both marginal and saturated iff $F = A(\var(F))$.
\end{lem}

We obtain that precisely all saturated clause-sets except the $A(V)$ are obtained as non-trivial saturations of some minimally unsatisfiable clause-set:
\begin{corol}\label{cor:charakmarsat}
  Consider $F \in \Smusat$.
  \begin{enumerate}
  \item\label{cor:charakmarsat1} $F$ is trivially the saturation of itself.
  \item\label{cor:charakmarsat2} If $F = A(\var(F))$, then this is also the only possibility for $F$ being a saturation, that is, if $F$ is the saturation of some $F' \in \Musat$, then we have $F' = F$.
  \item\label{cor:charakmarsat3} Otherwise $F$ is a saturation of some clause-set other than itself, that is, if $F \not= A(\var(F))$, then there is some $F' \in \Musat$ with $F' \not= F$ such that $F$ is a saturation of $F'$.
  \end{enumerate}
\end{corol}
\begin{prf}
Part \ref{cor:charakmarsat1} is trivial. For Part \ref{cor:charakmarsat2} assume that $F = A(\var(F))$, and we have $F = \saturate(F',C,x)$ for some $F' \in \Musat$: But since $F$ is marginal, $F'$ is not minimally unsatisfiable. Finally for Part \ref{cor:charakmarsat3} note, that if $F \ne A(\var(F))$, then by Lemma \ref{lem:charakmarsat} $F$ is not marginal, and thus there is $C \in F$ and $x \in C$ such that for $C' := C \sm \set{x}$ and $F' := (F \sm \set{C}) \addcup \set{C'}$ we have $F' \in \Musat$. Now $F = \saturate(F',C',x)$. \Qed
\end{prf}

By Lemma \ref{lem:charakmarsat} we know that $F \in \Smusat$ is marginal iff $F = A(\var(F))$; so if $F \in \Smusat$ is not full, then there is a literal occurrence which can be removed without destroying minimal unsatisfiability, that is, there is $C \in F$ and $x \in C$ such that $F' := (F \sm \set{C}) \addcup \set{C \sm \set{x}} \in \Musat$ (note that $F' \in \Usat$ in any case, but minimality in general is not maintained). So for inputs $F \in \Smusat$ the \emph{existence} of such $C, x$ is decidable in linear time (namely they exist iff $F$ is not full). But finding such $C, x$ should be hard in general, and the decision problem, whether a concrete literal can be removed, even for inputs $F \in \Smusat$ should be NP-complete:
\begin{quest}\label{que:elimlit}
  Is the promise problem for inputs $F \in \Smusat$, $C \in F$, $x \in C$, whether ``$F' := (F \sm \set{C}) \addcup \set{C \sm \set{x}} \in \Musat$ ?'', NP-complete? (That is, is there a polytime computation $G \in \Cls \leadsto (F,C,x) \in \Smusat \times \Cl \times \Lit$, with $x \in C \in F$, such that $G \in \Sat \Lra F' = (F \sm \set{C}) \addcup \set{C \sm \set{x}} \in \Musat$ ?) Note that a proof of $F' \in \Musat$ consists in providing for each $D \in F \sm \set{C}$ a satisfying assignment for $F' \sm \set{D}$ (unsatisfiability of $F'$ is trivial). And is the promise problem for input $F \in \Musat$, whether $F$ is marginal, coNP-complete? (That is, is there a polytime computation $G \in \Cls \leadsto F \in \Musat$, such that $G \in \Usat$ iff $F$ is marginal?)
\end{quest}

The decision ``$F' \in \Musat$ ?'' in Question \ref{que:elimlit} is easy for $F \in \Uclash$, namely iff no ``subsumption resolution'' with another clause containing $\ol{x}$ can be performed, i.e., there is no $D \in F$ with $\ol{x} \in D$ and $C \sm \set{x} \sse D$ (obviously this is necessary):
\begin{lem}\label{lem:uhit1l}
  Consider $F \in \Uclash$, $C \in F$ and $x \in C$. Let $C' := C \sm \set{x}$, and let $F' := (F \sm \set{C}) \addcup \set{C'}$. Then $F' \in \Musat$ iff there is no $D \in F \sm \set{C}$ with $C' \subset D$.
\end{lem}
\begin{prf}
If there is $D \in F \sm \set{C}$ with $C' \subset D$, then $F' \notin \Musat$. So assume there is no such $D$. Assume $F' \notin \Musat$. Thus there is $E \in F'$ with $F' \sm \set{E}$ unsatisfiable. We must have $E \not= C'$, since otherwise $F \sm \set{C}$ would be unsatisfiable. Since $F$ is hitting, $E$ clashes with every clause of $F' \sm \set{C'}$. It follows that $C' \subset E$ must hold (since the falsifying assignments for $E$ are disjoint with those for any clause in $F' \sm \set{C'}$), contradicting the minimal unsatisfiability of $F$. \Qed
\end{prf}

Some examples on removable literal occurrences illustrate Lemma \ref{lem:uhit1l}:
\begin{examp}\label{exp:addlitcluhit}
  For $F := \set{\set{a,b},\set{\ol{a},b},\set{\ol{b}}} \in \Uclashi{\delta=1}$ we can exactly remove one of the two literal-occurrences of $b$ and still obtain a clause-set in $\Musat$ (of course not in $\Uclash$ anymore; the resulting clause-sets are in fact marginally minimally unsatisfiable). For $\Dt{2} = A_2 = \set{\set{1,2},\set{-1,-2},\set{-1,2},\set{-2,1}} \in \Uclashnsi{\delta=2}$ we can not remove any literal occurrence without leaving $\Musat$ (i.e., $\Dt{2}$ is marginal).
\end{examp}

\subsection{Splitting}
\label{sec:Splitting}

An important characterisation of saturation for $F \in \Cls$, shown in \cite{Kullmann2007ClausalFormZII} (extending \cite[Lemma C.1]{Ku99dKo}), is that splitting on a variable $v$ yields minimally unsatisfiable clause-sets $\pao v0 * F$, $\pao v1 * F$. This enables induction on the number of variables, a central method for this \Schrift{}; see Lemma \ref{lem:auxminvardeg} for the basic example.
\begin{lem}\label{lem:auxSMUSAT}
  Consider $F \in \Cls$ not containing $C \subset D$ with $\abs{C} + 1 = \abs{D}$.
  \begin{enumerate}
  \item\label{lem:auxSMUSAT3} If there is $v \in \Va$ with $\pao{v}{0} * F, \pao{v}{1} * F \in \Musat$, then $F \in \Musat$.
  \item\label{lem:auxSMUSAT2} If there is $v \in \Va$ with $\pao{v}{0} * F, \pao{v}{1} * F \in \Smusat$, then $F \in \Smusat$.
  \item\label{lem:auxSMUSAT1} $F \in \Smusat$ iff $F \ne \top$ and $\fa\, v \in \var(F) \; \fa\, \ve \in \set{0,1}: \pao{v}{\ve} * F \in \Musat$.
  \end{enumerate}
\end{lem}
\begin{prf}
 For Part \ref{lem:auxSMUSAT3} assume $F \notin \Musat$; thus there is $C \in F$ with $F' := F \sm \set{C} \in \Usat$.  We consider three cases:
\begin{enumerate}
\item $v \notin \var(C)$: Due to the assumption on subsumption-freeness we have $C \addcup \set{v} \notin F'$. Now $C \in \pao v0 * F$, while $(\pao v0 * F) \sm \set{C} = \pao v0 * F' \in \Usat$, contradicting $\pao v0 * F \in \Musat$.
\item $v \in C$: By assumption holds $C' := C \sm \set{v} \notin F'$. Now $C' \in \pao v0 * F$, while $(\pao v0 * F) \sm \set{C'} = \pao v0 * F' \in \Usat$, contradicting $\pao v0 * F \in \Musat$.
\item $\ol{v} \in C$: By assumption holds $C' := C \sm \set{\ol{v}} \notin F'$. Now $C' \in \pao v1 * F$, while $(\pao v1 * F) \sm \set{C'} = \pao v1 * F' \in \Usat$, contradicting $\pao v1 * F \in \Musat$.
\end{enumerate}
Now consider Part \ref{lem:auxSMUSAT2}. By Part \ref{lem:auxSMUSAT3} we already know that $F \in \Musat$ holds. Assume that $F \notin \Smusat$; thus there is $C \in F$ and a literal $x$ with $F' := \saturate(F,C,x) \in \Usat$. So by Lemma \ref{lem:basicharacpsat} we have $F \sm \set{C} \models C' := C \addcup \set{\ol{x}}$. There exists at least one $\ve \in \set{0,1}$ with $\pao v{\ve} * \set{C'} \ne \top$, and then $\pao{v}{\ve} * (F \sm \set{C}) \models \pao{v}{\ve} * C'$. If $\var(x) = v$, then this contradicts minimal unsatisfiability of $\pao{v}{\ve} * F$. And if $\var(x) \ne v$, then $\pao{v}{\ve} * F \sm \pao{v}{\ve} * \set{C} \models (\pao{v}{\ve} * C) \addcup \set{\ol{x}}$, contradicting saturatedness of $\pao{v}{\ve} * F$ by Lemma \ref{lem:basicharacpsat}.

Part \ref{lem:auxSMUSAT1} is Corollary 5.3 in \cite{Kullmann2007ClausalFormZII}: we will re-prove the direction from left to right in Subsection \ref{sec:Controllingdeficiency} of this \Schrift, while for the direction from right to left the additional assumption is missing in \cite{Kullmann2007ClausalFormZII}, and so we give the (easy) proof:  By Part \ref{lem:auxSMUSAT3} we obtain $F \in \Musat$. If $F \notin \Smusat$, then there is a clause $C \in F$ and a literal $x \in \lit(F) \sm \lit(C)$, such that replacing $C$ by $C \addcup \set{x}$ is still unsatisfiable, but then $\pao x1$ deletes a further clause from $\pao x1 * F$, which is minimally unsatisfiable. \Qed
\end{prf}

The essence of the assumption on special subsumption-freeness in Lemma \ref{lem:auxSMUSAT} (automatically fulfilled if $F \in \Musat$) is to make sure that no contraction takes place when applying $\pao v0, \pao v1$. Alternatively we could use multi-clause-sets, since then no contractions would be performed, and the doubled clauses would destroy minimal unsatisfiability. In \cite[Lemma 1]{KullmannZhao2011Bounds} and the underlying report \cite[Lemma 2.1]{KullmannZhao2010Bounds} that additional assumption is missing by mistake for Parts \ref{lem:auxSMUSAT3}, \ref{lem:auxSMUSAT2}:
\begin{examp}\label{exp:mistake}
  Counter-examples for Parts \ref{lem:auxSMUSAT3}, \ref{lem:auxSMUSAT2} (without the additional assumption) are obtained by taking $F_0 \in \Musat$ resp.\ $F_0 \in \Smusat$ and $C \in F_0$ together with $v \in \Va \sm \var(F)$ and letting $F := F \addcup \set{C \cup \set{v}} \notin \Musat$. The counter-examples for Part \ref{lem:auxSMUSAT1} are the clause-sets $\set{\bot, \set{x}}$ and $\set{\bot,\set{x},\set{\ol{x}}}$ for some literal $x$ (there are no other counter-examples). For both $\ve \in \set{0,1}$ we have $\pao x{\ve} * F \in \Uclash$, but $F \notin \Musat$. Note that for a multi-clause-set $F$ the contraction would not occur, but we had e.g.\ $\pao{x}{\ve} * F = \set{2 * \bot}$ in the second case.
\end{examp}

\section{Variable-minimal unsatisfiability}
\label{sec:vmusat}

In \cite{ChenDing2006VMU} the generalisation of minimal unsatisfiability to ``variable-minimal unsatisfiability'' has been introduced, and the class of all such clause-sets is denoted by $\Vmusat$, the set of clause-sets $F \in \Usat$ such that for every $F' \sse F$ with $F' \in \Usat$ holds $\var(F') = \var(F)$; see \cite{BelovIvriiMatsliahMarquesSilva2012VMUSs} for related algorithms for computing $\var(F')$ for some $F' \sse F$ with $F' \in \Musat$ for input $F \in \Usat$. The corresponding class $\ul{\Vmusat}$ of multi-clause-sets is invariant under multiplicities. Thus, as with $\Lean$ (and different from $\Musat$), regarding variable-minimal unsatisfiability w.l.o.g.\ multi-clause-sets can be cast down to clause-sets. The main application of $\Vmusat$ in this \Schrift{} is obtained in Corollary \ref{cor:auxminvardegsigma}, where we will see that unsatisfiable clause-sets with equal surplus and deficiency are in $\Vmusat$. In Corollary \ref{cor:mleanvmu} we conclude from that, that for lean clause-sets the sub-instances with minimal deficiency obtained via restriction are in $\Vmusat$, which will allow us to lift the general upper bound on min-var-degrees from $\Musat$ to $\Lean$ (and sharpening deficiency by surplus). We now develop the basic theory of $\Vmusat$ from scratch, correcting some errors from the literature. The basic (trivial) characterisation of $\Vmusat$ is:

\begin{lem}\label{lem:characVMUtrivial}
  For $F \in \Cls$ holds $F \in \Vmusat$ if and only if $F \in \Usat$ and for all $v \in \var(F)$ holds $\set{C \in F : v \notin \var(C)} \in \Sat$.
\end{lem}

By definition we have $\Musat \subset \Vmusat$, moreover, as shown in \cite[Lemma 6]{ChenDing2006VMU}, for every deficiency $k \ge 2$ we have $\Musati{\delta=k} \subset \Vmusati{\delta=k}$ (for example, for every $F \in \Musati{\delta=k}$, $k \in \NN$, and every non-full clause $C \in F$, i.e., $\var(C) \subset \var(F)$, we can add to $F$ a full clause subsumed by $C$, obtaining $F' \in \Vmusati{\delta=k+1} \sm \Musati{\delta=k+1}$).

In \cite{ChenDing2006VMU} there is the false statement ``$\Vmusat \not\sse \Lean$'', based on the following erroneous example:
\begin{examp}\label{exp:falseexampvmusatlean}
  \cite[Page 266]{ChenDing2006VMU} gives the example $F_4 := \set{\set{a},\set{\ol{b}},\set{\ol{a},b},\set{a,\ol{b}}}$ with the assertion ``$F_4 \in \Vmusat \sm \Lean$''. Obviously we have $F_4 \in \Vmusat$, but we also have $F_4 \in \Lean$. Using the characterisation from \cite{Ku98e} (which is the only characterisation used in \cite{ChenDing2006VMU}), that $F \in \Lean$ holds iff every clause of $F$ can be used in a tree-resolution refutation of $F$, we see this as follows: the sole subset of $F_4$ in $\Musat$ is $\set{\set{a},\set{\ol{b}},\set{\ol{a},b}}$, while the clause $\set{a,\ol{b}}$ can(!) also be used in a tree-resolution refutation --- it is obviously superfluous, but nevertheless there is a tree-resolution refutation using it, namely via $(\set{a} \res \set{\ol{a},b}) \res \set{a,\ol{b}} = \set{a}$.
\end{examp}

Based on the characterisation of lean clause-sets via autarkies, it is easy to show that $\Vmusat$ consists of special lean clause-sets (thus Figure 1 in \cite{ChenDing2006VMU} needs to be corrected, showing instead that $\Lean$ is indeed a superclass of $\Vmusat$):
\begin{lem}\label{lem:vmusatlean}
  $\Vmusat \subset \Lean \sm \set{\top}$.
\end{lem}
\begin{prf}
In \cite{ChenDing2006VMU} the characterisation of $\Lean$ via variables usable in resolution refutation was (only) used. Here we use the equivalent characterisation via autarkies, shown in \cite[Theorem 3.16]{Ku98e}, and used as our definition in Subsection \ref{sec:prelimAut}, namely that for $F \in \Cls$ holds $F \in \Lean$ iff there is no autarky $\vp$ for $F$ with $\var(\vp) \cap \var(F) \ne \es$: If we had such an autarky for $F \in \Vmusat$, then $\vp * F \in \Usat$ with $\vp * F \subset F$ and $\var(\vp * F) \sse \var(F) \sm \var(\vp)$, contradicting $F \in \Vmusat$. That $\Vmusat$ is a strict subset of $\Lean \sm \set{\top}$, can for example be seen by \cite[Lemma 3.2]{Ku98e}, showing that if we extend a minimally unsatisfiable clause-set via Extended Resolution, then we stay in $\Lean$, while adding new clauses with new variables, and thus leaving $\Vmusat$; another example is clause-set $F_3$ from \cite[Page 266]{ChenDing2006VMU}. \Qed
\end{prf}

Thus it follows $\Vmusati{\delta=1} = \Musati{\delta=1}$ (shown in \cite[Lemma 6]{ChenDing2006VMU}), since by \cite[Corollary 5.7]{Ku98e} holds $\Leani{\delta=1} \cap \Usat = \Musati{\delta=1}$.

For $F \in \Vmusat$ obviously there is some $F' \sse F$ with $\var(F') = \var(F)$ and $F' \in \Musat$; \cite[Lemma 5]{ChenDing2006VMU} asserts the converse, but this is false, as the following simple example shows:
\begin{examp}\label{exp:falseexampvmu}
  Consider $F := \set{\bot,\set{v},\set{\ol{v}}}$ and $F' := \set{\set{v},\set{\ol{v}}}$; we have $F' \in \Musat$ and $\var(F') = \var(F)$, but $F \in \Lean \sm \Vmusat$, since $\set{\bot} \in \Usat$. If we don't want to use the empty clause, then we can consider any $F' \in \Musat$ with $v \in \var(F')$ and $\set{\set{v},\set{\ol{v}}} \cap F' = \es$, and let $F := F' \addcup \set{\set{v},\set{\ol{v}}}$ --- again we have $F' \in \Musat$ and $\var(F') = \var(F)$, but $F \in \Lean \sm \Vmusat$ (note $\var(F') \supset \set{v}$).
\end{examp}

The corrected version of \cite[Lemma 5]{ChenDing2006VMU} is as follows:
\begin{lem}\label{lem:crtivmusat}
  For $F \in \Cls$ let $\UU_F$ be the set of $F' \sse F$ with $\var(F') = \var(F)$, $\delta(F') \ge 1$, and $F' \in \Usat$. Then $F \in \Vmusat$ if and only if $F \in \UU_F$ and all minimal elements of $\UU_F$ w.r.t.\ the subset-relation are minimally unsatisfiable.
\end{lem}
\begin{prf}
The condition is necessary, since if $F \in \Vmusat$, then on the one hand we have $F \in \Lean \sm \set{\top}$, and thus $\delta(F) \ge 1$ by \cite{Ku98e} (or use \cite[Lemma 3]{ChenDing2006VMU}); and on the other hand if there would be a minimal element $F' \in \UU_F$ which wouldn't be minimally unsatisfiable, then there would be some $F'' \subset F'$ with $F'' \in \Musat$, whence by definition of $\UU_F$ we get $\var(F'') \subset \var(F)$ contradicting $F \in \Vmusat$.

For the other direction assume, that we have $\UU_F$ as specified, and we have to show $F \in \Vmusat$. Since $F \in \UU_F$, we have $F \in \Usat$. Consider now some $F' \sse F$ with $F' \in \Usat$, and assume $\var(F') \subset \var(F)$. Consider some minimal $F'' \in \UU_F$ (regarding inclusion) with $F' \subset F'' \sse F$. Furthermore consider a minimal element $G \in \UU_F$ with $G \sse F''$; by assumption $G \in \Musat$, and since $F' \subset F''$, we have $G \subset F''$. If for $C \in F''$ we had $\var(F'' \sm \set{C}) \subset \var(F'')$, then there would be $x \in C$ such that $x$ or $\ol{x}$ is pure in $F''$, thus also pure in $G$, whence $C \notin G$ (since $G \in \Musat$), contradicting $\var(G) = \var(F)$. Now choose some $C \in F'' \sm F'$ (we have $\var(F'' \sm \set{C}) = \var(F'')$); by minimality of $F''$ we now have $\delta(F'' \sm \set{C}) \le 0$ (otherwise all conditions for $\UU_F$ are fulfilled for $F'' \sm \set{C}$), whence $\delta(F'') = 1$. Due to $\var(F'') = \var(G)$ and $G \subset F''$ it follows $\delta(G) \le 0$, contradicting $G \in \Musat$. \Qed
\end{prf}

The following examples show applications of Lemma \ref{lem:crtivmusat}:
\begin{examp}\label{exp:characVMU}
  Consider the two (non-)examples from Example \ref{exp:falseexampvmu}:
  \begin{enumerate}
  \item For $F = \set{\bot,\set{v},\set{\ol{v}}}$ we have the minimal element  $\set{\bot,\set{v}}$ of $\UU_F$ which is not minimally unsatisfiable.
  \item For $F = F' \addcup \set{\set{v},\set{\ol{v}}}$ consider a minimal $F'' \sse F'$ with $\var(G) = \var(F)$ and $\delta(G) \ge 1$ for $G := F'' \addcup \set{\set{v},\set{\ol{v}}}$ (note $\top \subset F'' \subset F$): now $G$ is a minimal element of $\UU_F$ which is not minimally unsatisfiable.
  \end{enumerate}
\end{examp}

Based on \cite[Lemma 5]{ChenDing2006VMU}, also the proof of Theorem 3 in \cite{ChenDing2006VMU} is false (the procedure goes astray on the clause-sets of Example \ref{exp:falseexampvmu}). Fortunately we can give a simple proof of the assertion, which even shows fixed-parameter tractability (fpt) of the decision problem ``$F \in \Vmusati{\delta=k}$ ?'' in the parameter $k$:
\begin{thm}\label{thm:vmusatk}
  Membership decision ``$F \in \Vmusati{\delta=k}$ ?'' for input $F \in \Cls$ is fpt in the parameter $k \in \ZZ$.
\end{thm}
\begin{prf}
If $F \notin \Mlean$, then $F \notin \Vmusat$ (by Lemma \ref{lem:vmusatlean}). So we can assume now $F \in \Mlean$, and thus we have $\delta(F') < \delta(F)$ for all $F' \subset F$. Now the decisions of Lemma \ref{lem:characVMUtrivial}, as discussed in Example \ref{exp:MLEANfpt}, are fpt in $k$. \Qed
\end{prf}

\cite[Theorem 1]{ChenDing2006VMU} shows that decision of ``$F \in \Vmusat$ ?'' is $D^P$-complete. We now turn to the two basic reduction processes of this \Schrift{}, elimination of singular variables (Section \ref{sec:elimcreatesing}), and full subsumption resolution (Section \ref{sec:2subr}).

\section{Eliminating and creating singularity}
\label{sec:elimcreatesing}

In this section we continue the study of singular variables in minimally unsatisfiable clause-sets, as initiated in \cite{KullmannZhao2012ConfluenceC,KullmannZhao2012ConfluenceJ}. In Section \ref{sec:singDPred} we look at the reduction process, eliminating singular variables. A main insight is Lemma \ref{lem:DPminvdeg}, showing that the elimination is harmless concerning the minimum variable-degree. In Subsection \ref{sec:singDPext} we introduce the inverse elimination (``extension''); the main point here is the precise statement of the various conditions. Finally in Subsection \ref{sec:MUunit} we consider a special case of singularity, namely unit-clauses.

\subsection{Singular DP-reduction}
\label{sec:singDPred}

In \cite{KullmannZhao2012ConfluenceJ} the process of ``singular DP-reduction'' has been studied for minimally unsatisfiable clause-sets. By it we can reduce the case of arbitrary $F \in \Musat$ to (nonsingular) $F' \in \Musatns$ (that is, for every $v \in \var(F')$ we have $\ldeg_{F'}(v), \ldeg_{F'}(\ol{v}) \ge 2$). The definition is as follows (see \cite[Definition 8]{KullmannZhao2012ConfluenceJ}):
\begin{defi}[\cite{KullmannZhao2012ConfluenceJ}]\label{def:sdp}
  The relation \bmm{F \tsdp F'} (singular DP-reduction) holds for clause-sets $F, F' \in \Cls$, if there is a singular variable $v$ in $F$, such that $F'$ is obtained from $F$ by DP-reduction on $v$, that is, $F' = \dpi{v}(F)$. The reflexive-transitive closure of this relation is denoted by \bmm{F \tsdps F'}.
\end{defi}
By $\sdp(F) \subset \Musat$ for $F \in \Musat$ the set of nonsingular $F' \in \Musat$ with $F \tsdps F'$ is denoted. For us the main property of $\sdp(F)$ is that it is not empty. In \cite{KullmannZhao2012ConfluenceJ} it is shown that for $S \in \Smusat$ we have $\abs{\sdp(F)} = 1$, and that for arbitrary $F \in \Musat$ and $F', F'' \in \sdp(F)$ we have $n(F') = n(F'')$.

\begin{examp}\label{exp:sdp}
  In \cite{KullmannZhao2012ConfluenceJ} the following is shown for $F \in \Musat$:
  \begin{enumerate}
  \item For $\delta(F)=1$ we have $\sdp(F) = \set{\set{\bot}}$.
  \item For $\delta(F)=2$ all elements of $\sdp(F)$ are pairwise isomorphic.
  \item For $\delta(F) \ge 3$ in general there are non-isomorphic elements in $\sdp(F)$.
  \end{enumerate}
\end{examp}

By Sections 3.1, 3.2 in \cite{KullmannZhao2012ConfluenceJ} we have the following basic preservation properties:
\begin{lem}[\cite{KullmannZhao2012ConfluenceJ}]\label{lem:sDPrMU}
  For $F, F' \in \Musat$ with $F \tsdps F'$ we have:
  \begin{enumerate}
  \item\label{lem:sDPrMU1} $\delta(F') = \delta(F)$.
  \item\label{lem:sDPrMU20} $F \in \Musat \Ra F' \in \Musat$.
  \item\label{lem:sDPrMU2} $F \in \Smusat \Ra F' \in \Smusat$.
  \item\label{lem:sDPrMU3} $F \in \Uclash \Ra F' \in \Uclash$.
  \end{enumerate}
\end{lem}

Although singular DP-reduction can reduce the variable-degree of some variables, it can not decrease the \emph{minimum} variable-degree:
\begin{lem}\label{lem:DPminvdeg}
  For $F, F' \in \Musat$ with $F \tsdps F'$ we have $\minvdeg(F') \ge \minvdeg(F)$.
\end{lem}
\begin{prf}
It is sufficient to consider the case $F' = \dpi{v}(F)$ for a singular variable $v$. Assume $\minvdeg(F') < \minvdeg(F)$;  thus $\var(F') \not= \es$ (otherwise we have $\minvdeg(F') = +\infty$), and we consider $w \in \var(F')$ with $\vdeg_{F'}(w) = \minvdeg(F')$. So we have $\vdeg_{F'}(w) < \vdeg_F(w)$, and thus by \cite[Lemma 24]{KullmannZhao2012ConfluenceJ}, for all clauses $C \in F$ with $v \in \var(C)$ we have $w \in \var(C)$. But then $\minvdeg(F') = \vdeg_{F'}(w) \ge \vdeg_F(v) \ge \minvdeg(F) > \minvdeg(F')$, a contradiction. \Qed
\end{prf}

Thus in order to determine the minimum variable-degree for minimally unsatisfiable clause-sets in dependency on the deficiency, w.l.o.g.\ one can restrict attention to saturated and nonsingular instances:
\begin{corol}\label{cor:wlognonsingsat}
  For all $k \in \NN$, $k \ge 2$, holds:
  \begin{enumerate}
  \item\label{cor:wlognonsingsat1} $\minnonmer(k) = \minvdeg(\Musati{\delta=k}) = \minvdeg(\Smusatnsi{\delta=k})$.
  \item\label{cor:wlognonsingsat2} $\minnonmerh(k) = \minvdeg(\Uclashi{\delta=k}) = \minvdeg(\Uclashnsi{\delta=k})$.
  \item\label{cor:wlognonsingsat3} $\minnonmer(1) = \minnonmerh(1) = \minvdeg(\set{\set{1},\set{-1}}) = 2$.
  \end{enumerate}
\end{corol}
\begin{prf}
For Part \ref{cor:wlognonsingsat1} we note that by Lemma \ref{lem:DPminvdeg} for every $F \in \Musati{\delta=k}$ we can find $F' \in \Smusatnsi{\delta=k}$ with $\minvdeg(F') \ge \minvdeg(F)$, and thus $\minvdeg(\Musati{\delta=k}) \le \minvdeg(\Smusatnsi{\delta=k})$, while $\minvdeg(\Musati{\delta=k}) \ge \minvdeg(\Smusatnsi{\delta=k})$ holds due to $\Smusatnsi{\delta=k} \sse \Musati{\delta=k}$. The same reasoning applies for Part \ref{cor:wlognonsingsat2}. We used that by singular DP-reduction not all variables can vanish here, since $\delta(\set{\bot}) = 1$. For $F \in \Musati{\delta=1}$ however we can reach $\set{\bot}$ (recall Example \ref{exp:MU1}), and thus we stop one step before, which proves Part \ref{cor:wlognonsingsat3}. \Qed
\end{prf}

See Lemma \ref{lem:conjimplcomp} for some conditions under which the maps $k \in \NN \mapsto \minnonmerh(k)$ and $k \in \NN \mapsto \minnonmer(k)$ would be computable (computability of $\minnonmer(k)$ is part of Conjecture \ref{con:minnonmernotcomp}, while Conjecture \ref{con:minvdeghit} states $\minnonmerh = \minnonmer$).

\subsection{Singular DP-extensions}
\label{sec:singDPext}

We consider now the reverse direction of singular DP-reduction, from $\dpi{v}(F)$ to $F$, as a \emph{singular extension}, and also generalise it to arbitrary clause-sets. This process was mentioned in \cite[Examples 15,19,54]{KullmannZhao2012ConfluenceJ} for minimally unsatisfiable $\dpi{v}(F)$, called ``inverse singular DP-reduction'' there:
\begin{defi}\label{def:singext}
  Consider a clause-set $G \in \Cls$, a variable $v \in \Va \sm \var(G)$, and $m \in \NN$ with $m \le c(G)$. A \textbf{singular $m$-extension of $G$ with $v$} is a clause-set $F \in \Cls$ obtained as follows (employing four choice steps):
  \begin{enumerate}
  \item[(i)] $m$ different clauses $D_1, \dots, D_m \in G$ are chosen.
  \item[(ii)] A subset $C \sse \bca_{i=1}^m D_i$ is chosen.
  \item[(iii)] Clauses $D_i' \in \Cl$ for $i \in \tb 1m$ with $(D_i \sm C) \sse D_i' \sse D_i$ are chosen.
  \item[(iv)] A literal $x$ with $\var(x) = v$ is chosen.
  \end{enumerate}
  Let $C' := C \addcup \set{x}$ and $D_i'' := D_i' \addcup \set{\ol{x}}$ for $i \in \tb 1m$. $F$ is obtained by adding $C'$ and replacing $D_i$ with $D_i''$, i.e., $F := (G \sm \set{D_1,\dots,D_m}) \addcup \set{C', D_1'',\dots,D_m''}$.
\end{defi}
Note that the order of the clauses $D_i$ does not matter at all (Step (i)), so in reality an $m$-element subset of $G$ is chosen, and that the choice of $x$ (Step (iv)), either $x=v$ or $x=\ol{v}$, is inessential (the resulting $F$ just differ by flipping variable $v$). We also note that the clauses $D_1', \dots, D_m'$ are pairwise different. Obviously the four choices are always possible (only $1 \le m \le c(G)$ is needed for that).

\begin{examp}\label{exp:isdp}
  Consider $G :=  \set{\set{a,b,c},\set{a,b,\ol{c}}}$, $m := 2$, and so $\set{D_1,D_2} = G$ and $\bca_{i=1}^m D_i = \set{a,b}$, and the choices $C := \set{a}$, $D_1' := \set{b,c}$, $D_2' := \set{a,b,\ol{c}}$, and $x := v$. Then the $2$-extension $F$ of $G$ is $F = \set{\set{v,a},\set{\ol{v},b,c},\set{\ol{v},a,b,\ol{c}}}$.
\end{examp}

By definition we have for an $m$-extension $F$ of $G \in \Cls$ with $v$ the following simple properties: $c(F) = c(G) + 1$, $n(F) = n(G) + 1$, $\delta(F) = \delta(G)$, $v$ is singular for $F$, $\vdeg_F(v) = m+1$, $\dpi{v}(F) = G$. Indeed Definition \ref{def:singext} captures the inversion of singular DP-reduction:
\begin{lem}\label{lem:simpcharmext}
  Consider $m \in \NN$, $G, F \in \Cls$ and $v \in \Va$. Then $F$ is an $m$-extension of $G$ by $v$ iff $v$ is singular for $F$, $\vdeg_F(v) = m+1$, $\dpi{v}(F) = G$, $c(G) = c(F) - 1$.
\end{lem}
\begin{prf}
If $F$ is an $m$-extension of $G$ by $v$, then the four properties hold, as we have already mentioned. Now assume these four properties hold. Let the $m$ clauses $D_1,\dots,D_m$ be the result of singular DP-reduction on $v$ for $F$; they must be pairwise different, and all $m$ resolutions must be possible, otherwise $c(G) < c(F) - 1$. And let $C$ be the singular occurrence of $v$ minus the variable $v$. Now all properties of a singular $m$-extension (Definition \ref{def:singext}) are easily checked. \Qed
\end{prf}

Singular extensions behave well regarding minimal unsatisfiability:
\begin{lem}\label{lem:charmext}
  Consider $m \in \NN$, $G \in \Cls$ and an $m$-extension $F$ of $G$ by $v \in \Va$. Then $F \in \Musat \Lra G \in \Musat$.
\end{lem}
\begin{prf}
This follows by Lemma \ref{lem:simpcharmext} together with \cite[Lemma 9, Parts 1, 2]{KullmannZhao2012ConfluenceJ}. \Qed
\end{prf}

In the situation of Lemma \ref{lem:charmext}, regarding saturatedness we only have the direction $F \in \Smusat \Ra G \in \Smusat$, while for the other direction the conditions of \cite[Lemma 12]{KullmannZhao2012ConfluenceJ} need to be observed (this would yield ``saturated extensions'', which however we do not need here).

\subsection{Unit clauses}
\label{sec:MUunit}

We conclude this section by considering unit-clauses in minimally unsatisfiable clause-sets. The following (fundamental, simple) lemma is \cite[Lemma 14]{KullmannZhao2012ConfluenceJ}; in \cite[Subsection 3.3]{KullmannZhao2012ConfluenceJ} one finds further information.
\begin{lem}[\cite{KullmannZhao2012ConfluenceJ}]\label{lem:fullunit}
  Consider $F \in \Musat$.
  \begin{enumerate}
  \item\label{lem:fullunit1} If $v$ is full and singular in $F$, then we have $\set{v} \in F$ or $\set{\ol{v}} \in F$.
  \item\label{lem:fullunit2} If $\set{x} \in F$, then $v := \var(x)$ is singular in $F$ (with $\ldeg_F(x) = 1$). If here $F$ is saturated, then $v$ is also full in $F$.
  \end{enumerate}
\end{lem}

So unit-clauses in minimally unsatisfiable clause-sets are strong cases of singular variables. They can obviously be removed by singular DP-reduction, while on the other hand singular $\ge 2$-extensions can \emph{not} remove all unit-clauses:
\begin{lem}\label{lem:singextunit}
  Consider a clause-set $F \in \Musat$ containing at least one unit-clause, and obtain $F'$ from $F$ by a singular $m$-extension, where $m \ge 2$. Then also $F'$ must contain at least one unit-clause.
\end{lem}
\begin{prf}
For a unit-clause $\set{x} \in F$ to be removed in $F'$, it needs to be one of the $D_i$ (using the terminology of Definition \ref{def:singext}). Then the intersection $C$ must be empty (otherwise any other $D_i$ needed to contain $x$, and since $m \ge 2$ this would mean a subsumption in $F$). Thus the extension introduces the new unit-clause $C'$. \Qed
\end{prf}

The following examples show that the assumptions $F \in \Musat$ and $m \ge 2$ in Lemma \ref{lem:singextunit} are needed:
\begin{examp}\label{exp:singextunit}
  First consider $F := \set{\set{a},\set{\ol{a},b},\set{\ol{a},\ol{b}}} \in \Uclashi{\delta=1}$. To see necessity of the condition ``$m \ge 2$'', consider a 1-singular extension, obtaining $F' := \set{\set{v,a},\set{\ol{v},a},\set{\ol{a},b},\set{\ol{a},\ol{b}}} \in \Uclashi{\delta=1}$, which has no unit-clauses. On the other hand, using the notation from Definition \ref{def:singext}, any $\ge 2$-singular extension of $F$ which touches $\set{a}$ (i.e., $\set{a}$ is one of the $D_i$) has intersection $C = \bot$, and thus $C'$ is a new unit-clause, while if $\set{a}$ is not touched, then this unit-clause is simply maintained.

  For the condition ``$F \in \Musat$'' consider $F := \set{\set{a},\set{a,b}} \in \Cls \sm \Musat$, where $F' := \set{\set{v,a},\set{\ol{v},a},\set{\ol{v},b}}$ is a 2-extension without unit-clauses.
\end{examp}

For some $F$ the existence of a unit-clause is indeed necessary for singularity:
\begin{lem}\label{lem:existunitmu2}
  Consider $F \in \Musati{\delta=2}$ with $\minvdeg(F) \ge 4$. Then $F$ is singular if and only if $F$ contains a unit-clause.
\end{lem}
\begin{prf}
That if $F$ contains a unit-clause, then $F$ must be singular, follows by Lemma \ref{lem:fullunit}, Part \ref{lem:fullunit2}. So assume now that $F$ is singular, and we have to show that $F$ contains a unit-clause. Consider a reduction sequence $F = F_0 \tsdp F_1 \tsdp \dots \tsdp F_m$, where $F_m$ is nonsingular (note $m \ge 1$). So there exists $n \ge 2$ such that $F_m$ is isomorphic to $\Dt{n}$ (recall Example \ref{exp:MU2}), and thus every variable of $F_m$ has degree $4$. So by Lemma \ref{lem:DPminvdeg} we know $\minvdeg(F_i) = 4$ for $i \in \tb 0m$. We show by induction on $m$ that $F$ contains a unit-clause. If $m=1$, then in order to obtain the min-var-degree of at least $4$, at least $3$ side-clauses $D_1,\dots,D_3 \in \Dt{n}$ for the singular extension have to be chosen (using Definition \ref{def:singext}), but every literal occurs precisely twice in $\Dt{n}$ (because of variable-degree $4$ and nonsingularity), and thus the intersection $C$ has to be empty, and the new clause introduced by the singular extension is a unit-clause, whence $F$ contains a unit-clause. Finally assume $m > 1$. So by induction hypothesis, $F_1$ contains a unit-clause, and thus by Lemma \ref{lem:singextunit} also $F_0$ contains a unit-clause. \Qed
\end{prf}

We will see (Theorem \ref{thm:MUminvdegdef}) that the condition $\minvdeg(F) \ge 4$ in Lemma \ref{lem:existunitmu2} is equivalent to $\minvdeg(F) = 4$; the following examples show that this bound is sharp:
\begin{examp}\label{exp:minvdef2}
  $F_1 := \set{\set{1,2,3},\set{1,2,-3},\set{-1,2},\set{1,-2},\set{-1,-2}} \in \Musati{\delta=2}$ is a $1$-singular extension of $A_2$,
  where $F_1$ has no unit-clause and $\minvdeg(F_1) = 2$. While a $2$-singular extension of $A_2$ without unit-clauses and with min-var-degree $3$ is $F_2 := \set{\set{1,3},\set{1,2,-3},\set{1,-2,-3},\set{-1,2},\set{-1,-2}} \in \Musati{\delta=2}$.
\end{examp}

We conclude with a simple form of adding a new variable, by adding it in one sign as unit-clause, and adding it in the other sign to all given clauses (so that we obtain a full variable, occurring in all clauses):
\begin{defi}\label{def:extucl}
  A \textbf{full singular unit-extension} of a clause-set $F \in \Cls$ (by unit-clause $\set{x}$) is a clause-set $F' \in \Cls$ obtained from $F$ by adding a unit-clause $\set{x}$ with $\var(x) \notin \var(F)$, and by adding literal $\ol{x}$ to all clauses of $F$, i.e., $F' := \set{\set{x}} \addcup \set{C \addcup \set{\ol{x}} : C \in F}$ for some $x \in \Lit \sm \lit(F)$.
\end{defi}
A full singular unit-extension $F'$ of $F \ne \top$ by $\set{x}$ is a special case of a singular $c(F)$-extension of $F$ with $\var(x)$ (recall Definition \ref{def:singext}), and thus $F' \tsdp F$.\footnote{The case $m=0$ is excluded in Definition \ref{def:singext}, since it is not needed, and would only complicate the formulation.}

\begin{examp}\label{exp:fsue}
  Starting with $\set{\bot}$, the first full singular unit-extension is $\set{\set{v},\set{\ol{v}}}$ (up to the choice of the new literal), the second is $\set{\set{w},\set{v,\ol{w}},\set{\ol{v},\ol{w}}}$. In this way we get special examples of $\Smusati{\delta=1}$ (since we started with $\set{\bot} \in \Smusati{\delta=1}$).

  If we start with $\top$ instead, then first we get $\set{\set{v}}$, and then $\set{\set{w},\set{v,\ol{w}}}$.

  \cite[Example 15, Part 1]{KullmannZhao2012ConfluenceJ} contains two example of ``inverse unit elimination'', where Example (a) there is an example of a full singular unit-extension, while Example (b) there would be a non-full singular unit-extension (where the new variable is not full; this is not used in the present \Schrift{}).
\end{examp}

The process of full singular unit-extension of a clause-set $F$ maintains many properties of $F$, and we list here those we use:
\begin{lem}\label{lem:fullsingue}
  Consider a full singular unit-extension $F'$ of $F$ (by $\set{v}$):
  \begin{enumerate}
  \item\label{lem:fullsingue1} $n(F') = n(F) + 1$ and $c(F') = c(F) + 1$.
  \item\label{lem:fullsingue2} $\delta(F') = \delta(F)$.
  \item\label{lem:fullsingue3} $\surp(F') = \surp(F)$ for $F \ne \set{\bot}$.
  \item\label{lem:fullsingue4} $\minvdeg(F') = \minvdeg(F)$ for $n(F) > 0$.
  \item\label{lem:fullsingue5} $F'$ is satisfiable iff $F$ is satisfiable.
  \item\label{lem:fullsingue6} For $F \not= \top$: $F'$ is lean iff $F$ is lean.
  \item\label{lem:fullsingue7} $F'$ is (saturated) minimally unsatisfiable iff $F$ is (saturated) minimally unsatisfiable.
  \item\label{lem:fullsingue8} $F'$ is hitting iff $F$ is hitting.
  \end{enumerate}
\end{lem}
\begin{prf}
Parts \ref{lem:fullsingue1}, \ref{lem:fullsingue2} follow directly by definition. For Part \ref{lem:fullsingue3} we notice that for $F = \top$ we have $\surp(F') = \surp(F) = 0$, while for $n(F) > 0$ consider $\es \subset V \sse \var(F')$: if $v \notin V$, then $F'[V] = F[V]$, and thus the minimisation for $\surp(F)$ is included in $\surp(F')$, and if $v \in V$, then $\delta(F'[V]) = c(F') - \abs{V} \ge \delta(F') = \delta(F) \ge \surp(F)$, and thus these $V$ do not contribute to the minimisation.

For Part \ref{lem:fullsingue4} we just note that the variables of $F$ keep their degrees in $F'$, while the new variable has degree $\vdeg_{F'}(v) = c(F') > c(F)$, and thus does not contribute to the min-var-degree. Part \ref{lem:fullsingue5} is trivial, and follows also by the satisfiability-equivalence of $\dpi{v}(F)$ and $F$. For Part \ref{lem:fullsingue6} we note, that an autarky for $F'$ involving $v$ must be a satisfying assignment for $F'$, while the autarkies for $F'$ not involving $v$ are the same as the autarkies for $F$. Part \ref{lem:fullsingue7} concerning (just) minimal unsatisfiability follows with \cite[Lemma 9]{KullmannZhao2012ConfluenceJ}, while regarding saturatedness we can use \cite[Lemma 12]{KullmannZhao2012ConfluenceJ} (both assertions also follow easily by direct reasoning). Part \ref{lem:fullsingue8} is trivial. \Qed
\end{prf}

So our fundamental classes are respected by full singular unit-extension:
\begin{corol}\label{cor:fullsingue}
  If $F \in \Musati{\delta=k}$ ($k \in \NN$), then every full singular unit-extension is also in $\Musati{\delta=k}$. If furthermore $F$ is saturated resp.\ hitting, then every full singular unit-extension is also saturated resp.\ hitting.
\end{corol}

Obviously, full singular unit-extension is unique up to isomorphism:
\begin{lem}\label{lem:fullsingueiso}
  Consider a clause-set $F \in \Cls$ and clause-sets $F', F'' \in \Cls$ obtained from $F$ by repeated full singular unit-extensions. Then $F', F''$ are isomorphic if and only if $n(F') = n(F'')$.
\end{lem}
\begin{prf}
The number of repeated full singular unit-extensions leading to $F'$ resp.\ $F''$ is the number of variables in these clause-sets with degree strictly greater than $c(F)$, and sorting these variables by increasing degree yields the sequence of extensions. Thus just from knowing the number of variables in $F', F''$ we can reconstruct them up to isomorphism (using that a full singular unit-extension of $F$ by $\set{x}$ is isomorphic to one by $\set{y}$, for arbitrary literals $x,y$ with new variables). \Qed
\end{prf}

\section{Full subsumption resolution / extension}
\label{sec:2subr}

In this section we investigate the second reduction concept for this \Schrift, ``full subsumption resolution''. As with singular DP-reduction from Section \ref{sec:elimcreatesing}, in general this reduction uncovers hidden structure, while the inverse process, ``full subsumption extension'', serves as a generator for minimally unsatisfiable clause-sets with various properties. For this \Schrift{}, full subsumption resolution starting from some $A(V)$ is of special importance, while a more general use will be important for \cite{KullmannZhao2015FullClauses,KullmannZhao2016UHitSAT}. Subsection \ref{sec:2subdef} discusses the basic definitions (there are various technicalities one needs to be aware of), and first applications are given in Subsection \ref{sec:2subext}.

The basic idea is, for a clause-set $F$ containing two clauses $R \addcup \set{v}, R \addcup \set{\ol{v}} \in F$, to replace these two clauses by the clause $R$, i.e., we consider the case where the resolvent $R$ of parent clauses $C, D$ subsumes both parent clauses (thus the name). This is a very old procedure, based on the trivial observation that $(R \oder v) \und (R \oder \neg v)$ is logically equivalent to $R$. If we perform this in the inverse direction, as an ``extension'', then every clause-set $F \in \Cls$ can be transformed into its (equivalent) ``distinguished'' or ``canonical'' CNF $F' \sse A(\var(F))$ (just expand every non-full clause), which is uniquely determined, namely $F'$ is the set of $C \in A(\var(F))$ such that there is $D \in F$ with $D \sse C$.

We however have to be more careful about deficiency and membership in $\Musat$, and thus will consider only ``full subsumption resolution'', where the resolvent must not be present already, while for the ``strict'' form additionally the resolution variable $v$ must occur also in other clauses. For the inverse forms we have to be equally carefully, making sure that none of the two parent clauses is already present (this prevents the above expansion of arbitrary $F \in \Cls$ to $A(\var(F))$) --- from $A(V)$ by strict full subsumption resolution we can obtain precisely the $F \in \Uclash$ with $\var(F) = V$ (Lemma \ref{lem:inv2subrescan}). A main tool is Lemma \ref{lem:aux2subsumpR}, where especially Part \ref{lem:aux2subsumpR3ab} is somewhat subtle, and can be easily overlooked. Via this tool we have a controlled way of transforming $F \in \Musat$ resp.\ $F \in \Uclash$ into $A(\var(F))$. In Theorem \ref{thm:minnumvarmu} we obtain the determination of the possible numbers of variables and clauses in minimally unsatisfiable clause-sets of a given deficiency.

The (more general) well-known ``subsumption resolution'' is the reduction $F \leadsto (F \sm \set{C}) \cup \set{C \sm \set{x}}$ for $F \in \Cls$, that is the removal of a literal $x \in C$ from a clause $C \in F$, in case there exists $D \in F$ with $\ol{x} \in D$ and $D \sm \set{\ol{x}} \sse C$ (note that $C \res D = C \sm \set{x}$ subsumes $C$). An early use is in \cite{Ro65}, under the name ``replacement principle'', while the terminology ``subsumption resolution'' is used in \cite{GeTs96} (for SAT solving). The earliest sources with a systematic treatment appear to be \cite[Section 7]{KuLu97} and \cite[Section 7]{KuLu98}. An experimental study of the practical importance of subsumption resolution in connection with DP-reductions $F \leadsto \dpi{v}(F)$ (under suitable additional conditions to make DP-reduction feasible; see \cite[Subsection 1.3]{KullmannZhao2012ConfluenceJ} for an overview on such restrictions) is performed in \cite{EenBiere2005Satelite}, under the name of ``self-subsuming resolution'', and continued in \cite{HanSomenzi2009}. A theoretic (similar) use one finds in \cite[Section 4]{OrdyniakPaulusmaSzeider2013Acyclic}, where a variable $v$ is called ``DP-simplicial'' for $F \in \Cls$ iff all resolutions performed by the reduction $F \leadsto \dpi{v}(F)$ are subsumption resolutions.

\subsection{Basic definitions}
\label{sec:2subdef}

Before defining ``full subsumption resolution'' $F \leadsto (F \sm \set{R \addcup \set{v}, R \addcup \set{\ol{v}}}) \addcup \set{R}$ in Definition \ref{def:gen2subres} (so $R$ is new and the two clauses $R \addcup \set{v}, R \addcup \set{\ol{v}}$ vanish), we introduce the ``strict'' form, which is more important to us, and which has the additional condition that $v$ must still occur (in other clauses of $F$; the ``non-strict'' form on the other hand guarantees that $v$ vanishes (see Definition \ref{def:gen2subres})):
\begin{defi}\label{def:2subsumptionres}
  For clause-sets $F,F' \in \Cls$ by \bmm{F \tsubres F'} we denote that $F'$ is obtained from $F$ by one step of \textbf{strict full subsumption resolution}, that is,
  \begin{itemize}
  \item there is a clause $R \in F'$
  \item and a literal $x$ with $\var(x) \notin R$
  \item such that for the clauses $C := R \addcup \set{x}$ and $D := R \addcup \set{\ol{x}}$
  \item we have $F = (F' \sm \set{R}) \addcup \set{C,D}$;
  \item we furthermore require $\var(x) \in \var(F')$.
  \item As usual, the literals $x, \ol{x}$ are the resolution literals, $\var(x)$ is the resolution variable, $C, D$ are the parent clauses, and $R$ is the resolvent.
  \end{itemize}
  We write $F \tsubresk{k} F'$ for $k \in \NNZ$ if exactly $k$ steps have been performed, while we write $F \tsubress F'$ for an arbitrary number of steps (including zero).
\end{defi}

We require $R \notin F$, that is, the (full subsumption) resolvent is not already present in the original clause-set. This is of course satisfied if $F \in \Musat$. We also require that the variable $v$ does not vanish, for the sake of keeping control on the deficiency.

\begin{examp}\label{exp:s2subsr}
  Some simple examples are:
  \begin{enumerate}
  \item $\Dt{2} = \set{\set{1,2},\set{-1,-2},\set{-1,2},\set{-2,1}} \tsubres \set{\set{2}, \set{-1,-2},\set{-2,1}}$, and no further reduction is possible (note that the only possibility is blocked, since variable $1$ would vanish).
  \item $\set{\set{v},\set{\ol{v}}} \not\tsubres \set{\bot}$, as $v$ vanishes, while $\set{\set{v},\set{\ol{v}}, \set{v,x}} \tsubres \set{\bot, \set{v,x}}$.
  \item $\set{\set{v,w},\set{\ol{v},w}, \set{v,\ol{w}}} \not\tsubres \set{\set{v,w},\set{\set{w}, \set{v,\ol{w}}}}$, as one parent clause would be kept, while $\set{\set{v,w},\set{\ol{v},w}, \set{v,\ol{w}}} \tsubres \set{\set{\set{w}, \set{v,\ol{w}}}}$.
  \item $\set{\set{v,w},\set{\ol{v},w}, \set{v,\ol{w}}, \set{v}, \set{w}}$ can not be reduced by strict-full-subsumption resolution, since all possible resolvents are already there.
  \end{enumerate}
\end{examp}

The expansion of a clause $R$ to two clauses $R \addcup \set{v}, R \addcup \set{\ol{v}}$ under the above requirements is called ``extension'':
\begin{defi}\label{def:invss2subres}
  For clause-sets $F, F' \in \Cls$ we say that \textbf{$F$ is obtained from $F'$ by strict full subsumption extension} if $F \tsubres F'$. And for $k \in \NNZ$ we say that \textbf{$F$ is obtained from $F'$ by strict full subsumption extension with $k$ steps} if $F \tsubresk{k} F'$.
\end{defi}
So one step of strict full subsumption extension for a clause-set $F$ uses a non-full clause $R \in F$ and a variable $v \in \var(F) \sm \var(R)$, and replaces $R$ by the two clauses $R \addcup \set{v}, R \addcup \set{\ol{v}}$, where none of them is already present.

\begin{examp}\label{exp:is2subsr}
  From $\set{\set{a},\set{b}}$ by one step of strict full subsumption extension we obtain $\set{\set{a,b},\set{a,\ol{b}},\set{b}}$ and $\set{\set{a},\set{a,b},\set{\ol{a},b}}$; note that no new variable is introduced, that the original clause ($\set{a}$ resp.\ $\set{b}$) vanished, and that the replacement clauses were not already present. For $\set{\set{a,b},\set{a}}$ no strict full subsumption extension is possible. Further examples are obtained by ``reading Example \ref{exp:s2subsr} backwards''.
\end{examp}

The basic properties of strict full subsumption resolution are as follows:
\begin{lem}\label{lem:aux2subsumpR}
  For clause-sets $F, F' \in \Cls$ with $F \tsubresk{k} F'$ ($k \in \NNZ$) we have:
  \begin{enumerate}
  \item\label{lem:aux2subsumpR2n} $\var(F') = \var(F)$.
  \item\label{lem:aux2subsumpR1} $c(F') = c(F) - k$, $\delta(F') = \delta(F) - k$.
  \item\label{lem:aux2subsumpR6} $\minvdeg(F) \ge \minvdeg(F')$.
  \item\label{lem:aux2subsumpR2} $F'$ is logically equivalent to $F$.
  \item\label{lem:aux2subsumpR3a} $F \in \Musat \Ra F' \in \Musat$.
  \item\label{lem:aux2subsumpR4} $F \in \Smusat \Ra F' \in \Smusat$.
  \item\label{lem:aux2subsumpR3ab} Assume $k=1$ with resolution variable $v$ and resolvent $R$, and assume $F' \in \Musat$. Then exactly one of the following three possibilities holds:
    \begin{enumerate}
    \item $\saturate(F',R,v)$ is a partial saturation of $F'$ (recall Definition \ref{def:partsaturation}).
    \item $\saturate(F',R,\ol{v})$ is a partial saturation of $F'$.
    \item $F \in \Musat$.
    \end{enumerate}
  \item\label{lem:aux2subsumpR4a} $F' \in \Smusat \Ra F \in \Musat$.
  \item\label{lem:aux2subsumpR5} $F \in \Clash \Lra F' \in \Clash$, $F \in \Uclash \Lra F' \in \Uclash$.
  \end{enumerate}
\end{lem}
\begin{prf}
Parts \ref{lem:aux2subsumpR2n}, \ref{lem:aux2subsumpR1}, \ref{lem:aux2subsumpR6}, \ref{lem:aux2subsumpR2} follow directly from the definition (using for the first three parts that no variable vanishes). Part \ref{lem:aux2subsumpR3a} holds since we strengthen two clauses into one, which is logically equivalent to its parent clauses, and for Part \ref{lem:aux2subsumpR4} additionally note that a saturation of $F'$ by one literal can also be done on $F$ (if the resolvent is involved, then the resolution variable can not be $v$, since $F \in \Musat$). Now consider Part \ref{lem:aux2subsumpR3ab}. That the two possibilities for partial saturation exclude each other follows by the characterisation of partial saturations in Lemma \ref{lem:basicharacpsat} (and $F' \sm \set{R} \not\models R$). And that each possibility for partial saturation excludes $F \in \Musat$ follows by definition. While that the negation of the two partial saturation possibilities implies $F \in \Musat$ follows again by Lemma \ref{lem:basicharacpsat}. Finally Part \ref{lem:aux2subsumpR4a} follows by Part \ref{lem:aux2subsumpR3ab}, while Part \ref{lem:aux2subsumpR5} follows by trivial combinatorics. 
\Qed
\end{prf}

Part \ref{lem:aux2subsumpR3ab} of Lemma \ref{lem:aux2subsumpR} handles a subtle source for errors: One could easily think that for $F' \in \Musat$ a strict full subsumption extension always yields another $F \in \Musat$, but this is not so, as there are three possible cases to be considered here:

\begin{examp}\label{exp:props2subsr}
  Consider $F := \set{\set{v,a},\set{\ol{v},a},\set{\ol{v}},\set{v,\ol{a}}}$. So $F \tsubresk{1} F'$ for $F' := \set{\set{a},\set{\ol{v}},\set{v,\ol{a}}}$. We have $F' \in \Musat$, but $F \notin \Musat$, and indeed $\saturate(F',R,v) = \set{\set{a,v},\set{\ol{v}},\set{v,\ol{a}}}$ is a partial saturation of $F'$ (while $\saturate(F',R,\ol{v})$ isn't one).

  An example that Part \ref{lem:aux2subsumpR4a} can not be strengthened to ``$F' \in \Smusat \Ra F \in \Smusat$'' is obtained from $F' := \Dt{4} \in \Smusatnsi{\delta=2}$ (recall Example \ref{exp:MU2}) by one strict full subsumption extension on $\set{-1,2}$ with resolution variable $3$, obtaining $F = \set{\set{1,2,3,4},\set{-1,-2,-3,-4},\set{-1,2,3},\set{-1,2,-3},\set{-2,3},\set{-3,4},\set{-4,1}} \in\\ \Musatnsi{\delta=3} \sm \Smusat$, where $\pao 40 * F = \set{\set{1,2,3},\set{-1,2,3},\set{-1,2,-3},\set{-2,3},\set{-3}} \\ \notin \Musat$ (clause $\set{-1,2,-3}$ is superfluous; recall Lemma \ref{lem:auxSMUSAT}, Part \ref{lem:auxSMUSAT1}).
\end{examp}

The condition on the resolution variable for strict full subsumption resolution (that it must not vanish) is exactly needed for Parts \ref{lem:aux2subsumpR2n}, \ref{lem:aux2subsumpR1}, \ref{lem:aux2subsumpR6} of Lemma \ref{lem:aux2subsumpR}. If this condition is dropped, then we speak of \emph{full subsumption resolution}:

\begin{defi}\label{def:gen2subres}
  \textbf{Full subsumption resolution} is defined as strict full subsumption resolution, but now the resolution variable is allowed to vanish. If the resolution variable definitely vanishes, then we speak of if \textbf{non-strict full subsumption resolution}. In the other direction we speak of \textbf{full subsumption extension} resp.\ \textbf{non-strict full subsumption extension}.
\end{defi}
So if $F'$ is obtained from $F$ by one step of non-strict full subsumption extension, then we have $c(F') = c(F) + 1$, $n(F') = n(F) + 1$ and $\delta(F') = \delta(F)$. As mentioned, Lemma \ref{lem:aux2subsumpR} holds for full subsumption resolution except of  Parts \ref{lem:aux2subsumpR2n}, \ref{lem:aux2subsumpR1}, \ref{lem:aux2subsumpR6}.

\begin{examp}\label{exp:g2subsr}
  Considering the non-examples from Example \ref{exp:s2subsr}:
  \begin{enumerate}
  \item $\set{\set{v},\set{\ol{v}}} \not\tsubres \set{\bot}$, but by full subsumption resolution we obtain $\set{\bot}$.
  \item $\set{\set{v,w},\set{\ol{v},w}, \set{v,\ol{w}}} \not\tsubres \set{\set{v,w},\set{\set{w}, \set{v,\ol{w}}}}$, and the transition is also not possible by full subsumption resolution.
  \item $\set{\set{v,w},\set{\ol{v},w}, \set{v,\ol{w}}, \set{v}, \set{w}}$ is irreducible by full subsumption resolution.
  \end{enumerate}
  As follows from the characterisation of $\Smusati{\delta=1} = \Uclashi{\delta=1}$ in \cite{Ku99dKo}, a clause-set $F \in \Cls$ can be reduced by a series of non-strict full subsumption resolutions to $\set{\bot}$ iff $F \in \Smusati{\delta=1} = \Uclashi{\delta=1}$.
\end{examp}

\subsection{Extensions to full clause-sets}
\label{sec:2subext}

If we start with the full clause-sets $A(V)$, then by strict full subsumption resolution we obtain exactly all unsatisfiable hitting clause-sets:
\begin{lem}\label{lem:inv2subrescan}
  If for some finite $V \subset \Va$ we have $A(V) \tsubress F$, then $F \in \Uclash$ holds. And for $F \in \Uclash$ we have $A(\var(F)) \tsubress F$.
\end{lem}
\begin{prf}
The first part follows by Lemma \ref{lem:aux2subsumpR}, Part \ref{lem:aux2subsumpR5} (and $A(V) \in \Uclash$). And for the second part note, that if $F \in \Uclash$ has a non-full clause, then a strict full subsumption extension step can be applied, where the result is still in $\Uclash$ (again by Lemma \ref{lem:aux2subsumpR}, Part \ref{lem:aux2subsumpR5}; if $F$ has only full clauses, then $F = A(\var(F))$). \Qed
\end{prf}

Recall that in Example \ref{exp:props2subsr} we have seen, that strict full subsumption extension does not maintain minimal unsatisfiability in general. Now we show that from arbitrary minimally unsatisfiable $F$ we can indeed go all the way up to $A(\var(F))$, while staying in $\Musat$, when we additionally allow partial saturations:
\begin{lem}\label{lem:inv2satusubrescan}
  For $F \in \Musat$, which is not full, we can always perform a strict full subsumption extensions or a partial saturation. Performing these operations, for any order, any choice, as long as possible, at least one strict full subsumption extension will be performed, and the final result is $A(\var(F))$.
\end{lem}
\begin{prf}
If $F \in \Musat$ has a non-full clause, and if strict full subsumption extension can not be applied in order to obtain $F' \in \Musat$, then by Lemma \ref{lem:aux2subsumpR}, Part \ref{lem:aux2subsumpR3ab}, a partial saturation is possible. By Corollary \ref{cor:charakmarsat}, Part \ref{cor:charakmarsat2}, we can not reach a full-clause-set from a non-full one just by saturation. \Qed
\end{prf}

We obtain sharp upper bounds on deficiency and number of clauses in terms of the number of variables, showing that for $F \in \Musat$ with $n$ variables the worst-case concerning deficiency and number of clauses is reached by $A_n \in \Uclash$:
\begin{corol}\label{cor:upbdefn}
  For $F \in \Musat$ holds $\delta(F) \le 2^{n(F)} - n(F)$ (equivalently, $c(F) \le 2^{n(F)}$). We have equality iff $F$ is full (i.e., $F = A(\var(F))$).
\end{corol}
\begin{prf}
By Lemma \ref{lem:inv2satusubrescan} we can transform $F$ into $A(\var(F))$ by a series of steps not decreasing the number of clauses. Thus $c(F) \le c(A(\var(F))) = 2^{n(F)}$. For non-full $F$, at least one strict full subsumption extension was performed in the transformation (Lemma \ref{lem:inv2satusubrescan}), and so here $c(F) < 2^{n(F)}$. \Qed
\end{prf}

The upper bound $c(F) \le 2^{n(F)}$ for $F \in \Musat$ follows also by the observation, that every clause of $F$ must cover uniquely at least one of the $2^{n(F)}$ total assignments for $F$. We explicitly state the instructive reformulation, that the $A_n$ are the minimally unsatisfiable clause-sets of maximal deficiency for given number $m$ of variables:
\begin{corol}\label{cor:maxdefmun}
  Consider $m \in \NNZ$ and $F \in \Musati{n=m}$ such that $\delta(F)$ is maximal.\footnote{That is, $F \in \Musat$, $n(F) = m$, and for all $F' \in \Musat$ with $n(F') = m$ we have $\delta(F') \le \delta(F)$.} Then $F = A(\var(F))$. Thus the maximal deficiency for $F \in \Musati{n=m}$ is $2^m - m$ (realised by $A_m \in \Uclashi{n=m} \cap \Uclashi{\delta=2^m-m}$).
\end{corol}

So for $m=0,1,2,3,4,5,6$ variables the maximal deficiency of minimally unsatisfiable clause-sets is $1,1,2,5,12,27,58$; in general the deficiencies of the form $2^m - m$ are central for our investigations (note that the function $m \in \NNZ \mapsto 2^m - m \in \NN$ is monotonically increasing). We now determine the numbers of variables and numbers of clauses possible for minimally unsatisfiable clause-sets with a given deficiency. For $k \in \NN$ let $\bmm{\odef(k)} \in \NNZ$ be the smallest $n \in \NNZ$ with $2^n - n \ge k$. So $\odef(k)$ by definition is the smallest $n \ge 0$ with $\delta(A_n) \ge k$ (thus the notation ``nA''). We have $\odef(1) = 0$, $\odef(2) = 2$, $\odef(3) = \dots = \odef(5) = 3$, $\odef(6) = \dots = \odef(12) = 4$ and $\odef(13) = \dots = \odef(27) = 5$. Excluding the first term (note the anomaly that both $A_0$ and $A_1$ have deficiency $1$), the sequence $(\odef(k))_{k \in \NN}$ is sequence \url{http://oeis.org/A103586} in the ``On-Line Encyclopedia of Integer Sequences'' (\cite{Sloane2008OEIS}); as noted there we have $\odef(k) = 1 + \fld(k + \fld(k))$ for $k \ge 2$ (recall $\fld(k) = \floor{\ld(k)}$, that is, $\fld(k)$ is the largest $n \in \NNZ$ with $2^n \le k$).
\begin{thm}\label{thm:minnumvarmu}
  For $\mc{C} \sse \Cls$ we use $n(\mc{C}) := \set{n(F) : F \in \mc{C}}$. For $k \in \NN$ holds $n(\Uclashi{\delta=k}) = n(\Musati{\delta=k}) = \set{n \in \NNZ : n \ge \odef(k)}$, where $\odef(k) \in n(\Uclashnsi{\delta=k})$.
\end{thm}
\begin{prf}
By Corollary \ref{cor:upbdefn} we see $n(\Musati{\delta=k}) \sse \set{n \in \NNZ : n \ge \odef(k)}$. Increasing the number of variables by one while keeping the deficiency constant is achieved by one non-strict full subsumption extension step, which maintains the hitting property, and so it remains to show the existence of $F \in \Uclashnsi{\delta=k}$ with $n(F) = \odef(k)$.

 For $k=1$ we have $F = \set{\bot}$, so assume $k > 1$ (thus $\odef(k) \ge 2$). Let $F_0 := A_{\odef(k)-1}$ (so $\delta(F_0) < k$). Add a variable by one step of non-strict full subsumption extension, obtaining $F_1 \in \Uclash$ with one new variable and $\delta(F_1) = \delta(F_0)$, and then take a clause in $F_1$ without that new variable and perform one step of strict full subsumption extension (on that new variable), obtaining $F_2 \in \Uclashns$ with $n(F_2) = n(F_1) = \odef(k)$ and $\delta(F_2) = \delta(F_1) + 1 \le k$. Now by Lemma \ref{lem:inv2subrescan}, further strict full subsumption extensions yield $F$ as desired. \Qed
\end{prf}

In the context of QMA (recall Subsection \ref{sec:prelimQCA}), Theorem \ref{thm:minnumvarmu} without the assertion on non-singularity corresponds to \cite[Corollary 3.6]{HKKZ2003Tight}, but non-singularity is not obtained by their examples (which are also using non-strict full subsumption extension, in the form of \cite[Lemma 3.4]{HKKZ2003Tight}).

\section{Non-Mersenne numbers}
\label{sec:nonmer}

In this section we study the function $\nonmer: \NN \ra \NN$ via a recursive definition (Definition \ref{def:minvdegdef}), where numerical values are given in Table \ref{tab:valuesminvdegdef}.

\begin{table}[h]
  \centering
  \begin{tabular}{c||c|*{3}{c}|*{3}{c}|*{3}{c}|*{3}{c}|c}
    $k$ & 1 & 2& 3 & 4 & 5 & $\cdots$ & 11 & 12 & $\cdots$ & $26$ & 27 & $\cdots$ & 57 & 58\\
    \hline
    $\nonmer(k)$ & 2 & 4 & 5 & 6 & 8 & $\cdots$ & 14 & 16 & $\cdots$ & 30 & 32 & $\cdots$ & 62 & 64
  \end{tabular}
  \caption{Values for $\nonmer(k)$, $k \in \tb 1 {58}$}
  \label{tab:valuesminvdegdef}
\end{table}

The understanding of this recursion is the underlying topic of this section. This recursion is naturally obtained from splitting on variables with minimum occurrence in minimally unsatisfiable clause-sets, and will be used in Theorem \ref{thm:MUminvdegdef} later to prove the upper bound on the minimum var-degree. In Subsection \ref{sec:nMBackground} we provide background, motivation and the definition, in Subsection \ref{sec:nmbasic} we show simple, basic properties, in Subsection \ref{sec:nonmerjump} we obtain the central combinatorial characterisation of the recursion, and finally in Subsection \ref{sec:nonmerappl} we obtain various closed formulas.

\subsection{Background}
\label{sec:nMBackground}

The sequence $\nonmer$ is sequence \url{http://oeis.org/A062289} in the ``On-Line Encyclopedia of Integer Sequences'':
\begin{itemize}
\item It can be defined as the enumeration of those natural numbers containing ``10'' in their binary representation; in other words, exactly the numbers whose binary representation contain only $1$'s are skipped.
\item Thus the sequence leaves out exactly the numbers of the form $2^n-1$ for $n \in \NN$ (that is, $1,3,7,15,31,\dots$; sequence \url{http://oeis.org/A000225}), which are often called \href{http://mathworld.wolfram.com/MersenneNumber.html}{``Mersenne numbers''}.
\item $\nonmer$ consists of arithmetic progressions of slope $1$ and length the Mersenne numbers, each such progression separated by an additional step of $+1$: these \textbf{blocks} of length $1,3,7,\dots$ are shown in Table \ref{tab:valuesminvdegdef} via the vertical bars.
\end{itemize}

The key deficiencies (values of $k$) in Table \ref{tab:valuesminvdegdef} are the following two classes:
\begin{enumerate}
\item The $k$-values $k = 1,2,5,12,27,58, \dots$ (sequence \url{http://oeis.org/A000325}) are the deficiencies $k = 2^n-n$ of the clause-sets $A_n$, $n \in \NN$, while the corresponding values $\nonmer(k) = 2^k$ are the minimum variable-degree of the clause-sets $A_n$ (see Lemma \ref{lem:defAn}), as explained in Subsection \ref{sec:introbasicint}.
\item The $k$-values $1,4,11,26,57, \dots$ (\url{http://oeis.org/A000295}) are the positions just before these deficiencies, as also discussed in Subsection \ref{sec:introbasicint}; we call them ``jump positions'', since precisely at these positions the function value increases by $2$ for the next argument (compare Definition \ref{def:jump}).
\end{enumerate}

The recursion in Definition \ref{def:minvdegdef} is new, and so we can not use these characterisations, but must directly prove the basic properties; indeed we give a complete self-contained account:

\begin{defi}\label{def:minvdegdef}
  For $k \in \NN$ let $\nonmer(k) := 2$ if $k = 1$, while else
  \begin{displaymath}
    \nonmer(k) := \max_{i \in \tb{2}{k}} \min \big (2 \cdot i, \, \nonmer(k-i+1) + i \big ).
  \end{displaymath}
\end{defi}
The intuition underlying Definition \ref{def:minvdegdef} of $\nonmer(k)$, as later unfolded in Theorem \ref{thm:MUminvdegdef}, is that we want to get an upper bound on the min-var-degree of an $F \in \Musati{\delta=k}$ (recall Subsection \ref{sec:introfour}). We consider a variable $v \in \var(F)$ of minimum var-degree, consider $\ve \in \set{0, 1}$ such that the literal-degree $i$ of $\ol v$ resp.\ $v$ is maximal, and infer an upper bound on $\vdeg_F(v)$ from the two splitting results as follows. The index $i$ runs over the possible literal-degrees of $\ol v$ resp.\ $v$, and thus we have to maximise over it. Since $i$ is the maximum degree over both signs, we can take the minimum with $i + i = 2i$. In the splitting result $\pao v{\ve} * F$ the deficiency is reduced by $i-1$, since $i$ occurrences (i.e., clauses) and one variable are lost, and we apply recursively the lower bound $\nonmer(k-(i-1))$, where the $i$ cancelled occurrences have to be re-added.
\begin{examp}\label{exp:nM} We have $\nonmer(2) = \min(2 \cdot 2, \nonmer(2-2+1)+2) = \min(4,4) = 4$ and $\nonmer(3) = \max( \min(2 \cdot 2, \nonmer(3-2+1)+2), \min(2 \cdot 3, \nonmer(3-3+1)+3) ) = \max(\min(4,6), \min(6,5)) = 5$.
\end{examp}

An outline of our analysis of $\nonmer(k)$ is as follows: A basic insight is that we always have $\nonmer(k+1) - \nonmer(k) \in \set{1,2}$ (Lemma \ref{lem:stepNM}). In order to show that $\nonmer(k)$ is indeed the sequence as described above, it suffices thus to show that the blocks (given by the contiguous intervals of enumerated values) indeed have length $2^m-1$, $m=1,2,\dots$. First, to gain control over the index $i$ in Definition \ref{def:minvdegdef}, we introduce the index function $\inonmer(k)$, for which we have $\nonmer(k) = \nonmer(k-\inonmer(k)+1)+\inonmer(k)$ (Lemma \ref{lem:simprecnm}). Via the helper functions $i'(k) := k - \inonmer(k) + 1$ and $h(k) := \nonmer(i'(k))$ thus holds $\nonmer(k) = h(k) + \inonmer(k)$; this ``canonical partition'' of $\nonmer(k)$ into two nearly equal parts is of special importance for our applications. We then turn to the determination of the jump positions, the set $J = \set{k \in \NN: \nonmer(k+1) = \nonmer(k) + 2}$; so the blocks are the left-open right-closed intervals from one jump position to the next. The main combinatorial characterisation is given in Lemma \ref{lem:characjumpb}, which shows that simple local patterns characterise these jumps. In Corollary \ref{cor:undip} we then understand, that $i'(k)$ in general moves in the pattern ``repeat, increment, repeat, increment, ...'', while at a jump position we have a double repetition; see Table \ref{tab:auxfunc} for numerical values. From that we conclude in Theorem \ref{thm:characjumpc} first, that $i'$ maps a block to the previous block, and second, that indeed block $m$ has length $2^m-1$, since the preimage of $k$ under $i'$ has size $2$, except at a jump, where the preimage has size $3$, and so the length of block $m+1$, using the length $2^m-1$ of block $m$, is $2 \cdot (2^m-1) + 1 = 2^{m+1}-1$. We furthermore obtain $J = \set{2^m - m - 1 : m \in \NN_{\ge 2}}$. Closed formulas and special cases for $\nonmer(k)$ are derived in Subsection \ref{sec:nonmerappl}, and we conclude by two corollaries on $\inonmer(k)$. Later we will obtain two further alternative characterisations of $\nonmer$:
\begin{itemize}
\item A combinatorial characterisation is obtained in Corollary \ref{cor:leansharp}, where we will see that $\nonmer(k)$ for $k \in \NN$ is the maximal min-var-degree for lean clause-sets or variable-minimal unsatisfiable clause-sets with deficiency $k$.
\item In Subsection \ref{sec:recpdp} we will develop a general recursion scheme, which has the function $\nonmer$ ``built-in'', as shown in Theorem \ref{thm:altrecnonmersenne}, exploiting the ``canonical partition'' $\nonmer(k) = h(k) + \inonmer(k)$.
\end{itemize}
The importance of the partition $\nonmer(k) = h(k) + \inonmer(k)$ comes from its meaning in the proof of Theorem \ref{thm:MUminvdegdef}, as explained above: $\inonmer(k)$ stands for the degree of the literal set to true, while $i'(k)$ stands for the deficiency of the clause-set after setting this literal to true, and thus $h(k)$ is the upper bound on the minimum var-degree from the induction hypothesis.

\subsection{Basic properties}
\label{sec:nmbasic}

We begin our investigations into $\nonmer(k)$ by some simple bounds:
\begin{lem}\label{lem:nmsimpleb}
  Consider $k \in \NN$.
  \begin{enumerate}
  \item\label{lem:nmsimpleb1} $k+1 \le \nonmer(k) \le 2 \cdot k$ for $k \in \NN$.
  \item\label{lem:nmsimpleb2} For $k \ge 2$ we have $\nonmer(k) \ge 4$.
  \end{enumerate}
\end{lem}
\begin{prf}
The upper bound of Part \ref{lem:nmsimpleb1} follows directly from the definition (by the min-component $2 i$). The lower bounds follows by induction: $\nonmer(1) = 2 \ge 1 + 1$, while for $k > 1$ we have $\nonmer(k) \ge \min(2 k, \nonmer(k-k+1) + k) = \min(2 k, 2+k) = k+2$. Part \ref{lem:nmsimpleb1} follows by Part \ref{lem:nmsimpleb1} and $\nonmer(2) = 4$. \Qed
\end{prf}

A basic tool for investigating sequences is the Delta-operator, which measures the differences in values between two neighbouring arguments:
\begin{defi}\label{def:delta}
  For a sequence $a: I \ra \RR$, where $I \sse \ZZ$ is stable under increment ($n \in I \Ra n+1 \in I$), the sequence $\Delta a: I \ra \RR$ is defined for $k \in I$ by $\bmm{\Delta a(k)} := a(k+1) - a(k)$ (i.e., the step in the value of the sequence from $k$ to $k+1$).
\end{defi}
A few obvious properties of the Delta-operator are as follows:
\begin{enumerate}
\item $\Delta: \RR^{I} \ra \RR^{I}$ is linear: $\Delta(\lambda \cdot a + \mu \cdot b) = \lambda \cdot \Delta(a) + \mu \cdot \Delta(b)$.
\item $a \in \RR^{I}$ is constant iff $\Delta a = (0)$.
\item $a$ is increasing iff $\Delta a \ge 0$, while $a$ is strictly increasing iff $\Delta a > 0$. Here for sequences $a, b: \RR^{I} \ra \RR^{I}$ of real numbers we use $a \le b :\Lra \fa\, n \in I : a_n \le b_n$, and $a < b :\Lra \fa\, n \in I : a_n < b_n$.
\end{enumerate}
The first key insight is, that the next number in the sequence of non-Mersenne numbers is obtained by adding $1$ or $2$ to the previous number:
\begin{lem}\label{lem:stepNM}
  For $k \in \NN$ holds $\Delta \nonmer(k) \in \set{1,2}$.
\end{lem}
\begin{prf}
  For $k=1$ we get $\Delta \nonmer(1) = 2$. Now consider $k \ge 2$. We have
  \begin{multline*}
    \nonmer(k+1) =
    \max(\min(4,\nonmer(k)+2), \max_{i \in \tb 3{k+1}} \min(2 i, \nonmer(k-i+2)+i)) =\\
    \max_{i \in \tb 3{k+1}} \min(2 i, \nonmer(k-i+2)+i) =\\
    \max_{i \in \tb 2k} \min(2 (i+1), \nonmer(k-(i+1)+2)+(i+1)) =\\
  \max_{i \in \tb 2k} \min(2 i + 2, \nonmer(k-i+1)+i+1) = 1 + \max_{i \in \tb 2k} \min(2 i + 1, \nonmer(k-i+1)+i).
  \end{multline*}
Thus on the one hand we have $\nonmer(k+1) \ge 1 + \max_{i \in \tb 2k} \min(2 i, \nonmer(k-i+1)+i) = 1 + \nonmer(k)$, and on the other hand $\nonmer(k+1) \le 1 + \max_{i \in \tb 2k} \min(2 i + 1, \nonmer(k-i+1)+i+1) = 2 + \nonmer(k)$. \Qed
\end{prf}

Thus increasing the deficiency $k$ by $1$ increases $\nonmer(k)$ at least by $1$:
\begin{corol}\label{cor:NMmon}
  $\nonmer: \NN \ra \NN$ is strictly increasing.
\end{corol}

And changing $\nonmer(a+b)$ to $\nonmer(a) + b$ can not increase the value:
\begin{corol}\label{cor:NMdiff}
  We have $\nonmer(a+b) \ge \nonmer(a) + b$ for $a \in \NN$ and $b \in \NNZ$, and thus $\nonmer(a - b) \le \nonmer(a) - b$ for $b < a$.
\end{corol}
\begin{prf}
We have $\nonmer(a+b) - \nonmer(a) = \sum_{i=0}^{b-1} \Delta \nonmer(a+i) \ge b \cdot 1$, whence the first inequality. Applying it yields $\nonmer(a-b) + b \le \nonmer(a-b+b) = \nonmer(a)$. \Qed
\end{prf}

Instead of considering the maximum over $k-1$ cases $i \in \tb 2k$ to compute $\nonmer(k)$ (according to Definition \ref{def:minvdegdef}), we can now simplify the recursion to only one case $\inonmer(k) \in \tb 2k$, and for that case also consideration of the minimum is dispensable. $\inonmer(k)$ is the first index $i$ in Definition \ref{def:minvdegdef}, where the minimum is attained by the $\nonmer$-term, that is, where $2 i \ge \nonmer(k-i+1) + i$:
\begin{defi}\label{def:i(k)}
  Let $\bmm{\inonmer}: \NN_{\ge 2} \ra \NN$ be defined for $k \in \NN$, $k \ge 2$, by $\inonmer(k) := i$ for the smallest $i \in \tb 2k$ with $i \ge \nonmer(k-i+1)$ (note that $k \ge \nonmer(k-k+1) = 2$, and thus $\inonmer(k)$ is well-defined).
\end{defi}

\begin{examp}\label{exp:i(k)}
  We have $\inonmer(2) = 2$ and $\inonmer(3) = 3$, since $\nonmer(3-2+1) = 4$, $\nonmer(3-3+1) = 2$.
\end{examp}

As promised, from $\inonmer(k)$ we can compute $\nonmer(k)$ by one recursive call of $\nonmer$:
\begin{lem}\label{lem:simprecnm}
  For $k \in \NN$, $k \ge 2$, we have:
  \begin{enumerate}
  \item\label{lem:simprecnm1} $0 \le \inonmer(k) - \nonmer(k-\inonmer(k)+1) \le 2$.
  \item\label{lem:simprecnm3} $\Delta \inonmer(k) \in \set{0,1}$.
  \item\label{lem:simprecnm2} $\nonmer(k) = \nonmer(k-\inonmer(k)+1)+\inonmer(k)$.
  \end{enumerate}
\end{lem}
\begin{prf}
For Part \ref{lem:simprecnm1} we consider the sequence $i \mapsto f_k(i) := i - \nonmer(k-i+1)$; this sequence starts with $f_k(2) = 2 - \nonmer(k-1) \le 0$, and finishes with $f_k(k) = k - \nonmer(1) \ge 2$, and $\inonmer(k)$ is the smallest $i$ with $f_k(i) \ge 0$. By Lemma \ref{lem:stepNM} we have $\Delta f_k(i) = \Delta i(i) - \Delta \nonmer(k-i+1)(i) \in \set{1+1,1+2} = \set{2,3}$. So for $\inonmer(k) - \nonmer(k-\inonmer(k)+1) = f_k(\inonmer(k))$ by definition we have $f_k(\inonmer(k)) \ge 0$, while $f_k(\inonmer(k)) \le 2$ due to $\Delta f_k(\inonmer(k)) \le 3$ (otherwise $\inonmer(k)$ wouldn't be minimal).

For Part \ref{lem:simprecnm3} we consider the sequence $k \mapsto g_i(k) := i - \nonmer(k-i+1)$. Again by Lemma \ref{lem:stepNM} we get $\Delta g_i(k) \in \set{-1,-2}$. It follows immediately $\Delta \inonmer(k) \ge 0$. Now assume $\Delta \inonmer(k) \ge 1$; thus $-2 \le g_{\inonmer(k)}(k+1) < 0$, whence, as shown before, $g_{\inonmer(k)+1}(k+1) \ge -2 + 2 = 0$, and thus $\Delta \inonmer(k) = 1$.

For Part \ref{lem:simprecnm2} we consider the sequence $i \mapsto h_k(i) := \nonmer(k-i+1) + i$; by Lemma \ref{lem:stepNM} we have $\Delta h_k(i) \in \set{-1+1,-2+1} = \set{0,-1}$. Thus, and by definition of $\inonmer(k)$, we get $\nonmer(k) = \max(2 \cdot 1, \dots, 2 \cdot (\inonmer(k)-1), h_k(\inonmer(k))) = \max(2 \inonmer(k)-2, h_k(\inonmer(k)))$. Finally $h_k(\inonmer(k)) \ge 2 \inonmer(k) - 2 \Lra \nonmer(k - \inonmer(k) + 1) + 2 \ge \inonmer(k)$, which holds by Part \ref{lem:simprecnm1}. \Qed
\end{prf}

We obtain an alternative, functional characterisation of $\inonmer(k)$:
\begin{corol}\label{cor:characinonmer}
  For $k \in \NN$, $k\ge 2$, the value $\inonmer(k)$ is the unique $i \in \tb 2k$ fulfilling the two inequalities $\nonmer(k-i+1) \le i \le \nonmer(k-i+2)$.
\end{corol}
\begin{prf}
As shown in the first part of the proof of Lemma \ref{lem:simprecnm}, the sequence $i \in \tb 1k \mapsto f_k(i) := i - \nonmer(k-i+1) \in \ZZ$ is strictly increasing, while the two inequalities are equivalent to $f_k(i) \ge 0$, $f_k(i-1) \le -1$, and so they determine $i$ as the smallest $i \in \tb 2k$ with $f_k(i) \ge 0$, which is the definition of $\inonmer(k)$. \Qed
\end{prf}

\begin{examp}
  For $k=3$ we have $2 = \nonmer(k-3+1) \le 3 \le \nonmer(k-3+2) = 4$, while for $k=4$ we have $2 = \nonmer(k-4+1) \le 4 \le \nonmer(k-4+2) = 4$.
\end{examp}

\subsection{Characterising the jumps}
\label{sec:nonmerjump}

After these preparations we are able to characterise the ``jump positions'', which are defined as those $k$ where the function $\nonmer$ increases by $2$:
\begin{defi}\label{def:jump}
  Let $\bmm{J} := \set{k \in \NN : \Delta \nonmer(k) = 2}$ be the set of \textbf{jump positions}.
\end{defi}
Thus $\Delta \nonmer(k) = 1$ iff $k \notin J$, and by Table \ref{tab:valuesminvdegdef} we see $J = \set{1,4,11,26,57,\dots}$. Note that $\nonmer(k) = 1 + k + \abs{\set{k' \in J : k' < k}}$. It is useful to define two auxiliary functions:
\begin{defi}\label{def:twoaux}
  Let $\bmm{i'}, \bmm{h}: \NN_{\ge 2} \ra \NN$ be defined by $i'(k) := k-\inonmer(k)+1 \in \NN$ for $k \in \NN$, $k \ge 2$, while $h(k) := \nonmer(i'(k)) \in \NN$.
\end{defi}
Some basic properties:
\begin{enumerate}
\item We have $\Delta i'(k) = 1 - \Delta \inonmer(k)$.
\item Thus by Lemma \ref{lem:simprecnm}, Part \ref{lem:simprecnm3}, holds $\Delta i'(k) \in \set{0,1}$.
\item By Lemma \ref{lem:simprecnm}, Part \ref{lem:simprecnm2}, we have $\nonmer(k) = h(k) + \inonmer(k)$.
\item Thus $\Delta h(k) = \Delta \nonmer(k) - \Delta \inonmer(k)$.
\item By Lemmas \ref{lem:stepNM} and \ref{lem:simprecnm}, Part \ref{lem:simprecnm3} we get $\Delta h(k) \in \set{0,1,2}$.
\item By Lemma \ref{lem:simprecnm}, Part \ref{lem:simprecnm1} we have $\inonmer(k) - h(k) \in \set{0,1,2}$.
\item By Corollary \ref{cor:characinonmer} we have $h(k) = \nonmer(i'(k)) \le \inonmer(k) \le \nonmer(i'(k)+1)$.
\end{enumerate}

It is instructive to consider initial values of the auxiliary functions in Table \ref{tab:auxfunc}.

\begin{table}[h]
  \centering
  \begin{tabular}[c]{c||c|c||c|c||c|c||c|c||c}
    $k$ & $\nonmer$ & $\Delta \nonmer$ & $\inonmer$ & $\Delta \inonmer$ & $i'$ & $\Delta i'$ & $h$ & $\Delta h$ & $\inonmer - h$\\
    \hline\hline
    1 & 2 & 2 & - & - & - & - & - & - & -\\
    \hline
    2 & 4 & 1 & 2 & 1 & 1 & 0 & 2 & 0 & 0\\
    3 & 5 & 1 & 3 & 1 & 1 & 0 & 2 & 0 & 1\\
    4 & 6 & 2 & 4 & 0 & 1 & 1 & 2 & 2 & 2\\
    \hline
    5 & 8 & 1 & 4 & 1 & 2 & 0 & 4 & 0 & 0\\
    6 & 9 & 1 & 5 & 0 & 2 & 1 & 4 & 1 & 1\\
    7 & 10 & 1 & 5 & 1 & 3 & 0 & 5 & 0 & 0\\
    8 & 11 & 1 & 6 & 0 & 3 & 1 & 5 & 1 & 1\\
    9 & 12 & 1 & 6 & 1 & 4 & 0 & 6 & 0 & 0\\
    10 & 13 & 1 & 7 & 1 & 4 & 0 & 6 & 0 & 1\\
    11 & 14 & 2 & 8 & 0 & 4 & 1 & 6 & 2 & 2\\
    \hline
    12 & 16 & 1 & 8 & 1 & 5 & 0 & 8 & 0 & 0
  \end{tabular}
  \caption{Values of auxiliary functions; underlined the jump positions}
  \label{tab:auxfunc}
\end{table}

First we show some further simple properties of the auxiliary functions:
\begin{lem}\label{lem:characjumpa}
  Consider $k \ge 2$.
  \begin{enumerate}
  \item\label{lem:characjumpa1} If $\Delta \inonmer(k) = 0$, then:
    \begin{enumerate}
    \item\label{lem:characjumpa1a} $\Delta \inonmer(k+1) = 1$.
    \item\label{lem:characjumpa1b} $\inonmer(k) - h(k) \in \set{1,2}$.
    \item\label{lem:characjumpa1c} $\inonmer(k+1) = h(k+1)$.
    \end{enumerate}
  \item\label{lem:characjumpa2n} $\Delta \inonmer(k) = 1 \Lra \Delta i'(k) = 0 \Lra \Delta h(k) = 0$.
  \item\label{lem:characjumpa3} If $\Delta \inonmer(k) = 1$, then:
    \begin{enumerate}
    \item\label{lem:characjumpa3a} $k \notin J$.
    \item\label{lem:characjumpa3b} $\inonmer(k) - h(k) \in \set{0,1}$.
    \end{enumerate}
  \end{enumerate}
\end{lem}
\begin{prf}
For Part \ref{lem:characjumpa1a} assume $\Delta \inonmer(k+1) = 0$ (and thus $\Delta i'(k+1) = 1$ due to $\Delta i' = 1 - \Delta \inonmer$). Because of $\Delta h = \Delta \nonmer - \Delta \inonmer$ we obtain $\Delta h(k+1) \ge 1$. Thus $\inonmer(k) = \inonmer(k+2) \ge h(k+2) \ge h(k+1) + 1 = \nonmer(i'(k+1)) + 1 = \nonmer(i'(k)+1) + 1$, contradicting $\inonmer(k) \le \nonmer(i'(k)+1)$. For the remainder of Part \ref{lem:characjumpa1} note $\Delta h(k) = \Delta \nonmer(k) \ge 1$.

For Part \ref{lem:characjumpa1b} note $\inonmer(k) = \inonmer(k+1) \ge h(k+1) \ge h(k) + 1$.

For Part \ref{lem:characjumpa1c} assume $\inonmer(k+1) > h(k+1)$. Thus $\inonmer(k) = \inonmer(k+1) \ge h(k+1)+1 \ge h(k)+2$, whence $\inonmer(k) = h(k)+2$. If we would have $\Delta h(k) = 2$, then $\inonmer(k) = \inonmer(k+1) > h(k+1) = h(k) + 2$; thus $h(k+1) = h(k)+1$. Now $\inonmer(k)= h(k) + 2 = h(k+1) + 1 = \nonmer(i'(k+1)) + 1 = \nonmer(i'(k)+1) + 1$, a contradiction.

Part \ref{lem:characjumpa2n} is obvious, and Part \ref{lem:characjumpa3a} follows. Finally, Part \ref{lem:characjumpa3b} follows by $\inonmer(k+1) \le h(k+1) + 2$ and $\inonmer(k+1) = \inonmer(k)+1$, while $h(k+1) = h(k)$ due to Part \ref{lem:characjumpa2n}, whence $\inonmer(k) \le h(k) + 1$. \Qed
\end{prf}

We obtain characterisations of the jump positions via the auxiliary functions:
\begin{lem}\label{lem:characjumpb}
  For $k \ge 2$ the following conditions are equivalent:
  \begin{enumerate}
  \item\label{lem:characjump3a} $k \in J$
  \item\label{lem:characjump3c} $\Delta h(k) = 2$
  \item\label{lem:characjump3d} $\inonmer(k) = h(k) + 2$
  \item\label{lem:characjump3g} $\Delta \inonmer(k-1) = 1$ and $\inonmer(k-1) = h(k-1) + 1$
  \item\label{lem:characjump3e} $\Delta \inonmer(k-2) = \Delta \inonmer(k-1) = 1$ (yielding various equivalent forms via Lemma \ref{lem:characjumpa}, Part \ref{lem:characjumpa2n}).
  \end{enumerate}
\end{lem}
\begin{prf}
Condition \ref{lem:characjump3a} implies Condition \ref{lem:characjump3c} due to $\Delta \inonmer(k) = 0$ in case of $k \in J$ by Lemma \ref{lem:characjumpa}, Part \ref{lem:characjumpa3a}. Condition \ref{lem:characjump3c} implies Condition \ref{lem:characjump3d}, since $\Delta h(k) = 2$ implies $\Delta \inonmer(k) = 0$, and so by Lemma \ref{lem:characjumpa}, Part \ref{lem:characjumpa1c} we have $\inonmer(k) = \inonmer(k+1) = h(k+1)$, while the assumption says $h(k+1) = h(k) + 2$. In turn Condition \ref{lem:characjump3d} implies Condition \ref{lem:characjump3a}, since by Lemma \ref{lem:characjumpa}, Part \ref{lem:characjumpa3b} we get $\Delta \inonmer(k) = 0$, and thus $\Delta \nonmer(k) = \Delta h(k)$, where in case of $\Delta h(k) \le 1$ we would have $h(k) + 2 = \inonmer(k) \le \nonmer(i'(k)+1) = \nonmer(i'(k+1)) = h(k+1) \le h(k)+1$. So now we can freely use the equivalence of these three conditions.

Condition \ref{lem:characjump3d} implies Condition \ref{lem:characjump3g}, since we have $\Delta \inonmer(k) = 0$, and thus $\Delta \inonmer(k-1) = 1$ with Lemma \ref{lem:characjumpa}, Part \ref{lem:characjumpa1a}, from which we furthermore get $\inonmer(k) = \inonmer(k-1) + 1$ and $h(k-1) = h(k)$, and so $\inonmer(k-1) = \inonmer(k) - 1 = h(k) + 1 = h(k-1) + 1$. Condition \ref{lem:characjump3g} implies Condition \ref{lem:characjump3e}, since in case of $\Delta \inonmer(k-2) = 0$ we had $\inonmer(k-1) = h(k-1)$ with Lemma \ref{lem:characjumpa}, Part \ref{lem:characjumpa1c}. In turn Condition \ref{lem:characjump3e} implies Condition \ref{lem:characjump3d}, since $\inonmer(k) = \inonmer(k-1) + 1 = \inonmer(k-2) + 2$, while $h(k) = h(k-1) = h(k-2)$, where by definition $\inonmer(k-2) \ge h(k-2)$ holds, whence $\inonmer(k) \ge h(k) + 2$, which implies $\inonmer(k) = h(k) + 2$. \Qed
\end{prf}

We understand now the shape of the four $\Delta$-sequences:
\begin{corol}\label{cor:characDeltas}
  By definition the sequence $(\Delta \nonmer(k))_{k \in \NN}$ is $1$ except at the jump positions $k$, where it is $2$. The other three $\Delta$-sequences are shaped as follows:
  \begin{enumerate}
  \item\label{cor:characDeltas1} The sequence $(\Delta \inonmer(k))_{k \in \NN, k \ge 2}$ consists of alternating $0,1$'s except the two positions $k-2, k-1$ before a jump position $k \in J$, where we have two consecutive $1$'s (while at the jump position we have $0$).
  \item\label{cor:characDeltas2} The sequence $(\Delta i'(k))_{k \in \NN, k \ge 2}$ consists of alternating $0,1$'s except two positions before a jump position $k$, where we have two consecutive $0$'s.
  \item\label{cor:characDeltas3} The sequence $(\Delta h(k))_{k \in \NN, k \ge 2}$ consists of alternating $0,1$'s except two positions before a jump position $k$, where we have two consecutive $0$'s, followed by a $2$ at the jump position $k$, which is followed by $0$.
  \end{enumerate}
\end{corol}
\begin{prf}
Part \ref{cor:characDeltas1}: By Lemma \ref{lem:characjumpa}, \ref{lem:characjumpa1a} we have $\Delta \inonmer(k) = 0 \Ra \Delta \inonmer(k+1) = 1$, while by Lemma \ref{lem:characjumpb}, Part \ref{lem:characjump3e} we have $\Delta \inonmer(k) = \Delta \inonmer(k+1) = 1 \Ra k+2 \in J$, and by Lemma \ref{lem:characjumpa}, Part \ref{lem:characjumpa3a} we have $k \in J \Ra \Delta \inonmer(k) = 0$.

Part \ref{cor:characDeltas2} follows from Part \ref{cor:characDeltas1} by $\Delta i' = 1 - \Delta \inonmer$.

Part \ref{cor:characDeltas3}: By Lemma \ref{lem:characjumpa}, Part \ref{lem:characjumpa2n} the $0's$ in the sequence $\Delta h$ are precisely the $1$'s in the sequence $\Delta \inonmer$, while a $0$ of $\Delta \inonmer$ translates into a $2$ precisely at the jump positions by Lemma \ref{lem:characjumpb}, Part \ref{lem:characjump3c}. The assertion follows now by Part \ref{cor:characDeltas1}. \Qed
\end{prf}

Especially instructive is understanding of the $i'$-sequence:
\begin{corol}\label{cor:undip}
  The $i'$-sequence $(i'(k))_{k \in \NN, k \ge 2}$ consists of doublets $m,m$ for consecutive $m=1,2,\dots,$, except for $k \in J \sm \set{1}$, where we have at positions $k-2,k-1,k$ a triplet $m, m, m$. These triplet-values occur exactly when $m \in J$.
\end{corol}
\begin{prf}
The doublet/triplet structure follows by Corollary \ref{cor:characDeltas}, Part \ref{cor:characDeltas2}. Now consider a triplet $i'(k-2) = i'(k-1) = i'(k) = m$ for $k \in J \sm \set{1}$ , $m \in \NN$. By definition we have $\Delta \nonmer(m) = \Delta h(k)$ (due to $h(k) = h(i'(k)) = \nonmer(m)$, $h(k+1) = \nonmer(i'(k+1)) = \nonmer(i'(k)+1) = \nonmer(m+1)$). By Lemma \ref{lem:characjumpb}, Part \ref{lem:characjump3c} we have thus have $\Delta \nonmer(m) = 2$, i.e., $m \in J$. The triplets do not leave out some jump-value in $J$, since for $m \in J$ and for the last position $k$ with $i'(k) = m$ we have $\Delta \nonmer(m) = \Delta h(k)$. \Qed
\end{prf}

We can finish now our determination of $\nonmer$, by showing that the blocks (as mentioned in Subsection \ref{sec:nMBackground}) indeed have length $2^m - 1$:
\begin{thm}\label{thm:characjumpc}
  We partition $\NN$ into \textbf{blocks} via left-open right-closed intervals from consecutive elements of $J$; so the first block is $\set{1}$, the second is $\set{2,3,4}$, and so forth. More precisely, let $k_m$ for $m \in \NN$ be the $m$th element of $J$, plus $k_0 := 0$, so that the $m$-th block for $m \in \NN$ is $\tb{k_{m-1}+1}{k_m}$. Furthermore we extend $i'$ to $i': \NN \ra \NNZ$ by $i'(1) := 0$, and let block $0$ be $\set{0}$.
  \begin{enumerate}
  \item\label{lem:characjump1} $i'$ maps the $m$-th block for $m \in \NN$ surjectively to the $(m-1)$-th block.
  \item\label{lem:characjumpc2} So especially $i'(k_m) = k_{m-1}$ for $m \in \NN$.
  \item\label{lem:characjumpc3} The length of block $m$ for $m \in \NN$ is $2^m-1$.
  \item\label{lem:characjumpc4} $k_m = 2^{m+1} - (m+1) - 1$ for $m \in \NNZ$.
  \end{enumerate}
\end{thm}
\begin{prf}
Parts \ref{lem:characjump1}, \ref{lem:characjumpc2} follows by induction and Corollary \ref{cor:undip}. Part \ref{lem:characjumpc3} follows by induction from Part \ref{lem:characjump1} and Corollary \ref{cor:undip}: the length of block $1$ is $1 = 2^1 - 1$, while the length of block $m+1$ is obtained by duplicating every element of block $m$ and adding $1$ (for the triplicated last element): $2 \cdot (2^m - 1) + 1 = 2^{m+1} - 1$. Part \ref{lem:characjumpc4} follows from Part \ref{lem:characjumpc3} and the observation, that by definition $k_m$ is the sum of length of blocks $1, \dots, m$: $k_m = \sum_{i=1}^m (2^i - 1) = 2^{m+1} - 2 - m$. \Qed
\end{prf}

\subsection{Applications}
\label{sec:nonmerappl}

Now the closed formula for $\nonmer(k)$ can be proven (recall $\fld(x) = \floor{\log_2(x)}$); by Theorem \ref{thm:characjumpc} we have now established that $\nonmer(k)$ is the function as discussed in Subsection \ref{sec:nMBackground}, and thus the following formula is known, but there appears to be no proof in the literature, and so for completeness we give the simple proof:
\begin{thm}\label{thm:solveN}
  For $k \in \NN$ holds $\nonmer(k) = k + \fld(k+1 + \fld(k+1))$. 
\end{thm}
\begin{prf} Let $g(k) := \fld(k+1 + \fld(k+1))$ and $f(k) := k + g(k)$ (so $\nonmer(k) = f(k)$ is to be shown, for $k \ge 1$). We have $f(1) = 1 + \fld(2 + \fld(2)) = 1 + \fld(3) = 2 = \nonmer(1)$. We will now prove that the function $g(k)$ changes values exactly at the transitions $k \mapsto k+1$ for $k \in J$, that is, for indices $k = k_m := 2^{m+1} - m - 2$ (using Theorem \ref{thm:characjumpc}, Part \ref{lem:characjumpc4}) with $m \in \NN$ we have $\Delta g(k_m) = 1$, while otherwise we have $\Delta g(k_m) = 0$, from which the assertion follows (by the definition of $J$).

We have $g(1) = 1$ and $g(2) = 2$. Now consider $m \in \NN$ and $k_m + 1 \le k \le k_{m+1}$. We show $g(k) = m+1$, which proves the claim. Note that $g(k)$ is monotonically increasing. Now $g(k) \ge g(k_m+1) = \floor{\ld(2^{m+1}-m + \floor{\ld(2^{m+1}-m)})} = \floor{\ld(2^{m+1}-m + m)} = m+1$ and $g(k) \le g(k_{m+1}) = \floor{\ld(2^{m+2}-m-2 + \floor{\ld(2^{m+2}-m-2)})} \le \floor{\ld(2^{m+2}-m-2 + m+1)} = \floor{\ld(2^{m+2}-1)} = m+1$. \Qed
\end{prf}

Using $\odef: \NN \ra \NNZ$ as introduced in Theorem \ref{thm:minnumvarmu}, we get:
\begin{corol}\label{cor:nmo}
  For $k \in \NN$ holds $\nonmer(k) = k - 1 + \odef(k+1)$.
\end{corol}

From Theorem \ref{thm:solveN} we obtain very precise bounds for $\nonmer(k)$:
\begin{corol}\label{cor:upperboundnonmer}
  $k + \fld(k+1) \le \nonmer(k) \le k+1+\fld(k)$ holds for $k \in \NN$.
\end{corol}
\begin{prf} The lower bound follows trivially. The upper bound holds (with equality) for $k \le 2$, so assume $k \ge 3$. We have to show $g(k) = \fld(k+1+\fld(k+1)) \le 1 + \fld(k)$, which follows from $\ld(k+1+ \fld(k+1)) \le 1 + \ld(k)$. Now $\ld(k+1+ \fld(k+1)) \le \ld(k+1+ \ld(k+1)) \le \ld(k+k) = 1 + \ld(k)$. \Qed
\end{prf}

Note that $(k+1+\fld(k)) - (k + \fld(k+1)) \in \set{0,1}$, where this difference is zero iff $k+1$ is a power of $2$. Finally we can prove the already mentioned characterisation, which motivates the terminology of ``non-Mersenne numbers'', namely that $(\nonmer(k))_{n \in \NN}$ enumerates $\NN \sm \set{2^n-1 : n \in \NN}$.\footnote{Note that we are not speaking of ``non-Mersenne \emph{primes}''.} For that we consider the positions directly after the jump positions, which by Theorem \ref{thm:characjumpc}, Part \ref{lem:characjumpc4} are the positions $2^n - n$ for $n \ge 2$. From that position on until the next jump position, which is $2^{n+1} - n - 2$, the $\nonmer$-values increase constantly by $1$ per step. So we just need to understand the values of $\nonmer(2^n-n)$, to understand all of $\nonmer$, which is achieved as follows (note that $(2^{n+1} - n - 2) - (2^n - n) = 2^n - 2$):
\begin{corol}\label{cor:nonmerfp}
  Consider $n \in \NN$, $k := 2^n - n$, and $m \in \NNZ$ with $m \le 2^n-1$.
  \begin{enumerate}
  \item\label{cor:nonmerfp1} $\nonmer(k) = 2^n$.
  \item\label{cor:nonmerfp2} More generally for $m < 2^n-1$ holds $\nonmer(k+m) = 2^n + m$.
  \item\label{cor:nonmerfp3} For $m = 2^n-1$ we have $k+m = 2^{n+1} - (n+1)$, and thus $\nonmer(k + m) = 2^{n+1}$.
  \end{enumerate}
\end{corol}
\begin{prf}
By Theorem \ref{thm:solveN}, Part \ref{cor:nonmerfp1} follows with $\nonmer(2^n-n) = 2^n-n + \fld(2^n-n+1+\fld(2^n-n+1)) = 2^n-n+\fld(2^n-n+1+(n-1)) = 2^n-n + \fld(2^n) = 2^n$. Part \ref{cor:nonmerfp2} follows by Theorem \ref{thm:characjumpc}, Part \ref{lem:characjumpc4}, and Part \ref{cor:nonmerfp3} follows by Part \ref{cor:nonmerfp1}. \Qed
\end{prf}

Besides $\nonmer(2^n - n) = 2^n$, the value at the jumps is useful to note:
\begin{corol}\label{cor:nonmbeforejump}
  For $n \in \NN$, $n \ge 2$, we have $\nonmer(2^n-n-1) = 2^n-2$.
\end{corol}

It is also useful to have simple formulas for the $\inonmer(k)$-values around the jump positions:
\begin{corol}\label{cor:inonmeraroundjump}
  For $n \in \NN$, $n \ge 3$ the values of $\inonmer(2^n - n + m)$ are as follows, using $p := 2^{n-1}$ (where for $m = -4$ we need $n \ge 4$):\vspace{1ex}

  \begin{tabular}[c]{c|ccccccccc}
    $m$ & $-4$ & $-3$ & $-2$ & $-1$ & $0$ & $1$ & $2$ & $3$ & $4$\\
    \hline
    $\inonmer$ & $p - 2$ & $p - 2$ & $p - 1$ & $p$ & $p$ & $p + 1$ & $p + 1$ & $p + 2$ & $p + 2$
  \end{tabular}
\end{corol}
\begin{prf}
We have $\inonmer(2^n - n) = 2^{n-1}$ by Corollary \ref{cor:characinonmer}:
\begin{eqnarray*}
  2^{n-1} & \ge & \nonmer(2^n - n - 2^{n-1} + 1) = \nonmer(2^{n-1} - (n-1)) = 2^{n-1}\\
  2^{n-1} & \le & \nonmer(2^n - n - 2^{n-1} + 2) = \nonmer(2^{n-1} - (n-1) + 1) = 2^{n-1} + 1.
\end{eqnarray*}
The remaining values follow by Corollary \ref{cor:characDeltas}, Part \ref{cor:characDeltas1}. \Qed
\end{prf}

 We conclude with an alternative characterisation of the jump-set $J$:
\begin{corol}\label{cor:compnmi}
  For $k \in \NN$ the following conditions are equivalent:
  \begin{enumerate}
  \item If $k \ge 2$, then $\nonmer(k) < 2 \cdot \inonmer(k) - 1$.
  \item If $k \ge 2$, then $\nonmer(k) = 2 \cdot \inonmer(k) - 2$.
  \item $k \in J$, that is, $k = 2^{m+1}-m-2$ for some $m \in \NN$.
  \end{enumerate}
\end{corol}
\begin{prf}
Due to $1 \in J$ we assume $k \ge 2$. If $k \in J$, then by Lemma \ref{lem:characjumpb}, Part \ref{lem:characjump3d}, we have $\nonmer(k) = 2 \cdot \inonmer(k) - 2 < 2 \cdot \inonmer(k) - 1$. And if $k \notin J$, then by the same lemma we have $\inonmer(k) \le h(k) + 1$, and thus $\nonmer(k) = h(k) + \inonmer(k) \ge 2 \cdot \inonmer(k) - 1$. \Qed
\end{prf}

So for $k \in \NN \sm J$ and $e_0, e_1 \in \ZZ$, $e_1 \ge e_0$, with $e_0 + e_1 = \nonmer(k)$ we get $e_1 \ge \inonmer(k)$, which will be used in the proof of Lemma \ref{lem:twovarmvindeg}.

\section{The upper bound for $\minnonmer$}
\label{sec:specialcaseMU}

In this central section we prove in Subsection \ref{sec:Applyingtherecursion} the upper bound on the minimum var-degree for $\Musati{\delta=k}$, providing first in Subsection \ref{sec:Controllingdeficiency} the tools for splitting (in order to apply induction), in a self-contained way.

\subsection{Controlling deficiency}
\label{sec:Controllingdeficiency}

We take up again the argumentation of the first three case studies from Subsection \ref{sec:introfour}, adding now the missing details. In a sense the main auxiliary lemma of this \Schrift{} is the following Lemma \ref{lem:auxminvardeg} on the deficiencies obtained when splitting a saturated minimally unsatisfiable clause-set. This receives its importance from the fact that every minimally unsatisfiable clause-set can be saturated (recall Subsection \ref{sec:Saturation}; this method was first applied in this context in \cite{Ku99dKo}).

Consider $F \in \ul{\Cls}$, $v \in \var(F)$, $\ve \in \set{0,1}$, and let $m_{\ve}$ be the number of clauses satisfied by $v \ra \ve$, that is, the clauses containing $\ol{v}$ for $\ve=0$ resp.\ the clauses containing $v$ for $v=1$. By \cite[Lemma 11.4, Part 3]{Kullmann2007ClausalFormZII} we have $\delta^*(\pao v{\ve} * F) = \max_{F' \le \pao v{\ve} * F} \delta(F') \le \delta(F) - \min(m_{\ve}, \surp(F)) + 1$, but to be self-contained we do not use this here; we also need a stronger dependency on $m_{\ve}$, and thus we consider only special variables $v$, starting with the assumption that $v$ is not pure in $F$. Furthermore we consider here only $F \in \Cls$, and thus, as explained for Lemma \ref{lem:auxSMUSAT}, we assume that there are no $C, D \in F$ with $C \subset D$ and $\var(C) \cup \set{v} = \var(D)$ (to avoid contractions). Now
\begin{multline}
  \label{eq:deltapa}
  \delta(\pao v{\ve} * F) = c(\pao v{\ve} * F) - n(\pao v{\ve} * F) =\\
  (c(F) - m_{\ve}) - (n(F) - \abs{\var(F) \sm \var(\pao v{\ve} * F)} =\\
  \delta(F) - m_{\ve} + 1 + r_{\ve},
\end{multline}
where $r_{\ve}$ is the number of variables in $F$ which occur only in the vanishing clauses (as counted by $m_{\ve}$); we note that $\abs{\var(F) \sm \var(\pao v{\ve} * F)} = 1 + r_{\ve}$ uses that $v$ occurs in both signs. We want $r_{\ve} = 0$, and this is guaranteed in Lemma \ref{lem:auxminvardeg} by using $v$ of minimal degree.

If now $F \in \Smusat$ holds, then $\pao{v}{\ve} * F \in \Musat$ by Lemma \ref{lem:auxSMUSAT}, Part \ref{lem:auxSMUSAT1}. Since this is an essential insight (which the first author got from \cite{FlRe94}), we prove this here for the sake of completeness, where w.l.o.g.\ $\ve = 1$ and $F' := \pao v1 * F$. If there were $C \in F'$ with $F' \sm \set{C} \in \Usat$, then we have to consider two cases:
\begin{itemize}
\item If $C \in F$, then we had the partial saturation $\saturate(F,C,v) \in \Usat$ (recall Definitions \ref{def:saturation}, \ref{def:partsaturation}, and note that $C \cup \set{v} \notin F$, using $F \in \Musat$), due to $\pao v0 * \saturate(F,C,v) = \pao v0 * F \in \Usat$ and $\pao v1 * \saturate(F,C,v) = F' \sm \set{C} \in \Usat$, contradicting saturation of $F$.
\item So we must have that $C' := C \cup \set{\ol{v}} \in F$. But now $G := F \sm \set{C'} \in \Usat$, due to $\pao v0 * G = \pao v0 * F \in \Usat$ and $\pao v1 * G = F' \in \Usat$, contradicting the minimality of $F$, which concludes the proof.
\end{itemize}

We have all ideas together for the main auxiliary lemma. Since we use variables of minimal degree many times, it is useful to have a notation for them:
\begin{defi}\label{def:varsetmvd}
  For $F \in \Cls$ let $\bmm{\varmvd(F)} \sse \var(F)$ be the set of variables of minimal degree, that is, $\varmvd(F) := \set{v \in \var(F) : \vdeg_F(v) = \minvdeg(F)}$.
\end{defi}
$\varmvd(F) \not= \es$ iff $n(F) > 0$, and $\varmvd(F) = \var(F)$ holds iff $F$ is variable-regular.

\begin{lem}\label{lem:auxminvardeg}
  Consider $F \in \Cls$ and $v \in \varmvd(F)$, and let $m_0 := \ldeg_F(\ol{v})$ and $m_1 := \ldeg_F(v)$. Assume $v$ is not pure in $F$ (i.e., $m_0, m_1 \ge 1$). Consider $\ve \in \set{0,1}$.
  \begin{enumerate}
  \item\label{lem:auxminvardeg1} $\var(\pao v{\ve} * F) = \var(F) \sm \set{v}$ and $n(\pao v{\ve} * F) = n(F) - 1$.
  \item\label{lem:auxminvardeg2} If there are no clauses $C, D \in F$ with $C \subset D$ and $\var(C) \cup \set{v} = \var(D)$, then
    $\delta(\pao v{\ve} * F) = \delta(F) - m_{\ve} + 1$, where $\delta(F) - m_{\ve} + 1 < \delta(F)$ iff $m_{\ve} \ge 2$.
  \item\label{lem:auxminvardeg3} If $F \in \Smusati{\delta=k}$ for $k \in \NN$, then $\pao v{\ve} * F \in \Musati{\delta=k-m_{\ve}+1}$, where $m_{\ve} \le k$.
  \end{enumerate}
\end{lem}
\begin{prf}
Part \ref{lem:auxminvardeg1} follows by the fact, that due to degree-minimality and non-pureness of $v$, no variable can have all its occurrences only in clauses containing $v$ resp.\ $\ol{v}$. Part \ref{lem:auxminvardeg2} follows by Part \ref{lem:auxminvardeg1} (and the above \eqref{eq:deltapa}). Part \ref{lem:auxminvardeg3} follows by Lemma \ref{lem:auxSMUSAT}, Part \ref{lem:auxSMUSAT1} (as proven above) together with Part \ref{lem:auxminvardeg2}, where $m_{\ve} \le k$ follows by the general fact $\delta(F) \ge 1$ for $F \in \Musat$. \Qed
\end{prf}

If in the situation of Lemma \ref{lem:auxminvardeg}, Part \ref{lem:auxminvardeg3} (the essential part for us) the value of $m_{\ve}$ is minimal, i.e., $m_{\ve} = 1$, then we have $\delta(\pao v{\ve} * F) = \delta(F) = k$, while if $m_{\ve}$ is maximal, i.e., $m_{\ve} = k$, then we have $\delta(\pao v{\ve} * F) = 1$. The deficiency is strictly decreased for both splitting results $\pao v{\ve} * F$ iff $v$ is nonsingular.

\begin{examp}\label{exp:splitting}
  For $\Dt{n} \in \Smusatnsi{\delta=2}$ ($n \ge 2$; recall Example \ref{exp:MU2}) and each $v \in \var(\Dt{n})$, $\ve \in \set{0,1}$, holds $\pao v{\ve} * \Dt{n} \in \Musati{\delta=1}$.
\end{examp}

To demonstrate the usage of Lemma \ref{lem:auxminvardeg}, Part \ref{lem:auxminvardeg3}, we show the simple proof of the basic bound on the minimum var-degree (as outlined in Subsection \ref{sec:introbasicdegree}):
\begin{lem}[{{\cite[Lemma C.2]{Ku99dKo}}}]\label{lem:uppbdldg}
  For $F \in \Musat \sm \set{\set{\bot}}$ there is $v \in \var(F)$ with $\ldeg_F(v), \ldeg_F(\ol{v}) \le \delta(F)$.
\end{lem}
\begin{prf}
Let $F' \in \Smusat$ be a saturation of $F$ (recall Definition \ref{def:partsaturation}); note $\delta(F') = \delta(F)$. For $v \in \varmvd(F')$ by Lemma \ref{lem:auxminvardeg}, Part \ref{lem:auxminvardeg3} holds $\ldeg_{F'}(v), \ldeg_{F'}(\ol{v}) \le \delta(F')$, showing the assertion by $\var(F') = \var(F)$ and $\fa\, x \in \lit(F) : \lit_F(x) \le \lit_{F'}(x)$. \Qed
\end{prf}

We can find a variable realising Lemma \ref{lem:uppbdldg} by a simple search; we now give an example that choosing a variable of minimal degree however is not sufficient:
\begin{examp}\label{exp:uppbdldg}
  In the situation of Lemma \ref{lem:uppbdldg}, for $\delta(F) \ge 2$ the choice of $v \in \varmvd(F)$ can not guarantee the literal-degree bound, as we see by considering $A_2$ and performing a $3$-singular extension (recall Definition \ref{def:singext}), obtaining $F := \set{\set{3},\set{1,2,-3},\set{-1,2,-3},\set{1,-2,-3},\set{-1,-2}} \in \Musati{\delta=2}$; this can also be seen as a \emph{partial} singular unit-extension, similar to Definition \ref{def:extucl}, but not adding the new literal to all clauses. $F$ is variable-regular, every variable has degree $4$, but we have $\ldeg_F(-3) = 3 > \delta(F)$. In Example \ref{exp:uppbdldg2} we show that in general not even the existence of $v \in \varmvd(F)$ respecting the bounds is guaranteed.
\end{examp}

\subsection{Applying the recursion}
\label{sec:Applyingtherecursion}

The definition of $\nonmer(k)$ (recall Definition \ref{def:minvdegdef}) matches the recursion-structure of Lemma \ref{lem:auxminvardeg}, and we obtain an upper bound on the min-var-degree for minimally unsatisfiable clause-sets:
\begin{thm}\label{thm:MUminvdegdef}
  For all $k \in \NN$ and $F \in \Musati{\delta \le k} \sm \set{\set{\bot}}$ we have $\minvdeg(F) \le \nonmer(k)$.
\end{thm}
\begin{prf}
For $k=1$ the assertion (first shown in \cite{DDK98} for this case) follows by Lemma \ref{lem:uppbdldg}; so assume $k \ge 2$ (which yields $n(F) > 1$).\footnote{Indeed also for $k=2$ the assertion follows by Lemma \ref{lem:uppbdldg}.} We apply induction on $k$. Assume $\delta(F) = k$ (recall that $\nonmer$ is monotonically increasing). Saturate $F$ and obtain $F'$. Consider $v \in \varmvd(F')$. If $\vdeg_{F'}(v) = 2$ then we are done, so assume $\vdeg_{F'}(v) \ge 3$. Let $i := \max(\ldeg_{F'}(v),\ldeg_{F'}(\ol{v}))$; so $i \ge 2$ and $\vdeg_{F'}(v) \le 2 i$. W.l.o.g.\ assume that $i = \ldeg_{F'}(v)$. Consider $G := \pao v1 * F'$; by Lemma \ref{lem:auxminvardeg}, Part \ref{lem:auxminvardeg3} we have $G \in \Musat$ with $\delta(G) = k-i+1$ and $i \le k$. By the induction hypothesis we obtain a variable $w \in \var(G)$ with $\vdeg_G(w) \le \nonmer(k-i+1)$. By definition we have $\vdeg_{F'}(w) \le \vdeg_G(w) + \ldeg_{F'}(v)$. Altogether $\minvdeg(F) \le \min(2 i, \nonmer(k-i+1) + i) \le \nonmer(k)$. \Qed
\end{prf}

We can choose any $v \in \varmvd(F')$ of the saturation $F'$ in the proof of Theorem \ref{thm:MUminvdegdef} to realise $\vdeg_F(v) \le \nonmer(k)$, and thus by the proof of Lemma \ref{lem:uppbdldg}:
\begin{corol}\label{cor:MUminvdegdef}
  For $k \in \NN$ and $F \in \Musati{\delta \le k} \sm \set{\set{\bot}}$ there is $v \in \var(F)$ with $\vdeg_F(v) \le \nonmer(k)$ and $\ldeg_F(v), \ldeg_F(\ol{v}) \le k$.
\end{corol}

\begin{examp}\label{exp:uppbdldg2}
  Continuing Example \ref{exp:uppbdldg}, by definition we can realise the existence-statement in Theorem \ref{thm:MUminvdegdef} by choosing any $v \in \varmvd(F)$; we will now see that this guarantees $\ldeg_F(v), \ldeg_F(\ol{v}) \le k$ (as in Corollary \ref{cor:MUminvdegdef}) at least for some $v \in \varmvd(F)$ \emph{only for} $k \le 4$. If we have $\minvdeg(F) = \nonmer(k)$, then for every saturation $F'$ of $F$ we have $\minvdeg(F') = \minvdeg(F)$ (while in general only $\minvdeg(F') \ge \minvdeg(F)$ holds), and thus there is $v \in \varmvd(F)$ as required. So in order to find counter-examples, we have to consider the case $\minvdeg(F) < \nonmer(k)$. And indeed for $k \le 4$ all $v \in \varmvd(F)$ work as required, since if e.g.\ $\ldeg_F(v) > k$, then $\vdeg_F(v) \ge k+2 \ge \nonmer(k)$. However for $k=5$ we can consider $A_3$ and perform a singular $6$-extension via variable $4$, as a partial singular unit-extension, obtaining $F \in \Musat$ with $\delta(F) = 2^3 - 3 = 5$, $\minvdeg(F) = 1+6 = 7$, where $\varmvd(F) = \set{4}$ and $\ldeg_F(-4) = 6 > 5$.
\end{examp}

The upper bound on the minimum variable-degree of Theorem \ref{thm:MUminvdegdef} is not sharp, and will be further investigated from Section \ref{sec:strengtheningbound} on. However the bound is attained for infinitely many deficiencies, and we show in Lemma \ref{lem:sharpjump} that the jump positions (the set $J$; recall Definition \ref{def:jump}) are such deficiencies. To investigate the remaining deficiencies, we now show that they always have at least two variables realising the bound (if the bound is attained at all); this will be used to prove Theorem \ref{thm:supminvdegk6} (the first deficiency where Theorem \ref{thm:MUminvdegdef} is not sharp). So we consider ``extremal'' $F \in \Musati{\delta=k}$ with $\minvdeg(F) = \nonmer(k)$, and we show that such extremal clause-sets have at least two different variables of minimal degree, if $k \notin J$.

\begin{lem}\label{lem:twovarmvindeg}
  Consider $k \in \NN$.
  \begin{enumerate}
  \item\label{lem:twovarmvindeg1} For $k \notin J$ and $F \in \Musati{\delta=k}$ with $\minvdeg(F) = \nonmer(k)$ we have $\abs{\varmvd(F)} \ge 2$.
  \item\label{lem:twovarmvindeg2} For $k \in J$ there is $F \in \Uclashi{\delta=k}$ with $\minvdeg(F) = \nonmer(k)$ and $\abs{\varmvd(F)} = 1$.
  \end{enumerate}
\end{lem}
\begin{prf}
First assume $k \notin J$; we have to show the existence of different $v, w \in \varmvd(F)$. W.l.o.g.\ $F$ is saturated. Consider $v \in \minvdeg(F)$. By Corollary \ref{cor:compnmi} we have $\nonmer(k) \ge 2 \cdot \inonmer(k) - 1$. Because of $\ldeg_F(v) + \ldeg_F(\ol{v}) = \nonmer(k)$ thus w.l.o.g.\ $e_1 := \ldeg_F(v) \ge \inonmer(k)$. Let $F' := \pao v1 * F$. So $\delta(F') = k - e_1 + 1$. Recall $\nonmer(k) = \nonmer(k-\inonmer(k)+1) + \inonmer(k)$ (Lemma \ref{lem:simprecnm}, Part \ref{lem:simprecnm2}), and thus $\nonmer(k) \ge \nonmer(k-e_1+1)+e_1$. Since $n(F) \ge 2$, we can consider $w \in \varmvd(F')$. We have $\vdeg_{F'}(w) \le \nonmer(k-e_1+1)$ and $\vdeg_F(w) = \vdeg_{F'}(w) + e_1$. Thus $w \in \varmvd(F)$.

Now assume $k \in J$, i.e., $k = 2^{m+1}-m-2$ for $m \ge 1$. For $k=1$ we have the example $\set{{1},{-1}}$, so assume $k \ge 2$. Thus $\nonmer(k) = 2^{m+1}-2$. We obtain an example from $A_{m+1}$ by performing one strict full subsumption resolution: The resolution variable occurs $2^{m+1}-2$ times, the other $m-1$ variables occur $2^{m+1}-1$ times. \Qed
\end{prf}

\section{The min-var-degree upper bound for $\Lean$}
\label{sec:leansurp}

In this section we prove Theorem \ref{thm:leanminvardeg}, the upper bound $\nonmer(k)$ on the min-var-degree for lean clause-sets of surplus $k$, and the sharpness of this bound for any class between $\Vmusat$ and $\Lean$ in Theorem \ref{thm:vmusharp}. The proof consists in lifting Theorem \ref{thm:MUminvdegdef} to the general case in Subsection \ref{sec:proofgencase}, while sharpness of the upper bound is considered in Subsection \ref{sec:leansharp}. As a preparation, in Subsection \ref{sec:classSED} we introduce the class $\Sed$ of clause-sets, where deficiency and surplus coincide; this class is also interesting in its own right. Theorem \ref{thm:characSEDMSAT} characterises the elements of $\Sed$ as those clause-sets, where after removal of all clauses containing any given variable we obtain a matching-satisfiable clause-set. It follows (Corollary \ref{cor:auxminvardegsigma}) that unsatisfiable elements of $\Sed$ are variable-minimally unsatisfiable.

\subsection{Clause-sets with extremal surplus}
\label{sec:classSED}

We consider the task of generalising Theorem \ref{thm:MUminvdegdef} to $F \in \Lean$. Consider an arbitrary (multi-)clause-set $F$. Consider a set of variables $\es \not= V \sse \var(F)$ realising the surplus of $F$, i.e., such that $\delta(F[V])$ is minimal (recall Definition \ref{def:surp}). If $F[V]$ would be satisfiable, then a satisfying assignment would give a non-trivial autarky for $F$. Assuming that $F$ is lean thus yields that $F[V]$ must be unsatisfiable. So there exists a minimally unsatisfiable $F' \sse F[V]$. If now $\var(F') \not= \var(F[V]) = V$ would be the case, then we would loose control over the deficiency of $F'$. But this can not happen, as we will show in Corollary \ref{cor:mleanvmu}, namely $F[V]$ here is \emph{variable-minimally unsatisfiable}. We conclude this from the fact that $F[V]$ has a special structure, namely it belongs to the following class of clause-sets with maximal surplus (relative to the deficiency), whose basic theory we develop first:
\begin{defi}\label{def:sed}
  Let the class $\bmm{\Sed} \subset \Cls$ (``\ul{s}urplus \ul{e}qual \ul{d}eficiency'') consist of those clause-sets $F \in \Cls$ with $\surp(F) = \delta(F)$.
\end{defi}

So for $F \in \ul{\Cls}$ with $n(F) > 0$ we have $F \in \ul{\Sed} = \set{F \in \ul{\Cls} : \surp(F) = \delta(F)}$ iff for all $\es \subset V \sse \var(F)$ holds $\delta(F[V]) \ge \delta(F)$. It seems the class $\Sed$ crosses the classes considered in this \Schrift{} in interesting extremal cases.

\begin{examp}\label{exp:sed}
  Some basic examples:
  \begin{enumerate}
  \item We have $\top \in \Sed$ and $\set{\bot} \notin \Sed$, and for every $F \in \Sed$ we have $\bot \notin F$.
  \item For $F := \set{\set{1},\set{2}}$ we have $\surp(F) = \delta(F) = 0$, and thus $F \in \Sed$. On the other hand, for the multi-clause-set $F' := \set{2*\set{1},\set{2}}$ we have $\delta(F') = 1$, while still $\surp(F') = 0$, and thus $F' \notin \ul{\Sed}$.
  \item For every $F \in \ul{\Cls}$ with $c(F) \le 2$ and $\bot \notin F$ holds $F \in \ul{\Sed}$.
  \item $A_n \in \Sed$ for $n \ge 1$ (Example \ref{exp:surp}).
  \item $\Musati{\delta=1} \sm \set{\set{\bot}} \subset \Sed$ (since for $F \in \Musat \sm \set{\set{\bot}}$ holds $\surp(F) \ge 1$).
  \item $\Dt{n} \in \Sed$ for $n \ge 2$ (Example \ref{exp:MU2}).
  \item For $F := \set{\set{1,2,3},\set{1,2,-3},\set{1,-2},\set{-1,2},\set{-1,-2}}$ we have $F \in \Musat$ with $\delta(F) = 2$, but $\surp(F) = 1$, and thus $F \notin \Sed$.
  \item In Definition \ref{def:Mlcr} we introduce the class $\Mlcr \subset \Sed \cap \Sat \cap \Mlean$, and Example \ref{exp:mlcr} shows elements of this class.
  \item See also Corollary \ref{cor:sedextrcls}, Example \ref{exp:smallsurp}, and Question \ref{que:surpminv}.
  \end{enumerate}
  Finally we note that $F \in \ul{\Sed}$ iff $F' \in \ul{\Sed}$, where $F'$ is the multi-clause-set obtained from $F$ by forgetting all signs of the literals, i.e., replacing clauses $C \in F$ by $\var(C)$ (since $\delta(F') = \delta(F)$ and $\surp(F') = \surp(F)$).
\end{examp}

The class $\ul{\Sed}$ is not invariant under multiplicities; consider a multi-clause-set $F$ and the underlying clause-set $F'$:
\begin{enumerate}
\item If $F' \in \Sed$, then in general we do not have $F \in \ul{\Sed}$ (Example \ref{exp:sed}).
\item However in general holds $F \in \ul{\Sed} \Ra F' \in \Sed$, since if $F' \notin \Sed$, then $\surp(F') < \delta(F')$, and adding a duplicated clause to a multi-clause-set increases $\delta$ by $+1$, while $\surp$ is at most increased by $+1$ (it may also stay unchanged).
\end{enumerate}
We proceed by showing in the following Lemma \ref{lem:eqsed}, that $F \in \ul{\Sed}$ iff for every strict subset $V$ of variables of $F$, the number of clauses of $F$ with all variables contained in $V$ is at most $\abs{V}$. This basic but fundamental characterisation also yields a stronger Corollary \ref{cor:stabsed} than the above ``$F \in \ul{\Sed} \Ra F' \in \Sed$''. We use in this subsection for $F \in \ul{\Cls}$ and $V \sse \Va$ the notation $F^V \in \ul{\Cls}$, which is the sub-multi-clause-set of $F$ (i.e., $F^V \le F$) consisting of all $C \in F$ with $\var(C) \sse V$ (with the same multiplicities), that is, for $C \in \Cl$ with $\var(C) \sse V$ we set $F^V(C) := F(C)$, while otherwise $F^V(C) := 0$. So we have $F^{\es} = \set{F(\bot) * \bot}$ and $F^{\var(F)} = F$, and for $V \sse V'$ and $F \le F'$ holds $F^V \le F'^{V'}$; fundamental is
\begin{displaymath}
  c(F[V]) + c(F^{\var(F) \sm V}) = c(F[\var(F) \sm V] + c(F^V) = c(F).
\end{displaymath}
We always have $\var(F^V) \sse V$, but also for $V \sse \var(F)$ we might have $\var(F^V) \subset V$, for example for $F := \set{\set{1,2},\set{1,3}}$ we have $F^{\set{2,3}} = \top$ (so $\var(F^{\set{2,3}}) = \es$). Obviously we have $F^{\var(F^V)} = F^V$.
\begin{lem}\label{lem:eqsed}
  For a multi-clause-set $F$ the following three conditions are equivalent:
  \begin{enumerate}
  \item[(i)] $F \in \ul{\Sed}$.
  \item[(ii)] $\fa\, \es \sse V \subset \var(F) : c(F^V) \le \abs{V}$.
  \item[(iii)] $\fa\, \es \sse V \subset \var(F) : \delta(F^V) \le 0$.
  \end{enumerate}
\end{lem}
\begin{prf}
First we show the equivalence of Conditions (i) and (ii). For $\es \subset V' \sse \var(F)$ let $V := \var(F) \sm V'$, which runs through all $\es \sse V \subset \var(F)$. We have $c(F[V']) = c(F) - c(F^V)$, and thus we get:
\begin{multline*}
  \delta(F[V']) \ge \delta(F) \Lra c(F[V']) - \abs{V'} \ge c(F) - n(F) \Lra\\
  c(F) - c(F^V) - \abs{V'} \ge c(F) - n(F) \Lra c(F^V) \le n(F) - \abs{V'} = \abs{V},
\end{multline*}
and thus indeed (i) and (ii) are equivalent. Condition (iii) implies trivially (ii), since $\delta(F^V) = c(F^V) - n(F^V) \ge c(F^V) - \abs{V}$. Finally assume that (ii) holds, and consider $\es \sse V \subset \var(F)$. For $V' := \var(F^V)$ by (ii): $c(F^{V'}) \le \abs{V'}$, where $c(F^{V'}) = c(F^V)$, and thus $\delta(F^V) = c(F^V) - n(F^V) = c(F^{V'}) - \abs{V'} \le 0$.
\Qed
\end{prf}

We obtain, that decreasing multiplicities in $F \in \ul{\Sed}$ does not leave this class (even if the multiplicity drops to zero):
\begin{corol}\label{cor:stabsed}
  For $F \in \ul{\Sed}$ and $F' \le F$ we have $F' \in \ul{\Sed}$.
\end{corol}
\begin{prf}
For all $\es \sse V \subset \var(F')$ by Lemma \ref{lem:eqsed} we have $c(F'^V) \le c(F^V) \le \abs{V}$, and thus $F' \in \ul{\Sed}$. \Qed
\end{prf}

Even stronger, we can characterise the elements of $\ul{\Sed}$ as clause-sets close to matching-satisfiable clause-sets:
\begin{thm}\label{thm:characSEDMSAT}
  For $F \in \ul{\Cls}$ holds $F \in \ul{\Sed}$ if and only if for all $v \in \var(F)$ holds $F^{\var(F) \sm \set{v}} \in \ul{\Msat}$.
\end{thm}
\begin{prf}
Recall that by \cite[Lemma 7.2]{Kullmann2007ClausalFormZII} we have $G \in \ul{\Msat}$ iff $\delta^*(G) = 0$ for $G \in \ul{\Cls}$, where $\delta^*(G) = \max_{G' \le G} \delta(G)$. If $F \in \ul{\Sed}$, but $F^{\var(F) \sm \set{v}} \notin \ul{\Msat}$ for some $v \in \var(F)$, then there would be $G \le F^{\var(F) \sm \set{v}}$ with $\delta(G) > 0$, violating Lemma \ref{lem:eqsed} by choosing $V := \var(G)$ (namely $\delta(F^V) \le 0$, but $\delta(G) = c(G) - n(G) \le c(F^V) - n(G) = \delta(F^V)$). And if $F \notin \ul{\Sed}$, then there is $V \subset \var(F)$ with $\delta(F^V) \ge 1$, where $F^V \le F^{\var(F) \sm \set{v}}$ for $v \in \var(F) \sm V$. \Qed
\end{prf}

In Example \ref{exp:sed} we have seen satisfiable as well as unsatisfiable elements of $\ul{\Sed}$, and for the satisfiable elements we have matching-satisfiable as well as matching-unsatisfiable ones. Recalling Lemma \ref{lem:characVMUtrivial}, where we characterised the elements of $\Vmusat$ as those $F \in \Usat$, such that for all $v \in \var(F)$ holds $F^{\var(F) \sm \set{v}} \in \Sat$, we see that the unsatisfiable elements of $\Sed$ are variable-MU:
\begin{corol}\label{cor:auxminvardegsigma}
  $\Sed \cap \Usat \subset \Vmusat$ (and thus $\ul{\Sed} \cap \ul{\Usat} \subset \ul{\Vmusat}$).
\end{corol}
\begin{prf}
$\Sed \cap \Usat \subset \Vmusat$ follows from Theorem \ref{thm:characSEDMSAT} by Lemma \ref{lem:characVMUtrivial}. An example of $F \in \Musat \sm \Sed$ is given in Example \ref{exp:sed}. And for $F \in \ul{\Sed} \cap \ul{\Usat}$ and for the underlying clause-set $F'$ holds $F' \in \Sed \cap \Usat$, whence $F' \in \Vmusat$, and thus $F \in \ul{\Vmusat}$. \Qed
\end{prf}

By definition follows that for $\es \subset V \sse \var(F)$ with minimal $\delta(F[V])$ we have $F[V] \in \Sed$, which is a central insight for generalising the upper bound in the subsequent Subsection \ref{sec:proofgencase}:
\begin{corol}\label{cor:mleanvmu}
  Consider $F \in \ul{\Cls}$, $n(F) > 0$, and $\es \subset V \sse \var(F)$ with $\delta(F[V]) = \surp(F)$. Then we have:
  \begin{enumerate}
  \item\label{cor:mleanvmu1} $F[V] \in \ul{\Sed}$.
  \item\label{cor:mleanvmu2} If $F \in \Lean$, then $F[V] \in \Vmusat$.
  \end{enumerate}
\end{corol}
\begin{prf}
Part \ref{cor:mleanvmu1} holds since $F[V]$ has maximal deficiency amongst all $F[V']$ for $V' \ne \es$, and for $V' \sse V$ holds $F[V][V'] = F[V']$. Part \ref{cor:mleanvmu2} follows by Part \ref{cor:mleanvmu1}, Corollary \ref{cor:auxminvardegsigma}, and the unsatisfiability of $F[V]$ (recall that a satisfying assignment for $F[V]$, when restricted to $V$, is an autarky for $F$). \Qed
\end{prf}

In Part \ref{cor:mleanvmu2} of Corollary \ref{cor:mleanvmu} we use $\Lean$, $\Vmusat$ instead of $\ul{\Lean}$, $\ul{\Vmusat}$, since these classes are invariant under multiplicities (different from $\Sed$).

We conclude this subsection by considering the complexity of SAT decision for $F \in \Sed_{\delta=k}$ for fixed parameter $k \in \NN$. By Corollary \ref{cor:auxminvardegsigma} we could use Theorem \ref{thm:vmusatk}, however we can do better, namely we have $\Sed_{\delta=k} \subset \Mlean$ (due to $\surp(F) = k \ge 1$ for $F \in \Sed_{\delta=k}$), and so we can apply the fpt-result discussed in Example \ref{exp:MLEANfpt}, obtaining that SAT decision for inputs in $\Sed_{\delta=k}$ is fpt in $k$. The question is whether we actually have overall polytime SAT decision:
\begin{quest}\label{que:SED}
  Can SAT decision for $\Sed$ be done in polynomial time? If so, can we also find a satisfying assignment quickly?
\end{quest}
In Subsection \ref{sec:findaut} we will see a (potential) application of Question \ref{que:SED}.

\subsection{The generalised upper bound}
\label{sec:proofgencase}

Back to the main task, we first show that from $F[V] \in \Vmusat$ we obtain variables of low degree for $F$ itself:
\begin{lem}\label{lem:unsatrestr}
  Consider $F \in \ul{\Cls}$ and $\es \subset V \sse \var(F)$ such that $F[V] \in \Vmusat$. Then there exists  $v \in V$ with $\vdeg_F(v) \le \nonmer(\delta(F[V]))$ and $\ldeg_F(v), \ldeg_F(\ol{v}) \le \delta(F[V])$.
\end{lem}
\begin{prf}
Let $F' := F[V]$ and consider some minimally unsatisfiable $F'' \sse F'$; by assumption we have $\var(F'') = \var(F')$. So we get $\delta(F'') = \delta(F') - (c(F') - c(F''))$. By Corollary \ref{cor:MUminvdegdef} there is $v \in \var(F'')$ with $\vdeg_{F''}(v) \le \nonmer(\delta(F'')) = \nonmer(\delta(F') - (c(F') - c(F''))) \le \nonmer(\delta(F')) - (c(F') - c(F''))$ and $\ldeg_{F''}(v), \ldeg_{F''}(\ol{v}) \le \delta(F'') = \delta(F') - (c(F') - c(F''))$. Finally we have $\vdeg_F(v) \le \vdeg_{F''}(v) + (c(F') - c(F''))$ (note all occurrences of $v$ in $F$ are also in $F'$), and similarly for the literal-degrees. \Qed
\end{prf}

We obtain the generalisation and strengthening of Theorem \ref{thm:MUminvdegdef}:
\begin{thm}\label{thm:leanminvardeg}
  We have $\minvdeg(F) \le \nonmer(\surp(F))$ for a lean multi-clause-set $F$ with $n(F) > 0$. More precisely, there exists a variable $v \in \var(F)$ with $\vdeg_F(v) \le \nonmer(\surp(F))$ and $\ldeg_F(v), \ldeg_F(\ol{v}) \le \surp(F)$.
\end{thm}
\begin{prf}
By Corollary \ref{cor:mleanvmu}, Part \ref{cor:mleanvmu2}, and  Lemma \ref{lem:unsatrestr}. \Qed
\end{prf}

We recall $\surp(F) \le \minvdeg(F) - 1$ for arbitrary $F$, and thus we get:
\begin{corol}\label{cor:leanminvardegdef}
  For a lean multi-clause-set $F$, $n(F) > 0$, we have
  \begin{align*}
    \surp(F) + 1 \le \minvdeg(F) \le \nonmer(\surp(F)) \le \surp(F) + 1 + \fld(\surp(F))\\
    \minvdeg(F) \le \nonmer(\delta(F)) \le \delta(F) + 1 + \fld(\delta(F)).
  \end{align*}
\end{corol}

The bounds from Corollary \ref{cor:leanminvardegdef} are sharp in general:
\begin{examp}\label{exp:sharpbmvd}
  First we consider any lean clause-set $F \ne \top$, and perform a non-strict full subsumption extension $F \leadsto F'$. Obviously $F'$ is lean as well (with $\delta(F') = \delta(F)$). We have $\minvdeg(F') = 2$ and $\surp(F') = 1$, and thus $2 = \surp(F') + 1 = \minvdeg(F') = \nonmer(\surp(F')) = \surp(F) + 1 + \fld(\surp(F))$, while $\delta(F')$ is unbounded. This construction will be taken up again in Lemma \ref{lem:Lean1NPc}

  Now we turn to the $\delta$-upper bounds. For $n \ge 2$ consider $A_n$. We have $\surp(A_n) = \delta(A_n) = 2^n - n$ by Example \ref{exp:surp}. Thus here the inequalities of Corollary \ref{cor:leanminvardegdef} are $2^n -n + 1 = \surp(A_n) + 1 < 2^n = \minvdeg(A_n) = \nonmer(\delta(A_n)) = \delta(A_n) + 1 + \fld(\delta(A_n))$ (using Corollary \ref{cor:nonmerfp}).
\end{examp}

\begin{corol}\label{cor:sedextrcls}
  For $F \in \Lean \sm \set{\top}$ with $\minvdeg(F) = \nonmer(\delta(F))$ holds $F \in \Sed$.
\end{corol}
\begin{prf}
If $\surp(F) < \delta(F)$, then $\minvdeg(F) \le \nonmer(\surp(F)) < \nonmer(\delta(F))$ (Corollary \ref{cor:NMmon}). \Qed
\end{prf}

\subsection{Sharpness of the bound for $\Vmusat$}
\label{sec:leansharp}

We now show that for every deficiency $k$ there are variable-minimally unsatisfiable clause-sets where the min-var-degree is $\nonmer(k)$ (strengthening Example \ref{exp:sharpbmvd}). The examples are obtained as follows:
\begin{lem}\label{lem:leansharpaux}
  For a clause-set $F \in \Cls$, $\bot \notin F$ and $n(F) > 0$, with at least one full clause, we construct $F' \in \Cls$, also containing a full clause, as follows:
  \begin{enumerate}
  \item Let $C$ be a full clause of $F$.
  \item Let $F''$ be a full singular unit-extension of $F$ (recall Definition \ref{def:extucl}).
  \item Let $F' := F'' \addcup \set{C}$.
  \end{enumerate}
  We have the following properties:
  \begin{enumerate}
  \item $\surp(F') = \surp(F) + 1$, $\delta(F') = \delta(F) + 1$, $\minvdeg(F') = \minvdeg(F) + 1$.
  \item $F \in \Sed \Ra F' \in \Sed$.
  \item $F \in \Usat \Ra F' \in \Usat$.
  \end{enumerate}
\end{lem}
\begin{prf}
With Lemma \ref{lem:fullsingue} we get $\surp(F'') = \surp(F)$, $\delta(F'') = \delta(F)$, $\minvdeg(F'') = \minvdeg(F)$. Obviously $\delta(F') = \delta(F'') + 1$. Let $\var(F'') \sm \var(F) = \set{v}$. To see $\minvdeg(F') = \minvdeg(F'') + 1$, we note that for $w \in \var(F)$ we have $\vdeg_{F'}(w) = \vdeg_{F''}(w) + 1$, while $\vdeg_{F'}(v) = c(F') - 1 = c(F'') \ge \vdeg_{F''}(w) + 1$.

To prove $\surp(F') = \surp(F'') + 1$, we consider $\es \subset V \sse \var(F') = \var(F'')$. If $v \notin V$, then $\delta(F'[V]) = \delta(F''[V]) + 1$, since $C$ is full for $F$. If $V = \set{v}$, then $\delta(F'[V]) = \delta(F''[V]) = c(F'') \ge \surp(F'') + 1$. Finally, if $V \supset \set{v}$, then $\delta(F'[V]) = c(F') - \abs{V} \ge \delta(F') = \delta(F'') + 1 \ge \surp(F'') + 1$.

The implication $F \in \Sed \Ra F' \in \Sed$ follows now by definition of $\Sed$, and $F \in \Usat \Ra F' \in \Usat$ is trivial. \Qed
\end{prf}

Thus the upper bound on the min-var-degree of LEANs is tight for VMUs:
\begin{thm}\label{thm:vmusharp}
  For $\Vmusat \cap \Sed \sse \mc{C} \sse \Lean$ and $ k \in \NN$: $\minvdeg(\mc{C}_{\delta=k}) = \nonmer(k)$.
\end{thm}
\begin{prf}
By Theorem \ref{thm:leanminvardeg} it remains to show the lower bound $\minvdeg(\Vmusati{\delta=k} \cap \Sed) \ge \nonmer(k)$. For deficiencies $k = 2^n-n$, $n \in \NN$ we have $\nonmer(k) = 2^n$, and thus $A_n$ serves as lower bound example (as shown in Example \ref{exp:sharpbmvd}), while until the next jump position we can use Lemma \ref{lem:leansharpaux} together with Corollary \ref{cor:auxminvardegsigma}, where due to Corollary \ref{cor:nonmerfp} in this range also $\nonmer$ increases only by $1$ for $k \leadsto k+1$. \Qed
\end{prf}

Using Lemma \ref{lem:vmusatlean}, we can now determine the min-var-degrees for the classes $\Lean$ and $\Vmusat$, separated into layers via deficiency or surplus:
\begin{corol}\label{cor:leansharp}
  For $k \in \NN$ holds $\nonmer(k) = \minvdeg(\Leani{\delta=k}) = \minvdeg(\Leani{\surp=k}) = \minvdeg(\Vmusati{\delta=k}) = \minvdeg(\Vmusati{\surp=k})$.
\end{corol}

\begin{examp}\label{exp:leansharp}
  The example $F$ for $k=6$, using the full clause $\set{1,2,3} \in A_3$ and variable $4$ for the unit-extension, consists of the unit-clause $\set{4}$, the $8$ clauses of $A_3$, each with added literal $-4$, and the clause $\set{1,2,3}$. We have $F \in \Vmusat$ with $\delta(F) = 1+8+1 - 4 = 6$ and $\minvdeg(F) = 8+1 = 9 = \nonmer(6)$.
\end{examp}

\section{Algorithmic implications}
\label{sec:algoimpl}

In Subsections \ref{sec:autred}, \ref{sec:findaut} we consider the algorithmic implications of Theorem \ref{thm:leanminvardeg} (the upper bound on the minimum var-degree for lean clause-sets). First in Theorem \ref{thm:genautred} we show that via a poly-time autarky-reduction every clause-set $F \in \Cls$ can be reduced to some $F' \sse F$, where $F'$ fulfils the min-var-degree upper bound of Theorem \ref{thm:leanminvardeg}, although $F'$ might not be lean. For this autarky-reduction we do not know whether we can efficiently compute a certificate, the autarky, and we discuss Conjecture \ref{con:findauthard}, that efficient computation is possible, in Subsection \ref{sec:findaut}. We conclude with some remarks on the surplus in Subsection \ref{sec:remsurp}.

\subsection{Autarky reduction}
\label{sec:autred}

By Theorem \ref{thm:leanminvardeg} lean clause-sets fulfil a condition on the minimum variable-degree --- if that condition is not fulfilled, then there exists an autarky. In this section we try to pinpoint these autarkies. We consider a vast generalisation of lean clause-sets, namely matching-lean clause-sets (recall Subsection \ref{sec:prelimAut}, especially that a multi-clause-set $F$ with $n(F) > 0$ is matching-lean iff $\surp(F) \ge 1$). Recall the observation that a (multi-)clause-set $F$ has a non-trivial autarky (is not lean) iff there is $\es \subset V \sse \var(F)$ such that $F[V]$ is satisfiable, where the corresponding autarky reduction of $F$ removes all clauses containing some variable of $V$; note that to perform this autarky reduction the autarky itself is not needed, only its set $V$ of variables.
First we obtain a sufficient criterion for the existence of a non-trivial autarky by considering the converse of Theorem \ref{thm:leanminvardeg}:
\begin{lem}\label{lem:charakappcor}
  Consider $F \in \ul{\Mlean}$ with $n(F) > 0$ and $\minvdeg(F) > \nonmer(\surp(F))$. For all $F' := F[V]$ with $\es \subset V \sse \var(F)$ and $\delta(F[V]) = \surp(F)$ we have:
  \begin{enumerate}
  \item\label{lem:charakappcor1} $\delta(F') = \surp(F') = \surp(F)$, and so $F' \in \ul{\Sed} \cap \ul{\Mlean}$.
  \item\label{lem:charakappcor2} $\minvdeg(F') > \nonmer(\surp(F'))$.
  \item\label{lem:charakappcor3} $F' \in \Sat$ (yielding an autarky $\vp$ for $F$ with $\var(\vp) = V$).
  \end{enumerate}
\end{lem}
\begin{prf} Part \ref{lem:charakappcor1} follows by Corollary \ref{cor:mleanvmu}, Part \ref{cor:mleanvmu1}. For Part \ref{lem:charakappcor2} note that $\minvdeg(F') \le \nonmer(\surp(F'))$ implies $\minvdeg(F) \le \minvdeg(F') \le \nonmer(\surp(F))$ contradicting the assumption. Finally for Part \ref{lem:charakappcor3} assume $F' \in \Usat$. Then by Part \ref{lem:charakappcor1} and Corollary \ref{cor:auxminvardegsigma} we obtain $F' \in \Vmusat$, which contradicts Part \ref{lem:charakappcor2} by Lemma \ref{lem:unsatrestr}. \Qed
\end{prf}

Some simple examples for Lemma \ref{lem:charakappcor}:
\begin{examp}\label{exp:charakappcor}
  For $F := \set{3*\set{1,2}}$ we have $\delta(F) = 1$ and $F \in \ul{\Sed} \cap \ul{\Mlean} \cap \ul{\Sat}$ with $\minvdeg(F) = 3 > \nonmer(1) = 2$; for $F'$ as in Lemma \ref{lem:charakappcor} we have $F' = F$. The same holds for $F := A_2 \sm \set{\set{1,2}}$. These are examples of the class $\ul{\Mlcr}$, as investigated later in Subsection \ref{sec:findaut} (and special cases of Example \ref{exp:mlcr}).
\end{examp}

To better understand the background, we recall two fundamental facts regarding the surplus $\surp(F)$ for multi-clause-set $F$ with $n(F) > 0$:
\begin{enumerate}
\item $\surp(F)$ together with some $\es \subset V \sse \var(F)$ with $\surp(F) = \delta(F[V])$ can be computed in polynomial time (see \cite[Subsection 11.1]{Kullmann2007ClausalFormZI}).
\item If $\surp(F) \le 0$, then one can compute a non-trivial matching autarky for $F$ in polynomial time (see \cite[Section 7]{Ku00f} or \cite[Section 9]{Kullmann2007ClausalFormZI}).
\end{enumerate}

We see now that we can reach the conclusion of Theorem \ref{thm:leanminvardeg} for arbitrary inputs $F$ in polynomial time, via some autarky reduction (maintaining satisfiability-equivalence), where for simplicity of formulation we consider only clause-sets (for multi-clause-set first the reduction to the underlying clause-set is performed):
\begin{thm}\label{thm:genautred}
  Consider a clause-set $F$. We can find in polynomial time a sub-clause-set $F' \sse F$ such that:
    \begin{enumerate}
    \item There exists an autarky $\vp$ for $F$ with $F' = \vp * F$.
    \item If $n(F') \ge 1$, then $\surp(F') \ge 1$ and $\minvdeg(F') \le \nonmer(\surp(F'))$.
    \end{enumerate}
\end{thm}
\begin{prf}
The reduction process for $F$, yielding the final $F'$, is a loop with two steps:
\begin{enumerate}
\item Eliminate matching-autarkies, i.e., reduce $F$ to its the matching-lean kernel, such that we obtain $\surp(F) \ge 1$.
\item If $\minvdeg(F) \le \nonmer(\surp(F))$, then stop, otherwise apply the autarky-reduction of Part \ref{lem:charakappcor3} of Lemma \ref{lem:charakappcor}, removing all clauses containing a variable of $V$.
\end{enumerate}
This loop is aborted if $n(F)=0$ is reached at any step. All autarkies are composed together (as shown in \cite{Ku98e} for the general case, the composition of autarkies is again an autarky; here the autarkies are variable-disjoint, and thus we can just take their union), yielding the final $\vp$. \Qed
\end{prf}

In Theorem \ref{thm:genautred} we can only show the \emph{existence} of an autarky $\vp$ for $F$ with $F' = \vp * F$, however we currently do not know how to compute it efficiently. We conjecture that it can be found in polynomial time:
\begin{conj}\label{con:findauthard}
  For $F \in \Cls_{\surp \ge 1}$ there is a poly-time algorithm for computing a non-trivial autarky $\vp$ for $F$ in case of $\minvdeg(F) > \nonmer(\surp(F))$.
\end{conj}
Note that we ask only to find \emph{some} autarky $\vp$, not necessarily one given by Lemma \ref{lem:charakappcor} (i.e., with $\var(\vp) = V$ as in Part \ref{lem:charakappcor3} of Lemma \ref{lem:charakappcor}). This is enough since the number of variables is reduced by such a reduction, and this by some autarky:

\begin{lem}\label{lem:findauthard}
  If Conjecture \ref{con:findauthard} is true, then for the algorithm from Theorem \ref{thm:genautred}, which reduces a clause-set $F$ to some (satisfiability-equivalent) $F' \sse F$, we can also compute an autarky $\vp$ for $F$ with $F' = \vp * F$ in polynomial time.
\end{lem}
\begin{prf}
In the loop as given in the proof of Theorem \ref{thm:genautred}, we can replace the autarky-reduction according to Part \ref{lem:charakappcor3} of Lemma \ref{lem:charakappcor} by the reduction $F \leadsto \vp * F$ according to a (non-trivial) autarky as given by Conjecture \ref{con:findauthard}. \Qed
\end{prf}

We now discuss what we know about Conjecture \ref{con:findauthard}.

\subsection{On finding the autarky}
\label{sec:findaut}

Consider a matching-lean multi-clause-set $F$ with $n(F) > 0$ (thus $\surp(F) \ge 1$), where Lemma \ref{lem:charakappcor} is applicable, that is, we have $\minvdeg(F) > \nonmer(\surp(F))$. So we know that $F$ has a non-trivial autarky. Conjecture \ref{con:findauthard} states that finding such a non-trivial autarky in this case can be done in polynomial time (recall that finding a non-trivial autarky in general is NP-complete, which was shown in \cite{Ku00f}).

The task of actually finding the autarky can be considered as finding a satisfying assignment for the following class $\ul{\Mlcr} \subset \ul{\Sat} \cap \ul{\Mlean}$ of satisfiable(!) multi-clause-sets $F$, obtained by considering all $F[V]$ for minimal sets of variables $V$ with $\delta(F[V]) = \surp(F)$ (where ``CR'' stands for ``critical''):
\begin{defi}\label{def:Mlcr}
  Let \bmm{\ul{\Mlcr}} be the class of multi-clause-sets $F$ fulfilling the following three conditions:
  \begin{enumerate}
  \item\label{def:Mlcr1} $F \in \ul{\Mlean}$, $\bot \notin F$, $F \not= \top$.
  \item\label{def:Mlcr2} For all $\es \subset V \subset \var(F)$ holds $\delta(F[V]) > \surp(F)$.
  \item\label{def:Mlcr3} $\minvdeg(F) > \nonmer(\surp(F))$.
  \end{enumerate}
  The definition of $\Mlcr$ just uses ``$F \in \Mlean$'' instead.
\end{defi}
The basic properties of this class are collected in the following lemma:
\begin{lem}\label{lem:bpmlcr}
  For $F \in \ul{\Mlcr}$ holds:
  \begin{enumerate}
  \item $\delta(F) = \surp(F) \ge 1$ (whence $F \in \ul{\Sed}$).
  \item $F \in \ul{\Sat}$.
  \end{enumerate}
\end{lem}
\begin{prf}
Since $F \in \ul{\Mlean}$ and $n(F) > 0$, we have $\surp(F) \ge 1$. By $\bot \notin F$ we get $F = F[\var(F)]$, and thus $\surp(F) = \delta(F[\var(F)]) = \delta(F)$, while $F \in \ul{\Sat}$ follows by Lemma \ref{lem:charakappcor}. \Qed
\end{prf}

The examples we know for elements of $\Mlcr$ are as follows:
\begin{examp}\label{exp:mlcr}
  First we consider $F := \set{m*\set{1}}$ for $m \in \NN$: $\delta(F) = \surp(F) = m-1$, while $\minvdeg(F) = m$, so the first condition in Definition \ref{def:Mlcr} is fulfilled iff $m \ge 2$, while the second condition is trivially fulfilled, but the third condition is never fulfilled, and thus always $F \notin \ul{\Mlcr}$ holds. So consider $F := \set{m*\set{1,2}}$ for $m \in \NN$: $\delta(F) = \surp(F) = m-2$, while $\minvdeg(F) = m$, so $F \in \ul{\Mlcr} \Lra m = 3$. This example shows that $\ul{\Mlcr}$ is not invariant under multiplicities --- both increasing and decreasing multiplicities can lead outside of $\ul{\Mlcr}$.

  A simple example for $F \in \Mlcr_{\delta=1} \cap \Clash$ is given by
  \begin{displaymath}
    F := \set{\set{1,2}, \set{-1,2,-3},\set{-2,3},\set{1,-2,-3}}.
  \end{displaymath}
  We have $\delta(F) = 4-3 = 1$ and $\minvdeg(F) = 3$; for $\surp(F) = 1$ and Condition \ref{def:Mlcr2} of Definition \ref{def:Mlcr} notice, that any two variables cover all four clauses, and thus the minimum of $\delta(F[V])$ is only attained for $V = \var(F)$; finally by $\minvdeg(F) = 3 > \nonmer(1) = 2$ we get $F \in \Mlcr$, while $F \in \Clash$ by definition (any two clauses have a clash).

  Another class of example is obtained by full clause-sets. Let $F$ be a full clause-set and $n := n(F)$, $m := c(F)$. Then $F \in \Mlcr$ iff $n < m < 2^n$:
  \begin{enumerate}
  \item We have $\delta(F) = m - n$ and thus $\delta(F) \ge 1 \Lra m > n$.
  \item Furthermore $F \in \Sat \Lra m < 2^n$.
  \item For $\es \subset V \subset \var(F)$ we have $\delta(F[V]) = m - \abs{V}$. Thus $\surp(F) = \delta(F)$, and Condition \ref{def:Mlcr2} of Definition \ref{def:Mlcr} is fulfilled.
  \item It remains to show the condition on the min-var-degree: We have $\minvdeg(F) = m$, while $\nonmer(\surp(F)) = \nonmer(m-n)$. By Theorem \ref{thm:solveN} we obtain $\nonmer(m-n) = m - n + \fld(m - n + 1 + \fld(m - n + 1))$. We obtain for $n \ge 1$:
    \begin{multline*}
      \minvdeg(F) > \nonmer(\surp(F)) \Lra m >  m - n + \fld(m - n + 1 + \fld(m - n + 1)) \Lra\\
      \fld(m - n + 1 + \fld(m - n + 1)) < n \Lla\\
      \fld(2^n - 1 - n + 1 + \fld(2^n - 1 - n + 1)) < n \Lra\\
      \fld(2^n - n + \fld(2^n - n)) = \fld(2^n - n + (n-1)) = \fld(2^n - 1) < n.
    \end{multline*}
  \end{enumerate}
  Finally we note that the class $\ul{\Mlcr}$ is invariant against changes of polarities of literal occurrences (if clauses become equal in this way, then their multiplicities have to be added), and thus for example replacing all clauses $C \in F \in \ul{\Mlcr}$ by their positive forms, $\var(C)$, we obtain a positive (no complementations occur) multi-clause-set $F' \in \ul{\Mlcr}$ (with $c(F') = c(F)$ and $n(F') = n(F)$).
\end{examp}
The importance of $\ul{\Mlcr}$ is, that it is sufficient to find a non-trivial autarky for this class of satisfiable clause-sets. In order to show this, we need to strengthen the polytime computation of $\surp(F)$:
\begin{lem}\label{lem:surpmin}
  For a multi-clause-set $F$ with $n(F) > 0$ we can compute in polynomial time a minimal subset $\es \subset V \sse \var(F)$ with $\delta(F[V]) = \surp(F)$.
\end{lem}
\begin{prf}
Let $V := \var(F)$. Check whether there is $v \in \var(F)$ with $\surp(F[V \sm \set{v}]) = \surp(F)$ --- if yes, then $V := V \sm \set{v}$ and repeat, if not, then $V$ is the desired result. \Qed
\end{prf}

We are ready to show that $\ul{\Mlcr}$ is really the ``critical class'' for the problem of finding the witness-autarky underlying the reduction $F \leadsto F'$ of Theorem \ref{thm:genautred}:
\begin{thm}\label{thm:HardMlcr}
  Consider $F \in \ul{\Cls}$ with $\surp(F) \ge 1$ and $\minvdeg(F) > \nonmer(\surp(F))$.
  \begin{enumerate}
  \item\label{thm:HardMlcr1} For every minimal subset $\es \subset V \sse \var(F)$ with $\delta(F[V]) = \surp(F)$ we have $F[V] \in \ul{\Mlcr}$.
  \item\label{thm:HardMlcr2} We can compute in polytime some $\es \subset V \sse \var(F)$ with $F[V] \in \ul{\Mlcr}$.
  \end{enumerate}
  So Conjecture \ref{con:findauthard} is equivalent to the statement, that finding a non-trivial autarky for clause-sets in $\ul{\Mlcr}$ can be achieved in polynomial time.
\end{thm}
\begin{prf}
Part \ref{thm:HardMlcr1} follows with Lemma \ref{lem:charakappcor}. Part \ref{thm:HardMlcr2} follows from Part \ref{thm:HardMlcr1} with Lemma \ref{lem:surpmin}. The final statement follows with Part \ref{thm:HardMlcr2} (note that every autarky for some $F[V]$ yields an autarky for $F$). \Qed
\end{prf}

Since $\ul{\Mlcr} \subset \ul{\Sed}$, if both questions of Question \ref{que:SED} have a yes-answer, then this would prove Conjecture \ref{con:findauthard}. Before concluding this section with remarks on the surplus, a note on the role of multi-clause-sets in Theorem \ref{thm:HardMlcr}: Since for $F \in \ul{\Sed}$ also the underlying clause-set fulfils $\ulcls(F) \in \Sed$, if we can find in polytime satisfying assignments for satisfiable elements of $\Sed$ (Question \ref{que:SED}), we don't need to consider multi-clause-sets, and can consider $\ulcls(F[V])$ in Theorem \ref{thm:HardMlcr}.  Otherwise however it might be necessary to find satisfying assignments for $\ul{\Mlcr}$. Even if we start with $F \in \Cls$, the restriction $F[V]$ as considered in Theorem \ref{thm:HardMlcr} in general is a multi-clause-set, and although we have $F[V] \in \ul{\Mlcr}$, in general we don't have $\ulcls(F[V]) \in \Mlcr$.

\subsection{Final remarks on the surplus}
\label{sec:remsurp}

It is instructive to investigate the precise relationship between minimum variable-degree and surplus for lean clause-sets, which by Corollary \ref{cor:leanminvardegdef} are indeed very close. Small values behave as follows:
\begin{lem}\label{lem:charsurp1}
  Consider $F \in \Lean$, $n(F) > 0$ (so $\surp(F) \ge 1$ and $\minvdeg(F) \ge 2$).
  \begin{enumerate}
  \item\label{lem:charsurp1a} $\surp(F) = 1$ holds if and only if $\minvdeg(F) = 2$ holds.
  \item\label{lem:charsurp1b} $\minvdeg(F) = 3$ implies $\surp(F) = 2$.
  \item\label{lem:charsurp1c} $\surp(F) = 2$ implies $\minvdeg(F) \in \set{3,4}$.
  \item\label{lem:charsurp1d} $\minvdeg(F) = 4$ implies $\surp(F) \in \set{2,3}$.
  \end{enumerate}
\end{lem}
\begin{prf} First consider Part \ref{lem:charsurp1a}. If $\surp(F) = 1$ (so $n(F) > 0$), then by Theorem \ref{thm:leanminvardeg} we have $\minvdeg(F) \le \nonmer(1) = 2$, while in case of $\minvdeg(F) = 1$ there would be a matching autarky for $F$. If on the other hand $\minvdeg(F) = 2$ holds, then by definition $\surp(F) \le 2 - 1 = 1$, while $\surp(F) \ge 1$ holds since $F$ is matching lean. For Part \ref{lem:charsurp1b} note that due to $\surp(F)+1 \le \minvdeg(F)$ we have $\surp(F) \le 2$, and then the assertion follows by Part \ref{lem:charsurp1a}; Part \ref{lem:charsurp1d} follows in the same way. Finally Part \ref{lem:charsurp1c} follows by Part \ref{lem:charsurp1a} and $\nonmer(2) = 4$. \Qed
\end{prf}

So for $F \in \Musatns \sm \set{\set{\bot}}$ we have $\surp(F) \ge 2$; for a general polytime reduction of arbitrary (also non-boolean) clause-sets achieving surplus at least $2$ see \cite[Lemma 11.9]{Kullmann2007ClausalFormZI}. We illustrate Lemma \ref{lem:charsurp1}, using $F$ with $\delta(F) = \surp(F)$:
\begin{examp}\label{exp:smallsurp}
  Examples for cases $\surp(F) \in \set{2,3}$ in Lemma \ref{lem:charsurp1}:
  \begin{enumerate}
  \item An example for $\minvdeg(F) = 4$ in Part \ref{lem:charsurp1c} with $F \in \Uclash \cap \Sed$ is given by $A_2$.
  \item For $\setb{\set{a,b,c},\set{\ol{a},b,c},\set{a,\ol{b},c},\set{\ol{a},\ol{b},c},\set{a,\ol{c}},\set{\ol{a},\ol{c}}} \in \Uclash \cap \Sed$ we have $\minvdeg(F) = 4$ and $\surp(F)=3$ (Part \ref{lem:charsurp1d}).
  \end{enumerate}
\end{examp}

\begin{quest}\label{que:surpminv}
  is there for every $k \in \NN$ an $F \in \Uclash \cap \Sed$ with $\surp(F) = k$ and $\minvdeg(F) = k+1$?
\end{quest}

As we have for $\Musat$ the levels $\Musati{\delta=k}$ for $k = 1, 2, \dots$, we can consider for $\Lean$ the levels $\Leani{\surp=k}$ for $k = 1, 2, \dots$. However, while the levels $\Musati{\delta=k}$ as well as $\Leani{\delta=k}$ all are decidable in polynomial time, already the first level $\Leani{\surp=1}$ is NP-complete:
\begin{lem}\label{lem:Lean1NPc}
  Consider the map $E: \Cls \ra \Cls$, which has $E(\top) := \top$, while otherwise for $F \in \Cls \sm \set{\top}$ it chooses (by some rule --- it doesn't matter) a clause $C \in F$ and a variable $v \in \Va \sm \var(F)$, and replaces $C$ by $C \addcup \set{v}, C \addcup \set{\ol{v}}$; in other words, an non-strict full subsumption extension $F \leadsto E(F)$ is performed, as in Example \ref{exp:sharpbmvd}. Then we have for $F \in \Cls$:
  \begin{enumerate}
  \item $F \in \Lean$ iff $E(F) \in \Lean$.
  \item $F \in \Musat$ iff $E(F) \in \Musat$.
  \item $\surp(F) \le 1$.
  \end{enumerate}
  Thus $\Leani{\surp=1}$ is coNP-complete, while $\Musati{\surp=1}$ is $D^P$-complete.
\end{lem}
\begin{prf}
The properties of the map $E$ are trivial. The completeness-properties follow with the coNP-completeness of $\Lean$ (\cite{Ku00f}) and the $D^P$-completeness of $\Musat$ (\cite{PW88}). \Qed
\end{prf}

With Lemma \ref{lem:Lean1NPc} we also get easy examples for minimally unsatisfiable clause-sets of surplus $= 1$ and arbitrary deficiency $\ge 1$.

\section{Matching lean clause-sets}
\label{sec:genboundml}

In this section, which concludes our considerations on generalisations (beyond $\Musat$), we consider the question whether Theorem \ref{thm:leanminvardeg} can incorporate non-lean clause-sets. We consider the large class $\Mlean$ of matching lean clause-sets, which is natural, since a basic property of $F \in \Musat$ used in the proof of Theorem \ref{thm:leanminvardeg} is $\delta(F) \ge 1$ for $F \not= \top$, and this actually holds for all $F \in \Mlean$. We will construct for arbitrary deficiency $k \in \NN$ and $K \in \NN$ clause-sets $F \in \Mlean$ of deficiency $k$, where every variable occurs positively at least $K$ times. Thus neither the upper bound $\max(\ldeg_F(v), \ldeg_F(\ol{v})) \le f(\delta(F))$ nor $\ldeg_F(v) + \ldeg_F(\ol{v}) = \vdeg_F(v) \le f(\delta(F))$ for some chosen variable $v$ and for any function $f$ does hold for $\Mlean$. 

\newcommand{\Fsf}{F_{\mr{sf}}}
An example $\Fsf \in \Clashi{\delta=1}$, $n(F) > 0$, with $\minldeg(F) := \min_{x \in \lit(F)} \ldeg_F(x) \ge 2$ (the minimum literal-degree; and thus $\minvdeg(F) \ge 4$) is given in \cite[Section 5]{Ku2003e}, a ``star-free'' clause-set, as discussed in Subsection \ref{sec:introbicliques} (furthermore between any two different clauses of $\Fsf$ there is exactly one conflict). In \cite[Subsection 9.3]{Kullmann2007ClausalFormZI} it is shown that $\Fsf \in \Mlean$ holds. ``Star-freeness'' in our context means, that there are no singular variables (occurring in one sign only once). The simpler construction of this section pushes the number of positive occurrences arbitrary high, but there are variables with only one negative occurrence (i.e., there are singular variables).

For a finite set $V$ of variables let $M(V) \sse A(V)$ be the full clause-set over $V$ containing all full clauses with at most one complementation; e.g.\ $M(\set{1,2}) = \set{\set{1,2},\set{-1,2},\set{1,-2}}$:
\begin{enumerate}
\item Obviously $n(M(V)) = \abs{V}$, $c(M(V)) = \abs{V} + 1$ and $\delta(M(V)) = 1$ holds.
\item We have already seen that $M(V) \in \Mlean$ (indeed $M(V) \in \Mlcr$, as shown in Example \ref{exp:mlcr}).
\item By definition we have $\ldeg_{M(V)}(v) = \abs{V}$ and $\ldeg_{M(V)}(\ol{v}) = 1$ for all $v \in V$.
\end{enumerate}

\begin{lem}\label{lem:exmpmlean}
  For $k \in \NN$ and $K \in \NN$ there is $F \in \Mleani{\delta=k}$ such that for all variables $v \in \var(F)$ we have $\ldeg_F(v) \ge K$, and with $F \in \Usat$ for $k \ge 2$.
\end{lem}
\begin{prf} For $k = 1$ we can set $F := M(\tb {v_1}{v_K})$; so assume $k \ge 2$. Consider any clause-set $G \in \Musati{\delta=k-1}$ with $n := n(G) \ge K$, and let $V := \var(G)$. Consider a disjoint copy of $V$, that is a set $V'$ of variables with $V' \cap V = \es$ and $\abs{V'} = \abs{V}$, and consider two enumerations of the clauses $M(V) = \set{C_1, \dots, C_{n+1}}$, $M(V') = \set{C_1', \dots, C_{n+1}'}$. Now
\begin{displaymath}
  F := G \addcup \setb{ C_i \addcup C_i' : i \in \tb{1}{n+1}} \in \Usat
\end{displaymath}
(with $\var(F) = V \addcup V'$) has no matching autarky: If $\vp$ is a matching autarky for $F$, then $\var(\vp) \cap V = \es$, since $G$ is matching lean, thus $\var(\vp) \cap V' = \es$, since $M(V')$ is matching lean, and thus $\vp$ must be trivial. Furthermore we have $n(F) = 2 n$ and $c(F) = c(G) + n + 1$, and thus $\delta(F) = c(G) + n + 1 - 2 n = \delta(G)+1 = k$. By definition for all variables $v \in \var(F)$ we have $\ldeg_F(v) \ge n$. \Qed
\end{prf}

For $k=1$ the examples of Lemma \ref{lem:exmpmlean} for $K \ge 3$ are necessarily satisfiable, since $\Mleani{\delta=1} \cap \Usat = \Musati{\delta=1}$. It remains the questions whether the singular variables can be eliminated:
\begin{quest}\label{que:exmpmlean1}
  Are there examples for deficiency $k \in \NN$ of $F \in \Mleani{\delta=k}$, $n(F) > 0$, with $\minldeg(F) \ge k+1$ ? The above mentioned star-free $\Fsf$ shows that this is the case for $k=1$. In general, by Theorem \ref{thm:leanminvardeg} for such examples we have $F \notin \Lean$. What about the stronger condition $\minldeg(F) \ge K$ for arbitrary $K \in \NN$ ?
\end{quest}

\section{Lower bounds for $\minnonmer$}
\label{sec:strengtheningbound}

We now return to minimally unsatisfiable clause-sets. By Theorem \ref{thm:MUminvdegdef} we have $\minnonmer(k) = \minvdeg(\Musati{\delta=k}) \le \nonmer(k)$ for all $k \in \NN$. The task of precisely determining $\minnonmer(k)$ for all $k$ seems a deep question, and is the subject of the remainder of this \Schrift{}. All our examples yielding lower bounds on $\minnonmer(k)$ are actually (unsatisfiable) hitting clause-sets, and thus we believe (recall Definition \ref{def:nonmerstar})
\begin{conj}\label{con:minvdeghit}
  For all $k \in \NN$ holds $\minnonmer(k) = \minnonmerh(k)$.
\end{conj}

We will see in Theorem \ref{thm:supminvdegk6} that $\minnonmer \ne \nonmer$. We believe that $\minnonmer$ is a highly complicated function, but the true values deviate only at most by one from $\nonmer$:

\begin{conj}\label{con:sharpness}
  For all $k \in \NN$ we have $\minnonmer(k) \ge \nonmer(k)-1$.
\end{conj}

\begin{examp}\label{exp:consharpness}
  Consider $F \in \Musati{\delta=k}$ with $\minvdeg(F) = \minnonmer(k)$. If $\minnonmer(k) = \nonmer(k)$, then $F \in \Sed$ (recall Corollary \ref{cor:sedextrcls}), and by Lemma \ref{lem:twovarmvindeg} for $k \notin J$ holds $\abs{\varmvd(F)} \ge 2$. While in case of $\minnonmer(k) = \nonmer(k)-1$ we obtain by Theorem \ref{thm:leanminvardeg}, that if $F \notin \Sed$, then $\surp(F) = k-1$ and $k-1 \notin J$.
\end{examp}

\begin{quest}\label{que:extrSED}
  Does $F \in \Sed$ hold for all $F \in \Musati{\delta=k}$ with $\minvdeg(F) = \minnonmer(k)$?
\end{quest}

Concerning numerical bounds, by Corollary \ref{cor:upperboundnonmer} we get:
\begin{lem}\label{lem:preciseboundminnonmercon}
  If Conjecture \ref{con:sharpness} holds, then $k - 1 + \fld(k+1) \le \minnonmer(k) \le k+1+\fld(k)$ holds for $k \in \NN$.
\end{lem}

Later in Lemma \ref{lem:minvdegvalid} we will see that $\minnonmer: \NN \ra \NN$ is monotonically increasing. In Theorem \ref{thm:imprupperbound} we will (implicitly) construct a correction function $\gamma_1: \NN \ra \set{0,1}$ such that $\minnonmer \le \nonmer - \gamma_1$, where we remark in the Conclusion (Section \ref{sec:open}) that also $\minnonmer \ne \nonmer - \gamma_1$ holds. Note that Conjecture \ref{con:sharpness} says that there exists $\gamma: \NN \ra \set{0,1}$ with $\minnonmer = \nonmer - \gamma$, while for every $\gamma: \NN \ra \set{0,1}$ the function $\nonmer - \gamma$ is still monotonically increasing (by Lemma \ref{lem:stepNM}), and is thus a possible candidate.

In Subsection \ref{sec:preciseval} we provide a general method for obtaining \emph{lower bounds}, via considering full clauses (while in Section \ref{sec:improveupbdMU} we turn to improved upper bounds):
\begin{defi}\label{def:numfcl}
  For a clause-set $F \in \Cls$ let $\bmm{\nfc(F)} \in \NNZ$ be the number of full clauses, that is $\nfc(F) := \abs{\set{C \in F : \var(C) = \var(F)}}$. And for a class $\mc{C} \sse \Cls$ of clause-sets we define $\bmm{\nfc(\mc{C})} := \set{\nfc(F) : F \in \mc{C}} \sse \NNZ$ as the set of all possible numbers of full clauses, while $\bmm{\maxnfc(\mc{C})} \in \nni$ is the supremum of $\nfc(\mc{C})$.
\end{defi}
The maximum possible number of full clauses is an interesting quantity:
\begin{defi}\label{def:maxsmar}
  For $k \in \NN$ let
  \begin{eqnarray*}
    \bmm{\maxsmar(k)} & := & \maxnfc(\Musati{\delta=k}) \in \NN\\
    \bmm{\maxsmarh(k)} & := & \maxnfc(\Uclashi{\delta=k}) \in \NN
  \end{eqnarray*}
  (these numbers are finite due to $\maxsmarh(k) \le \maxsmar(k) \le \minnonmer(k)$).
\end{defi}
According to our numerical investigations, $\maxsmar$ is very close to $\minnonmer$:
\begin{conj}\label{con:maxnfcminvdeg}
  For all $k \in \NN$ we have $\maxsmar(k) \ge \minnonmer(k) - 1$.
\end{conj}

The smallest deficiency $k$ with $\maxsmar(k) < \minnonmer(k)$ is $k=3$, as shown in Lemma \ref{lem:nfcdef123} (together with Theorem \ref{thm:supminvdegk6}). The stronger form ``$\maxsmar(k) \ge \nonmer(k) - 1$'' is refuted by $\minnonmer(14) = \nonmer(14) - 1 = 17$, as shown in \cite{KullmannZhao2014Sharper} (discussed in Subsection \ref{sec:concmnM}), from which $\maxsmar(14) \le 17 - 1 = 16$ follows, since, as we will show in Corollary \ref{cor:evenfc}, if $\maxsmar(k) = \minnonmer(k)$ holds, then $\maxsmar(k)$ must be even.

Regarding $\maxsmarh(k)$, the maximal number of full clauses for unsatisfiable hitting clause-sets with deficiency $k$ (which indeed is always even, as shown in \cite{KullmannZhao2015FullClauses}), we conjecture the difference $\maxsmar(k) - \maxsmarh(k)$ is unbounded; this conjecture follows from the conjecture $\maxsmarh = S_2$ together with Conjecture \ref{con:maxnfcminvdeg} (since the difference $\nonmer - S_2$ is unbounded), as discussed in Subsection \ref{sec:concSmar}.

We show $\maxsmarh(k) = \nonmer(k)$ for two infinite classes of deficiencies $k$ (Lemmas \ref{lem:sharpnesssimpldef}, \ref{lem:sharpjump}). The main point here could be considered as (just) the equalities $\minnonmer(k) = \nonmer(k)$, for which in these two cases the proofs don't needed to consider full clauses, and so the general method for computing lower bounds on $\maxsmar(k)$, with the beginnings developed in Subsection \ref{sec:fullcl}, is not applied for these two lemmas. However in future work we will employ this method more fully (see Subsection \ref{sec:concSmar}), and, more important for the \Schrift{} at hand, we need for the proof of Theorem \ref{thm:supminvdegk6} (that $\minnonmer(6) = \nonmer(6) - 1 = 8$) the fact $\maxsmar(3) \le 4$, shown in Lemma \ref{lem:nfcdef123}.

\subsection{Some precise values for $\minnonmer$}
\label{sec:preciseval}

A general lower-bound method for $\minnonmer$ is provided by the number $\nfc(F)$ of full clauses in a clause-set $F$. The supremum $\maxnfc(\Musati{\delta=k}) = \maxsmar(k)$ of this number over all elements of $\Musati{\delta=k}$ for fixed $k$ is an interesting quantity in its own right, but in this \Schrift{} we only touch on this subject, providing the bare minimum of information needed in our context. See Subsection \ref{sec:concSmar} for an outlook on the interesting properties of this quantity. Some simple examples:
\begin{examp}\label{exp:numfcl}
  $\nfc(\top) = 0$, $\nfc(\set{\bot}) = 1$, and $\nfc(\set{\set{1},\set{-1,2}}) = 1$. While $\nfc(\es) = \es$, thus $\maxnfc(\es) = 0$, and $\nfc(\Cls) = \NNZ$, thus $\maxnfc(\Cls) = +\infty$.
\end{examp}
By definition we have:
\begin{lem}\label{lem:nfcminvdeg}
   $\nfc(F) \le \minvdeg(F)$ holds for every $F \in \Cls$ (since every variable in $F$ has degree at least $\nfc(F)$), and thus $\maxnfc(\mc{C}) \le \minvdeg(\mc{C})$ for every $\mc{C} \sse \Cls$.
\end{lem}
We obtain that for lean clause-sets (especially minimally unsatisfiable clause-sets) of fixed deficiency the number of full clauses is bounded:
\begin{corol}\label{cor:maxnumfullcl}
  A lean clause-set of deficiency $k$ can have at most $\nonmer(k)$ many full clauses; i.e., for all $k \in \NN$ we have $\maxnfc(\Leani{\delta=k}) \le \nonmer(k)$.
\end{corol}

\begin{quest}\label{que:maxnumfullcl}
  Is the upper bound of Corollary \ref{cor:maxnumfullcl} sharp?
\end{quest}
The first test-case for Question \ref{que:maxnumfullcl} is deficiency $3$, where we will see in Corollary \ref{cor:maxsmar}, that $\maxsmar(3) = 4 = \nonmer(3) - 1$, and where for sharpness a lean clause-set $F$ with $\delta(F) = 3$ and $\nfc(F) = 5$ needed to be demonstrated.

Precise values for $\maxnfc(\Musati{\delta=k}) = \maxsmar(k)$ we show for two infinite classes of deficiencies. The simplest class is given by the deficiencies directly after the jumps (recall Subsection \ref{sec:nonmerjump}), the deficiencies of the $A_n$:
\begin{lem}\label{lem:sharpnesssimpldef}
  For $n \in \NN$ and $k := 2^n - n$ holds $\maxsmarh(k) = \maxsmar(k) = \minnonmerh(k) = \minnonmer(k) = \nonmer(k) = 2^n$.
\end{lem}
\begin{prf}
  We have $\maxnfc(A_n) = 2^n$ (recall Lemma \ref{lem:defAn}), and thus $\maxnfc(\Uclashi{\delta=k}) \ge 2^n$, while by Corollary \ref{cor:nonmerfp} we have $\nonmer(k) = 2^n$. \Qed
\end{prf}

Also for the jumps themselves the same conclusions hold, namely by Lemma \ref{lem:twovarmvindeg}, Part \ref{lem:twovarmvindeg2} (and the proof) we have:
\begin{lem}\label{lem:sharpjump}
  For all $k \in J$ holds $\maxsmarh(k) = \maxsmar(k) = \minnonmerh(k) = \minnonmer(k) = \nonmer(k)$.
\end{lem}

Note that for $k \in J$ there is $n \in \NN$, $n \ge 2$, with $k = 2^n-n-1$ and $\nonmer(k) = 2^n - 2$. The underlying method of Lemmas \ref{lem:sharpnesssimpldef}, \ref{lem:sharpjump} is simple (as explained in Subsection \ref{sec:introbasicint}): start with $A_n$ and apply strict full subsumption resolution to full clauses. Zero steps have been used in Lemma \ref{lem:sharpnesssimpldef}, one step in Lemma \ref{lem:sharpjump}, and one example for two steps will be seen in the proof of Theorem \ref{thm:def15} (clause-set $F_3$ there). The further development of this method we have to leave for future work:
\begin{quest}\label{que:sfsAn}
  Explore the application of strict full subsumption resolution to $A_n$ in order to obtain lower bounds for $\minnonmerh(k)$ (and also $\maxsmarh(k)$). (By Lemma \ref{lem:inv2subrescan} we have $A_n \tsubress F$ iff $F \in \Uclash$ with $\var(F) = \var(A_n) = \tb 1n$; we obtain the lower bounds $\minnonmerh(\delta(F)) \ge \minvdeg(F)$ and $\maxsmarh(\delta(F)) \ge \nfc(F)$.)
\end{quest}

\subsection{On the number of full clauses}
\label{sec:fullcl}

We have a special interest in those $F \in \Musat$ where the lower bound $\nfc(F)$ meets the upper bound $\minvdeg(F)$. In this case this number must be even, and we obtain another $F' \in \Musat$ by resolving on any variable realising the minimum variable-degree (where any such variable only occurs in full clauses here):
\begin{lem}\label{lem:evenfc}
  Consider $F \in \Musat$ with $\nfc(F) = \minvdeg(F)$. Then $\nfc(F)$ is even, and moreover for each $v \in \varmvd(F)$ the set of full clauses of $F$ is partitioned into full-subsumption-resolvable pairs with resolution variable $v$.
\end{lem}
\begin{prf}
The occurrences of $v$ are exactly in the full clauses of $F$. Every full clause $C$ must be resolvable on $v$ with another full clause $D$, yielding $E := C \res D$, and thus the full clauses of $F$ can be partitioned into pairs $\set{v} \addcup E, \set{\ol{v}} \addcup E$ for $\frac{\nfc(F)}2$ many clauses $E$ (of length $n(F) - 1$; note that because of fullness, for a given $E$ the clauses $C, D$ are uniquely determined up to order). \Qed
\end{prf}

Thus, if lower and upper bound match, they must be even numbers:
\begin{corol}\label{cor:evenfc}
  If $\maxsmar(k) = \nonmer(k)$ or (weaker) $\maxsmar(k) = \minnonmer(k)$ for some $k \in \NN$, then $\maxsmar(k)$ is even.
\end{corol}

In \cite[Corollary 2.6]{KullmannZhao2015FullClauses} we show that $\maxsmarh(k)$ indeed is always even (based on \cite[Utterly Trivial Observation]{Zeilberger2001CoveringSystems}), while \cite[Lemma 7.2]{KullmannZhao2015FullClauses} shows that the first $k$ with odd $\maxsmar(k)$ is given by $\maxsmar(7) = 9 = \nonmer(7) - 1$. Later in Lemma \ref{lem:nfcdef123} we see that the first $k$ with $\maxsmar(k) = \nonmer(k) - 1$ is $k=3$.

Another property of $\nfc(\Musati{\delta=k})$ related to evenness is that if $m$ is a possible number of full clauses, then $2 m$ is a possible number for $\delta = k+m-1$:
\begin{lem}\label{lem:simpconsfc}
  $2 m \in \nfc(\Musati{\delta=k+m-1})$ for $k \in \NN$ and $m \in \tb 1{\maxnfc(\Musati{\delta=k})}$.
\end{lem}
\begin{prf}
Consider $F \in \Musati{\delta=k}$ with $\nfc(F) = \maxnfc(\Musati{\delta=k})$. Choose $m$ full clauses of $F$, and choose a new variable $v \notin \var(F)$. Replace each of the chosen full clauses $C \in F$ by two clauses $C \addcup \set{v}, C \addcup \set{\ol{v}}$ (one non-strict and $m-1$ strict full subsumption extensions), obtaining $F'$. We have $F' \in \Musati{\delta=k+m-1}$ and $\nfc(F') = 2 m$. \Qed
\end{prf}

As a special case we obtain that $2, 4$ are always possible for the number of full clauses (except for $k=1$):
\begin{corol}\label{cor:alwaysfc}
  For $k \in \NN$ holds $2 \in \nfc(\Musati{\delta=k})$, and if $k \ge 2$ then $4 \in \nfc(\Musati{\delta=k})$.
\end{corol}
\begin{prf}
We show the assertion by induction on $k$, using Lemma \ref{lem:simpconsfc}, as follows: We have $\set{\set{1},\set{-1}} \in \Musati{\delta=1}$, so consider $k \ge 2$. We know $2 \in \nfc(\Musati{\delta=k-1})$, thus $4 \in \nfc(\Musati{\delta=k-1+2-1 = k})$. And once we have any $F \in \Musati{\delta=k}$ with a full clause, we get $F' \in \Musati{\delta=k}$ with $\nfc(F') = 2$ by performing a non-strict full subsumption extension on that full clause. \Qed
\end{prf}

We now turn to the determination of $\maxnfc(\Musati{\delta=k})$ for $k = 1,2,3$.
\begin{lem}\label{lem:nfcdef123}
  We have:
  \begin{enumerate}
  \item\label{lem:nfcdef1231} $\maxsmarh(1) = \maxsmar(1) = 2 = \nonmer(1)$.
  \item\label{lem:nfcdef1232} $\maxsmarh(2) = \maxsmar(2) = 4 = \nonmer(2)$.
  \item\label{lem:nfcdef1233} $\maxsmar(3) = 4 = \nonmer(3)-1$.
  \end{enumerate}
\end{lem}
\begin{prf}
Part \ref{lem:nfcdef1231}: By Corollary \ref{cor:maxnumfullcl} we have $\maxnfc(\Musati{\delta=1}) \le 2$ (which can also be deduced from the fact, that between two clauses of $F \in \Musati{\delta=1}$ there is at most one conflict). Due to $A_1 \in \Uclashi{\delta=1}$ we have $\maxnfc(\Uclashi{\delta=1}) \ge 2$.

Part \ref{lem:nfcdef1232}: By $A_2 \in \Uclashi{\delta=2}$ holds $\maxnfc(\Uclashi{\delta=2}) \ge 4$, by Corollary \ref{cor:maxnumfullcl} we have $\maxnfc(\Musati{\delta=2}) \le 4$.

Part \ref{lem:nfcdef1233}: By Corollary \ref{cor:alwaysfc} we have $\maxnfc(\Musati{\delta=3}) \ge 4$, by Corollary \ref{cor:evenfc} we have $\maxnfc(\Musati{\delta=3}) \le 4$. \Qed
\end{prf}

In Corollary \ref{cor:maxsmar} together with Theorem \ref{thm:supminvdegk6} these determinations will be completed for $k \le 6$.

\section{A method for improving upper bounds for $\minnonmer$}
\label{sec:improveupbdMU}

We now present a framework for generalising the argumentation of Theorem \ref{thm:MUminvdegdef} together with the analysis of the underlying recursion from Section \ref{sec:nonmer}. The idea is as follows:
\begin{enumerate}
\item We start with upper bounds $\minnonmer(k) \le a_k$ for $k=1,\dots,p$, collected in a ``valid bounds-function'' $f$.
\item For deficiency $p+1$ and an envisaged min-var-degree $m$ we consider the set $\potp_f(p+1,m)$ of ``possible'' degree-pairs of variables (the degrees of the positive and negative literals) in an envisaged clause-set $F \in \Smusati{\delta=p+1}$ with $\minvdeg(F) = m$.
\item If $\potp_f(p+1,m) = \es$, then $m$ is ``inconsistent'', that is, impossible to realise, and thus $\minnonmer(p+1) < m$ (as shown in Theorem \ref{thm:basicproppp}).
\item While in case of $\potp_f(p+1,m) \ne \es$ there might exist such an $F$ or not (the formal reasoning underlying the definition of $\potp_f(p+1,m)$ is not complete).
\end{enumerate}
So here we generalise the approach of Section \ref{sec:nonmer} for describing the function $\nonmer$ to a general recursion scheme, obtaining a general method for improved \emph{upper bounds}. The applications in this \Schrift{} are as follows:
\begin{itemize}
\item In Theorem \ref{thm:altrecnonmersenne} we obtain an alternative description of $\nonmer(k)$.
\item In Section \ref{sec:strbMU} we will first show that the smallest $k$ with $\minnonmer(k) < \nonmer(k)$ is $k=6$, namely $\minnonmer(6) = 8 = \nonmer(6) - 1$ (Theorem \ref{thm:supminvdegk6}).
\item By the general recursion scheme then follows from this improvement, that for all $k=2^m-m+1$ for $m \ge 3$ we have $\minnonmer(k) \le \nonmer(k)-1$. This improved upper bound is denoted by by $\nonmer_1: \NN \ra \NN$ (Theorem \ref{thm:imprupperbound}).
\end{itemize}

\subsection{Analysing splitting-situations}
\label{sec:analyspsit}

``Valid bounds-functions'' shall be monotonically increasing --- we know that $\nonmer$ is (strictly) monotonically increasing, and we show that $\minvdeg$ is also monotonically increasing (not strictly, as we will later see in Theorem \ref{thm:supminvdegk6}):
\begin{lem}\label{lem:minvdegvalid}
  The map $\minnonmer$ is monotonically increasing ($\minnonmer(k) \le \minnonmer(k+1)$ for $k \in \NN$).
\end{lem}
\begin{prf}
  For $F \in \Musati{\delta=k}$, $n(F) \not= 0$, we can construct $F' \in \Musati{\delta=k+1}$ with $\minvdeg(F) \le \minvdeg(F')$ as follows:
  \begin{enumerate}
  \item If $F$ is full, then obtain a non-full $F'' \in \Musati{\delta=k}$ with $\minvdeg(F) = \minvdeg(F'')$ by a full singular unit-extension (Lemma \ref{lem:fullsingue}), and replace $F$ by $F''$.
  \item If $F$ is not full, then first saturate $F$, and then perform a strict full subsumption extension (Lemma \ref{lem:aux2subsumpR}), obtaining the desired $F'$. \Qed
  \end{enumerate}
\end{prf}

We define now ``valid bounds-functions'', which are sensible as upper bounds on $\minnonmer$, and we also define how to obtain such a function from initial upper bounds $\minnonmer(k) \le a_k$ for $k=1,\dots,p$:
\begin{defi}\label{def:validboundsf}
  A \textbf{valid bounds-function} is a function $f: \NN \ra \nnpi$  fulfilling the following three conditions:
  \begin{enumerate}
  \item $f(1) = 2$.
  \item $f$ is monotonically increasing (i.e., $\fa\, k, k' \in \NN: k \le k' \Ra f(k) \le f(k')$).
  \item $f(k)$ is an upper bound for the minimum variable-degree of minimally unsatisfiable clause-sets of deficiency $k$ (i.e., $\fa\, k \in \NN: \minnonmer(k) \le f(k)$).
  \end{enumerate}
  The set of all valid bounds-functions is denoted by $\bmm{\svbf} \subset (\nnpi)^{\NN} = \set{f: \NN \ra \nnpi}$. And by $\bmm{\svbfs} := \set{f \in \svbf : f \le \nonmer} \subset \NN^\NN$ we denote the set of valid bounds-functions (pointwise) less-or-equal than the non-Mersenne function.

  For $a_1,\dots,a_p \in \NN$, $p \in \NN$. such that $a_1 = 2$, $a_i \le a_j$ for $i \le j$, and $a_i \ge \minnonmer(i)$, we define \bmm{[a_1,\dots,a_p]} as that $f \in \svbf$ with $f(k) = a_k$ for $k \in \tb 1p$, while $f(k) = \infty$ for $k > p$.
\end{defi}
By Lemma \ref{lem:minvdegvalid} $\minnonmer$ is a valid bounds-function, namely the smallest possible one. By Theorem \ref{thm:MUminvdegdef} and Corollary \ref{cor:NMmon} also $\nonmer$ is a valid bounds-function. In Corollary \ref{cor:altrecnonmersenne} we will see, that the continuation $[a_1,\dots,a_p](k) = \infty$ for $k \ge p+1$ is harmless in the sense, that $\nonmer$ can automatically be taken into account, via the improvement of valid bounds-functions through the use of potential degree-pairs defined below.

\begin{lem}\label{lem:compllatvb}
  $\svbf$ as well as $\svbfs$, together with $\le$, is a complete lattice, where infima resp.\ suprema are given by pointwise minimum resp.\ pointwise supremum. The smallest elements of both lattices is $\minnonmer$, while the largest is $[2]$ resp.\ $\nonmer$.
\end{lem}
The following definition reflects the main method to analyse and improve a given upper bound $f$ for $\minnonmer(k)$, namely it determines the numerical possibilities compatible with $f$ (recall the proof of Theorem \ref{thm:MUminvdegdef}):
\begin{defi}\label{def:possiblepairs}
  Consider $k, m \in \NN$ with $k \ge 2$ and $m \ge 4$, together with a valid bounds-function $f$. The \textbf{set of potential degree-pairs w.r.t.\ $f$} for (deficiency) $k$ and (minimum variable-degree) $m$, denoted by \bmm{\potp_f(k,m)}, is the set of pairs $(e_0, e_1) \in \NN^2$ fulfilling the following conditions:
  \begin{enumerate}[(i)]
  \item\label{def:possiblepairs1} $e_0, e_1 \ge 2$
  \item\label{def:possiblepairs2} $e_0, e_1 \le k$
  \item\label{def:possiblepairs3} $e_0 + e_1 = m$
  \item\label{def:possiblepairs4} $e_0 \le e_1$
  \item\label{def:possiblepairs5} $\fa\, \ve \in \set{0,1} : f(k-e_{\ve}+1) +e_{\ve} \ge m$.
  \end{enumerate}
  We set $\bmm{\potp(k,m)} := \potp_{\minnonmer}(k,m)$.
\end{defi}
The motivation for Definition \ref{def:possiblepairs} is to assume $F \in \Smusati{\delta=k}$ with $\minvdeg(F) = m$ and $v \in \varmvd(F)$, and to determine the possible literal-degrees $e_0 = \ldeg_F(\ol{v})$, $e_1 = \ldeg_F(v)$, ``possible'' in a formal sense. ``e'' stands for ``eliminated clauses'', namely $e_{\ve}$ is the number of clauses eliminated by $\pao v{\ve}$. The ``high'' values of $m$ (for fixed $k$) are of real interest; compare Lemma \ref{lem:montonepotp}. The basic properties of $\potp_f$ are as follows:
\begin{enumerate}
\item\label{item:possiblepairs_ext} For every valid $f$ and $k \ge 2$ we have $\potp_f(k,4) = \set{(2,2)}$ and $\potp_f(k,m) = \es$ for $m > 2k$.
\item\label{item:possiblepairs2} Discussion of the five conditions (i) - (v) in Definition \ref{def:possiblepairs}:
  \begin{enumerate}[(i)]
  \item Only nonsingular variables are considered, since only in this way the deficiency strictly decreases.
  \item The deficiency of $F_{\ve} := \pao{v}{\ve} * F$ is $k_{\ve} := k-e_{\ve}+1 \ge 1$ (recall Lemma \ref{lem:auxminvardeg} with $m_{\ve} = e_{\ve}$), splitting on a variable with minimal degree.
  \item $e_0, e_1$ are the literal-degrees of $\ol{v}, v$, which sum up to the variable-degree $m$ of $v$.
  \item W.l.o.g.\ we can restrict attention to such degree-pairs, since $F$ plays a role only up to isomorphism, and thus one can flip the sign of $v$ in $F$.
  \item We have $F_{\ve} \in \Musati{\delta=k_{\ve}}$ (assuming $F$ is saturated). And for $w \in \var(F_{\ve})$ we have $\vdeg_F(w) \le \vdeg_{F_{\ve}}(w) + e_{\ve}$. If for some $\ve \in \set{0,1}$ we would have $\minvdeg(\Musati{\delta=k_{\ve}}) +e_{\ve} < m$, then for $w \in \varmvd(F_{\ve})$ we would have $\vdeg_F(w) \le \vdeg_{F_{\ve}}(w) + e_{\ve} \le \minvdeg(\Musati{\delta=k_{\ve}}) + e_{\ve} < m$, but by assumption on $w$ we have $\vdeg_F(w) \ge m$.
  \end{enumerate}
\item\label{item:possiblepairs3} An important special case of $\potp_f(k,m)$ is $\potp_{\nonmer}(k,m)$; we have $\potp(k,m) \sse \potp_{\nonmer}(k,m)$ (see Lemma \ref{lem:potpincl} for a generalisation). The main point in using functions $f$ is that the precise values of $\minnonmer(k)$ might not be known.
\item\label{item:possiblepairs6} To compute $\potp_f(k,m)$ according to the definition,  only the values $f(k')$ for $k' \in \tb 1{k-1}$ are needed.
\end{enumerate}

\begin{examp}\label{exp:potp}
  Consider $f := [\nonmer(1),\nonmer(2),\nonmer(3)] = [2,4,5]$. First we determine $\potp_f(4,7)$:
  \begin{enumerate}
  \item By Conditions (i) - (iv) only $\set{(3,4)}$ remains.
  \item Now $f(4-3+1) + 3 = f(2) + 3 = 7$, but $f(4-4+1) + 4 = f(1) + 4 = 6 < 7$.
  \item Thus $\potp_f(4,7) = \es$.
  \end{enumerate}
  We will see in Theorem \ref{thm:basicproppp} that we can conclude $\minnonmer(4) \le 6$ (there is no ``formal'' possibility to reach the min-var-degree of $7$ for deficiency $4$). Now we determine $\potp_f(4,6)$:
  \begin{enumerate}
  \item By Conditions (i) - (iv), $\set{(2,4),(3,3)}$ are the possibilities.
  \item Checking Condition (v) for $(2,4)$: $f(4-2+1) + 2 = f(3) + 2 = 7$, $f(4-4+1) + 4 = f(1) + 4 = 6$.
  \item Checking Condition (v) for $(3,3)$: $f(4-3+1) + 3 = f(2) + 3 = 7$.
  \item Thus $\potp_f(4,6) = \set{(2,4),(3,3)}$.
  \end{enumerate}
  The intuitive meaning of this is, that a min-var-degree of $6$ can not be excluded by this type of formal reasoning, and so $4 \mapsto 6$ is the refinement of $[2,4,5]$, as will be constructed in Lemma \ref{lem:newboundspotp}.
\end{examp}

We invite the reader to compute the following special case of what we show later (in the proof of Theorem \ref{thm:imprupperboundhit}; it might also be useful to consider Table \ref{tab:auxfunc}):
\begin{examp}\label{exp:pp13}
  $\potp_{\nonmer}(13,17) = \set{(8,9)}$, while for any valid bounds-function $f$ with $f(k) = \nonmer(k)$ for $k \in \tb 15$ and $f(6) = \nonmer(6) - 1 = 8$ holds $\potp_f(13,17) = \es$.
\end{examp}

If we have a potential degree-pair for $m$, then also for $m' \le m$:
\begin{lem}\label{lem:montonepotp}
  Consider $k, m, m' \in \NN$ with $k \ge 2$ and $4 \le m' \le m$, and consider a valid bounds-function $f$. If $\potp_f(k,m) \not= \es$, then also $\potp_f(k,m') \not= \es$.
\end{lem}
\begin{prf}
  Consider $(e_0, e_1) \in \potp_f(k,m)$. Consider any $2 \le e_0' \le e_0$ and $2 \le e_1' \le e_1$ with $e_0' \le e_1'$ and $e_0' + e_1' = m'$. Now $f(k - e'_{\ve} + 1) + e'_{\ve} \ge f(k - e_{\ve} + 1) + e'_{\ve} = f(k - e_{\ve} + 1) + e_{\ve} - e_{\ve} + e'_{\ve} \ge m - e_{\ve} + e'_{\ve} = m - (m - e_{\ol{\ve}}) + (m' - e'_{\ol{\ve}}) = e_{\ol{\ve}} + m' - e'_{\ol{\ve}} \ge m'$ for $\ve \in \set{0,1}$, and thus $(e_0', e_1') \in \potp_f(k,m') \not= \es$. \Qed
\end{prf}

Using a smaller bounds-function can not yield more potential degree-pairs, as is obvious from Definition \ref{def:possiblepairs}:
\begin{lem}\label{lem:potpincl}
  Consider $k, m \in \NN$ with $k \ge 2$, $m \ge 4$, and valid bounds-functions $f_1, f_2$ with $f_1 \le f_2$ (pointwise). Then $\potp_{f_1}(k,m) \sse \potp_{f_2}(k,m)$. Especially for any valid bounds-function $f$ holds $\potp(k,m) \sse \potp_f(k,m)$.
\end{lem}

Again directly by definition (using monotonicity of valid bounds functions) we get that increasing $k$ while keeping $m$ can not remove potential degree-pairs:
\begin{lem}\label{lem:potpincdef}
  Consider $k, m \in \NN$ with $k \ge 2$, $m \ge 4$, and a valid bounds-function $f$. Then $\potp_f(k,m) \sse \potp_f(k+1,m)$.
\end{lem}

The main use of potential degree-pairs is to yield upper bounds on $\minnonmer(k)$:
\begin{thm}\label{thm:basicproppp}
  Consider $k, m \in \NN$ with $k \ge 2$, $m \ge 4$, and a valid bounds-function $f$. If $\potp_f(k,m) = \es$, then $\minnonmer(k) < m$.
\end{thm}
\begin{prf}
Assume $\minvdeg(\Musati{\delta=k}) \ge m$. Then there is $F \in \Smusatnsi{\delta=k}$ with $\minvdeg(F) \ge m$ (using Corollary \ref{cor:wlognonsingsat}). Consider $v \in \varmvd(F)$; if $\ldeg_F(\ol{v}) \le \ldeg_F(v)$ holds, then let $e := (\ldeg_F(\ol{v}), \ldeg_F(v))$, while otherwise flip the components of this pair. Now we have $e \in \potp(k,\minvdeg(F))$ (using Remark \ref{item:possiblepairs2} to Definition \ref{def:possiblepairs}), and thus $\potp_f(k,m) \not= \es$ by Lemmas \ref{lem:montonepotp}, \ref{lem:potpincl}, contradicting the assumption. \Qed
\end{prf}

\subsection{Recursion on potential degree-pairs}
\label{sec:recpdp}

Theorem \ref{thm:basicproppp} is used to improve valid bounds-functions $f$, taking $f$ as providing \emph{additional} upper bounds, besides what reasoning via potential degree-pairs yields:
\begin{lem}\label{lem:newboundspotp}
  Consider $f \in \svbf$. We obtain $f' \in \svbf$ recursively as follows:
  \begin{enumerate}
  \item $f'(1) := 2$.
  \item For $k > 1$ consider the largest $4 \le m \le 2 k$ such that $\potp_{f'}(k,m) \not= \es$, using Remark \ref{item:possiblepairs6} to Definition \ref{def:possiblepairs} (that we only need $f'(k')$ for $k' < k$).
  \item Now $f'(k) := \min(m, f(k))$.
  \end{enumerate}
\end{lem}
\begin{prf}
$f'(k)$ is well-defined for $k > p$ due to Remark \ref{item:possiblepairs_ext} to Definition \ref{def:possiblepairs}. That $f'$ is valid follows by induction as follows. We have to show $f'(k) \le f'(k+1)$ and $\minnonmer(k) \le f'(k)$ for all $k \in \NN$. For $k=1$ both properties are true by definition. And the induction step follows for monotonicity by Lemma \ref{lem:potpincdef}, and for the upper-bound-condition by Theorem \ref{thm:basicproppp}. \Qed
\end{prf}

The mapping $f \in \svbf \mapsto f' \in \svbf$ we call the ``non-Mersenne operator'':
\begin{defi}\label{def:potprec}
  For $f \in \svbf$ let the $f' \in \svbf$ according to Lemma \ref{lem:newboundspotp} be denoted by $\bmm{\potprec(f)} := f'$ (defined via ``recursion on potential degree-pairs''); we call $\potprec: \svbf \ra \svbf$ the ``non-Mersenne operator''.
\end{defi}

The basic properties of the non-Mersenne operator are that of a \href{http://en.wikipedia.org/wiki/Closure_operator#Closure_operators_on_partially_ordered_sets}{``kernel operator''}, which are order-theoretic properties as follows:
\begin{lem}\label{lem:kernelop}
  The map $\potprec: \svbf \ra \svbf$ is a kernel operator of the complete lattice $\svbf$, that is, for all $f, g \in \svbf$ holds:
  \begin{enumerate}
  \item $\potprec(f) \le f$ (intensive)
  \item $\potprec(\potprec(f)) = \potprec(f)$ (idempotent)
  \item $f \le g \Ra \potprec(f) \le \potprec(g)$ (monotonically increasing).
  \end{enumerate}
\end{lem}
\begin{prf}
Intensitivity follows by definition of $\potprec$ (note that in Lemma \ref{lem:newboundspotp} we have defined $f'(k)$ such that $f'(k) \le f(k)$ holds). Also idempotence follows directly from the definition in Lemma \ref{lem:newboundspotp}, namely that $f'(k)$ for $ k > 1$ already uses the improved values $f'(k')$ for $k' < k$. Monotonicity follows by Lemma \ref{lem:potpincl}. \Qed
\end{prf}

By Lemma \ref{lem:kernelop} we get that $\potprec(f)$ for $f \in \svbf$ is the supremum of the set of $f' \le f$ with $\potprec(f') = f'$. By Theorem \ref{thm:basicproppp} we get $\potprec(\minnonmer) = \minnonmer$. In order to show that the non-Mersenne operator at most reproduces $\nonmer$, that is, for all $f \in \svbf$ holds $\potprec(f) \le \nonmer$, we need to provide potential degree-pairs for $\nonmer$:
\begin{lem}\label{lem:nonmerpotp}
  For $k \ge 2$ we have (recall Definition \ref{def:twoaux}):
  \begin{enumerate}
  \item\label{lem:nonmerpotp1} $(h(k), \inonmer(k)) \in \potp_{\nonmer}(k, \nonmer(k))$.
  \item\label{lem:nonmerpotp2} $\potp_{\nonmer}(k, \nonmer(k) + 1) = \es$.
  \end{enumerate}
\end{lem}
\begin{prf}
Let $m := \nonmer(k)$. For Part \ref{lem:nonmerpotp1} let $e_0 := h(k)$, $e_1 := \inonmer(k)$; so we have to show $(e_0, e_1) \in \potp_{\nonmer}(k, m)$. Consider the conditions (i) - (v) in Definition \ref{def:possiblepairs}. We have $e_0 \ge 2$, since $\nonmer \ge 2$ in general, and $e_1 \ge 2$ by Definition \ref{def:i(k)}. As shown in  Corollary \ref{cor:characinonmer} we have $e_0 \le e_1$, where $e_1 \le k$ by definition. Furthermore we have $e_0 + e_1 = m$ by Lemma \ref{lem:simprecnm}, Part \ref{lem:simprecnm2}. Altogether we have now shown conditions (i) - (iv), and it remains to show that $\nonmer(k - e_{\ve} + 1) + e_{\ve} \ge m$ holds for both $\ve \in \set{0,1}$; for $\ve = 1$ we have equality, as already remarked, and it remains to show $\nonmer(k - e_0 + 1) + e_0 \ge m$, which is equivalent to
\begin{displaymath}
  \nonmer(k - e_0 + 1) \ge \inonmer(k).
\end{displaymath}
By Definition \ref{def:i(k)} of $\inonmer(k)$ (as the smallest $i$) this is implied by $\nonmer(k - e_0 + 1) \ge \nonmer(k - \nonmer(k - e_0 + 1) + 1)$. By the monotonicity of $\nonmer$ this is implied by $e_0 \le \nonmer(k - e_0 + 1)$, i.e., $\nonmer(k-\inonmer(k)+1) \le \nonmer(k - e_0 + 1)$. Again by monotonicity, this is implied by $\inonmer(k) \ge e_0$, i.e., $e_1 \ge e_0$, which we have already shown.

For Part \ref{lem:nonmerpotp2} we have to show $\potp_{\nonmer}(k, m + 1) = \es$. Assume $(e_0, e_1) \in \potp_{\nonmer}(k, m + 1)$ according to Definition \ref{def:possiblepairs}. Thus $\nonmer(k - e_1 + 1) + e_1 \ge m+1$, where $2 \le e_1 \le k$. Because of $e_0 + e_1 = m+1$ and $e_0 \le e_1$ holds $e_1 \ge \frac 12 (m+1)$, so $\min(2 e_1, \nonmer(k - e_1 + 1) + e_1) \ge m+1$, and thus $\nonmer(k) \ge m+1$ (Definition \ref{def:minvdegdef}). \Qed
\end{prf}

We obtain an alternative recursion for $\nonmer(k)$ (recall Definition \ref{def:minvdegdef}):
\begin{thm}\label{thm:altrecnonmersenne}
  $\potprec([2]) = \potprec(\nonmer) = \nonmer$.
\end{thm}
\begin{prf}
By Definition \ref{def:potprec} and Lemma \ref{lem:nonmerpotp} we get $\potprec([2]) = \nonmer$. Since $\potprec$ is idempotent, we also get $\potprec(\nonmer) = \nonmer$. \Qed
\end{prf}

So the non-Mersenne operator yields $\nonmer$ in the worst-case:
\begin{corol}\label{cor:altrecnonmersenne}
  $\potprec: \svbf \ra \svbfs$, that is, for every $f \in \svbf$ holds $\potprec(f) \le \nonmer$.
\end{corol}

\section{Strengthening the upper bound for $\minnonmer$}
\label{sec:strbMU}

In this final section many techniques introduced in this \Schrift{} come together, and we give some initial sharpness results (considering small deficiencies), and some non-sharpness results in the form of improved bounds (improving $\nonmer$ for infinitely many deficiencies). In Subsection \ref{sec:delta6} we determine $\minnonmer(k)$ for $1 \le k \le 6$ as values $2,4,5,6,8,8$, where the main achievement is Theorem \ref{thm:supminvdegk6}, showing $\minnonmer(6) = 8 = \nonmer(6) - 1$ (see \cite{KullmannZhao2015FullClauses} for a list of values up to $k=13$). Applying the non-Mersenne operator, we obtain the improved upper bound $\minnonmer(k) \le \nonmer_1(k)$ in Subsection \ref{sec:sharpeningbound}, where $\nonmer_1$ is like $\nonmer$, but with a duplication after the jump positions, that is, $\Delta \nonmer(k) = \Delta \nonmer_1 (k) = 2$ is followed by $\Delta \nonmer_1(k+1) = \Delta \nonmer(k+1) - 1 = 0$.

\subsection{Deficiencies $1,\dots,6$}
\label{sec:delta6}

We show that the first deficiency $k$, for which the bound $\minnonmer(k) \le \nonmer(k)$ is not sharp, is $k=6$. First we prove sharpness for the first five values:
\begin{thm}\label{thm:def15}
  For $k \in \tb 15$ we have $\minnonmer(k) = \minnonmerh(k) = \nonmer(k)$.
\end{thm}
\begin{prf}
We have to give examples showing that the upper bound $\nonmer(k)$ is attained for examples in $\Uclashi{\delta=k})$. Lemma \ref{lem:sharpnesssimpldef} covers deficiencies $k = 1,2,5$, namely 
\begin{enumerate}
\item\label{thm:def15k1} $A_1 \in \Uclashi{\delta=1}$ has $\nfc(A_1) = 2 = \nonmer(1)$ (recall Example \ref{exp:MU1}).
\item\label{thm:def15k2} $A_2 \in \Uclashi{\delta=2}$ has $\nfc(A_2) = 4 = \nonmer(2)$ (recall Example \ref{exp:MU2}).
\item\label{thm:def15k5} $A_3 \in \Uclashi{\delta=5}$ has $\nfc(A_3) = 8 = \nonmer(5)$.
\end{enumerate}
Deficiency $k=4$ is a jump position, and thus covered by Lemma \ref{lem:sharpjump}, where the example is as follows:
\begin{enumerate}
\item[4.]\label{thm:def15k4} For $F_4 := \set{\set{1,2},\set{-1,2,3},\set{1,-2,3},\set{-1,-2,3},\set{-1,2,-3},\set{1,-2,-3},\\ \set{-1,-2,-3}}$ we have $F_4 \in \Uclashi{\delta=4}$ with $\nfc(F_4) = 6 = \nonmer(4)$.
\end{enumerate}
The remaining case $k=3$ we obtain via strict full subsumption resolution from $F_4$:
\begin{enumerate}
\item[5.]\label{thm:def15k3} For $F_3 := \set{\set{1,2},\set{-1,3},\set{1,-2,3},\set{-1,2,-3},\set{1,-2,-3}, \set{-1,-2,-3}}$ we have $F_3 \in \Uclashi{\delta=3}$ with $\minvdeg(F_3) = 5 = \nonmer(3)$. \Qed
\end{enumerate}
\end{prf}

The examples in the proof of Theorem \ref{thm:def15} together with Lemma \ref{lem:nfcdef123} yield
\begin{corol}\label{cor:maxsmar}
  For $k \in \tb 15 \sm \set{3}$ we have $\maxsmarh(k) = \maxsmar(k) = \nonmer(k)$, while $\maxsmarh(3) = \maxsmar(3) = 4 = \minnonmer(3) - 1 = \nonmer(3) - 1$.
\end{corol}
In the sequel of this subsection we consider $k=6$. A computation shows that there is only one potential degree-pair for the min-var-degree as given by $\nonmer(6) = 9$:
\begin{lem}\label{lem:potp69}
  $\potp_{\nonmer}(6,9) = \set{(4,5)}$.
\end{lem}
\begin{prf}
  Conditions (i) - (iv) of Definition \ref{def:possiblepairs} yield $\potp_{\nonmer}(6,9) \sse \set{(3,6), (4,5)}$. Condition (v) excludes $(3,6)$, since we have $\nonmer(6-6+1) + 6 = 8 \not\ge 9$, while $(4,5)$ fulfils this condition due to $\nonmer(6-4+1) + 4 = 5 + 4 \ge 9$ and $\nonmer(6-5+1) + 5 = 4 + 5 \ge 9$. \Qed
\end{prf}

However, the potential degree-pair of Lemma \ref{lem:potp69} actually are not realisable, and thus $\minnonmer(6) < \nonmer(6)$, as we will show now. The proof works by assuming there is $F \in \Musati{\delta=6}$ with $\minvdeg(F) = 9$, and considering a splitting of $F$ into $F_{\ve} = \pao{v}{\ve} * F$, $\ve = 0,1$, on some $v \in \varmvd(F)$, as in the proof of Theorem \ref{thm:MUminvdegdef}. By Lemma \ref{lem:potp69} we know that w.l.o.g.\ the deficiencies of the two splitting results are $\delta(F_0) = 3$ and $\delta(F_1) = 2$. We can argue that $F_1$ is isomorphic to some $\Dt{m}$ (recall Example \ref{exp:MU2}), and thus every variable in $F_1$ has degree $4$. It follows that every variable of $F_1$ must be in each of the five clauses satisfied by the assignment (otherwise $\minvdeg(F) < 9$), which in turn means that $F_0$ has at least five full clauses, contradicting $\maxsmar(3) = 4$.

\begin{thm}\label{thm:supminvdegk6}
  $\maxsmarh(6) = \maxsmar(6) = \minnonmerh(6) = \minnonmer(6) = 8 = \nonmer(6) - 1$.
\end{thm}
\begin{prf}
$\maxsmarh(6) \ge 8$ is confirmed by the variable-clause matrix
\begin{displaymath}
  \begin{pmatrix}
    + & + & + & - & + & - & - & - & + & -\\
    + & + & + & + & + & + & + & + & - & -\\
    + & + & - & + & - & + & - & - & 0 & 0\\
    + & - & + & + & - & - & + & - & 0 & 0
  \end{pmatrix}
\end{displaymath}
($4$ variables, $10$ clauses, $8$ full clauses; unsatisfiability is given by $8 \cdot 2^{-4} + 2 \cdot 2^{-2} = 1$). We note this clause-set is in $\Sed$ by Theorem \ref{thm:characSEDMSAT}. Assume now that there exists $F \in \Musati{\delta=6}$ with $\minvdeg(F) = 9$. By Lemmas \ref{lem:sDPrMU}, \ref{lem:DPminvdeg} w.l.o.g.\ we can assume that $F$ is saturated and nonsingular. By Theorem \ref{thm:minnumvarmu} we know $n(F) \ge 4$. Consider $v \in \var(F)$ with $\vdeg_F(v) = 9$. W.l.o.g.\ we assume $\ldeg_F(v) \ge \ldeg_F(\ol{v})$. By Lemma \ref{lem:potp69} we have $\ldeg_F(v) = 5$, $\ldeg_F(\ol{v}) = 4$. For $F_{\ve} := \pao {v}{\ve} * F$, $\ve \in \set{0,1}$, we have $\delta(F_0) = 6 - 4 + 1 = 3$, $\delta(F_1) = 6 - 5 + 1 = 2$.

Let the $5$ occurrences of $v$ in $F$ be $C_1, \dots, C_5 \in F$, and let $C_i' := C_i \sm \set{v}$. And let the $4$ occurrences of $\ol{v}$ in $F$ be $D_1, \dots, D_4 \in F$, and let $D_i' := D_i \sm \set{\ol{v}}$. Using $G := \set{C \in F : v \notin \var(C)} = F \sm \set{C_1,\dots,C_5,D_1,\dots,D_4}$ we get
\begin{align*}
  F_0 = & \ \set{C_1', \dots, C_5'} \hspace{2pt} \addcup G & \qquad \text{($D_1,\dots,D_4$ vanish)}\\
  F_1 = & \ \set{D_1', \dots, D_4'} \addcup G & \qquad \text{($C_1,\dots,C_5$ vanish)},
\end{align*}
where $c(F_0) = 5 + c(G) = c(F) - 4$ and $c(F_1) = 4 + c(G) = c(F) - 5$.

Consider first $F_0 \in \Musati{\delta=3}$. We have $\minvdeg(F_0) \ge 9 - 4 = 5$, and thus $\minvdeg(F_0) = 5$ (due to $\minvdeg(F_0) \le \nonmer(3) = 5$). Every $w \in \varmvd(F_0)$ has at least $9$ occurrences in $F$ (since $\minvdeg(F) = 9$), from which at most $4$ are eliminated, and thus actually such variables have $\vdeg_F(w) = 9$, and furthermore $w \in \var(D_i)$ for all $i \in \tb 14$. By Lemma \ref{lem:twovarmvindeg}, Part \ref{lem:twovarmvindeg1}, we have $\abs{\varmvd(F_0)} \ge 2$, and so we have $\abs{D_i} \ge 3$ for all $i \in \tb 14$.

Now consider $F_1 \in \Musati{\delta=2}$. We have $\minvdeg(F_1) \ge 9 - 5 = 4$, thus $\minvdeg(F_1) = 4$ (due to $\minvdeg(F_1) \le \nonmer(2) = 4$), and thus by Lemma \ref{lem:existunitmu2} $F_1$ is nonsingular iff $F_1$ does not contain unit-clauses. If $F_1$ would contain a unit-clause, then there would be a binary clause $\set{\ol{v},x} \in F$, contradicting that all $D_i$ contain at least three literals. So $F_1$ is nonsingular, and thus $F_1$ is isomorphic to some $\Dt{m}$ for some $m \ge 2$. It follows that $F_1$ is $4$-variable-regular, where all the variables of $F_1$ have at least $9$ occurrences in $F$, and thus we have $\var(F_1) \sse \var(C_i')$ for all $i \in \tb 15$ (since five occurrences of every variable in $F_1$ must vanish). By $\var(F_0) = \var(F_1) = \var(F) \sm \set{v}$ we get $\var(C_i') = \var(F_0)$ for all $i \in \tb 15$.

Coming back to the structure of $F_0$, we now know that $F_0$ has five full clauses $C_1', \dots, C_5'$, which contradicts Lemma \ref{lem:nfcdef123}, Part \ref{lem:nfcdef1233}. \Qed
\end{prf}

\subsection{Sharpening the bound}
\label{sec:sharpeningbound}

Based on recursion on potential degree-pairs, we can improve the upper bound $\nonmer(k)$ for $\minvdeg(\Musati{\delta=k})$ for $k \ge 6$ (generalising Example \ref{exp:pp13}):
\begin{defi}\label{def:nm1}
  Let $\nonmer_1: \NN \ra \NN$ be defined as $\nonmer_1 := \potprec([2,4,5,6,8,8])$ (recall Definition \ref{def:potprec}).
\end{defi}

By Lemma \ref{lem:newboundspotp} together with Theorem \ref{thm:supminvdegk6} we get:
\begin{thm}\label{thm:imprupperbound}
  For all $k \in \NN$ we have $\minnonmer(k) \le \nonmer_1(k)$.
\end{thm}

It remains to determine $\nonmer_1$ numerically:
\begin{thm}\label{thm:imprupperboundhit}
  In Table \ref{tab:valuesminvdegdefhit} we find the values of $\nonmer_1(k)$ for $k \le 30$. We have $\nonmer_1(k) = \nonmer(k)$ for $k \notin \set{2^m-m+1 : m \in \NN, m \ge 3}$, while for $k = 2^m-m+1$ we have $\nonmer_1(k) = \nonmer(k) - 1 = 2^m$.
\end{thm}

\begin{table}[h]
  \centering
  \setlength{\tabcolsep}{5pt}
  \renewcommand{\arraystretch}{1.1}
  \begin{tabular}{c||*{15}{c}}
    $k$ & \bmm{1} & 2& 3 & \bmm{4} & 5 & 6 & 7 & 8 & 9 & 10 & \bmm{11} & 12 & 13 & 14 & 15\\
    $\nonmer_1(k)$ & 2 & 4 & 5 & 6 & 8 &\ul{8} & 10 & 11 & 12 & 13 & 14 & 16 & \ul{16} & 18 & 19\\
    \hline
    $k$ & 16 & 17 & 18 & 19 & 20 & 21 & 22 & 23 & 24 & 25 & \bmm{26} & 27 & 28 & 29 & 30\\
    $\nonmer_1(k)$ & 20 & 21 & 22 & 23 & 24 & 25 & 26 & 27 & 28 & 29 & 30 & 32 & \ul{32} & 34 & 35
  \end{tabular}
  \caption{Values of $\nonmer_1(k)$ for $k \in \tb 1{30}$, in bold the jump-values (i.e., $k \in J$), and underlined the changed values compared to $\nonmer(k)$; we see that directly after the jump we have stagnation, followed by a second jump.}
  \label{tab:valuesminvdegdefhit}
\end{table}

\begin{prf}
We use induction on $k$. Due to $2^3-3+1=6$ the assertion holds for $k \le 6$, which is the induction basis, and in the sequel we assume $k \ge 7$. We show the following, which implies the theorem (using $\nonmer_1 \le \nonmer$ by Lemma \ref{lem:kernelop}):
\begin{enumerate}
\item For $k = 2^m-m+1$, $m \ge 4$, we have
  \begin{enumerate}
  \item $\potp_{\nonmer_1}(k,\nonmer(k)) = \es$ and
  \item $\potp_{\nonmer_1}(k,\nonmer(k)-1) \ne \es$.
  \end{enumerate}
\item Otherwise $\potp_{\nonmer_1}(k,\nonmer(k)) \ne \es$.
\end{enumerate}

\paragraph{\textbf{Part 1}} We consider $k = 2^m-m+1$, $m \ge 4$. We have $\nonmer(k) = 2^m+1$.

\textbf{Part (a).} To show $\potp_{\nonmer_1}(2^m-m+1,2^m+1) = \es$, we assume $(e_0, e_1) \in \potp_{\nonmer_1}(2^m-m+1,2^m+1)$. Thus we know $e_0, e_1 \ge 2$, $e_0, e_1 \le 2^m-m+1$, $e_0 + e_1 = 2^m + 1$, $e_0 \le e_1$, whence $e_0 \le 2^{m-1}$, and
\begin{equation}\label{eq:proof1}
  \nonmer_1(2^m - m + 1 - e_{\ve} + 1) + e_{\ve} \ge 2^m + 1
\end{equation}
for both $\ve \in \set{0,1}$.

\textbf{Case (a.1).} Assume $e_0 \le 2^{m-1} - 1$, and thus $e_1 \ge 2^{m-1} + 2$.

From \eqref{eq:proof1} we get $\nonmer(2^m - m + 1 - e_1 + 1) + e_1 \ge 2^m + 1$, where (using Corollary \ref{cor:NMmon}):
\begin{multline*}
  \nonmer(2^m - m + 1 - e_1 + 1) + e_1 \ge 2^m + 1 \Ra\\
  \nonmer(2^m - m + 1 - (2^{m-1} + 2) + 1) + 2^{m-1} + 2 \ge 2^m + 1 \Lra\\
  \nonmer(2^{m-1} - m) \ge 2^{m-1} - 1,
\end{multline*}
where by Corollary \ref{cor:nonmbeforejump} we have $\nonmer(2^{m-1} - m) = \nonmer(2^{m-1} - (m-1) - 1) = 2^{m-1} - 2$, and we obtained a contradiction, finishing Case (a.1). $\surd$

\textbf{Case (a.2).} It remains $e_0 = 2^{m-1}$. From \eqref{eq:proof1} we get $\nonmer_1(2^m - m + 1 - e_0 + 1) + e_0 \ge 2^m + 1$, where $2^m - m + 1 - e_0 + 1 = 2^m - m + 1 - 2^{m-1} + 1 = 2^{m-1} - (m - 1) + 1$, and thus by induction hypothesis we get $\nonmer_1(2^m - m + 1 - e_0 + 1) + e_0 = 2^{m-1} + 2^{m-1} = 2^m$, a contradiction. $\surd$ This concludes Part (a).

\textbf{Part (b).} We show $(2^{m-1}, 2^{m-1}) \in \potp_{\nonmer_1}(k,\nonmer(k)-1)$.\footnote{We have  $(2^{m-1}, 2^{m-1}) = (h(k), \inonmer(k)-1)$, but we don't need this here.} For this it remains to show $\nonmer_1(2^m - m + 1 - 2^{m-1} + 1) + 2^{m-1} \ge 2^m$, and indeed $\nonmer_1(2^m - m + 1 - 2^{m-1} + 1) = \nonmer_1(2^{m-1} - (m-1) + 1) = 2^{m-1}$ by induction hypothesis. $\surd$ This concludes Part 1.

\paragraph{\textbf{Part 2}} $k \ne 2^m-m+1$ for any $m \ge 4$. We have to show $\potp_{\nonmer_1}(k, \nonmer(k)) \ne \es$.

\textbf{Case (a).} $k = 2^m-m+2$; thus $\nonmer(k) = 2^m + 2$.

We have $(2^{m-1}, 2^{m-1}+2) \in \potp_{\nonmer_1}(k, 2^m+2)$\footnote{We have  $(2^{m-1}, 2^{m-1}+2) = (h(k)-1, \inonmer(k)+1)$.}, due to $\nonmer_1(2^m-m+2 - 2^{m-1} + 1) + 2^{m-1} = \nonmer_1(2^{m-1}-(m-1) + 2) + 2^{m-1} = 2^{m-1} + 2 + 2^{m-1}$ and $\nonmer_1(2^m-m+2 - (2^{m-1}+2) + 1) + 2^{m-1} + 2 = \nonmer_1(2^{m-1}-(m-1)) + 2^{m-1} + 2 = 2^{m-1} + 2^{m-1} + 2$. $\surd$

\textbf{Case (b).} $k = 2^m-m+3$; thus $\nonmer(k) = 2^m + 3$.

We have $(2^{m-1}, 2^{m-1}+3) \in \potp_{\nonmer_1}(k, 2^m+3)$\footnote{We have  $(2^{m-1}, 2^{m-1}+3) = (h(k)-1, \inonmer(k)+1)$.}, due to $\nonmer_1(2^m-m+3 - 2^{m-1} + 1) + 2^{m-1} = \nonmer_1(2^{m-1}-(m-1) + 3) + 2^{m-1} = 2^{m-1} + 3 + 2^{m-1}$ and $\nonmer_1(2^m-m+3 - (2^{m-1}+3) + 1) + 2^{m-1} + 3 = \nonmer_1(2^{m-1}-(m-1)) + 2^{m-1} + 3 = 2^{m-1} + 2^{m-1} + 3 = 2^m + 3$. $\surd$

For all remaining cases (c) -- (e) we have
\begin{equation}
  \label{eq:nm1bounds}
  2^m-m+4 \le k \le 2^{m+1} - (m+1);
\end{equation}
we show
\begin{displaymath}
  (e_0, e_1) := \big ( h(k), \inonmer(k) \big ) \in \potp_{\nonmer_1}(k, \nonmer(k)),
\end{displaymath}
which is the same pair as we showed in Lemma \ref{lem:nonmerpotp}, Part \ref{lem:nonmerpotp1} to be in $\potp_{\nonmer}(k, \nonmer(k))$. Recall Corollary \ref{cor:inonmeraroundjump} for the computation of special values of $\inonmer(k)$.

Since Conditions (i) - (iv) from Definition \ref{def:possiblepairs} do not depend on the bounds-function, by Lemma \ref{lem:nonmerpotp}, Part \ref{lem:nonmerpotp1} we have already shown these conditions, and it remains to show
\begin{displaymath}
  \nonmer_1(k - e_{\ve} + 1) + e_{\ve} \ge \nonmer(k)
\end{displaymath}
for both $\ve \in \set{0,1}$; we call this condition ``(C$\ve$)''. If the argument $k - e_{\ve} + 1$ for $\nonmer_1$ in (C$\ve$), which we know is strictly less than $k$, is not of the special form $2^m - m + 1$ for some $m \ge 4$, then (C$\ve$) holds by  Lemma \ref{lem:nonmerpotp}, Part \ref{lem:nonmerpotp1} (and induction hypothesis); we call this sufficient condition ``NSF'' (``not special form'').

First we show (C1) via NSF. So $e_1 = \inonmer(k)$; we use $k - \inonmer(k) + 1 = i'(k)$ (recall Definition \ref{def:twoaux}), where $i'(k)$ is monotonically increasing, and thus we need to check only the lower and the upper bound on $k$ in \eqref{eq:nm1bounds}:
\begin{itemize}
\item For $k = 2^m-m+4$ holds $\inonmer(k) = 2^{m-1} + 2$, thus $i'(k) = 2^m-m+4 - (2^{m-1} + 2) + 1 = 2^{m-1} - m + 3 > 2^{m-1} - (m-1) + 1$.
\item For $k = 2^{m+1} - (m+1)$ holds $\inonmer(k) = 2^m$, thus $i'(k) = 2^{m+1} - (m+1) - 2^m + 1 = 2^m - m < 2^m - m + 1$.
\end{itemize}

So (C1) holds. $\surd$ It remains to show (C0); so $e_0 = h(k)$, and we use $h'(k) := k - h(k) + 1$. We do not have a jump-position within the $k$-range we consider, and thus by Lemma \ref{lem:characjumpb} we get, that $h'(k)$ is monotonically increasing for the $k$-range we consider, and thus, as above, we only need to consider the lower and upper bound on the $k$-range.

\textbf{Case (c).} $2^m-m+4 \le k \le 2^{m+1} - (m+1) - 2$.

NSF holds here:
\begin{itemize}
\item The lower bound $k = 2^m - m + 4$: $\inonmer(k) = 2^{m-1} + 2$, thus $h(k) = \nonmer(k) - \inonmer(k) = 2^m + 4 - 2^{m-1} - 2 = 2^{m-1} + 2$, which is the same as $\inonmer(k)$ above, and thus $h'(k) > 2^{m-1} - (m-1) + 1$.
\item The upper bound $k = 2^{m+1} - (m+1) - 2$: $\inonmer(k) = 2^m-1$, and thus $h(k) = \nonmer(k) - \inonmer(k) = 2^{m+1} - 3 - 2^m + 1 = 2^m - 2$, and so $h'(k) = 2^{m+1} - (m+1) - 2 - (2^m - 2) + 1 = 2^m - m < 2^m - m + 1$. $\surd$
\end{itemize}

\textbf{Case (d).} $k = 2^{m+1} - (m+1) - 1$; thus $\nonmer(k) = 2^{m+1} - 2$.

Now $\inonmer(k) = 2^m$, and so $h(k) = \nonmer(k) - \inonmer(k) = 2^m - 2$, whence $h'(k) = 2^{m+1} - (m+1) - 1 - (2^m - 2) + 1 = 2^m - m + 1$. So here NSF does not hold, but (C0) holds nevertheless: $\nonmer_1(2^m - m + 1) + 2^m - 2 = 2^m + 2^m - 2 = 2^{m+1} - 2$. $\surd$

\textbf{Case (e).} $k = 2^{m+1} - (m+1)$; thus $\nonmer(k) = 2^{m+1}$.

Now $h(k) = 2^m = \inonmer(k)$, and this is a special case of (C1). \Qed
\end{prf}

It is instructive to note the new $\Delta$-values explicitly:
\begin{corol}\label{cor:deltanm1}
  For $k \in \NN$ holds $\Delta \nonmer_1(k) \in \set{0,1,2}$, with
  \begin{enumerate}
  \item $\Delta \nonmer_1(k) = 0 \iff k = 2^m - m$ for some $m \in \NN$, $m \ge 3$.
  \item $\Delta \nonmer_1(k) = 2 \iff k = 2^m - m \pm 1$ for some $m \in \NN$, $m \ge 3$.
  \end{enumerate}
\end{corol}

\section{Conclusion and open problems}
\label{sec:open}

The main subject of this \Schrift{} can be seen in the study of $\minvdeg(\mc{C}_{\delta=k})$ for classes $\Uclash \sse \mc{C} \sse \Mlean$ and $k \in \NN$, that is, the study of the maximal minimum variable-degree of classes of matching-lean clause-sets containing all unsatisfiable hitting clause-sets, parameterised by the deficiency. If $\mc{C} \sse \Lean$, then this quantity is bounded, and indeed we have shown $\minvdeg(\Leani{\delta=k}) = \nonmer(k)$ (more generally this holds for every subclass of $\Lean$ containing $\Vmusat$). While for $\mc{C} = \Mlean$ this quantity is unbounded. For $\mc{C} = \Musat$ we have shown the improved bound $\minvdeg(\Musati{\delta=k}) = \minnonmer(k) \le \nonmer_1(k)$, where indeed also this bound is not sharp (as will be shown in \cite{KullmannZhao2014Sharper}; see Subsection \ref{sec:concmnM}) --- the question about the determination of $\minnonmer(k)$ is a major open research question for us.
For lean clause-sets we have shown the strengthened upper bound $\minvdeg(F) \le \nonmer(\surp(F))$ for the surplus $\surp(F) \le \delta(F)$, and indeed for every clause-set $F$ we can satisfiability-equivalently remove some clauses in polytime such that this upper bound holds.

\subsection{Conjectures and questions}
\label{sec:concmatchlean}

We made the following four conjectures:
\begin{enumerate}\setlength{\itemsep}{0pt}
\item Conjecture \ref{con:findauthard}: If a clause-set violates the upper bound on the min-var-degree for lean clause-sets, then it must have a non-trivial autarky. As we have seen, we can determine the set of variables involved, but the determination of the autarky itself is open --- the conjecture states that there is a poly-time algorithm for computing such an autarky. See Subsection \ref{sec:findaut} for more information on this topic.
\item Conjecture \ref{con:minvdeghit}: the maximum min-var-degree for unsatisfiable hitting clause-sets is the same as for the larger class of minimally unsatisfiable clause-sets. In Conjecture \ref{con:minvardeghit} we generalise this to non-boolean clause-sets.
\item Conjecture \ref{con:sharpness}: $\minnonmer$ is close to $\nonmer$, more precisely, $\nonmer -1 \le \minnonmer \le \nonmer$.
\item Conjecture \ref{con:maxnfcminvdeg}: $\maxsmar$ is close to $\minnonmer$, that is, $\minnonmer - 1 \le \maxsmar$. In Lemma \ref{lem:asympdet} we will state a weaker, but proven (in future work) lower bound.
\end{enumerate}
Five more conjectures are in this final section. We also asked the following questions:
\begin{enumerate}\setlength{\itemsep}{0pt}
\item Question \ref{que:maxvdegns} asks about the max-var-degree for nonsingular MUs.
\item Question \ref{que:uniformns} raises the possibility, that there are only finitely many possibilities for uniform and nonsingular MUs of a given deficiency.
\item Question \ref{que:vertdeghyp} asks, whether the minimum vertex degree is bounded for minimally non-2-colourable hypergraphs of bounded deficiency.
\item Question \ref{que:elimlit} is about some complexity problems around the elimination of literal occurrences in minimally unsatisfiable clause-sets.
\item  Question \ref{que:SED} is about the complexity of SAT decision for $\Sed$. At first sight it might seem easy to translate every $F \in \Cls$ into some sat-equivalent element of $\Sed$, and in fact to manipulate deficiency and surplus alone is rather easy, but we do not know how to handle them together.
\item Question \ref{que:surpminv} concerns the existence of unsatisfiable hitting clause-sets of arbitrary surplus equal deficiency (i.e., in $\Sed$) and a min-var-degree as low as possible. An underlying question is to better understand the surplus.
\item Question \ref{que:exmpmlean1} is about strengthening the construction of Lemma \ref{lem:exmpmlean}, for finding matching-lean clause-sets of a given deficiency with high minimum literal-degree (perhaps completely different constructions are needed).
\item Question \ref{que:extrSED} is about the structure of MUs with maximal min-var-degree.
\item Question \ref{que:maxnumfullcl} asks for the determination of $\maxnfc(\Leani{\delta=k})$.
\item Question \ref{que:sfsAn} is about lower bounds for $\minnonmerh(k), \maxsmarh(k)$ via strict full subsumption resolution starting with $A_n$.
\end{enumerate}

In the remainder we outline main research areas related to this \Schrift.

\subsection{Improved upper bounds for $\minnonmer$}
\label{sec:concmnM}

We know $\minnonmer \le \nonmer_1$, and we know $\minnonmer(k)$ precisely for $k \in \tb 16$ (extended in \cite{KullmannZhao2015FullClauses} to $\tb 1{13}$) and for $k \in J, J+1$. Also of high relevance here is to determine $\minnonmerh(k)$, which by Conjecture \ref{con:minvdeghit} is the same as $\minnonmer(k)$. Another major conjecture is Conjecture \ref{con:sharpness}, which says that $\minnonmer$ deviates at most by $1$ from $\nonmer$. Beyond the current \Schrift, we know the following improvements of the upper bound $\nonmer_1$, as developed in \cite{KullmannZhao2014Sharper}:
\begin{itemize}\setlength{\itemsep}{0pt}
\item Generalising the ideas of Theorem \ref{thm:supminvdegk6}, which is based on the improved upper bound for deficiency $2^3 - 3 + 1 = 6$, we can show also for deficiency $k=2^4-4+2 = 14$ that we have $\minnonmer(k) = \nonmer(k)-1$. Via the non-Mersenne operator, this yields the improved upper bound $\nonmer_2$.
\item More generally, we obtain a sequence of improved upper bounds $\nonmer_{m-2}$ for $m \ge 3$, improving the upper bound at deficiency $k=2^m-2$ from $\nonmer_{m-3}$ (which is unchanged from $\nonmer$ at this deficiency) and applying the non-Mersenne operator.
\item The infimum of $\nonmer_1, \nonmer_2, \dots$ is $\nonmer_{\omega}$.
\item However, this is not the end of it --- also for deficiency $k=15$ we have $\minnonmer(k) = \nonmer(k)-1$, obtaining $\nonmer_{\omega+1}$. This new improvement depends on new ideas --- will there be an infinite chain of ever-increasing complexity of such improvements?
\end{itemize}

We believe that a closed ``nice'' formula for $\minnonmer(k)$ is impossible:
\begin{conj}\label{con:minnonmernotcomp}
  The function $\minnonmer: \NN \ra \NN$ is ``complex'', and for no finite tuple $\vec{a}$ holds $\potprec(\vec{a}) = \minnonmer$, but $\minnonmer$ is computable in doubly-exponential time.
\end{conj}
See Lemma \ref{lem:conjimplcomp} for conditions implying the computability-part of Conjecture \ref{con:minnonmernotcomp}.

\subsection{Determining $\maxsmar$}
\label{sec:concSmar}

While Subsection \ref{sec:concmnM} was about improving the upper bound, here now we turn to the lower bound.
In Subsection \ref{sec:fullcl} we provided only the minimum needed in this \Schrift{} for the measure $\nfc(F)$ of full clauses. In the forthcoming \cite{KullmannZhao2015FullClauses} we show the following lower bound, using $S_2: \NN \ra \NN$, the function for the ``Smarandache Primitive Numbers'' introduced in \cite[Unsolved Problem 47]{Smarandache1993OPNS}, which for $k \in \NN$ is defined as the minimal natural number $s \in \NN$ such that $2^k$ divides $s!$. 
\begin{lem}[\cite{KullmannZhao2015FullClauses}]\label{lem:smarandache}
  For all $k \in \NN$ holds $\maxsmarh(k) \ge S_2(k)$.
\end{lem}
The basic structure of the proof of Lemma \ref{lem:smarandache} is as follows: Full subsumption extension (strict and non-strict) yields a non-deterministic process to create elements of $\Uclash$ with ``many'' full clauses. Maximising the number of full clauses for a given deficiency, a recursion is obtained, which is roughly similar to our recursion for $\nonmer(k)$ (Definition \ref{def:minvdegdef}). The analysis yields that the solution for this recursion is a \emph{meta-Fibonacci sequence} as introduced in \cite[Page 145]{Hofstadter1979GEB}, of a special form as analysed in \cite{Conolly1989MetaFibonacci}. We are then able to identify this special form as identical with $S_2(k)$ (as conjectured on the OEIS \cite{Sloane2008OEIS}).

Lemma \ref{lem:smarandache} yields the interesting inequality $S_2 \le \minnonmer \le \nonmer$ (but recall that in Lemma \ref{lem:preciseboundminnonmercon} we obtained a much sharper lower bound for $\minnonmer(k)$ from Conjecture \ref{con:sharpness}). This is relevant as the upper bound $\nonmer$ on $S_2$ as well as the lower bound $S_2$ on $\minnonmer$. From \cite{WenpengLiu2002Smarandache} we get that $k + 1 \le S_2(k)$ and thus by Corollary \ref{cor:upperboundnonmer} we get
\begin{lem}[\cite{KullmannZhao2015FullClauses}]\label{lem:asympdet}
  $k + 1 \le S_2(k) \le \nonmer(k) \le k + 1 + \fld(k)$ for $k \in \NN$.
\end{lem}
For sequences $a, b: \NN \ra \RR$ let \emph{asymptotic equality} be denoted by $a \sim b :\Lra \lim_{n \ra \infty} \frac{a_n}{b_n} = 1$. It is known that $S_2(k)$ and $\nonmer(k)$ are asymptotically equal to $(k)_{k \in \NN}$, and thus

\begin{corol}[\cite{KullmannZhao2015FullClauses}]\label{cor:asympdet}\sloppy
  The six sequences $S_2$, $\maxsmarh$, $\maxsmar$, $\minnonmerh$, $\minnonmer$, $\nonmer$ are asymptotically equal to $(k)_{k \in \NN}$.
\end{corol}

\begin{figure}[h]
  \centering
 $\xymatrix {
      & {\nonmer}\\
      & {\minnonmer} \ar[u]_{1}\\
      {\maxsmar} \ar[ur]_{1} \ar[uur]^{2} && {\minnonmerh} \ar[ul]^{=}\\
      & {\maxsmarh} \ar[ul]^{\infty} \ar[ur]_{\infty}\\
      & S_2 \ar[u]_{=}
    }$
  \caption{The four main combinatorial quantities, and the two numerical functions}
  \label{fig:sandwich}
\end{figure}
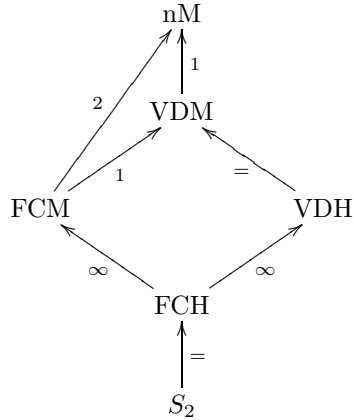

In Figure \ref{fig:sandwich} we show the six quantities from Corollary \ref{cor:asympdet} and the relations between them. An arrow means a (proven) $\le$-relation. If the arrow is labelled with ``$m$'', then we conjecture the difference is at most this number (where in all three cases here we know examples where this difference is attained), while the label ``$=$'' means that we conjecture equality, and the label ``$\infty$'' means that we conjecture that the difference is unbounded.\vspace{-1ex}

For a more precise asymptotic determination of these six quantities from Corollary \ref{cor:asympdet}, calling them $a_k$, we need to consider the six sequences $(a_k - k)_{k \in \NN}$. Currently we only know $\nonmer(k) - k \sim \ld(k)$. The place of $\nonmer_1$ and its refinements, as discussed in Subsection \ref{sec:concmnM}, in Figure \ref{fig:sandwich} is between $\nonmer$ and $\minnonmer$, as more and more refined approximations of $\minnonmer$ from above.

\subsection{Generalisation to non-boolean clause-sets}
\label{sec:concgennb}

It is interesting to generalise Theorem \ref{thm:MUminvdegdef} for generalised clause-sets; see \cite{Kullmann2007ClausalFormZI,Kullmann2007ClausalFormZII} for a systematic study, while the most general notion of generalised clause-sets, ``signed clause-sets'' are discussed in \cite{BeckertHaehnleManya2000SignedCNF}. Generalised clause-sets $F$ have literals $(v,\ve)$, meaning ``$v \not= \ve$'', for variables $v$ with non-empty finite domains $D_v$ and values $\ve \in D_v$. The deficiency is generalised by giving every variable a weight $\abs{D_v} - 1 \in \NNZ$ (which is $1$ in the boolean case), i.e., $\delta(F) = c(F) - \sum_{v \in \var(F)} (\abs{D_v} - 1) = c(F) + n(F) - \sum_{v \in \var(F)} \abs{D_v}$; see \cite[Subsection 7.2]{Kullmann2007ClausalFormZI}. A partial assignment is a map $\vp$ with some finite set of variables as domain $\dom(\vp) =: \var(\vp)$, which maps $v \in \var(\vp)$ to $\vp(v) \in D_v$. A partial assignment $\vp$ satisfies a clause-set $F$ iff for every $C \in F$ there is $(v,\ve) \in C$ with $v \in \var(\vp)$ and $\vp(v) \ne \ve$. Minimally unsatisfiable (generalised) clause-sets are defined as usual (they are unsatisfiable, while every strict subset is satisfiable). In \cite[Corollary 9.9]{Kullmann2007ClausalFormZI} it is shown that also all minimally unsatisfiable generalised clause-sets $F$ fulfil $\delta(F) \ge 1$ (based, like in the boolean case, on matching autarkies).

The degree $\vdeg_F(v)$ of a variable $v$ in a clause-set $F$ is the sum of the degrees of the literals $(v,\ve)$ for $\ve \in D_v$, and thus $\vdeg_F(v) = \abs{\set{C \in F : C \cap (\set{v} \times D_v) \ne \es}}$. For a given deficiency $k \in \NN$, the basic question is to determine the supremum of $\minvdeg(F)$ over all minimally unsatisfiable $F$ with $\delta(F) = k$. The base case of deficiency $k=1$ is handled in \cite[Lemma 5.4]{Kullmann2007ClausalFormZII}, showing that for generalised minimally unsatisfiable clause-sets of deficiency $1$ we have $\minvdeg(F) \le \max_{v \in \var(F)} \abs{D_v}$; actually all structural knowledge from \cite{AhLi86,DDK98,Ku99dKo} has been completely generalised in \cite[Subsection 5.2]{Kullmann2007ClausalFormZII}.

But $k \ge 2$ requires more work, since here the basic method of saturation is not available for generalised clause-sets, as discussed in \cite[Subsection 5.1]{Kullmann2007ClausalFormZII}: Saturated generalised clause-sets  (i.e., unsatisfiable clause-sets, where no literal occurrence can be added without rendering the clause-set satisfiable), with deficiency at least $2$, after splitting do not necessarily yield minimally unsatisfiable (generalised) clause-sets. Thus the proofs for the boolean case seem not to be generalisable for arbitrary minimally unsatisfiable (generalised) clause-sets.

In order to repair this, the ``substitution stability parameter regarding irredundancy'' $\stabpar(F) \in \ZZ_{\ge -1} \addcup \set{+\infty}$ is introduced in \cite[Subsection 5.3]{Kullmann2007ClausalFormZII}), defined as the supremum of $k \in \ZZ_{\ge -1}$ such that for every partial assignment with $n(\vp) := \abs{\var(\vp)} \le k$ the clause-set $\vp * F$, obtained as usual by application of $\vp$ to $F$, is minimally unsatisfiable. So $\stabpar(F) \ge 0$ iff $F$ is minimally unsatisfiable, and as shown in \cite[Corollary 4.8]{Kullmann2007ClausalFormZII}, $\stabpar(F) = +\infty$ iff $F$ is a hitting clause-set (i.e., for all $C, D \in F$, $C \ne D$, there are $x \in C$, $y \in D$ with $x = (v,\ve)$ and $y = (v,\ve')$ for some variable $v$ and $\ve, \ve' \in D_v$ with $\ve \ne \ve'$). And $\stabpar(F) \ge 1$ iff splitting on any variable yields always a minimally unsatisfiable clause-set. So for a boolean clause-sets $F$ holds $\stabpar(F) \ge 1$ iff $F$ is saturated, but for generalised clause-sets we only have that $\stabpar(F) \ge 1$ implies saturatedness (\cite[Corollary 5.3]{Kullmann2007ClausalFormZII}).

In \cite[Corollary 5.10]{Kullmann2007ClausalFormZII} one finds a generalisation of the basic bound $\minvdeg(F) \le 2 \delta(F)$ for the boolean case. Namely $\minvdeg(F) \le \max_{v \in \var(F)} \abs{D_v} \cdot \delta(F)$ is shown for $F$ with $\stabpar(F) \ge 1$. Since for (generalised) saturated $F$ with $\delta(F) = 1$ we have $\stabpar(F) = \infty$ (\cite[Corollary 5.6]{Kullmann2007ClausalFormZII}), this covers the above mentioned result $\minvdeg(F) \le \max_{v \in \var(F)} \abs{D_v}$ for (arbitrary) minimally unsatisfiable $F$ with $\delta(F) = 1$ (note that here saturation works as in the boolean case).

In \cite{KullmannZhao2014Smarandache} we concentrate on unsatisfiable hitting (generalised) clause-sets, and via generalised non-Mersenne numbers $\nonmer^d(k)$ we are able to generalise Theorem \ref{thm:MUminvdegdef} to generalised clause-sets. We believe that in general the minimum variable-degree of minimally unsatisfiable clause-sets $F$ with $\stabpar(F) \ge 1$ for a given deficiency is always obtained by unsatisfiable hitting clause-sets (generalising Conjecture \ref{con:sharpness}):
\begin{conj}\label{con:minvardeghit}
  Let $\Uclashi{\delta=k}^{\le d}$ denote the set of generalised unsatisfiable hitting clause-sets of deficiency $k \in \NN$ and with domain-sizes at most $d \in \NN$, and let $\Musati{\delta=k,\stabpar \ge 1}^{\le d}$ be defined in the same way. Then we have for all $k, d \in \NN$ that $\minvdeg(\Uclashi{\delta=k}^{\le d}) = \minvdeg(\Musati{\delta=k,\stabpar \ge 1}^{\le d})$.
\end{conj}

Furthermore, the ``$2$'' in $S_2(k)$ in Lemma \ref{lem:smarandache} comes from the boolean domain, and generalising the results of this \Schrift{} and from \cite{KullmannZhao2015FullClauses} in \cite{KullmannZhao2014Smarandache} to the non-boolean domain sheds light on $S_d(k)$ (the minimal $s \in \NN$ such that $d^k$ divides $s!$) for arbitrary prime numbers $d \in \NN$, as introduced in \cite[Unsolved Problem 49]{Smarandache1993OPNS} (while for non-prime-numbers $d$ the definition of $S_d$ has to be generalised). See \cite[Subsection III.1]{Ibstedt1998ComputerNumber} for basic properties of $S_p(k)$.

\subsection{Classification of $\Musat$}
\label{sec:concclassmu}

As mentioned in the introduction, a major motivation for us is the project of the classification of minimally unsatisfiable clause-sets in the deficiency (recall Examples \ref{exp:MU1}, \ref{exp:MU2}), where the main conjecture is:
\begin{conj}\label{con:classmu}
  For every deficiency $k \in \NN$ there are finitely many ``patterns'' which determine the nonsingular elements of $\Musati{\delta=k}$, as well as the saturated and hitting cases amongst them. Especially for every $k$ the isomorphism types of $\Musatnsi{\delta=k}$ can be efficiently enumerated (without repetitions), and for any given $F \in \Musatnsi{\delta=k}$ its isomorphism type can be determined in polynomial time.
\end{conj}
Conjecture \ref{con:classmu} has been shown for $k \le 2$ (recall Examples \ref{exp:MU1}, \ref{exp:MU2}). So, as a small example for the applications of Conjecture \ref{con:classmu}, if we have given some $F \in \Musat$ with $\delta(F) \le 2$, and want to know what is the ``cause'' of inconsistency, then in case of $\delta(F) = 1$ we can declare this as a ``trivial direct contradiction'', just $\bot$ concealed by non-strict subsumption extensions. While in case of $\delta(F) = 2$ we have precisely one basic pattern, a cycle establishing equivalences (of unique length), plus the requirement of some non-equivalences along the cycle. To see what ``patterns'' in Conjecture \ref{con:classmu} mean in general, the next step, classification of $\Musatnsi{\delta=3}$, is of great importance (a non-trivial task, as it seems). We remark that \cite{KleineBueningXu2005HomomorphismsMu} shows the necessity to consider nonsingular elements, since for every deficiency $k \ge 1$ the isomorphism-problem in $\Musati{\delta=k}$ is polytime-equivalent to (full) graph isomorphism.

As we discussed in Subsection \ref{sec:prelimcolouring}, the translation $e: \Cls \ra \Hyp$ has the property $F \in \Musati{\delta=k} \Lra e(F) \in \Mnc_{\defh=k-1}$ for $k \in \NN$, and so the classification of $\Musati{\delta=k}$ is a subtask of the classification of $\Mnc_{\defh=k-1}$. The possibility of a characterisation of $\Mnc_{\defh=0}$ was already raised in \cite{AhLi86} (where concentration on the special case of saturated (``strong'' there) minimally non-2-colourable hypergraphs was recommended), but is indeed still outstanding, which is understandable, given that polytime decision of $\Musati{\delta=1}$ is easy when compared with polytime decision of $\Mnc_{\defh=0}$. In the other direction, going more special than more general, the classification of $\Musati{\hdef=k} \subset \Musati{\delta \le k}$ (recall the hermitian defect $\hdef \ge \delta$ from Subsection \ref{sec:introbicliques}) could be a stepping stone (recall $\Musati{\hdef=1} = \Uclashi{\delta=1}$).

A major step towards Conjecture \ref{con:classmu} should be the classification of unsatisfiable (nonsingular) hitting clause-sets in dependency on the deficiency, i.e., determining the isomorphism types of $\Uclashnsi{\delta=k}$. We remark here that unsatisfiable hitting clause-sets do not seem to have a close correspondence in hypergraph colouring, due to the lack of complementation in hypergraphs. The main conjecture is the Finiteness Conjecture (a special case of Conjecture \ref{con:classmu}):
\begin{conj}\label{con:finhit}
  For every deficiency $k \in \NN$ there are only finitely many isomorphism types of nonsingular unsatisfiable hitting clause-sets, or equivalently, the number of variables of elements of $\Uclashnsi{\delta=k}$ is bounded.
\end{conj}
For $k \le 2$ finiteness has been established (Examples \ref{exp:MU1}, \ref{exp:MU2}), while recently we proved it for $k=3$ (\cite{KullmannZhao2016UHitSAT}), exploiting ``clause-factorisations'', which generalise singular DP-reduction and full subsumption resolution, and where the irreducible case corresponds to ``irreducible covering systems'' as investigated in \cite{Korec1984Covers,BergerFelzenbaumFraenkel1990Covers}. Assuming Conjecture \ref{con:finhit}, the question arises about the computability of the function, which maps $k \in \NN$ to the set of isomorphism types. Equivalently one can consider the computability of any function, which maps $k \in \NN$ to an upper bound on the number of variables of elements of $\Uclashnsi{\delta=k}$. It is conceivable that such functions grow so quickly that they are not computable, we however believe that a small bound holds, and we conjecture the following strengthened form of Conjecture \ref{con:finhit}:
\begin{conj}\label{con:uppboundfinhit}
  For every $k \ge 2$ holds $\max \set{n(F) : F \in \Uclashnsi{\delta=k}} =  4 k - 5$.
\end{conj}
In terms of Theorem \ref{thm:minnumvarmu} that means $n(\Uclashnsi{\delta=k}) \sse \tb{\odef(k)}{4k-5}$, and we furthermore conjecture equality here. By Example \ref{exp:MU2} indeed for $k=2$ the maximal number of variables is $4 \cdot 2 - 5 = 3$, and in \cite{KullmannZhao2016UHitSAT} we prove Conjecture \ref{con:uppboundfinhit} for $k=3$. We obtain computability of $\minnonmer$ as follows (using  Corollary \ref{cor:wlognonsingsat}):
\begin{lem}\label{lem:conjimplcomp}
  Assume Conjecture \ref{con:uppboundfinhit} holds. Then $k \in \NN \mapsto \minnonmerh(k)$ is computable, by enumerating all possible clause-sets $F$ with at most $4k-5$ variables, checking whether they are in $\Uclashnsi{\delta=k}$, and if so, including $\minvdeg(F)$ into the maximum-computation. If also Conjecture \ref{con:minvdeghit} holds, then also $\minnonmer$ is computable.
\end{lem}
Conjecture \ref{con:minnonmernotcomp} says additionally, that $\minnonmer$ should be ``complex''.

\bibliographystyle{elsarticle-num}

\newcommand{\noopsort}[1]{}

\newpage

\appendix
\addtocontents{toc}{\protect\setcounter{tocdepth}{2}}

\section{Overview on notations}
\label{sec:appOverview}

We recap here notations, and give links to their definitions, in four sections:
\begin{enumerate}
\item \ref{sec:appSets} is about sets (like $\Musat$).
\item \ref{sec:appMeasures} is about measures (like $\delta(F)$).
\item \ref{sec:appNumerical} is about numerical quantities (like $\minnonmer(k)$).
\item \ref{sec:appOperations} is about operations (like $\var(F)$).
\end{enumerate}
Each section is subdivided into three subsections, reviewing first the notations already discussed in the Introduction (Section \ref{sec:intro}), either in preliminary form or as Introduction-only, and then considering the main text (starting with Section \ref{sec:prelim}).

\subsection{Sets}
\label{sec:appSets}

\subsubsection{Preliminary definitions in Introduction}

$\Cls, \Sat, \Usat, \Musat$ (clause-sets, especially satisfiable, unsatisfiable, and minimally unsatisfiable ones) are defined in Subsection \ref{sec:introdef}. $\Lean$ (lean clause-sets) and $\Vmusat$ (variable-minimally unsatisfiable clause-sets) are defined in Subsection \ref{sec:introsurp}, $\Smusat$ (saturated minimally unsatisfiable clause-sets) in Subsection \ref{sec:prelimMU}, $\Musatns$ (nonsingular MUs) in Subsection \ref{sec:prelimTovey's}, while $\Mlean$, $\Msat$ (matching lean/satisfiable clause-sets) are introduced in Subsection \ref{sec:introautgen}. $\Clash$ (hitting clause-sets) and $\Uclash$ (unsatisfiable hitting clause-sets) are mentioned in Subsections \ref{sec:prelimQCA}, \ref{sec:introbicliques}.

\subsubsection{Definitions only for Introduction}

$\Pcls{p}$ (clause-sets with clauses of length at most $p$), $\Unicls$ (uniform clause-sets, i.e., all clauses have the same length), $\Punicls{p}$ ($p$-uniform clause-sets, i.e., all clauses have length $p$), $\Punimusat{p}$ ($p$-uniform minimally unsatisfiable clause-sets) are used in Subsection \ref{sec:prelimTovey's}.

$\Llean$ (linearly lean clause-sets) and $\Lsat$ (linearly satisfiable clause-sets) are used in Subsection \ref{sec:introautgen}.

$\Hyp$ (hypergraphs), $\Poscls$ (positive clause-sets), $\Mnc[k]$ (minimally non-$k$-colourable hypergraphs) are discussed in Subsection \ref{sec:prelimcolouring}. $\Ihyp$ (intersecting hypergraphs) are discussed in Subsection \ref{sec:introtrans}, together with \emph{bihitting clause-sets}.

\subsubsection{Main text}

$\Va$ (variables), $\Lit$ (literals), $\Cl$ (clauses) and $\Cls$ (clause-sets) are introduced in Subsection \ref{sec:prelimcls}, as well as $\Clash$ (hitting clause-sets). $\Pass$ (partial assignments), $\Sat$ (satisfiable clause-sets), $\Usat$ (unsatisfiable clause-sets), $\Uclash$ (unsatisfiable hitting clause-sets) are introduced in Subsection \ref{sec:prelimsem}.

The underlined versions of sets of clause-sets in general (e.g., $\ul{\Cls}$) as the corresponding sets of multi-clause-sets are discussed in Subsection \ref{sec:prelimmulti}.

$\Lean$ (lean clause-sets), $\Mlean$ (matching-lean clause-sets), $\Msat$ (matching satisfiable clause-sets) are introduced in Subsection \ref{sec:prelimAut}.

$\Musat$ (minimally unsatisfiable clause-sets), $\Smusat$ (saturated minimally unsatisfiable clause-sets), $\Musatns$ (nonsingular minimally unsatisfiable clause-sets), $\Smusatns$ (nonsingular saturated minimally unsatisfiable clause-sets) and $\Uclashns$ (nonsingular unsatisfiable hitting clause-sets) are introduced in Subsection \ref{sec:musubcl}. $\Vmusat$ (variable-minimally unsatisfiable clause-sets) is introduced in Section \ref{sec:vmusat}.

$J$ (jump positions for $\nonmer(k)$) is defined in Definition \ref{def:jump}.

$\Sed$ (clause-sets where the surplus equals the deficiency) is introduced in Subsection \ref{sec:classSED} (Definition \ref{def:sed}). $\Mlcr$ (matching-lean clause-sets which are ``critical'' concerning the open task of ``finding the autarky'') is introduced in Subsection \ref{sec:findaut} (Definition \ref{def:Mlcr}).

$\svbf$ (valid bounds-functions) and $\svbfs$ (valid bounds-functions bounded by the non-Mersenne numbers) are introduced in Subsection \ref{sec:analyspsit} (Definition \ref{def:validboundsf}).

For index-notations like $\Cls_{\delta=5}$ (all clause-sets $F \in \Cls$ with $\delta(F) = 5$) see Definition \ref{def:classescls}.

\subsection{Measures for clause-sets}
\label{sec:appMeasures}

\subsubsection{Preliminary definitions in Introduction}

$c(F), n(F), \delta(F)$ (number of clauses/variables and deficiency) are introduced in Subsection \ref{sec:introdef}. Also $\ldeg_F(x)$ (literal-degree), $\vdeg_F(v)$ (variable-degree) and $\minvdeg(F)$ (minimum variable-degree) are defined there. The surplus $\surp(F)$ is defined in Subsection \ref{sec:introsurp}.

\subsubsection{Definitions only for Introduction}

$\delta^*(F)$ (maximal deficiency) is introduced in Subsection \ref{sec:introdef}.

$\maxvdeg(F)$ (maximum variable-degree) is discussed in Subsection \ref{sec:prelimTovey's}.

$\defh(G)$ (deficiency of hypergraphs) is introduced in Subsection \ref{sec:prelimcolouring}.

$\bcp(A), \bcp(F)$ (the biclique partition number of matrix $A$ resp.\ clause-set $F$) is considered in Subsection \ref{sec:introbicliques}, together with $h(A), h(F)$ (hermitian rank of matrices/clause-sets) and $\hdef(F)$ (hermitian defect of clause-sets).

\subsubsection{Main text}

$c(F), n(F), \ell(F), \delta(F)$ (number of clauses/variables/literal occurrences and deficiency) are defined in Subsection \ref{sec:prelimcls} (and for multi-clause-sets in Subsection \ref{sec:prelimmulti}).

$\ldeg_F(x), \vdeg_F(v)$ (literal/variable-degree) are defined in Subsection \ref{sec:prelimDegrees}, together with $\minvdeg(F)$ (minimum variable-degree) and $\minvdeg(\mc{C})$ for a clause-set $F$ resp.\ a class $\mc{C}$ of clause-sets (Definition \ref{def:minvdeg}). $\surp(F)$ (surplus) is introduced in Subsection \ref{exp:prelimsurplus} (Definition \ref{def:surp}).

$\minldeg(F)$ (minimum literal-degree) is mentioned in Section \ref{sec:genboundml}. $\nfc(F)$ (number of full clauses) and $\maxnfc(\mc{C})$ (maximum number of full clauses for a class of clause-sets) is introduced in Section \ref{sec:strengtheningbound} (Definition \ref{def:numfcl}).

\subsection{Numerical quantities}
\label{sec:appNumerical}

\subsubsection{Preliminary definitions in Introduction}

$\minnonmer(k)$ (maximum of minimum variable-degrees for minimally unsatisfiable clause-sets of deficiency $k$) is introduced in Subsection \ref{sec:introdef}. The function $\nonmer(k)$ is discussed in Subsections \ref{sec:introdef}, \ref{sec:introsurp}, \ref{sec:introbasicint}.

\subsubsection{Definitions only for Introduction}

$\maxvdeg(\Punimusat{p}) = f(p)+1$ (minimum of maximum variable-degrees for $p$-uniform minimally unsatisfiable clause-sets) is discussed in Subsection \ref{sec:prelimTovey's}.

$m(p)$ (minimum number of hyperedges in $p$-uniform minimally non-2-colourable hypergraphs) and generalisations ($m(p,k), m^*(p,k)$) are mentioned in Subsection \ref{sec:prelimcolouring}.

\subsubsection{Main text}

$\ld(x)$ (binary logarithm), $\fld(x)$ (truncated binary logarithm) are defined at the beginning of Section \ref{sec:prelim}.

$\minnonmer(k), \minnonmerh(k)$ (maximum min-var-degrees for MUs of deficiency $k$) are defined in Subsection \ref{sec:musubcl} (Definition \ref{def:nonmerstar}).

$\odef(k)$ (minimum number of variables for MUs which can reach deficiency $k$)  is defined in Theorem \ref{thm:minnumvarmu}. Then $\nonmer(k)$ is the topic of Section \ref{sec:nonmer} (Definition \ref{def:minvdegdef}), with helper functions $\inonmer(k)$ (critical index) defined in Definition \ref{def:i(k)}, and helper functions $i'(k)$ (critical index as it appears in the recursion) and $h(k)$ (the recursive application) defined in Definition \ref{def:twoaux}.

$\maxsmar(k), \maxsmarh(k)$ (maximal number of full clauses for MUs/UHITs of deficiency $k$) are introduced in Section \ref{sec:strengtheningbound} (Definition \ref{def:maxsmar}).

$\nonmer_1(k)$ (strengthened non-Mersenne numbers) are introduced in Subsection \ref{sec:sharpeningbound} (Definition \ref{def:nm1}).

$S_2(k)$ (number-theoretic function) is discussed in Subsection \ref{sec:concSmar}.

\subsection{Operations}
\label{sec:appOperations}

\subsubsection{Preliminary definitions in Introduction}

Partial assignments $\pao v0, \pao v1$ and their application $\pao v{\ve} * F$ to clause-sets $F$ are introduced at the beginning of the introduction. Partial assignments $\vp$ in general (but as clauses) and their application $\vp * F$ are discussed in Subsection \ref{sec:prelimMU}. 

The set $\var(F)$ of variables in a clause-set is introduced in Subsection \ref{sec:introdef}. Also used is $\var(x)$, the variable of literal $x$, $\var(C)$, the set of variables in clause $C$, and $\var(\vp)$, the variables assigned by partial assignment $\vp$.

The unsatisfiable full clause-sets $A_n$ are explained in Subsection \ref{sec:introbasicint}.

Union $A \addcup B$ in case of disjointness is mentioned in Subsection \ref{sec:prelimcolouring}

\subsubsection{Definitions only for Introduction}

For a hypergraph $G$, by $V(G)$ the vertex-set and by $E(G)$ the hyperedge-set is denoted (Subsection \ref{sec:prelimcolouring}).

The translations of hypergraphs $G$ to clause-sets $F_2(G)$, and of clause-sets $F$ to hypergraphs $e(F)$, are discussed in Subsection \ref{sec:prelimcolouring}, where also the $k$-core of a hypergraph is mentioned.

For a hypergraph $G$ by $\Tr(G)$ the transversal hypergraph is denoted (Subsection \ref{sec:introtrans}).

For a matrix $M$ the qualitative class of $M$ is $\Q(M)$ (Subsection \ref{sec:prelimQCA}).

$n_+(M), n_-(M)$ are the numbers of positive/negative eigenvalues of matrix $M$ (Subsection \ref{sec:introbicliques}). Also discussed there is $\cmg(F)$, the conflict multigraph of clause-set $F$, while $\scf(F)$ is its conflict matrix.

\subsubsection{Main text}

In Subsection \ref{sec:prelimcls} complementation $\ol{x} \in \Lit,\ol{L} \sse \Lit$ of literals and sets of literals is introduced (where for literals $x \in \ZZ \sm \set{0}$ we use $-x = \ol{x}$). Furthermore we have $\var(x) \in \Va$, $\var(C) \subset \Va$, $\var(F) \subset \Va$ for variables of literals, clauses and clause-sets, while $\lit(F) \subset \Lit$ is the set of literals related to clause-set $F$.

Special constructions are $\top \in \Cls$ (empty clause-set) and $\bot \in \Cl$ (empty clause). Furthermore the full clause-set over $V \subset \Va$ is $A(V)$, with the special case $A_n$ for $V = \tb 1n$.

In Subsection \ref{sec:prelimsem} then $\var(\vp) \subset \Va$ and $\lit(\vp) \subset \Lit$ are introduced, the variables and literals of partial assignments $\vp$, together with $\vp^{-1}(\ve)$, the literals set to $\ve$ by $\vp$. Special constructions are $\epa \in \Pass$ (empty partial assignment) and $\pao{v}{\ve} \in \Pass$ (partial assignment $v \mapsto \ve$). The application of partial assignments to clause-sets is $\vp * F \in \Cls$. And semantical implication is denoted by $F \models C$.

In Subsection \ref{sec:prelimResolution} the (partial) resolution operation $C \res D \in \Cl$ is introduced, followed by the DP-operator $\dpi{v}(F) \in \Cls$.

The restriction $F[V] \in \Cls$ of a (multi-)clause-set $F$ to variable-set $V$ is defined in Definition \ref{def:restriction}.

The clause-sets $\Dt{n} \in \Musatnsi{\delta=2}$ are given in Example \ref{exp:MU2}.

A single potential saturation step $F \in \Cls \leadsto \saturate(F,C,x) \in \Cls$ is given in Definition \ref{def:saturation}.

In Section \ref{sec:elimcreatesing} singular DP-reduction in one step $F \tsdp F'$ and in many steps $F \tsdps F'$ is introduced (Definition \ref{def:sdp}). Furthermore the results of full reduction are collected in $\sdp(F) \subset \Musatns$.

In Section \ref{sec:2subr} we find strict full subsumption resolution in one step $F \tsubres F'$, in $k$ steps $F \tsubresk{k} F'$, and in many steps $F \tsubress F'$ (Definition \ref{def:2subsumptionres}).

$\Delta a$ (Definition \ref{def:delta}) for a sequence $a$ of numbers is the difference between neighbour terms.

$\varmvd(F) \sse \var(F)$ is the set of variables in $F$ realising the minimum degree $\minvdeg(F)$ (Definition \ref{def:varsetmvd}).

$F^V \le F$ for a multi-clause-set $F$ and a set of variables $V$ is the sub-multi-clause-set consisting of all clauses $C \in F$ with $\var(C) \sse V$ (Subsection \ref{sec:classSED}).

The full clause-set $M(V) \in \Mlean$ is used in Section \ref{sec:genboundml}.

$[a_1,\dots,a_p] \in \svbf$ (Definition \ref{def:validboundsf}) is the valid bounds-function given by these initial values. The set $\potp_f(k,m)$ of potential degree-pairs is introduced in Definition \ref{def:possiblepairs}. The improvement of a valid bounds-function $f \in \svbf$ to $\potprec(f) \in \svbf$ is defined in Definition \ref{def:potprec}.

\end{document}